\documentclass[11pt,a4paper]{article}

\usepackage[colorlinks]{hyperref}

\usepackage[table]{xcolor}
\usepackage{slashed}
\usepackage{xkeyval}
\usepackage{url,amsmath,ifthen,amssymb,latexsym,pstricks,mathrsfs,comment,amsthm,graphicx,tikz,tikz-cd,enumerate,accents,pgffor,cite,wrapfig,multicol,float,mathtools}

\usepackage{bibspacing}

\usepackage{enumitem}
\usepackage[margin=10pt,font=small,labelfont=bf, labelsep=period]{caption}

\usepackage[T1]{fontenc}
\usepackage{ wasysym }
\usepackage{stmaryrd}

\usepackage{geometry} 

\newlength{\normalparindent}
\AtBeginDocument{\setlength{\normalparindent}{\parindent}}
\newcommand{\longintertext}[1]{%
  \intertext{%
    \parbox{\linewidth}{%
      \setlength{\parindent}{\normalparindent}
      \noindent#1%
    }%
  }%
}

\begin{document}

\newcommand{\nc}{\newcommand}
\nc{\rnc}{\renewcommand}

\nc\oijn{1\leq i<j\leq n}
\nc\vt\vartheta
\nc\ew\iota
\nc\C{\mathcal C}
\nc\CA{\C(A)}

\nc\wh\hat
\nc\U{\mathcal U}
\nc\UC{\mathcal U_C}
\nc\phiC{\phi_C}
\nc\phiL{\phi_L}
\nc\UL{\mathcal U_L}
\nc\V{\mathcal V}
\nc\M{\mathcal M}
\nc\cS{\mathcal S}

\nc\chap\section
\nc\sect\subsection

\nc\da\downarrow
\nc\LR{\mathcal L}
\nc\Ef{E^\flat}
\nc\PB{\mathcal{PB}}
\nc\Pf{P^\flat}

\nc{\custpartn}[3]{{\lower1.4 ex\hbox{
\begin{tikzpicture}[scale=.3]
\foreach \x in {#1}
{ \uvert{\x}  }
\foreach \x in {#2}
{ \lvert{\x}  }
#3 \end{tikzpicture}
}}}
\nc{\uvert}[1]{\fill (#1,2)circle(.2);}
\rnc{\lvert}[1]{\fill (#1,0)circle(.2);}
\newcommand{\darcx}[3]{\draw(#1,0)arc(180:90:#3) (#1+#3,#3)--(#2-#3,#3) (#2-#3,#3) arc(90:0:#3);}
\newcommand{\uarcx}[3]{\draw(#1,2)arc(180:270:#3) (#1+#3,2-#3)--(#2-#3,2-#3) (#2-#3,2-#3) arc(270:360:#3);}
\nc{\darc}[2]{\darcx{#1}{#2}{.4}}
\nc{\uarc}[2]{\uarcx{#1}{#2}{.4}}
\nc\udotted[2]{\draw[dotted] (#1,2)--(#2,2);}
\nc\ddotted[2]{\draw[dotted] (#1,0)--(#2,0);}

\nc{\lv}[1]{\fill (#1,0)circle(.1);}
\nc{\lvs}[1]{{\foreach \x in {#1}{\lv{\x}}}}

\nc{\uv}[1]{\fill (#1,2)circle(.1);}
\nc{\uvs}[1]{{\foreach \x in {#1}{\uv{\x}}}}

\nc{\stline}[2]{\draw(#1,2)--(#2,0);}
\nc{\stlines}[1]{{\foreach \x/\y in {#1}{\stline{\x}{\y} }}}

\nc\Th\Theta
\nc\bS{{\bf S}}
\nc\bone{{\bf 1}}
\nc\bzero{{\bf 0}}
\nc\supp{\operatorname{supp}}

\nc{\noinclude}[1]{#1} 

\rnc\S{S}
\nc\F{\mathcal F}

\nc\Max{\operatorname{Max}}
\nc\FXG{F_X(G)}
\nc\FYG{F_Y(G)}
\nc\FXM{F_X(M)}
\nc\FXI{F_X(I)}
\nc\FYM{F_Y(M)}

\nc\ben{\begin{enumerate}[label=\textup{(\roman*)},leftmargin=7mm]}
\nc\bena{\begin{enumerate}[label=\textup{(\alph*)},leftmargin=7mm]}
\nc\een{\end{enumerate}}

\let\oldproofname=\proofname
\rnc{\proofname}{\rm\bf{\oldproofname}}

\nc\fin[1]{#1^{\operatorname{f}}}

\nc\Aut{\operatorname{Aut}}
\nc\FAut{\fin{\operatorname{Aut}}}
\nc\FPAut{\fin{\operatorname{PAut}}}
\nc\FEnd{\fin{\operatorname{End}}}
\nc\FSub{\fin{\operatorname{Sub}}}
\nc\FPEnd{\fin{\operatorname{PEnd}}}
\nc\End{\operatorname{End}}
\nc\PAut{\operatorname{PAut}}
\nc\PEnd{\operatorname{PEnd}}
\nc\Sing{\operatorname{Sing}}
\nc\codim{\operatorname{codim}}
\nc\Si\Sigma
\nc\mr\mathrel
\nc\lam\lambda
\nc\bn{{\bf n}}
\nc\bu{{\bf u}}
\nc\bv{{\bf v}}
\nc{\T}{\mathcal T}
\nc{\PT}{\mathcal P\mathcal T}
\nc\PTnSn{\PT_n\sm\G_n}
\nc\PTnTn{\PT_n\sm\T_n}
\nc{\ze}{\zeta}
\nc{\al}{\alpha}
\nc{\be}{\beta}
\nc{\ga}{\gamma}
\nc{\de}{\delta}
\nc\ka\kappa
\nc\ve\varepsilon
\rnc\th\theta
\nc{\rank}{\operatorname{rank}}
\nc{\corank}{\operatorname{corank}}
\nc{\Sub}{\operatorname{Sub}}
\nc{\dom}{\operatorname{dom}}
\nc{\im}{\operatorname{im}}
\nc{\R}{\mathrel{\mathscr R}}
\rnc{\L}{\mathrel{\mathscr L}}
\rnc{\H}{\mathrel{\mathscr H}}
\nc{\D}{\mathrel{\mathscr D}}
\nc{\J}{\mathrel{\mathscr J}}
\rnc{\implies}{\ \Rightarrow\ }
\nc{\sub}{\subseteq}
\nc\bl{{\bf l}}
\nc\br{{\bf r}}
\nc\bk{{\bf k}}
\nc\bp{{\bf p}}
\nc\bq{{\bf q}}
\nc\bs{{\bf s}}
\nc\bt{{\bf t}}
\nc\bz{{\bf z}}
\nc\trans[1]{\left(\begin{smallmatrix}#1\end{smallmatrix}\right)}

\nc\G{\mathcal G}
\nc\FG{\fin{\mathcal G}}
\nc\bm{{\bf m}}
\nc\Lam\Lambda
\nc\Om\Omega
\nc\De\Delta
\nc\E{\mathcal E}
\nc\I{\mathcal I}
\nc\Ga\Gamma
\newcommand{\restr}{\mathord{\restriction}}

\nc{\pre}{\preceq}
\rnc{\iff}{\ \Leftrightarrow\ }
\nc{\pfcase}[1]{\medskip \noindent {\bf Case #1.}}
\nc{\pfcasens}[1]{\noindent {\bf Case #1.}}
\nc{\pfstep}[1]{\medskip \noindent {\bf Step #1.}}
\nc{\es}{\varnothing}

\rnc{\P}{{\mathcal P}}

\nc{\cofix}{\operatorname{cofix}}
\nc{\fix}{\operatorname{fix}}
\nc{\Fix}{\operatorname{Fix}}
\nc{\la}{\langle}
\nc{\ra}{\rangle}
\nc{\si}{\sigma}
\nc\vs\varsigma
\nc{\sm}{\setminus}
\nc{\set}[2]{\{ {#1} : {#2} \}} 
\nc{\pres}[2]{\la {#1} : {#2} \ra} 
\nc{\Spres}[2]{\textup{\textsf{Sgp}}\la {#1} : {#2} \ra} 
\nc{\Mpres}[2]{\textup{\textsf{Mon}}\la {#1} : {#2} \ra} 
\nc\ol\overline
\nc\ul\underline
\nc{\bigset}[2]{\big\{ {#1} : {#2} \big\}} 
\nc{\pf}{\begin{proof}}
\nc{\epf}{\end{proof}}
\nc{\epfres}{\hfill\qed}
\nc{\epfreseq}{\tag*{\qed}}
\nc{\epfeq}{\tag*{\qed}}
\nc{\N}{\mathbb N}
\nc{\Z}{{\mathbb Z}}
\nc{\mt}{\mapsto}
\nc{\id}{\mathrm{id}}
\nc{\COMMA}{,\qquad}
\nc{\COMMa}{,\ \ \ }
\nc{\OR}{\qquad\text{or}\qquad}
\nc{\AND}{\qquad\text{and}\qquad}
\nc{\ANDSIM}{\qquad\text{and similarly}\qquad}
\nc{\ANd}{\quad\text{and}\quad}
\nc\WHERE{\qquad\text{where}\qquad}
\nc{\ba}{{\bf a}}
\nc{\bb}{{\bf b}}
\nc{\bc}{{\bf c}}
\nc{\bd}{{\bf d}}
\nc{\bx}{{\bf x}}

\nc{\bit}{\begin{itemize}}
\nc{\eit}{\end{itemize}}
\nc{\bmc}{\begin{multicols}}
\nc{\emc}{\end{multicols}}
\nc{\itemit}[1]{\item[\emph{(#1)}]}
\nc{\itemnit}[1]{\item[(#1)]}
\nc{\pfitem}[1]{\medskip \noindent #1.}
\nc{\firstpfitem}[1]{#1.}

\numberwithin{equation}{section}

\newtheorem{thm}[equation]{Theorem}
\newtheorem{lemma}[equation]{Lemma}
\newtheorem{cor}[equation]{Corollary}
\newtheorem{prop}[equation]{Proposition}

\theoremstyle{definition}

\newtheorem{defn}[equation]{Definition}
\newtheorem{rem}[equation]{Remark}
\newtheorem{eg}[equation]{Example}
\newtheorem{ass}[equation]{Assumption}

\title{Product decompositions of semigroups induced by action pairs~\vspace{-5ex}}

\date{}
\author{}

\maketitle
\begin{center}
{\large 
Scott Carson,%
\hspace{-.2em}\footnote{\label{fn:SC}Department of Mathematics, University of York, York YO10 5DD, U.K. {\it Emails:} {\tt sc1485@york.ac.uk}, {\tt victoria.gould@york.ac.uk}}
Igor Dolinka,%
\hspace{-.2em}\footnote{Department of Mathematics and Informatics, University of Novi Sad, Trg Dositeja Obradovi\'ca 4, 21101 Novi Sad, Serbia. {\it Email:} {\tt dockie@dmi.uns.ac.rs}}
James East,%
\hspace{-.2em}\footnote{Centre for Research in Mathematics and Data Science, Western Sydney University, Locked Bag 1797, Penrith NSW 2751, Australia. {\it Email:} {\tt j.east@westernsydney.edu.au}}
Victoria Gould,\hspace{-.2em}\textsuperscript{\ref{fn:SC}}
Rida-e Zenab%
\footnote{Department of Mathematics, Sukkur I.B.A.~University, Pakistan.  {\it Email:} {\tt ridaezenab@iba-suk.edu.pk}}
}
\end{center}

\begin{abstract}
This paper concerns a class of semigroups that arise as products $US$, associated to what we call `action pairs'.  Here $U$ and $S$ are subsemigroups of a common monoid and, roughly speaking, $S$ has an action on the monoid completion $U^1$ that is suitably compatible with the product in the over-monoid.

The semigroups encapsulated by the action pair construction include many natural classes such as inverse semigroups and (left) restriction semigroups, as well as many important concrete examples such as transformational wreath products, linear monoids, (partial) endomorphism monoids of independence algebras, and the singular ideals of many of these.  Action pairs provide a unified framework for systematically studying such semigroups, within which we build a suite of tools to ensure a comprehensive understanding of them. We then apply our abstract results to many special cases of interest.

The first part of the paper constitutes a detailed structural analysis of semigroups arising from action pairs.  We show that any such semigroup $US$ is a quotient of a semidirect product $U\rtimes S$, and we classify all congruences on semidirect products that correspond to action pairs.  We also prove several covering and embedding theorems, each of which naturally extends celebrated results of McAlister on proper (a.k.a.~$E$-unitary) inverse semigroups.  

The second part of the paper concerns presentations by generators and relations for semigroups arising from action pairs.  We develop a substantial body of general results and techniques that allow us to build presentations for $US$ out of presentations for the constituents $U$ and $S$ in many cases, and then apply these to several examples, including those listed above.  Due to the broad applicability of the action pair construction, many results in the literature are special cases of our more general ones.



\end{abstract}

\tableofcontents

\chap{Introduction}\label{chap:intro}

The motivation for the current work comes from many different directions.  The basic philosophical underpinning is the desire to understand mathematical structures by sets of data that are sufficient to determine the entire structure, yet are intrinsically simpler.  In particular, we draw inspiration from classical techniques and constructions in group theory, and seek to understand the (much) more complicated and wider context of these techniques and constructions within the realm of semigroups and monoids.

One approach is to use \emph{presentations}. Here the idea is to define an algebraic structure---often, and originally, a group---in terms of generators and relations between them. It is possible in this way to capture some infinite groups using a finite amount of information: i.e., via a finite number of generators and a finite number of relations, resulting in a finitely presented group.  Important examples include the Artin braid groups \cite{Artin1947,Artin1925} and Thompson's groups $F$, $T$ and $V$  \cite{CFP1996}.  On the other hand, it is also useful to determine finite groups via (finite) presentations.  For example, the well-known Coxeter presentation of a finite symmetric group \cite{Moore1897} tells us how to calculate all products by dint of the products between certain transpositions; cf.~\cite{Humphreys1990}.  Group presentations were introduced as far back as the mid-nineteenth century \cite{Hamilton1856,Dyck1882}, and have since been central in approaches to algebraic and differential topology, geometric group theory, algebraic combinatorics, representation theory and many other branches of mathematics. Semigroup presentations go back to the mid-twentieth century, and again have proven a crucial tool in the above fields.  For some early studies, see for example \cite{Aizenstat1962,Aizenstat1958,Popova1961,Popova1962}, and for early papers with connections to logic see \cite{Turing1950,Markov1947,Post1947}; some historical information can be found in \cite{NB2021}.

A second approach to the above philosophy involves the concept of \emph{products}, which covers a wide range of mathematical constructions prominently featuring in the areas of algebra, topology, category theory, mathematical logic and graph theory, just to name a few.  The idea behind the introduction and use of such constructions is two-fold. Firstly, one of the most basic tasks of almost any coherent mathematical field or theory is to describe the many ways in which the objects of study are constructed; in many cases, the philosophy is to go from smaller, simpler structures towards larger, more complicated ones. This is exactly how various types of products permeate mathematics: think of direct products in general algebra, semidirect, wreath, and free (amalgamated) products in group theory, product topologies, tensor products in linear algebra, ultraproducts in model theory, cap and cup products in algebraic topology, various types of graph products, pullbacks in category theory, and so on. However, there is also the reverse task: to understand the structure of a given mathematical object in as much detail as possible, one attempts to break it down into more elementary parts and describe the mutual relationships between the parts. Quite often, this happens by recognising that the given structure is isomorphic (or in some other way `equivalent' or `similar') to a product of some kind. Smaller structures obtained as a result of such a decomposition are often expected or required to be \emph{sub}structures of the original, but not always.

A flagship instance of the phenomenon just described is the decomposition of finite abelian groups into the direct product of cyclic (sub)groups of prime power order. For groups, and indeed monoids, there is an exact correlation between internal and external direct products.  However, as soon as one moves away from monoids, even just as far as semigroups, this correlation is lost.  In particular,  an external direct product $S\times T$ of semigroups $S$ and $T$ need not contain subsemigroups isomorphic to either $S$ or $T$, yet clearly $S\times T$ is determined by $S$ and $T$.  

Another natural class of examples illustrating the same phenomenon---but where the smaller parts are not necessarily substructures---are the \emph{almost-factorisable inverse semigroups} of Lawson~\cite{Lawson1994}.  A semigroup $S$ from this class decomposes as a product~$S=EG$, where~$G$ is the group of units of some \emph{different} inverse monoid $M$, and where~$E$ is the semilattice of singular (i.e., non-identity) idempotents of $M$.  Here we have $E\sub S$, yet~$G\not\sub S$; in fact, we even have $G\cap S=\es$.  Nevertheless, the structure of $S$ is closely governed by that of the constituents $E$ and~$G$, and in fact $S$ is an idempotent-separating homomorphic image of a semidirect product~$E\rtimes G$.  

Yet another natural example is the semigroup $S$ of singular (i.e., non-bijective) partial transformations of a finite set $X$.  This semigroup can be decomposed as $S=ET$, where $E$ is the (monoid) semilattice of all partial identities, and $T$ is the semigroup of singular (full) transformations of~$X$.  This time we have $T\sub S$ but $E\not\sub S$, although we do not have the extreme $E\cap S=\es$; rather,~$E\sm S$ consists only of the identity map.  \emph{Many} more such examples exist.

The main goal of the current study is to provide a unified framework for working with a large class of structures including the kind just described; presentations will form an important set of tools.  The central notion throughout is that of an \emph{action pair} $(U,S)$.  Such a pair consists of subsemigroups $U$ and $S$ of a common monoid $M$, such that $S$ acts on the monoid completion $U^1=U\cup\{1\}$ in a way that is suitably compatible with the product in $M$.  The formal definitions are given later in the paper, including various other conditions that are required.  But it is worth noting that `compatibility' here means that
\[
s\cdot u = {}^su\cdot s \qquad\text{for all $s\in S$ and $u\in U^1$.}
\]
(In the above expression, $\cdot$ is the product in $M$, and the action of $S$ on $U^1$ is denoted by $(s,u)\mt{}^su$.)
It follows quickly from this that the set product
\[
US=\set{us}{u\in U,\ s\in S}
\]
is a subsemigroup of $M$.  It is important to note that $US$ might contain neither $U$ nor $S$, nor even isomorphic copies of either.  On the other hand, if $U$ and $S$ are both submonoids of $M$, then $US$ contains both $U$ and $S$ as submonoids.

An important natural class of examples extending those discussed above come from (subsemigroups of) inverse semigroups \cite{Lawson1998}, and more generally left restriction semigroups \cite{Gould_notes} .  If $S$ is an inverse semigroup with semilattice of idempotents $E$, then $(E,S)$ is an action pair in the monoid completion $S^1$, with respect to the conjugation action: ${}^se=ses^{-1}$ for $e\in E^1$ and $s\in S$.  More generally, $(F,T)$ is an action pair for any subsemigroup $T\leq S$ and any subsemilattice~$F\leq E$ for which $F^1$ is closed under conjugation by elements of~$T$.  (Almost-)factorisable inverse monoids and semigroups \cite{CH1974,Lawson1994} are (very) special cases of this.  Other similar examples involve products of idempotents and one- or two-sided units of arbitrary monoids \cite{JE2020a}.  Many other examples are considered throughout the paper, some with very different behaviour to those discussed so far.  The class of semigroups arising from action pairs is extremely rich, and part of our motivation stems from the desire to obtain a comprehensive understanding of such semigroups, analogous to that of more specialised classes such as inverse and (left) restriction semigroups.  

As well as this introduction, and the preliminary Chapter \ref{chap:prelim}, the paper consists of two main parts:
\bit
\item Part \ref{part:I} constitutes a detailed structural analysis of semigroups arising from action pairs,~and
\item Part \ref{part:II} concerns presentations by generators and relations for semigroups in this class.
\eit
We now give a very brief overview of the main highlights of the paper.  More detailed summaries can be found in the opening remarks of both parts, and in the introductions to the chapters within.

As we have just indicated, the first part of the paper is devoted to analysing the structure of a semigroup $US$ arising from an action pair $(U,S)$.  Chapter \ref{chap:ASD} contains a number of preliminary results on actions and semidirect products that will be used extensively throughout the rest of the paper.  Our notion of semidirect product is slightly more general than the usual one, which is essential to capture the widest possible collection of applications.
Chapter~\ref{chap:AP} introduces action pairs, provides an extensive collection of examples (and non-examples for contrast), and contains the first main structural result:
\bit
\item Theorem \ref{thm:FT} shows that $US$ is a homomorphic image of a semidirect product $U\rtimes S$, and classifies all congruences on semidirect products that correspond to action pairs.
\eit
Chapter \ref{chap:PAP} concerns the special class of \emph{proper} action pairs; these extend the notion of proper (a.k.a.~$E$-unitary) inverse semigroups, and more generally proper left restriction semigroups.  The main results of this chapter are broad generalisations of McAlister's celebrated Covering Theorem \cite{McAlister1974} and P-Theorem \cite{McAlister1974b} for inverse semigroups.  These new results encapsulate many extensions of McAlister's theorems to other more specialised classes:
\bit
\item Theorem \ref{thm:cover} shows that $US$ has a \emph{proper cover}, by which we mean that there is a proper action pair $(U',S')$ in a suitable monoid, with ${U'\cong U}$ and $S'\cong S$, and a natural surmorphism~${U'S'\to US}$.  Further, if $S$ or $U$ is a submonoid of the monoid containing $US$, then the surmorphism~${U'S'\to US}$ restricts to an isomorphism $U'\to U$ or $S'\to S$, respectively.
\item Theorem \ref{thm:embed2} shows that a proper \emph{monoid} $US$ embeds in a semidirect product~${\U\rtimes(S/\si)}$, where $\U$ is a monoid containing $U$, and where~$\si$ is a special congruence on~$S$ akin to the least group congruence on an inverse semigroup.  Theorems \ref{thm:embed1}, \ref{thm:embed3} and \ref{thm:embed4} are variations on this result, and show that when $U$ has certain additional structural properties (such as being commutative, a semilattice or a left-regular band), then an embedding $US\to\U\rtimes(S/\si)$ of the above kind exists where~$\U$ also has the relevant structural property.
\eit
As applications, we give new proofs of the generalisations of McAlister's above-mentioned theorems to the class of left restriction semigroups \cite{OCarroll1976,Lawson1986,FGG2009,BGG2010}.

In the second part of the paper we turn our attention to presentations by generators and relations.  As discussed above, presentations are extremely important tools when working with any kind of algebraic structure, and several results exist for building presentations for algebras that arise from others via natural constructions.  See for example \cite{RRW1998,Lavers1998,EEF2005,FADEG2019,RRT2003,HR1994,DR2009,Ruskuc1999}; the introduction to \cite{JEptc} contains a fuller discussion and many more references.
Chapter \ref{chap:pres} contains many general results on presentations for a semigroup $US$ arising from an action pair $(U,S)$ in several important cases:
\bit
\item Theorems \ref{thm:ESmon0}, \ref{thm:ESmon}, \ref{thm:Msimp}, \ref{thm:M} and \ref{thm:ESmon00} concern the case in which $U$ and $S$ are both submonoids of the over-monoid.  As we will see, one complication that must be overcome is the fact that the semidirect product $U\rtimes S$ need not be a monoid in this case.
\item Theorems \ref{thm:ES}, \ref{thm:simp} and \ref{thm:US_sd} only assume that one of $U$ or $S$ is a submonoid.  Other conditions must be satisfied in these cases in order to obtain usable results.
\eit
The general machinery developed in Chapter \ref{chap:pres} is then applied to a number of important examples in the remaining chapters:
\bit
\item Free left restriction monoids are treated in Chapter \ref{chap:LR}; see Theorem \ref{thm:LX}.
\item Chapter \ref{chap:IA} concerns several monoids and semigroups of (partial) endomorphisms of an independence algebra.  See especially Theorem \ref{thm:SubA}, which gives a presentation for the $\cap$-semilattice of finite-codimensional subalgebras of an arbitrary strong independence algebra.  
\item A number of transformational wreath products are covered in Chapter \ref{chap:wreath}, the main results here being Theorems \ref{thm:MwrSingPTn}, \ref{thm:MwrPTn} and \ref{thm:MwrGn}--\ref{thm:MwrIn}.
\eit
The potential for further applications is vast.

The article also includes several other results of independent interest.  For example:
\bit
\item In Section \ref{sect:SD} we prove a number of results on (our more general) semidirect products, including the fact that a semidirect product $U\rtimes S$ is a monoid if and only if $U$ and $S$ are both monoids, with $S$ acting monoidally on $U$ by monoid morphisms; see Corollary \ref{cor:USmon}.
\item Theorem \ref{thm:US_sd} gives a presentation for a semidirect product $U\rtimes S$ in the case that $S$ is a monoid and $U$ an arbitrary semigroup.  This complements \cite[Theorem 3.1]{FADEG2019}, which treats the reverse case, where $U$ is a monoid and $S$ a semigroup.  By contrast, \cite[Corollary~2]{Lavers1998} gives a presentation for $U\rtimes S$ when it is a monoid (which occurs when the conditions discussed in the previous point hold).  Our Theorem~\ref{thm:M} extends this result, by relaxing the assumption that $S$ acts by \emph{monoid} morphisms, and applies to the largest submonoid of $U\rtimes S$ with identity $(1,1)$.
\item In Section \ref{sect:LRpres} we show that the free left restriction monoid over any set has a unique minimum (monoid) generating set; see Proposition \ref{prop:GaX}, Theorem \ref{thm:LX} and Remark \ref{rem:LX}.
\item In Section \ref{sect:MaxA}, we prove several results on independence algebras with suitably large subalgebras, and show that these are precisely the algebras displaying the most `freedom' in (partial) automorphisms; see especially Propositions \ref{prop:sub1}--\ref{prop:sub3}.
\item In Section \ref{sect:AutA}, we give a number of results on generating sets for groups of (finitary) automorphisms of independence algebras, which provide a common generalisation of the following two facts: (i) finite symmetric groups are generated by transpositions, and (ii) finite-dimensional general linear groups are generated by `elementary row operation' matrices.  See Theorems \ref{thm:n-2} and \ref{thm:n-1}.
\eit
Throughout the text we pose a number of open problems and enticing directions for future studies.

\subsection*{Acknowledgements}

The first author was supported by an EPSRC doctoral fellowship.  The second author is partially supported by the Ministry of Education, Science, and Technological Development of the Republic of Serbia. The third author is supported by Future Fellowship FT190100632 of the Australian Research Council.

\chap{Preliminaries}\label{chap:prelim}

This article is intended to be largely self-contained, and accessible to a broad audience of algebraists.  Accordingly, we now gather the basic definitions and facts concerning semigroups and presentations we need, and fix most of the notation we will be using.  
For more background on semigroups, and for proofs of the various assertions below, see for example \cite{Howie1995,Higgins1992,CPbook}.  For inverse semigroups in particular, see \cite{Lawson1998} or \cite[Chapter 5]{Howie1995}; for (left) restriction semigroups, see \cite{Gould_notes}.

\sect{Semigroups}\label{sect:sgp}

A \emph{semigroup} is a set with an associative binary operation, typically denoted by juxtaposition.  A \emph{monoid} is a semigroup with an identity element.  Unless otherwise specified, the identity of any monoid is denoted by $1$.  

A \emph{semigroup morphism} is a map $\phi:S\to T$, where $S$ and $T$ are semigroups, and where $(xy)\phi=(x\phi)(y\phi)$ for all $x,y\in S$.  A monoid morphism is a semigroup morphism between monoids that additionally maps the identity to the identity.

If $S$ is a semigroup, then $S^1$ denotes the \emph{monoid completion} of $S$.  Formally:
\bit
\item $S^1=S$ if $S$ is already a monoid; 
\item otherwise, $S^1=S\sqcup\{1\}$, where $1$ is a symbol not belonging to $S$, and acting as an identity for $S^1$.  (Here and elsewhere $\sqcup$ denotes disjoint union.)
\eit
If $S$ happens to be a subsemigroup of a specified monoid $M$, then we often assume the identity of $S^1$ is the same as that of $M$; we will always be clear about this.

An element $x$ of a semigroup $S$ is \emph{(von Neumann) regular} if $x=xax$ for some $a\in S$.  For $b=axa$, we then have $x=xbx$ and $b=bxb$.  A semigroup is \emph{(von Neumann) regular} if all its elements are regular.  An element $x\in S$ is an \emph{idempotent} if $x=x^2$.  We generally write $E(S)$ for the set of all idempotents of $S$.  Clearly any idempotent is regular.  A semigroup $S$ is a \emph{band} if $S=E(S)$: i.e., if every element of $S$ is an idempotent.  A \emph{semilattice} is a commutative band.  A band $S$ is \emph{left-regular} if $xyx=xy$ for all $x,y\in S$; right-regular bands are defined symmetrically, but they will not arise in our investigations.  Left-regular bands are of course generalisations of semilattices, but they are important for many reasons, one of which is that they include many semigroups of geometric significance \cite{MSS2015,MSS2014}.

A semigroup $S$ is \emph{inverse} if for every $x\in S$ there is a \emph{unique} element $a\in S$ such that $x=xax$ and $a=axa$; this element is then denoted $a=x^{-1}$.  Equivalently, a semigroup $S$ is inverse if it is regular and its idempotents commute; the set $E(S)$ is therefore a semilattice.  Inverse semigroups can also be defined as a variety of \emph{unary semigroups} (i.e., semigroups equipped with an additional basic unary operation); they are precisely the class of semigroups with a unary operation $x\mt x^{-1}$ satisfying the identities
\[
(x^{-1})^{-1} = x = xx^{-1} x \COMMA (xy)^{-1} = y^{-1}x^{-1} \AND xx^{-1} yy^{-1} = yy^{-1} xx^{-1}.
\]

A relation $\si$ on a semigroup $S$ is \emph{left-compatible} if $(x,y)\in\si \implies (ax,ay)\in\si$ for all $a\in S$.  \emph{Right-compatibility} is defined symmetrically, and a relation is \emph{compatible} if it is both left- and right-compatible.  A \emph{(left or right) congruence} is an equivalence relation that is (left- or right-) compatible.  We will need the following basic result a number of times; it is surely well known, but we include a short proof for convenience.

\begin{lemma}\label{lem:RC}
If $\si$ is a right congruence on a monoid $M$, then the $\si$-class of $1$ is a submonoid of $M$.
\end{lemma}

\pf
Let $T$ be the $\si$-class of $1$.  Clearly $1\in T$.  Now let $a,b\in T$.  Since $a,b\mr\si1$, and since $\si$ is a right congruence, we have $ab\mr\si1b=b\mr\si1$.
\epf

As in \cite{Ault1974,Yamada1955}, a \emph{mid-identity} (a.k.a.~\emph{mid-unit} or \emph{middle unit}) of a semigroup $S$ is an element $a\in S$ such that $xay=xy$ for all $x,y\in S$.  

A \emph{retraction} is a morphism $\phi:S\to T$, where~$S$ is a semigroup, $T\leq S$ a subsemigroup, and $\phi\restr_T=\id_T$.  Obviously retractions are surjective.

\newpage

\begin{lemma}\label{lem:MIS}
Suppose $a$ is a mid-identity of a semigroup $S$, and let $T\leq S$.
\ben
\item The following diagram commutes, with all sets subsemigroups of $S$, and all maps surmorphisms:
\[
\begin{tikzcd} 
~ & T \arrow[swap]{dl}{f_1:x\mt xa} \arrow{dr}{f_2:x\mt ax} & \\
Ta \arrow[swap]{dr}{g_1:x\mt ax} & & aT \arrow{dl}{g_2:x\mt xa}\\
& aTa & 
\end{tikzcd}
\]
\item If $Ta\leq T$, then $aTa\leq aT$, and $f_1$ and $g_2$ are retractions.
\item If $aT\leq T$, then $aTa\leq Ta$, and $f_2$ and $g_1$ are retractions.
\een
\end{lemma}

\pf
This is all easily checked, and we just give a couple of sample calculations.  To show that~$f_1$ is a morphism, let $x,y\in T$.  Then
\[
(xf_1)\cdot(yf_1) = xa\cdot ya = xay \cdot a = xy\cdot a = (xy)f_1.  
\]
If $Ta\leq T$, then for any $t\in T$ we have $(ta)f_1=taa=ta$, so that $f_1$ is a retraction.
\epf

\begin{rem}
Although $Ta$, $aT$ and $aTa$ are subsemigroups of $S$, they might not be subsemigroups of $T$.  See Remark \ref{rem:M2} for a specific example of this, where $aT\cap T=\es$; in fact, $aT=S\sm T$.
\end{rem}

\begin{rem}
Lemma \ref{lem:MIS} also follows from \cite[Theorem 2.15]{DE2018}, which concerns \emph{sandwich semigroups} in (locally small) categories.  The simplest case of a sandwich semigroup is a \emph{semigroup variant} \cite{Hickey1986,Hickey1983}.  The variant of a semigroup $S$ with respect to an element $a\in S$ is the semigroup $S^a=(S,\star_a)$, where the `sandwich operation' $\star_a$ is defined by $x\star_ay=xay$ for all $x,y\in S$.  Of course if $a$ is a mid-identity of $S$, then $\star_a$ is precisely the original operation in $S$.
\end{rem}

We also need the following simple result.

\begin{lemma}\label{lem:sl}
Let $S$ be a semigroup, and suppose $P\leq S$ is a semilattice.  If $x,y\in S$ are such that $[x\in Py$ and $y\in Px]$ or $[x\in yP$ and $y\in xP]$, then $x=y$.
\end{lemma}

\pf
By symmetry, we assume that $x\in Py$ and $y\in Px$, so that $x=py$ and $y=qx$ for some $p,q\in P$.  But then
\[
x = py = pqx = qpx = (qq)px = q(qpx) = qx = y.  \qedhere
\]
\epf

Green's $\R$, $\L$, $\J$, $\H$ and $\D$ relations on a semigroup $S$ \cite{Green1951} will not play an explicit role in our study, but since we do mention them in passing from time to time it is convenient to define them here.  The first three of these relations are defined, for $x,y\in S$, by
\[
x\R y \iff xS^1=yS^1 \COMMA x\L y \iff S^1x=S^1y \AND x\J y \iff S^1xS^1=S^1yS^1.
\]
These can also be characterised in terms of divisibility.  For example, $x\R y$ if and only if $x=ya$ and $y=xb$ for some $a,b\in S^1$; this is in turn equivalent to having either $x=y$ or else $x=ya$ and $y=xb$ for some $a,b\in S$.  The final two of Green's relations are defined by
\[
{\H} = {\R}\cap{\L} \AND {\D} = {\R}\vee{\L},
\]
where the latter denotes the join of $\R$ and $\L$ in the lattice of all equivalence relations on $S$: i.e.,~$\D$ is the least equivalence containing ${\R}\cup{\L}$.  It is well known that in fact ${\D}={\R}\circ{\L}={\L}\circ{\R}$.  For some classes of semigroups (some of) Green's relations can be characterised \emph{equationally}.  For example, if $S$ is inverse, then $x\R y\iff xx^{-1}=yy^{-1}$.

\sect{Transformation semigroups}\label{sect:ts}

Certain semigroups of (partial) transformations will play an important role throughout, so we revise their definitions here.  Let $X$ be an arbitrary set.  A \emph{partial transformation} of $X$ is a function $A\to X$ for some $A\sub X$.  The set $\PT_X$ of all such partial transformations of $X$ is a monoid under ordinary relational composition, called the \emph{partial transformation monoid}.  

As usual, we denote by $\dom(\al)$ and $\im(\al)$ the \emph{domain} and \emph{image} of a partial transformation $\al\in\PT_X$, with their standard meanings.  
For $x\in\dom(\al)$ we write $x\al$ for the image of $x$ under~$\al$, and we compose partial transformations from left to right.  
Note, for example, that
\[
\dom(\al\be) = \dom(\be)\al^{-1} = \bigset{x\in\dom(\al)}{x\al\in\dom(\be)} \qquad\text{for $\al,\be\in\PT_X$.}
\]
The \emph{kernel} of $\al\in\PT_X$ is the equivalence
\[
\ker(\al) = \bigset{(x,y)\in\dom(\al)\times\dom(\al)}{x\al=y\al},
\]
and the \emph{rank} of $\al$ is the cardinal
\[
\rank(\al)=\big|{\im(\al)}\big|=\big|{\dom(\al)/\ker(\al)}\big|.
\]
Here $\dom(\al)/\ker(\al)$ is the quotient of the set $\dom(\al)$ by the equivalence $\ker(\al)$: i.e., the set of all $\ker(\al)$-classes.
Important submonoids of $\PT_X$ include:
\bit
\item $\T_X=\set{\al\in\PT_X}{\dom(\al)=X}$, the \emph{full transformation monoid},
\item $\I_X=\set{\al\in\PT_X}{\al \text{ is injective}}$, the \emph{symmetric inverse monoid},
\item $\G_X=\set{\al\in\T_X}{\al \text{ is bijective}}$, the \emph{symmetric group}.
\eit
As the name suggests, $\I_X$ is an inverse monoid.  The inverse of $\al\in\I_X$ is the ordinary inverse mapping $\al^{-1}$, and the idempotents of $\I_X$ are precisely the partial identity maps $\id_A$ ($A\sub X$).  These are not the only idempotents of $\PT_X$ when $|X|\geq2$.

When $X$ is finite, the sets
\[
\Sing(\PT_X) = \PT_X\sm\G_X \COMMA
\Sing(\T_X) = \T_X\sm\G_X \AND
\Sing(\I_X) = \I_X\sm\G_X 
\]
are subsemigroups, indeed two-sided ideals, of $\PT_X$, $\T_X$ and $\I_X$, respectively.  For arbitrary $X$, we also have the subsemigroup $\PT_X\sm\T_X$ of all \emph{strictly} partial transformations of $X$, which is in fact a right ideal of $\PT_X$.  When $X=\{1,2,\ldots,n\}$ for some integer $n$, we typically denote $\PT_X$ by $\PT_n$, and similarly for $\T_n$, $\I_n$ and so on.  We often make use of the standard two-line notation for partial transformations.  For example, $\al = \trans{1&2&3&4&5&6\\2&-&3&2&6&6} \in \PT_6$ has domain $\{1,3,4,5,6\}$, and maps $1\mt2$, $3\mt3$, and so on.

\sect{Left restriction semigroups}\label{sect:lrs}

\emph{Left restriction semigroups} arise in many different contexts and have many different names.
Essentially, one can define them in three different ways:
\bit
\item by a representation: up to isomorphism, they are precisely subsemigroups of partial transformation semigroups closed under the map $\al\mt\id_{\dom(\al)}$;  
\item as a generalisation of inverse semigroups:  they possess a distinguished semilattice of idempotents $E$ (not necessarily \emph{all} the idempotents), such that every $\widetilde{\R}_E$-class contains an idempotent of $E$, and the so-called `ample condition' holds; here $\widetilde{\R}_E$ is a relation containing Green's $\R$-relation;
\item they form a variety of unary semigroups (see below).
\eit
For more details and background we refer the reader to \cite{Gould_notes}.  It is convenient here to take the third option:

\begin{defn}
A \emph{left restriction semigroup} is a unary semigroup $S$ satisfying the following identities, where we write the unary operation as $s\mt s^+$:
\begin{enumerate}[label=\textup{(L\arabic*)}]
\bmc2
\item \label{L1} $x^+x=x$,
\item \label{L2} $x^+y^+=y^+x^+$,
\item \label{L3} $(x^+y)^+=x^+y^+$,
\item \label{L4} $xy^+=(xy)^+x$.
\emc
\een
\end{defn}

\emph{Right} restriction semigroups are defined dually, and there is also a notion of a \emph{(two-sided) restriction semigroup}; the latter have two interacting unary operations, but it is not necessary to give the definitions here.  Archetypal examples of left restriction semigroups include:
\bit
\item any inverse semigroup $S$, with $s^+=ss^{-1}$,  
\item the partial transformation monoid $\PT_X$ for any set $X$, with $\al^+=\id_{\dom(\al)}$.
\eit
If the left restriction semigroup $S$ is a monoid with identity $1$, it follows immediately from \ref{L1} that $1^+=1$.  Any monoid trivially becomes a left restriction monoid, upon defining $s^+=1$ for all $s$, but such structures are generally not of interest to us.

The identities \ref{L1}--\ref{L4} are rather compact, and some important properties are not immediately apparent.  For example (writing $=_1$ to indicate an application of \ref{L1}, etc.), we have
\begin{enumerate}[label=\textup{(L\arabic*)}]
\addtocounter{enumi}{4}
\item \label{L5} $x^+x^+=_3(x^+x)^+=_1x^+$, 
\een
so that each $x^+$ is an idempotent.  We then also have
\begin{enumerate}[label=\textup{(L\arabic*)}]
\addtocounter{enumi}{5}
\item \label{L6} $(x^+)^+ =_5 (x^+x^+)^+ =_3 x^+(x^+)^+ =_2 (x^+)^+x^+ =_1 x^+$,
\een
so that the operation ${}^+$ is itself idempotent.  

It follows from \ref{L3} that the set
\[
P(S) = \set{s^+}{s\in S}
\]
is a subsemigroup of $S$.  By \ref{L2} and \ref{L5}, $P(S)$ is in fact a semilattice.  By \ref{L6}, $P(S)$ is closed (indeed, fixed pointwise) under ${}^+$.  The elements of $P(S)$ are called \emph{projections}.

\sect{Presentations}\label{sect:pres}

We conclude this preliminary chapter by establishing the notation we use for presentations.  Presentations exist for any variety of universal algebras; here we explain how they work for semigroups and monoids, both for the sake of concreteness, and since these are our intended applications. 

Let~$X$ be a set, and~$X^+$ the \emph{free semigroup} on $X$, which consists of all non-empty words over~$X$ under the operation of concatenation.  Let $R\sub X^+\times X^+$ be a set of pairs of words, and write $R^\sharp$ for the congruence on $X^+$ generated by $R$.  We say a semigroup $S$ has \emph{(semigroup) presentation}~${\Spres XR}$ if $S\cong X^+/R^\sharp$, or equivalently if there is a surmorphism $X^+\to S$ with kernel~$R^\sharp$; if~$\phi$ is such a surmorphism, we say $S$ has \emph{presentation $\Spres XR$ via $\phi$}.  The elements of $X$ and $R$ are called \emph{generators} and \emph{relations}, respectively, and a relation $(u,v)\in R$ is sometimes displayed as an equation, $u=v$.  A \emph{normal form function} is a function $N:S\to X^+$ such that $N(s)\phi=s$ for all $s\in S$.  If we write ${\sim}=R^\sharp=\ker(\phi)$, then it follows from the definitions that
\begin{equation}\label{eq:N}
N(u\phi)\sim u \AND N(st)\sim N(s)N(t) \qquad\text{for all $u\in X^+$ and $s,t\in S$.}
\end{equation}

There is also a corresponding notion of monoid presentations.  The \emph{free monoid} on the set~$X$ is denoted by $X^*$, and is defined by $X^*=X^+\cup\{\ew\}$.  Here $\ew$ is the empty word, which is the identity of $X^*$.  We say a monoid~$S$ has \emph{(monoid) presentation} $\Mpres XR$, where $R\sub X^*\times X^*$, if $S\cong X^*/R^\sharp$.  We also speak of monoid presentations via surmorphisms $X^*\to S$ (with kernel~$R^\sharp$), and normal form functions~$N:S\to X^*$; typically we assume that $N(1)=\ew$.

On occasions, we will identify $\Spres XR$ with the semigroup $X^+/R^\sharp$ itself, and similarly for monoid presentations.  Consider such a semigroup $\Spres XR$, and another semigroup $T$.  Any function $f:X\to T$ extends uniquely to a morphism $\phi:X^+\to T$, which is defined by $(x_1\cdots x_k)\phi=(x_1f)\cdots(x_kf)$ for any $x_1,\ldots,x_k\in X$.  Standard notions of universal algebra tell us that if $\phi$ \emph{preserves $R$}, in the sense that $u\phi=v\phi$ for every $(u,v)\in R$, then $R^\sharp\sub\ker(\phi)$.  In this case, $\phi$ induces a morphism $\Phi:\Spres XR\to T$, defined by $[w]\Phi=w\phi$ for $w\in X^+$, where $[w]$ denotes the $R^\sharp$-class of $w$.

The next result has a simple proof, and is essentially folklore.  Although we state it for monoids (as that is the context in which we will most often apply it), it holds more generally for arbitrary classes of (universal) algebras admitting presentations, so in particular there is a semigroup version as well, whose explicit statement we omit.

\begin{lemma}\label{lem:ST}
Suppose a monoid $S$ has presentation $\Mpres XR$ via $\phi:X^*\to S$, and let ${N:S\to X^*}$ be a normal form function.  Suppose $\pi:S\to T$ is a surmorphism onto a monoid~$T$, and suppose $\ker(\pi)=\Om^\sharp$, where $\Om\sub S\times S$.  Then $T$  has presentation $\Mpres X{R\cup R_\Om}$ via $\phi\pi:X^*\to T$, where $R_\Om = \bigset{(N(s),N(t))}{(s,t)\in\Om}$.  \epfres
\end{lemma}

The next result also has a semigroup version, which we omit.  For the statement, recall that left-regular bands were defined in Section \ref{sect:sgp}.

\begin{lemma}\label{lem:pres}
Let $M=\Mpres XR$, where $X$ is an alphabet and $R\sub X^*\times X^*$.  Also, let
\[
R_1 = \bigset{(x^2,x)}{x\in X} \COMMa R_2 = \bigset{(xy,yx)}{x,y\in X} \ANd R_3 = \bigset{(xyx,xy)}{x,y\in X}.
\]
\ben
\item \label{pres1} If $R_2\sub R^\sharp$, then $M$ is commutative.
\item \label{pres2} If $R_1\cup R_2\sub R^\sharp$, then $M$ is a semilattice.
\item \label{pres3} If $R_1\cup R_3\sub R^\sharp$, then $M$ is a left-regular band.
\een
\end{lemma}

\pf
\firstpfitem{\ref{pres1}}  This is clear.

\pfitem{\ref{pres2}}  Since $R_2\sub R^\sharp$, $M$ is commutative.  Since $R_1\sub R^\sharp$, $M$ is idempotent-generated.  Since any idempotent-generated commutative semigroup is a semilattice, the result follows.

\pfitem{\ref{pres3}}  Write ${\sim}=R^\sharp$, and for $w\in X^*$ write $\ol w$ for the $\sim$-class of $w$.  We must show that
\bena\bmc2
\item \label{pres3a} $u^2\sim u$ for all $u\in X^*$, and
\item \label{pres3b} $uvu\sim uv$ for all $u,v\in X^*$.
\emc\een
In fact, it suffices to prove \ref{pres3b}, as \ref{pres3a} follows upon taking $v=\ew$.  We first show that
\bena\addtocounter{enumi}{2}
\item \label{pres3c} $xvx\sim xv$ for all $x\in X$ and $v\in X^*$.
\een
We prove this by induction on $k=\ell(v)$, the length (number of letters) of $v$.  If $k=0$ then $v=\ew$, and \ref{pres3c} says $x^2\sim x$, which holds because $R_1\sub R^\sharp$.  So now suppose $k\geq1$, and write $v=v'y$ where $v'\in X^*$ and $y\in X$.  Since $\ell(v')=k-1$, we have $xv'\sim xv'x$ by induction.  Since $R_3\sub R^\sharp$, we have $xy\sim xyx$.  It then follows that
\[
xv = xv'y \sim xv'xy \sim xv'xyx \sim xv'yx = xvx,
\]
completing the proof of \ref{pres3c}.

We now prove \ref{pres3b} by induction on $l = \ell(u)$.  If $l=0$ then $u=\ew$, and \ref{pres3b} says $v\sim v$, which is obviously true.  So now suppose $l\geq1$, and write $u=zu'$, where $u'\in X^*$ and $z\in X$.  Since $\ell(u')=l-1$, we have $u'v\sim u'vu'$ by induction.  By \ref{pres3c} we have $zu'v\sim zu'vz$.  It then follows that
\[
uv = zu'v \sim zu'vu' \sim zu'vzu' = uvu,
\]
completing the proof of \ref{pres3b}, and hence of the lemma.
\epf

\begin{rem}
If $R_1\sub R^\sharp$, then of course $M$ is idempotent-generated, but it need not be a band.
\end{rem}

\newpage

\part{Structure}\label{part:I}

This paper concerns a certain class of semigroups arising from what we will call \emph{action pairs}.  This first part of the paper is devoted to a structural analysis of such semigroups and pairs.

Roughly speaking, an action pair $(U,S)$ consists of two subsemigroups of a common over-monoid $M$, with an action of $S$ on $U^1$ that is suitably compatible with the operation of $M$.  Associated to such a pair is another subsemigroup of $M$, namely the (internal) product
\[
US=\set{us}{u\in U,\ s\in S}.
\]
The semigroup $US$ does not always contain $U$ or $S$ as subsemigroups, and does not even need to contain isomorphic copies of either.  Nevertheless, we will see that the structure of $US$ can be concisely described in terms of the structure of $U$ and $S$, the action of $S$ on $U^1$, and some other data.  We will also see that action pairs provide a natural context in which to extend a number of important classical results on inverse semigroups, left restriction semigroups, and more general structures.

There are three chapters in this part of the paper.  Chapter \ref{chap:ASD} contains preliminary material on actions and semidirect products, which will underpin all that follows.

In Chapter \ref{chap:AP}, we introduce action pairs (and more generally weak action pairs), discuss a number of examples (including the class of left restriction semigroups), and prove a structure theorem that will be used extensively in the second part of the paper when we study presentations.  This result, Theorem \ref{thm:FT}, states that a semigroup $US$ arising from an action pair $(U,S)$ is a homomorphic image of a semidirect product $U\rtimes S$, so that $US\cong(U\rtimes S)/\th$ for a suitable congruence $\th$.  The theorem also classifies the congruences on semidirect products leading to action pairs, showing how they are built from natural families of right congruences on $S$.

Chapter \ref{chap:PAP} contains a number of further structure theorems, all of which involve so-called \emph{proper} action pairs.  The resulting `product semigroups' $US$ include all proper left restriction semigroups, and in particular all proper (a.k.a.~$E$-unitary) inverse semigroups.  A very brief summary of the main results is as follows:
\ben
\item A semigroup $US$ arising from an arbitrary (weak) action pair $(U,S)$ can be \emph{covered} by a proper semigroup $U'S'$ with $U'\cong U$ and $S'\cong S$.  See Theorem \ref{thm:cover}.
\item Any proper semigroup $US$ can be naturally embedded in a semidirect product ${\U\rtimes(S/\si)}$, where $\U$ contains $U$, and where $\si$ is a conguence on $S$ akin to the least group congruence on an inverse semigroup.  See Theorem \ref{thm:embed2}.  When $U$ has additional structural properties (such as being commutative, a semilattice or a left-regular band), we can adjust our construction to ensure that $\U$ has the same structural property as well; see Theorems \ref{thm:embed1}, \ref{thm:embed3} and \ref{thm:embed4}.
\een
As applications of these results, we provide short proofs of the above-mentioned classical results on left restriction semigroups; this class includes inverse semigroups, but with a different interpretation of the unary operation.

\chap{Actions and semidirect products}\label{chap:ASD}

Actions and semidirect products are important tools in practically every part of algebra, and they will play a crucial role throughout the current work.  In this chapter we gather various facts that we need in the rest of the paper.

\sect{Actions}\label{sect:A}

A \emph{(left) action} of a semigroup $S$ on a set $M$ is a map $S\times M\to M:(s,x)\mt {}^sx$ such that
\[
{}^s({}^tx)={}^{st}x \qquad\text{for all $s,t\in S$ and $x\in M$.}
\]
All actions are assumed to be left actions.  Our choice of the (left) superscript notation for actions is for readability, particularly as many expressions below involve both actions and products; it is convenient to keep the former `vertical' and the latter `horizontal'.  If $S$ happens to be a monoid with identity $1$, then we say the action is \emph{monoidal} if 
\[
{}^1x=x \qquad\text{for all $x\in M$.}
\]
The actions that arise in later chapters will almost always be monoidal when $S$ is a monoid.  (A notable exception is Theorem \ref{thm:FT}.)  We will also generally be interested in the case that the action preserves an algebraic structure on~$M$.  If $M$ happens to be a semigroup, then the action is \emph{by semigroup morphisms} if   
\[
{}^s(xy)={}^sx\cdot{}^sy \qquad\text{for all for all $s\in S$ and $x,y\in M$.}
\]
If $M$ happens to be a monoid with identity $1$, then the action is \emph{by monoid morphisms} if 
\[
{}^s(xy)={}^sx\cdot{}^sy \AND {}^s1=1 \qquad\text{for all for all $s\in S$ and $x,y\in M$.}
\]
On many occasions in the paper, we will be concerned with the case in which a semigroup $S$ acts on a monoid $M$ by semigroup morphisms, but \emph{not} by monoid morphisms; see Example \ref{eg:act}\ref{act2} for an important/motivating case.  In such cases, the elements ${}^s1$ ($s\in S$) play a crucial role, so it is convenient to record some of their properties.  (All of these properties are trivial when the action is by monoid morphisms.)

\begin{lemma}\label{lem:+}
Suppose a semigroup $S$ acts on a monoid $M$ by semigroup morphisms, and for $s\in S$ write $s^+={}^s1$.  Then for all $s,t\in S$ and $x\in M$ we have
\ben
\bmc2
\item \label{lem+1} $s^+s^+=s^+$,
\item \label{lem+2} ${}^s(t^+)=(st)^+$,
\item \label{lem+3} $s^+(st)^+=(st)^+=(st)^+s^+$,
\item \label{lem+4} $s^+\cdot{}^sx = {}^sx = {}^sx\cdot s^+$.
\emc
\een
\end{lemma}

\pf
These are all easily established, and there are some interdependencies; for example \ref{lem+1} follows from \ref{lem+4} with $x=1$.  As an example calculation, $s^+\cdot{}^sx = {}^s1\cdot{}^sx = {}^s(1\cdot x)={}^sx$, giving the first half of \ref{lem+4}.
\epf

\begin{rem}
The reader will note that we have used the notation $s^+$ for the unary operation in a left restriction semigroup (cf.~Section \ref{sect:lrs}), and also for the elements ${}^s1$ arising from a semigroup action on a monoid.  This choice is deliberate, and it will transpire that there are deep connections between the two ideas; see especially Section \ref{sect:LRAP}.
\end{rem}

We conclude this section with some useful examples that we will often revisit.

\begin{eg}\label{eg:act}
\ben

\item \label{act1}  Consider arbitrary semigroups $S$ and $M$, such that $M$ contains an idempotent $e$.  Then ${}^sx=e$ ($s\in S$, $x\in M$) defines an action of $S$ on $M$ by semigroup morphisms, which we call the \emph{constant action with image $e$}.  If $S$ is a monoid, this action is monoidal if and only if $M=\{e\}$.  If $M$ is a monoid, the action is by monoid morphisms if and only if~$e=1$.

\item \label{act2}  Consider an inverse semigroup $S$, and let $E=E(S)$ be the semilattice of idempotents.  Then ${}^se=ses^{-1}$ defines an action of $S$ on $E^1$ by semigroup morphisms, as follows from standard facts about inverse semigroups; cf.~Section \ref{sect:sgp}.  For example,
\[
{}^s(ef) = sefs^{-1} = ss^{-1}sefs^{-1} = ses^{-1}sfs^{-1} = {}^se\cdot{}^sf \qquad\text{for all $s\in S$ and $e,f\in E^1$,}
\]
where we used the fact that the idempotents $e$ and $s^{-1}s$ commute.  (Taking $f=e=e^2$ in the above shows that indeed ${}^se\in E^1$ for $s\in S$ and $e\in E^1$.)
With respect to this action, we have $s^+={}^s1=ss^{-1}$ for all $s\in S$.  Thus, the action is by monoid morphisms if and only if $ss^{-1}=1$ for all $s\in S$, which is equivalent to $S$ being a group.  If $S$ happens to be a monoid, then the action is monoidal.

\item \label{act3}  Let $X$ be a non-empty set, and consider the power set $\P(X)=\set{A}{A\sub X}$ as a monoid (indeed, semilattice) under $\cap$, with identity $X$.  The partial transformation monoid $\PT_X$ (cf.~Section \ref{sect:ts}) acts monoidally on $\P(X)$ by semigroup morphisms, via
\[
{}^\al\!A=A\al^{-1} = \set{x\in\dom(\al)}{x\al\in A} \qquad\text{for $A\sub X$ and $\al\in\PT_X$.}
\]
This action is not by \emph{monoid} morphisms, however, as ${}^\al\!X=\dom(\al)$ for all $\al\in\PT_X$.  On the other hand, the induced action of the \emph{full} transformation monoid $\T_X\leq\PT_X$ \emph{is} by monoid morphisms, as $\dom(\al)=X$ for all $\al\in\T_X$.
\een
\end{eg}

\sect{Semidirect products}\label{sect:SD}

Typically, to define a semidirect product $U\rtimes S$, for semigroups $U$ and $S$, one begins with an action of $S$ on $U$, or equivalently a representation of $S$ by endomorphisms of $U$; see for example \cite{Zenab2018,Preston1986}.  However, in order to encompass all of our desired applications, we need a slightly more general notion:

\begin{defn}\label{defn:SD}
Suppose $U$ and $S$ are semigroups, and suppose $S$ has a left action on $U^1$ (the monoid completion of $U$) by semigroup morphisms, denoted
\[
\varphi:S\times U^1\to U^1:(s,u)\mt{}^su.
\]
The \emph{semidirect product} $U\rtimes S=U\rtimes_\varphi S$ is the semigroup
\begin{equation}\label{eq:sd}
U\rtimes S = \bigset{(u,s)}{u\in U,\ s\in S} \qquad\text{with operation}\qquad (u,s)\cdot(v,t) = (u\cdot{}^sv,st).
\end{equation}
Throughout the paper, we omit the $\varphi$ subscript, and write $U\rtimes S=U\rtimes_\varphi S$.  Context will always make it clear which action is in play.
\end{defn}

\begin{rem}
In \eqref{eq:sd}, it is possible to have ${}^sv=1\not\in U$, but we always have $u\cdot{}^sv\in U$.  

If $S$ acts on $U$ itself, then it also acts on $U^1$; if $U\not=U^1$, then we additionally define ${}^s1=1$ for all $s\in S$.  In this case, $U\rtimes S$ as in Definition~\ref{defn:SD} agrees with the ordinary semidirect product.  Throughout the paper, we will frequently be interested in actions on $U^1$ that do \emph{not} arise in this way.  (An action of a semigroup on a monoid $M$ might not restrict to an action on $M\sm\{1\}$, even if the latter is a subsemigroup of $M$.)

As an example, let $S=\T_X$ for an arbitrary set $X$ with $|X|\geq2$, and let $U=\P(X)\sm\{X\}$ be the $\cap$-semilattice of all \emph{proper} subsets of $X$.  Then the action of $S$ on $U^1=\P(X)$ given in Example~\ref{eg:act}\ref{act3} does not arise from an action on $U$.  (If $\al\in\T_X$ has image $A\subsetneq X$, then $A\in U$, but ${}^\al\!A=X\not\in U$.)  In this case, we will see that the semidirect product $U\rtimes S$ maps naturally onto $\PT_X\sm\T_X$, the semigroup of all \emph{strictly} partial transformations of $X$; cf.~Example~\ref{eg:Tn} and Proposition \ref{prop:A1}.
\end{rem}

For the rest of this section, we fix semigroups $S$ and $U$, and assume that $S$ acts on $U^1$ by semigroup morphisms.  We begin the following observation.

\begin{lemma}\label{lem:USmon0}
If $U\rtimes S$ is a monoid, then $U$ and $S$ are both monoids, and the identity of $U\rtimes S$ is~$(1,1)$, where as usual we write $1$ for the identity of both $U$ and $S$.
\end{lemma}

\pf
Suppose $U\rtimes S$ is a monoid with identity $(e,f)$.  So for all $u\in U$ and $s\in S$,
\begin{equation}\label{eq:usef}
(u,s) = (u,s)\cdot(e,f) = (u\cdot{}^se,sf) \ANDSIM (u,s) = (e\cdot{}^fu,fs).
\end{equation}
Looking at the second coordinates of \eqref{eq:usef}, it follows that $S$ is a monoid with identity $f=1$.  Looking at the first coordinates, and taking $s=1\in S$ (and remembering $f=1$), we have
\begin{equation}\label{eq:uusee1u}
u = u\cdot{}^1e = e\cdot{}^1u \qquad\text{for all $u\in U$.}
\end{equation}
But then for any $u\in U$, several applications of \eqref{eq:uusee1u} gives
\[
ue = e\cdot{}^1(ue) = e \cdot {}^1u \cdot {}^1e = u \cdot {}^1e = u,
\]
so that $e\in U$ is a right identity for $U$.  In particular, $e$ is an idempotent; combining this with~\eqref{eq:uusee1u}, it quickly follows that $e$ is also a left identity for $U$.
\epf

So $U\rtimes S$ being a monoid implies that $U$ and $S$ are both monoids, but the next result shows that the converse is not true in general.

Recall that we have a fixed action of $S$ on $U^1$.  This can be extended to an action of $S^1$ on~$U^1$ in an obvious way: if $S\not=S^1$, then we additionally define~${}^1u=u$ for all $u\in U^1$.  This action of $S^1$ is of course monoidal if $S\not=S^1$, but need not be monoidal if $S=S^1$ is a monoid.  In any case, we also have the semidirect product $U^1\rtimes S^1$, which contains $U\rtimes S$ as a subsemigroup.  

We now define $\bone=(1,1)\in U^1\rtimes S^1$.  One might initially expect $U^1\rtimes S^1$ to be a monoid with identity~$\bone$, but this is not always true.  (But if $U^1\rtimes S^1$ \emph{is} a monoid, then its identity is $\bone$, by Lemma \ref{lem:USmon0}.)  Rather, we have
\begin{equation}\label{eq:bone}
\bone\cdot(u,s)=({}^1u,s) \AND (u,s)\cdot\bone = (us^+,s) \qquad\text{for any $u\in U^1$ and $s\in S^1$,}
\end{equation}
where again we write $s^+={}^s1$ for $s\in S^1$.

\begin{lemma}\label{lem:USmon}
\ben
\item \label{USmon1} $\bone$ is a left identity for $U^1\rtimes S^1$ if and only if the action of $S^1$ on $U^1$ is monoidal,
\item \label{USmon2} $\bone$ is a right identity for $U^1\rtimes S^1$ if and only if the action of $S^1$ on $U^1$ is by monoid morphisms,
\item \label{USmon3} $U^1\rtimes S^1$ is a monoid with identity $\bone$ if and only if the action of $S^1$ on $U^1$ is monoidal and by monoid morphisms.  
\een
\end{lemma}

\pf
This all follows from \eqref{eq:bone}.  The only (slightly) non-obvious part is the forwards implication in \ref{USmon2}.  For this, we note that if $\bone$ is a right identity, then $(1,s)=(1,s)\cdot\bone=(s^+,s) \implies s^+=1$ for all $s\in S$.
\epf

\begin{rem}\label{rem:11id}
In general, $\bone$ need not be a left or right identity for $U^1\rtimes S^1$, and in fact it need not even be an idempotent, as $\bone\cdot\bone=({}^11,1)=(1^+,1)$.  It is quite possible to have $1^+\not=1$.  For example, this is the case for the constant action of $S=S^1$ on $U^1$ whose image is a non-identity idempotent; cf.~Example \ref{eg:act}\ref{act1}.
\end{rem}

Combining Lemmas \ref{lem:USmon0} and \ref{lem:USmon}, we immediately obtain the following neat characterisation of monoid semidirect products:

\begin{cor}\label{cor:USmon}
The following are equivalent:
\ben
\item $U\rtimes S$ is a monoid,
\item $U$ and $S$ are both monoids, and the action of $S(=S^1)$ on $U(=U^1)$ is monoidal and by monoid morphisms.  \epfres
\een
\end{cor}

Even though $\bone$ need not be a left or right identity for $U^1\rtimes S^1$ (cf.~Remark \ref{rem:11id}), it does nevertheless have the following very important property.  

\begin{lemma}\label{lem:MI}
The element $\bone$ is a mid-identity for $U^1\rtimes S^1$.
\end{lemma}

\pf
For any $u,v\in U^1$ and $s,t\in S^1$, we have
\[
(u,s)\cdot\bone\cdot(v,t) = (us^+,s)\cdot(v,t) = (us^+\cdot{}^sv,st) = (u\cdot{}^sv,st) = (u,s)\cdot(v,t),
\]
where we used Lemma \ref{lem:+}\ref{lem+4} in the third step.
\epf

\begin{rem}
We noted in Remark \ref{rem:11id} that $\bone$ might not be an idempotent of $U^1\rtimes S^1$, but now we have seen that it is always a mid-identity.  It is well known, and easy to see, that any \emph{regular} mid-identity must be an idempotent.
Thus, $U^1\rtimes S^1$ gives a neat example of a semigroup with a non-regular mid-identity (when $1^+\not=1$).  The square of any mid-identity is always an idempotent mid-identity, and here this is $\bone\cdot\bone=(1^+,1)$.
\end{rem}

The next result is not essential for our purposes, but it seems worth recording; it of course leads to an alternative proof of Lemma \ref{lem:MI}.

\begin{prop}
The mid-identities of $U^1\rtimes S^1$ are precisely the elements of the form $(e,1)$ for which ${}^1e$ is a left identity for the subsemigroup $\set{{}^1u}{u\in U^1}\leq U^1$.
\end{prop}

\pf
Throughout the proof we write $V=\set{{}^1u}{u\in U^1}$.  (This is a subsemigroup because ${}^1u\cdot{}^1v={}^1(uv)$ for all $u,v\in U^1$.)

Suppose first that $(e,f)$ is a mid-identity for $U^1\rtimes S^1$, and let $u\in U^1$ be arbitrary.  Then
\[
({}^1u,1) = \bone\cdot(u,1) = \bone\cdot(e,f)\cdot(u,1) = ({}^1e\cdot{}^fu,f).  
\]
It follows that $f=1$, and then also that ${}^1u = {}^1e\cdot{}^fu = {}^1e\cdot{}^1u$, so that ${}^1e$ is a left identity for $V$.

Conversely, suppose $e\in U^1$ is such that ${}^1e$ is a left identity for $V$, and let $(u,s),(v,t)\in U^1\rtimes S^1$ be arbitrary.  Then
\[
(u,s)\cdot(v,t) = (u\cdot{}^sv,st) \AND (u,s)\cdot(e,1)\cdot(v,t) = (u\cdot{}^se\cdot{}^sv,st) = (u\cdot{}^s(ev),st),
\]
so we must show that $u\cdot{}^s(ev)=u\cdot{}^sv$.  But this follows quickly from
\[
{}^s(ev) = {}^{s\cdot1}(ev) = {}^s({}^1(ev)) = {}^s({}^1e\cdot{}^1v) = {}^s({}^1v) = {}^{s\cdot1}v = {}^sv.  \qedhere
\]
\epf

The remaining results of this section concern certain natural subsemigroups of $U\rtimes S$.  We begin by defining
\[
M_1 = \bigset{(u,s)\in U\rtimes S}{u=us^+} \AND 
M_2 = \bigset{(u,s)\in U\rtimes S}{u={}^1u} .
\]
(We will soon see that these are indeed subsemigroups.  In $M_2$, note that ${}^1u$ is defined even if $1\not\in S$.)
We also define
\[
M = M_1\cap M_2 = \bigset{(u,s)\in U\rtimes S}{{}^1u=u=us^+}.
\]

\begin{lemma}\label{lem:1us+}
We have $M=\bigset{(u,s)\in U\rtimes S}{u={}^1us^+}$.
\end{lemma}

\pf
We must show that for any $u\in U$ and $s\in S$, we have ${}^1u=u=us^+ \iff u={}^1us^+$.  The forwards implication being clear, suppose $u={}^1us^+$.  Then
\begin{align*}
{}^1u &= {}^1({}^1us^+) = {}^1({}^1u) \cdot {}^1(s^+) = {}^{1\cdot1}u\cdot(1\cdot s)^+ = {}^1u\cdot s^+ = u\\
\AND
us^+ &= {}^1us^+\cdot s^+ = {}^1us^+ = u,
\end{align*}
where we used parts \ref{lem+1} and \ref{lem+2} of Lemma \ref{lem:+}.
\epf

The statements below refer to the element $\bone=(1,1)\in U^1\rtimes S^1$, which might not belong to $U\rtimes S$, and in particular might not belong to any of the subsets $M_1$, $M_2$ or $M$.  When $\bone$ does belong to $U\rtimes S$ (i.e., when $U$ and $S$ are monoids), it belongs to any of these subsets if and only if it belongs to all of them, and this is equivalent to having $1=1^+$.

\begin{lemma}\label{lem:M1}
\ben
\item \label{M11} $M_1=(U\rtimes S)\cdot\bone$ is a left ideal (and hence a subsemigroup) of $U\rtimes S$, and $\bone$ is a right identity for $M_1$ (but not necessarily an element of $M_1$).
\item \label{M12} $M_1=U\rtimes S$ if and only if
\begin{equation}\label{eq:M12}
u=us^+ \qquad\text{for all $u\in U$ and $s\in S$,}
\end{equation}
in which case $M=M_2$.
\een
\end{lemma}

\pf
\firstpfitem{\ref{M11}}  First let $(u,s)\in M_1$.  Then $(u,s)\cdot\bone=(us^+,s)=(u,s)$.  This shows that ${M_1\sub(U\rtimes S)\cdot\bone}$, and also that $\bone$ is a right identity for $M_1$.  

Conversely, suppose $(u,s)\in (U\rtimes S)\cdot\bone$, so that $(u,s)=(v,t)\cdot\bone=(vt^+,t)$ for some $v\in U$ and $t\in S$.  It follows that $u=vt^+\in U$ (as $v\in U$) and $s=t\in S$, so that $(u,s)\in U\rtimes S$.  Using Lemma~\ref{lem:+}\ref{lem+1}, it also follows that $us^+ = vt^+\cdot t^+ = vt^+  = u$, so that $(u,s)\in M_1$.

This completes the proof that $M_1=(U\rtimes S)\cdot\bone$.  Since $(U\rtimes S)\cdot\bone$ is clearly closed under left multiplication by elements of $U\rtimes S$, and since $M_1\sub U\rtimes S$ (by definition), it follows that $M_1$ is a left ideal of $U\rtimes S$.

\pfitem{\ref{M12}}  This is immediate from the definitions.
\epf

\begin{rem}\label{rem:M1}
Condition \eqref{eq:M12} is satisfied, for example, when the action of $S$ (or equivalently of $S^1$) on $U^1$ is by monoid morphisms (i.e., $s^+=1$ for all $s\in S$), but this is not necessary.  For example, \eqref{eq:M12} still holds if $U$ has a right identity $e$ that is not a left identity, and the action of~$S$ on $U^1(\not=U)$ is constant with image $e$, as then $us^+=ue=u$ for all $u\in U$ and $s\in S$.
\end{rem}

The situation for $M_2$ is slightly more complicated, as we do not always have $M_2=\bone\cdot(U\rtimes S)$:

\begin{lemma}\label{lem:M2}
\ben
\item \label{M21} $M_2\sub\bone\cdot(U\rtimes S)$, with equality if and only if
\begin{equation}\label{eq:M21}
{}^1u\in U \qquad\text{for all $u\in U$.}
\end{equation}
\item \label{M22} $M_2$ is a right ideal (and hence a subsemigroup) of $U\rtimes S$, and $\bone$ is a left identity for $M_2$ (but not necessarily an element of $M_2$).
\item \label{M23} $M_2=U\rtimes S$ if and only if 
\begin{equation}\label{eq:M23}
u={}^1u \qquad\text{for all $u\in U$,}
\end{equation}
in which case~$M=M_1$.  Moreover, \eqref{eq:M23} is equivalent to the action of $S^1$ on $U^1$ being monoidal.
\een
\end{lemma}

\pf
\firstpfitem{\ref{M21}}  If $(u,s)\in M_2 $, then $\bone\cdot(u,s) = ({}^1u,s) = (u,s)$.  This shows that ${M_2\sub\bone\cdot(U\rtimes S)}$, and also that $\bone$ is a left identity for $M_2$.  

Next we assume that \eqref{eq:M21} holds, aiming to prove the reverse inclusion in this case.  To do so, let ${(u,s)\in\bone\cdot(U\rtimes S)}$, so that $(u,s)=\bone\cdot(v,t) = ({}^1v,t)$ for some $v\in U$ and $t\in S$, so $(u,s)\in U\rtimes S$ by \eqref{eq:M21}.  We also have ${}^1u = {}^1({}^1v) = {}^{1\cdot1}v = {}^1v = u$, so that $(u,s)\in M_2$, as required.  

If, on the other hand, \eqref{eq:M21} does not hold, then ${}^1u\not\in U$ for some $u\in U$.  In this case, for any $s\in S$ we have ${\bone\cdot(u,s)=({}^1u,s)\not\in U\rtimes S}$, so certainly $\bone\cdot(u,s)\not\in M_2$.

\pfitem{\ref{M22}}  Since the reverse inclusion in part \ref{M21} does not hold in general, we must prove this directly, so let $(u,s)\in M_2$ and $(v,t)\in U\rtimes S$.  So ${}^1u=u$, and we have $(u,s)\cdot(v,t)=(u\cdot{}^sv,st)\in U\rtimes S$, with
\[
{}^1(u\cdot{}^sv) = {}^1u \cdot{}^1({}^sv) = u\cdot{}^{1\cdot s}v = u\cdot{}^sv.
\]
We have already seen that $\bone$ is a left identity for $M_2$.

\pfitem{\ref{M23}}  The first claim is again immediate.  

For the second, we first observe that monoidality clearly implies \eqref{eq:M23}.  For the converse, suppose \eqref{eq:M23} holds.  We must show that ${}^11=1$.  This is clear if ${}^11\not\in U$, as then ${{}^11\in U^1\sm U=\{1\}}$, so we now assume that ${}^11\in U$.  For any $u\in U$, we use \eqref{eq:M23} to calculate
\[
{}^11 \cdot u = {}^11 \cdot {}^1u = {}^1(1\cdot u) = {}^1u = u \ANDSIM u\cdot{}^11 = u.
\]
Thus, ${}^11\in U$ is an identity element of $U$, so that $U^1=U$, and ${}^11=1$ by definition.
\epf

\begin{rem}\label{rem:M2}
Condition \eqref{eq:M21} does not hold in general.  For example, when $S=S^1$ and $U\not=U^1$, and the action of $S$ on $U^1$ is constant with image~$1(\not\in U)$, we have 
\[
M_2=\es \AND \bone\cdot(U\rtimes S)=\{1\}\rtimes S = (U^1\rtimes S) \sm (U\rtimes S).
\]

Clearly \eqref{eq:M23} implies \eqref{eq:M21}; that is, \eqref{eq:M21} holds when $S^1$ acts monoidally on $U^1$ (i.e., ${}^1u=u$ for all $u\in U^1$).  However, the reverse implication \eqref{eq:M21}$\implies$\eqref{eq:M23} need not hold.  For example, if $e$ is a non-identity idempotent of $U$, then the constant action with image~$e$ satisfies \eqref{eq:M21} but not~\eqref{eq:M23}.
\end{rem}

We now give the corresponding result for $M=M_1\cap M_2$.

\begin{lemma}\label{lem:M}
\ben
\item \label{M1} $M\sub\bone\cdot(U\rtimes S)\cdot\bone$, with equality if and only if
\begin{equation}\label{eq:M1}
{}^1us^+\in U \qquad\text{for all $u\in U$ and~$s\in S$.}
\end{equation}
\item \label{M2} $M$ is a subsemigroup of $U\rtimes S$, and $\bone$ is a two-sided identity for $M$ (but not necessarily an element of $M$).
\item \label{M3} $M$ is a monoid with identity $\bone$ if and only if $U$ and $S$ are monoids and $1=1^+$.
\item \label{M4} $M=U\rtimes S$ if and only if
\begin{equation}\label{eq:M4}
u={}^1us^+ \qquad\text{for all $u\in U$ and $s\in S$,}
\end{equation}
and this is equivalent to \eqref{eq:M12} and \eqref{eq:M23} both holding.
\een
\end{lemma}

\pf
Throughout the proof we use Lemmas \ref{lem:+}, \ref{lem:M1} and \ref{lem:M2} without explicit mention.  We begin with the second part.

\pfitem{\ref{M2}}  Since $M_1$ and $M_2$ are subsemigroups of $U\rtimes S$, so too is $M_1\cap M_2=M$.  Since $\bone$ is a right and left identity for $M_1$ and $M_2$, respectively, it is a two-sided identity for $M_1\cap M_2=M$.

\pfitem{\ref{M1}}  By part \ref{M2} we have $M=\bone\cdot M\cdot\bone\sub\bone\cdot(U\rtimes S)\cdot\bone$.

Next, we assume that \eqref{eq:M1} holds, and let $(u,s)\in\bone\cdot(U\rtimes S)\cdot\bone$.  So
\[
(u,s)=\bone\cdot(v,t)\cdot\bone=({}^1vt^+,t) \qquad\text{for some $v\in U$ and $t\in S$.}
\]
So $(u,s)\in U\rtimes S$ by \eqref{eq:M1}.  As in the proof of Lemma \ref{lem:1us+}, one can check that
\[
{}^1us^+={}^1({}^1vt^+)\cdot t^+ = {}^1vt^+ = u,
\]
so that $(u,s)\in M$, by the same lemma.

If, on the other hand, \eqref{eq:M1} does not hold, then ${}^1us^+\not\in U$ for some $u\in U$ and $s\in S$.  In this case, $\bone\cdot(u,s)\cdot\bone=({}^1us^+,s)\not\in U\rtimes S$, so certainly $\bone\cdot(u,s)\cdot\bone\not\in M$.

\pfitem{\ref{M3}}  By part \ref{M2}, $M$ is a monoid with identity $\bone$ if and only if $\bone\in M$, and we noted just before Lemma \ref{lem:M1} that this is equivalent to the stated conditions.

\pfitem{\ref{M4}}  The equivalence of $M=U\rtimes S$ with \eqref{eq:M4} follows from Lemma \ref{lem:1us+}.  The equivalence of~$M=U\rtimes S$ with~\eqref{eq:M12} and \eqref{eq:M23} follows immediately from the original definition of $M$.
\epf

\begin{rem}\label{rem:M}
As in Remark \ref{rem:M2}, condition \eqref{eq:M1} is not always satisfied, but it is (for example) when $S^1$ acts monoidally on $U^1$.

As in Remark \ref{rem:M1}, condition \eqref{eq:M4} is satisfied (for example) when the action of $S^1$ on $U^1$ is monoidal and by monoid morphisms, but this is again not necessary.  Indeed, consider the case in which $U$ has a right identity $e$ that is not a left identity, and suppose $S\not=S^1$.  Consider the constant action of $S$ on $U^1(\not=U)$ with image $e$, and extend this to a monoidal action of $S^1(\not=S)$ on~$U^1$ in the usual way.  Then for any $u\in U$ and $s\in S$, we have ${}^1us^+ = ue = u$.
\end{rem}

By Lemmas \ref{lem:MIS} and \ref{lem:MI}, we have the following commutative diagram, with all sets subsemigroups of $U^1\rtimes S^1$, and all maps surmorphisms:
\[
\begin{tikzcd} 
~ & U\rtimes S \arrow[swap]{dl}{f_1\colon\bx\ \mt\ \bx\cdot\bone} \arrow{dr}{f_2\colon\bx\ \mt\ \bone\cdot\bx} & \\
(U\rtimes S)\cdot\bone \arrow[swap]{dr}{g_1\colon\bx\ \mt\ \bone\cdot\bx} & & \bone\cdot(U\rtimes S) \arrow{dl}{g_2\colon\bx\ \mt\ \bx\cdot\bone}\\
& \bone\cdot(U\rtimes S)\cdot\bone & 
\end{tikzcd}
\]
By Lemma \ref{lem:M1}\ref{M11}, the semigroup on the western side of this diagram is equal to $M_1$.  By Lemma~\ref{lem:M2}\ref{M21} the eastern semigroup is $M_2$ precisely when \eqref{eq:M21} holds.  Since \eqref{eq:M21} clearly implies \eqref{eq:M1}, it follows from Lemma \ref{lem:M}\ref{M1} that the southern semigroup is $M$ in this case.  Thus, we have the following:

\begin{prop}\label{prop:M}
If ${}^1u\in U$ for all $u\in U$, then the following diagram commutes, with all sets subsemigroups of $U\rtimes S$, and all maps retractions:
\[
\begin{tikzcd} [column sep=large]
~ & U\rtimes S \arrow[swap]{dl}{f_1\colon\bx\ \mt\ \bx\cdot\bone} \arrow{dr}{f_2\colon\bx\ \mt\ \bone\cdot\bx} & \\
M_1 \arrow[swap]{dr}{g_1\colon\bx\ \mt\ \bone\cdot\bx} & & M_2 \arrow{dl}{g_2\colon\bx\ \mt\ \bx\cdot\bone}\\
& M & 
\end{tikzcd}
\]
\epfres
\end{prop}

\begin{rem}
It follows that under the ${}^1u\in U$ assumption from Proposition \ref{prop:M}, the subsemigroups $M_1$, $M_2$ and~$M$ of $U\rtimes S$ can also be identified with quotients of $U\rtimes S$:
\[
M_1\cong(U\rtimes S)/\ker(f_1) \COMMA
M_2\cong(U\rtimes S)/\ker(f_2) \AND
M\cong(U\rtimes S)/\ker(F) ,
\]
where $F=f_1g_1=f_2g_2$ is the retraction
\[
F:(U\rtimes S)\to M:\bx\mt\bone\cdot\bx\cdot\bone \qquad\text{given by}\qquad (u,s)F=({}^1us^+,s) \qquad\text{for $u\in U$ and $s\in S$.}
\]
The kernels of these maps can be described easily.  For example:
\bit
\item $(u,s)f_1 = (v,t)f_1 \iff us^+=vt^+ \text{ and } s=t$, 
\item $(u,s)f_2 = (v,t)f_2 \iff {}^1u={}^1v \text{ and } s=t$, and
\item $(u,s)F = (v,t)F \iff {}^1us^+={}^1vt^+ \text{ and } s=t$.
\eit
\end{rem}

\begin{rem}\label{rem:f1g2}
The ${}^1u\in U$ assumption from Proposition \ref{prop:M} is precisely condition \eqref{eq:M21}.  When \eqref{eq:M12} additionally holds,
\[
f_1:U\rtimes S\to M_1=U\rtimes S \AND g_2:M=M_2\to M
\]
are identity maps; cf.~Lemma~\ref{lem:M1}\ref{M12}.  The analogous statement holds for $f_2$ and $g_1$ when \eqref{eq:M23} holds.
\end{rem}

As we have already mentioned, the actions considered in the rest of the paper will almost always be monoidal.  It is therefore convenient to record the next result, which summarises the main points above in the case that the action of $S^1$ on $U^1$ is monoidal.  
For the statement, recall (see for example \cite{RSbook}) that the \emph{local monoid} at an idempotent $e$ of a semigroup $T$ is the set
\[
eTe = eT\cap Te = \set{ete}{t\in T} = \set{t\in T}{et=t=te},
\]
and is the largest subsemigroup of $T$ that happens to be a monoid with identity $e$.

\begin{prop}\label{prop:MM}
Suppose $S$ acts on $U^1$ by semigroup morphisms, and suppose this action is monoidal if $S$ is a monoid.  Then
\[
M=\bigset{(u,s)\in U\rtimes S}{u=us^+}
\]
is a subsemigroup of $U\rtimes S$, and the map
\[
F:U\rtimes S\to M:(u,s)\mt(us^+,s)
\]
is a retraction.  If $U$ and $S$ are monoids, then $M$ is a (local) monoid with identity $\bone=(1,1)$.
\end{prop}

\pf
The monoidality assumption in the statement is equivalent to the action of $S^1$ on $U^1$ being monoidal.  It follows from Lemma \ref{lem:M2}\ref{M23} that 
\[
M=M_1=\bigset{(u,s)\in U\rtimes S}{u=us^+}\leq U\rtimes S.
\]
Monoidality also means that Proposition \ref{prop:M} applies, and together with $M=M_1$, this gives the assertion regarding the retraction $F=f_1$; cf.~Remark \ref{rem:f1g2}.  The final assertion follows from parts \ref{M1} and \ref{M3} of Lemma~\ref{lem:M}, and the fact that monoidality implies \eqref{eq:M1} and $1=1^+$.
\epf

\chap{Action pairs}\label{chap:AP}

In this chapter we introduce the kinds of semigroups that are the focus of the paper.  These are certain products of the form $US$, where $U$ and $S$ are subsemigroups of a common over-monoid~$M$ satisfying natural axioms.  After giving the relevant definitions and basic properties in Section~\ref{sect:basic}, we discuss several examples in Section~\ref{sect:eg}.  Section \ref{sect:LRAP} treats the class of left restriction monoids, showing how these fit into our general framework, and providing the tools for studying many of our motivating examples in Part \ref{part:II} of the paper.  Section \ref{sect:structure} contains the main structural result of the chapter, Theorem~\ref{thm:FT}, which characterises our product semigroups as quotients of semidirect products by certain special congruences that we abstractly axiomatise.  Finally, Section \ref{sect:cong} contains several results on generating sets for these congruences, which then feed into presentations in Part~\ref{part:II}.

\sect{Basic definitions and properties}\label{sect:basic}

The general set-up we will be concerned with involves two subsemigroups of a common semigroup; since the latter embeds in its monoid completion, we may assume without loss of generality that it already is a monoid.  Thus, we have $U,S\leq M$, where $M$ is a monoid, subject to various natural assumptions to be detailed below.  We denote the identity of $M$ by $1$, and write $U^1=U\cup\{1\}$ and $S^1=S\cup\{1\}$.  Note that $U^1=U$ if $U$ happens to be a submonoid of $M$, and similarly for~$S^1$.  In any case, we also have $U^1,S^1\leq M$.

\begin{lemma}\label{lem:US}
If $U,S\leq M$, for some monoid $M$, and if $SU\sub U^1S^1$, then $US\leq M$.
\end{lemma}

\pf
This follows from $US\cdot US \sub UU^1\cdot S^1S=US$.
\epf

Note that the semigroups $U$ and $S$ need not be contained in the product $US$.  Indeed, this will be the case for several of our motivating examples.  Of course, if $U$ and $S$ are both submonoids of $M$, then $U,S\sub US$.

Here is the first key definition.

\begin{defn}\label{defn:AP}
A \emph{weak (left) action pair} in a monoid $M$ is a pair $(U,S)$ of subsemigroups of~$M$ for which the following condition holds:
\begin{enumerate}[label=\textup{(A\arabic*)}]
\item \label{A1} $S$ has a left action on $U^1$ by semigroup morphisms, written $(s,u)\mt {}^su$, such that
\[
s\cdot u={}^su\cdot s \qquad\text{for all $s\in S$ and $u\in U^1$.}
\]
\een
If the following condition also holds, we call $(U,S)$ a \emph{(left) action pair}:
\begin{enumerate}[label=\textup{(A\arabic*)}]\addtocounter{enumi}{1}
\item \label{A2} $us=vt \implies u\cdot{}^s1 = v\cdot{}^t1$ for all $u,v\in U^1$ and $s,t\in S$.
\een
\end{defn}

\begin{rem}
Strictly speaking, it would be more precise to speak of a \emph{(weak) action triple} $(U,S,\varphi)$, where $\varphi:S\times U^1\to U^1:(s,u)\mt {}^su$ is the specific action from \ref{A1}.  However, as with semidirect products in Definition \ref{defn:SD}, we omit the label $\varphi$ for brevity, and the action will always be clear from context.

For the time being, we have resisted using the notation $s^+$ for the elements ${}^s1$ ($s\in S$), for reasons that will soon become clear; see Proposition \ref{prop:+}.
\end{rem}

It of course follows from Lemma \ref{lem:US} and \ref{A1} that $US\leq M$ for any (weak) action pair $(U,S)$ in $M$.
Most of the results of this paper concern action pairs, but some hold more generally for weak action pairs; see for example Proposition \ref{prop:A1}, and Theorems \ref{thm:cover} and \ref{thm:ESmon0}.

The main role played by the monoid $M$ in Definition \ref{defn:AP} is to provide a common environment in which the semigroups $U$ and $S$ exist, or, more specifically, a context in which to form products of the form $us$ and $su$, for $u\in U$ and $s\in S$.  Consequently, $M$ could be replaced in the definition by any monoid containing $U$ and $S$, and we will often speak simply of `an action pair $(U,S)$', without explicitly specifying the common over-monoid $M$.

\begin{rem}\label{rem:AP2}
Note that we do not assume in \ref{A1} that $S$ acts by \emph{monoid} morphisms (i.e., that ${}^s1=1$ for all $s\in S$), but we will soon see that this special case is very important.  
On the other hand, if $S$ happens to be a submonoid of $M$, then it follows from \ref{A1} that $u = 1\cdot u = {}^1u\cdot 1 = {}^1u$ for all $u\in U^1$, meaning that $S$ acts \emph{monoidally} in this case.

In any case, is important to note that $s = s\cdot 1 = {}^s1\cdot s$ for all $s\in S$ in any weak action pair.  When we return to the $s^+={}^s1$ notation later, this says that
\[
s=s^+s \qquad\text{for all $s\in S$.}
\]
\end{rem}

Many of the action pairs we consider will satisfy stronger assumptions:

\begin{defn}\label{defn:SAP}
A \emph{strong (left) action pair} in a monoid $M$ is a pair $(U,S)$ of subsemigroups of~$M$ for which the following conditions both hold:
\begin{enumerate}[label=\textup{(SA\arabic*)},leftmargin=12mm]
\item \label{SA1} $sU^1\sub U^1s$ for all $s\in S$.
\item \label{SA2} $us=vt \implies u=v$ for all $u,v\in U^1$ and $s,t\in S$.
\een
We call \ref{SA2} the \emph{left-uniqueness} property.
\end{defn}

\begin{rem}\label{rem:AP4}
It is worth observing that $sU^1\sub U^1s$ in \ref{SA1} is equivalent to $sU\sub U^1s$.  The condition $sU\sub Us$ of course implies $sU^1\sub U^1s$, but not conversely in general; see Example~\ref{eg:Tn}.
\end{rem}

The terminology obviously suggests that any strong action pair is an action pair, even though no action was specified in Definition \ref{defn:SAP}.  This is indeed the case, as we will show in Lemma~\ref{lem:SAPAP}.  First, however, we give an alternative characterisation of action pairs, which shows that certain key properties of the elements ${}^s1$ ($s\in S$) determine the whole of the action.

\begin{prop}\label{prop:+}
Let $U$ and $S$ be subsemigroups of a monoid $M$.  Then $(U,S)$ is an action pair if and only if \ref{SA1} holds and there exists a map $S\to U^1:s\mt s^+$ satisfying the following conditions:  
\begin{align}
\label{+1}\tag*{($+$1)} s&=s^+s &&\text{for all $s\in S$,}\\
\label{+2}\tag*{($+$2)} st^+&=(st)^+s &&\text{for all $s,t\in S$,}\\
\label{+3}\tag*{($+$3)} (st)^+&=(st)^+s^+ &&\text{for all $s,t\in S$,}\\
\label{+4}\tag*{($+$4)} us=vt\ \ &\Rightarrow \ \ us^+=vt^+ &&\text{for all $u,v\in U^1$ and $s,t\in S$.}
\end{align}
\end{prop}

\pf
\firstpfitem{($\Rightarrow$)}  Suppose first that $(U,S)$ is an action pair.  Since \ref{A1}$\implies$\ref{SA1}, it remains to demonstrate the existence of a suitable map $s\mt s^+$.  For this, we define $s^+={}^s1$ for all $s\in S$.  Properties \ref{+3} and \ref{+4} then follow immediately from Lemma \ref{lem:+}\ref{lem+3} and \ref{A2}, respectively, and we noted in Remark \ref{rem:AP2} that \ref{+1} holds.  For \ref{+2}, let $s,t\in S$.  Then \ref{A1} and Lemma~\ref{lem:+}\ref{lem+2} give
\[
st^+={}^s(t^+)\cdot s = (st)^+s.
\]

\pfitem{($\Leftarrow$)}  Conversely, suppose \ref{SA1} holds, and also that $S\to U^1:s\mt s^+$ satisfies \ref{+1}--\ref{+4}.  

We begin by defining an action of $S$ on $U^1$.  To this end, let $s\in S$ and $u\in U^1$.  By~\ref{SA1} we have $su=vs$ for some $v\in U^1$, and we define ${}^su=vs^+$.  It follows quickly from \ref{+4} that this is well defined.  Thus, for any $s\in S$ and $u\in U^1$, we have
\begin{equation}\label{eq:su}
{}^su = vs^+ \qquad\text{for any $v\in U^1$ such that $su=vs$.}
\end{equation}
We will soon show that this does indeed define an action by semigroup morphisms.  But first we observe that for $s,u,v$ as above, we have 
\begin{equation}\label{eq:susus}
su=vs=vs^+s={}^su\cdot s,
\end{equation}
where we used \ref{+1} in the second step.  Next we note that
\begin{equation}\label{eq:s+}
{}^s1 = s^+ \qquad\text{for all $s\in S$.}
\end{equation}
Indeed, this follows immediately from \eqref{eq:su} with $u=v=1$.  Combining \eqref{eq:s+} with \ref{+4}, we obtain~\ref{A2}.  

In light of \eqref{eq:susus}, it remains to check that \eqref{eq:su} does indeed define an action of $S$ on $U^1$ by semigroup morphisms: i.e., that 
\begin{equation}\label{eq:SU}
{}^{st}u = {}^s({}^tu) \AND {}^s(uv) = {}^su\cdot{}^sv \qquad\text{for all $u,v\in U^1$ and $s,t\in S$.}
\end{equation}
To do this, we begin by showing that 
\begin{equation}\label{eq:sus+}
{}^su\cdot s^+={}^su \qquad\text{for all $u\in U^1$ and $s\in S$.}
\end{equation}
To see this, first note that for any $s\in S$, \ref{+1} and \ref{+4} give
\[
1\cdot s = s = s^+\cdot s \ \ \implies \ \ 1\cdot s^+=s^+\cdot s^+,
\]
so that each $s^+$ is an idempotent.  Combining this with \eqref{eq:su}, we quickly obtain \eqref{eq:sus+}.

Beginning with the second part of \eqref{eq:SU}, let $s\in S$ and $u,v\in U^1$.  Several applications of~\eqref{eq:susus} gives 
\[
{}^s(uv)\cdot s = s\cdot uv = {}^su\cdot s\cdot v = {}^su\cdot {}^sv\cdot s,
\]
and it then follows from \ref{+4} that ${}^s(uv)\cdot s^+ = {}^su\cdot{}^sv\cdot s^+$.  We then apply \eqref{eq:sus+} to both sides to deduce ${}^s(uv) = {}^su\cdot{}^sv$.

For the first part of \eqref{eq:SU}, let $s,t\in S$ and $u\in U^1$; we must show that ${}^{st}u={}^s({}^tu)$.  This time~\eqref{eq:susus} gives
\[
{}^{st}u\cdot st = st\cdot u = s\cdot{}^tu\cdot t = {}^s({}^tu)\cdot st.
\]
It then follows from \ref{+4} that ${{}^{st}u\cdot(st)^+ = {}^s({}^tu)\cdot(st)^+}$.  Since ${}^{st}u\cdot(st)^+={}^{st}u$, by \eqref{eq:sus+}, the proof will be complete if we can show that
\begin{equation}\label{eq:stust+}
{}^s({}^tu)\cdot(st)^+={}^s({}^tu).
\end{equation}
For this we use \eqref{eq:susus}, \eqref{eq:sus+}, \eqref{eq:susus} and \ref{+2} to calculate
\[
{}^s({}^tu)\cdot s = s\cdot {}^tu = s\cdot{}^tu\cdot t^+ = {}^s({}^tu)\cdot s\cdot t^+ = {}^s({}^tu)\cdot(st)^+\cdot s.
\]
It then follows from \ref{+4} that ${}^s({}^tu)\cdot s^+ = {}^s({}^tu)\cdot(st)^+\cdot s^+$.  Combining this with \eqref{eq:sus+} and~\ref{+3}, we finally deduce that
\[
{}^s({}^tu) = {}^s({}^tu)\cdot s^+ = {}^s({}^tu)\cdot(st)^+\cdot s^+ = {}^s({}^tu)\cdot(st)^+.
\]
This completes the proof of \eqref{eq:stust+}, and hence of the proposition.
\epf

\begin{rem}
Proposition \ref{prop:+} allows us to specify an action pair $(U,S)$ either by means of:
\bit
\item an action $(s,u)\mt{}^su$, as in Definition \ref{defn:AP}, or 
\item an appropriate map $S\to U^1:s\mt s^+$.
\eit
The proof of the proposition shows how the two formulations are interchangeable, but it is worth drawing this out explicitly here:
\bit
\item Given an action pair $(U,S)$, as in Definition \ref{defn:AP}, the map $S\to U^1 : s\mt s^+={}^s1$ satisfies \ref{+1}--\ref{+4}.  
\item Conversely, if $(U,S)$ satisfies \ref{SA1}, and if $S\to U^1:s\mt s^+$ satisfies \ref{+1}--\ref{+4}, then $(U,S)$ is an action pair (as in Definition \ref{defn:AP}) with respect to the action of $S$ on $U^1$ given in \eqref{eq:su}.
\eit
In all that follows, we will use either viewpoint, as convenient, typically without further comment.
\end{rem}

Here is the promised result justifying the terminology of strong action pairs; it also gives a criterion for distinguishing the strong ones.

\begin{lemma}\label{lem:SAPAP}
\ben
\item \label{SAPAP1} Any strong action pair is an action pair.  
\item \label{SAPAP2} An action pair $(U,S)$ is strong if and only if $s^+=1$ for all $s\in S$, meaning that the action of $S$ on $U^1$ is by monoid morphisms.
\een
\end{lemma}

\pf
\firstpfitem{\ref{SAPAP1}}  The trivial map ${S\to U^1:s\mt s^+=1}$ always satisfies conditions \ref{+1}--\ref{+3}.  If $(U,S)$ is a strong action pair, then this trivial map also satisfies condition \ref{+4}.

\pfitem{\ref{SAPAP2}}  Suppose first that the action pair $(U,S)$ is strong, and let $s\in S$.  Then from $s^+\cdot s = s = 1\cdot s$, it follows from \ref{SA2} that $s^+=1$.

Conversely, suppose $s^+=1$ for all $s\in S$.  Together with \ref{A2}, this clearly implies \ref{SA2}.  Since \ref{SA1} always follows from \ref{A1}, it follows that $(U,S)$ is strong.
\epf

It follows that strong action pairs have the following equivalent definition:

\begin{defn}\label{defn:SAP'}
A \emph{strong (left) action pair} in a monoid $M$ is a pair $(U,S)$ of subsemigroups of~$M$ for which the following conditions both hold:
\begin{enumerate}[label=\textup{(SA\arabic*)$'$},leftmargin=13mm]
\item \label{SA1'} $S$ has a left action on $U^1$ by monoid morphisms, written $(s,u)\mt {}^su$, such that
\[
s\cdot u={}^su\cdot s \qquad\text{for all $s\in S$ and $u\in U^1$.}
\]
\een
\begin{enumerate}[label=\textup{(SA\arabic*)$\textcolor{white}{'}$},leftmargin=13mm]\addtocounter{enumi}{1}
\item $us=vt \implies u = v$ for all $u,v\in U^1$ and $s,t\in S$.
\een
\end{defn}

\begin{rem}
We have seen that a pair $(U,S)$ of subsemigroups of a monoid $M$ is a strong action pair if and only if 
\bit
\item \ref{SA1} and \ref{SA2} both hold (cf.~Definition \ref{defn:SAP}), or equivalently 
\item \ref{SA1'} and \ref{SA2} both hold (cf.~Definition \ref{defn:SAP'}).
\eit
Since \ref{SA1'}$\implies$\ref{A1}$\implies$\ref{SA1}, $(U,S)$ is also strong if and only if
\bit
\item \ref{A1} and \ref{SA2} both hold.
\eit
In Section \ref{sect:2M1}, we will prove one result (Theorem \ref{thm:ESmon0}) concerning pairs $(U,S)$ satisfying \ref{SA1'}, but not necessarily \ref{SA2}.  This limited attention does not seem to warrant naming such pairs.
\end{rem}

The next lemma involves restricting action pairs to subsemigroups; it follows quickly from an examination of Definitions \ref{defn:AP} and \ref{defn:SAP'}.

\begin{lemma}\label{lem:US0}
If $(U,S)$ is a (strong) action pair, then so too is:
\ben
\item \label{US01} $(U,T)$ for any $T\leq S$,
\item \label{US02} $(V,S)$ for any $V\leq U$ for which $V^1$ is closed under the action of $S$.  \epfres
\een
\end{lemma}

For the next statement, we say that an action pair $(U,S)$ in $M$ \emph{extends to} $(V,T)$ if:
\bit
\item $(V,T)$ is an action pair in $M$,
\item $U\leq V$ and $S\leq T$, and 
\item the map $S\to U^1:s\mt s^+$ is the restriction of the map $T\to V^1:t\mt t^+$.
\eit
The third condition can be formulated equivalently in terms of actions of $S$ and $T$ on $U^1$ and~$V^1$.

\begin{lemma}\label{lem:US1}
\ben
\item \label{US11} Any action pair $(U,S)$ extends to $(U^1,S)$.  
\item \label{US12} An action pair $(U,S)$ extends to $(U,S^1)$ and $(U^1,S^1)$ if and only if the following condition holds:
\bit
\item $us=v \implies us^+=v$ for all $u,v\in U^1$ and $s\in S$.
\eit
\item \label{US13} Any strong action pair $(U,S)$ extends to $(U^1,S)$, $(U,S^1)$ and $(U^1,S^1)$, and these are all strong.  
\een
\end{lemma}

\pf
\firstpfitem{\ref{US11}}  This follows immediately from Definition \ref{defn:AP}, given that $(U^1)^1=U^1$.

\pfitem{\ref{US12}}  Given part \ref{US11}, we just need to prove the claim for the pair $(U,S^1)$.  We first observe that~\ref{SA1} clearly holds for this pair.  There is only one way to extend the map $S\to U^1:s\mt s^+$ to a map $S^1\to U^1$ satisfying \ref{+1}, since the latter implies $1^+=1$.  Items \ref{+1}--\ref{+3} all clearly hold for this extended map, while \ref{+4} holds if and only if the stated condition holds.

\pfitem{\ref{US13}}  Item \ref{SA1} is clear in each case.  For \ref{SA2}, suppose $us=vt$ for some $u,v\in U^1$ and $s,t\in S^1$.  Then for any $x\in S$ we have $u(sx)=v(tx)$, with $sx,tx\in S$, so it follows from left-uniqueness of $(U,S)$ that $u=v$.
\epf

The condition stated in part \ref{US12} of the previous lemma leads to the following useful fact.

\begin{lemma}\label{lem:submon}
If $(U,S)$ and $(U,S^1)$ are both action pairs in $M$, then $US$ is a submonoid of $M$ if and only if $U$ and $S$ are both submonoids.
\end{lemma}

\pf
With the backwards implication being clear, suppose $US$ is a submonoid, so that $us=1$ for some $u\in U$ and $s\in S$.  Since $us=1\cdot1$ with $u,1\in U^1$ and $s,1\in S^1$, we can apply \ref{A2} in the action pair $(U,S^1)$ to deduce that $us^+=1\cdot 1^+=1\cdot1=1$.  But then
\[
s = 1\cdot s = us^+\cdot s = us = 1,
\]
where we used \ref{+1} in the third step.  It follows that also $u=u\cdot1=us=1$.  Thus, ${1=u=s\in U\cap S}$.
\epf

The next result provides a necessary condition for an action pair to be strong.

\begin{lemma}\label{lem:dis}
If $(U,S)$ is a strong action pair, then 
\ben
\item \label{dis1} $U\cap S\sub\{1\}$,  
\item \label{dis2} $1\in US \iff 1\in U\cap S \iff U\cap S=\{1\} \iff U\cap S\not=\es$.
\een
\end{lemma}

\pf
\firstpfitem{\ref{dis1}} If $x\in U\cap S$, then from $x\cdot x = 1\cdot x^2$, left-uniqueness implies $x=1$.

\pfitem{\ref{dis2}}  The following implications are all obvious:
\[
U\cap S=\{1\} \implies 1\in U\cap S \implies U\cap S\not=\es \AND 1\in U\cap S\implies 1\in US.
\]
Part \ref{dis1} gives $U\cap S\not=\es \implies U\cap S=\{1\}$, so it remains to show that $1\in US\implies 1\in U\cap S$.  Since $(U,S)$ is strong, $(U,S^1)$ is a (strong) action pair, by Lemma \ref{lem:US1}\ref{US13}.  The required implication then follows from Lemma \ref{lem:submon}.
\epf

\begin{rem}
The assumption that $(U,S)$ is strong cannot be removed from Lemma \ref{lem:dis}.  For example, let $M$ be a left restriction monoid with semilattice of projections $P$; see Section~\ref{sect:lrs} for the definitions.  We will show in Proposition \ref{prop:LR1} that $(P,M)$ is an action pair, and we note that $P\cap M=P$ need not contain only the identity element.
\end{rem}

Now we provide a \emph{sufficient} condition.  Recall that a \emph{right unit} of a monoid $M$ is an element $s\in M$ such that $1=st$ for some $t\in M$.  We write $R(M)=\set{s\in M}{1\in sM}$ for the submonoid of all right units of $M$.  (The submonoid $R(M)$ is also Green's $\R$-class of $1\in M$.)

\begin{lemma}\label{lem:group}
Suppose $(U,S)$ is an action pair in a monoid $M$.
\ben
\item \label{group1} If $s\in S\cap R(M)$, then $s^+=1$.  
\item \label{group2} If $S\sub R(M)$, then $(U,S)$ is strong.
\een
\end{lemma}

\pf
By Lemma \ref{lem:SAPAP}\ref{SAPAP2}, it suffices to prove the first part.  For this, suppose ${s\in S\cap R(M)}$, so that $1=st$ for some $t\in M$.  Combining this with \ref{+1}, we have ${1 = st = s^+st = s^+}$.
\epf

\begin{rem}
It of course follows from Lemma \ref{lem:group}\ref{group2} that $(U,S)$ is strong if $S$ is contained in the group of (two-sided) units of $M$.

The condition $S\sub R(M)$ is certainly not necessary for $(U,S)$ to be strong; see Example \ref{eg:Tn} (out of \emph{many} other examples considered in the paper).
\end{rem}

The final result of this section shows that \ref{A2} is equivalent to an ostensibly stronger condition.

\begin{lemma}\label{lem:ap}
If $(U,S)$ is an action pair, then for any $u,v,w\in U^1$ and $s,t\in S$, we have
\[
us = vt \implies u\cdot{}^sw = v\cdot{}^tw.
\]
\end{lemma}

\pf
Combining $us=vt$ with \ref{A1}, we have
\[
u\cdot{}^sw \cdot s = us\cdot w = vt\cdot w = v\cdot {}^tw \cdot t.
\]
It then follows from \ref{A2} that $u\cdot{}^sw \cdot s^+ = v\cdot {}^tw \cdot t^+$, and Lemma \ref{lem:+}\ref{lem+4} completes the proof.
\epf

Before we move on, we note that there are obviously right-handed versions of our action pairs.  By symmetry/duality, results concerning \emph{right} pairs can be directly obtained from their left-handed counterparts, so we focus only on the left.  The left-handed theory is our chosen focus, as it applies directly to all of our motivating examples.  

There is also a notion of \emph{right}-uniqueness (for \emph{left} action pairs), where ${us=vt\implies s=t}$, in contrast to the left-uniqueness property \ref{SA2}.  Right-uniqueness is satisfied, for example, by free left restriction monoids, which are the topic of Chapter \ref{chap:LR}.  However, right-uniqueness is equivalent to certain right congruences on $S$ (introduced in Section \ref{sect:structure}) all being trivial, so there is no need to develop a parallel theory in this case.

\sect{Examples}\label{sect:eg}

The following collection of examples and non-examples highlight some subtleties in the definitions.  They also show that the semigroups $US$ arising from (weak/strong) action pairs $(U,S)$ include some well known classes.  Additional examples will be explored in greater depth in the next section, and in Part \ref{part:II} of the paper.  

Several of the examples discussed below use the transformation monoids $\PT_X$,~$\T_X$, $\Sing(\T_n)$, and so on, as defined in Section \ref{sect:ts}.  
Some involve inverse semigroups, as defined in Section~\ref{sect:sgp}.  For the calculations in such examples, we make use of the previously-mentioned fact that idempotents of an inverse semigroup commute.  It is also clear that every idempotent of an inverse semigroup is its own inverse.  On a number of occasions we make use of the easily checked fact that if $U$ and $S$ are submonoids of a monoid $M$, and if $SU\sub US$, then $US=\la U\cup S\ra$.

We begin with a somewhat degenerate construction:

\begin{eg}
Let $U$ and $S$ be arbitrary semigroups, and assume $U\cap S=\es$.  Let $0$ and $1$ be symbols not belonging to $U\cup S$, and let $M=U\cup S\cup\{0,1\}$.  Then $M$ is a monoid under the product that extends the original products in $U$ and $S$, and where additionally
\[
us=su=0 \COMMa x0=0x=0 \ANd x1=1x=x \qquad\text{for all $u\in U$, $s\in S$ and $x\in M$.}
\]
(So $M$ is the monoid completion of the \emph{0-direct union} of $U$ and $S$.)  For any action of $S$ on $U$ by morphisms, $(U,S)$ is a weak action pair in $M$ with respect to the extended action of $S$ on $U^1$.  We then of course have $US=SU=\{0\}$.
\end{eg}

The next series of examples involve idempotents and one- or two-sided units.

\begin{eg}\label{eg:FIM}
Let $M$ be an inverse monoid, and let
\[
E = E(M) = \set{e\in M}{e=e^2} \AND G = G(M) = \set{g\in M}{gg^{-1}=g^{-1}g=1}
\]
be the semilattice of idempotents, and group of units of $M$.  It is easy to see that ${(E,G)=(E^1,G^1)}$ is a strong action pair in $M$:
\begin{enumerate}[label=\textup{(SA\arabic*)},leftmargin=12mm]
\item For $g\in G$ and $e\in E(=E^1)$, we have $ge=(geg^{-1})g$, with $geg^{-1}\in E$.  This shows that $gE\sub Eg$.  In fact, it is easy to see that $gE=Eg$.  (The action of $G$ on $E(=E^1)$ from \ref{A1} or \ref{SA1'} is given by ${}^ge=geg^{-1}$.)
\item Suppose $eg=fh$ for some $e,f\in E$ and $g,h\in G$.  Then from $e=fhg^{-1}$ we obtain $e=fe$.  Similarly, $f=ef$, and since idempotents commute, it follows that $e=f$.
\een
(Alternatively, one could verify conditions \ref{A1} and \ref{A2}, and then apply Lemma \ref{lem:group}\ref{group2}.)
The resulting subsemigroup $EG=GE\leq M$ is the largest \emph{factorisable inverse submonoid} of $M$. See \cite{CH1974,FitzGerald2010} for more on factorisable inverse monoids, and \cite{EEF2005} for presentations.  By the symmetry afforded by the inversion map, $(E,G)$ is also a strong \emph{right} action pair.

If we write $\Ef=E\sm\{1\}$, then $(\Ef,G)$ is also a strong action pair in $M$ (cf.~Lemma \ref{lem:US0}\ref{US02}).  The subsemigroup $\Ef G=G\Ef=EG\sm G$ is an \emph{almost-factorisable inverse semigroup} in the sense of Lawson \cite{Lawson1994}; see also \cite{McAlister1976}.  Note that this class of inverse semigroups has an abstract definition, and it is a highly nontrivial result that any almost-factorisable inverse semigroup has the form $EG\sm G$ for suitable $E$ and $G$ (and $M$).
\end{eg}

Generalisations of Example \ref{eg:FIM} will be discussed in Section \ref{sect:LRAP}.

\begin{eg}\label{eg:FM}
Following on from the previous example, let $M$ be an arbitrary monoid, $E$ the set of idempotents of $M$, and $G$ the group of units.  Since $E$ need not be a submonoid, we also write $U=\la E\ra$ for the submonoid generated by the idempotents.  Since the conjugate of an idempotent by a unit is an idempotent, we have an action of $G$ on $U$ given by ${}^gu=gug^{-1}$, with respect to which $g\cdot u={}^gu\cdot g$.  Thus, \ref{A1} holds, meaning that $(U,G)$ is a weak action pair.  Of course it follows that \ref{SA1} holds as well.  We have $UG=GU=\la E\cup G\ra$ by \cite[Lemma 32]{EF2012}.

Since ${}^g1=1$ for all $g\in G$, conditions \ref{A2} and \ref{SA2} are equivalent for the pair $(U,G)$, and we note that these need not hold in general (although we have just seen in Example \ref{eg:FIM} that they do hold when~$M$ is inverse).  For a specific example of this, let $M=\T_2=\{\al,\be,\ga,\de\}$, where
\[
\al = \trans{1&2\\1&1} \COMMA
\be = \trans{1&2\\2&2} \COMMA
\ga = \trans{1&2\\1&2} \AND
\de = \trans{1&2\\2&1} .
\]
Then $U=E=\{\al,\be,\ga\}$ and $G=\{\ga,\de\}$, yet $\al \ga=\be \de$.
\end{eg}

\begin{eg}\label{eg:EL}
Again let $M$ be an arbitrary monoid, and let $G$, $E$ and $U=\la E\ra$ be as in the previous example.  This time let
\[
L = L(M) = \set{s\in M}{1\in Ms} \AND R = R(M) = \set{s\in M}{1\in sM}
\]
be the submonoids of left and right units.

Let $s\in L$, so that $1=s's$ for some $s'\in M$.  For any idempotent $e\in E$, we have $se=(ses')s$, and it is easy to check that $ses'$ is an idempotent.  This quickly leads to $sU\sub Us$, so that \ref{SA1} holds for the pair $(U,L)$.  We have $UL=\la E\cup L\ra$ by \cite[Lemma 2.5]{JE2019a}.  

Condition \ref{SA2} still does not hold for $(U,L)$ in general.  Indeed, this is again witnessed by $M=\T_2$ from Example \ref{eg:FM}, since when $M$ is finite every left unit is a (two-sided) unit so that $L=G$ is the group of units.  Even though $(U,L)$ is not left-unique in general, we still always have $U\cap L=\{1\}$; see \cite[Lemma 2.1]{JE2019a}, and cf.~Lemma \ref{lem:dis}.

Regarding condition \ref{A1} for the pair $(U,L)$, for each $s\in L$ we fix some $s'\in M$ such that $1=s's$.  Then for $s\in L$ and $u\in U(=U^1)$ we define ${}^su=sus'\in U$, and we have $su={}^su\cdot s$.  We also have ${}^s(uv)={}^su\cdot {}^sv$ for all $s\in L$ and $u,v\in U$.  However, for \ref{A1} to hold, we also require that ${}^{st}u={}^s({}^tu)$ for all $s,t\in L$ and $u\in U$: i.e., $(st)u(st)'=s(tut')s'$ for all such $s,t,u$.  This would be the case if we could choose the elements $s'$ in such a way that $(st)'=t's'$ for all $s,t\in L$, as is the case for example when $M$ is inverse (and $s'=s^{-1}$).  

The above discussion concentrated on the pair $(U,L)$.  We note that the pair $(U,R)$ might not even satisfy \ref{SA1}.  For a specific example of this, let $M=\PB_\N$ be the partial Brauer monoid over $\N=\{0,1,2,\ldots\}$.  The definition of this monoid is somewhat involved, so we refer the reader to \cite{JE2019a}.  
Consider the elements
\[
\al = \custpartn{1,...,5}{1,...,7}{\stline13\stline24\stline35\stline46\stline57\darc12\udotted6{10}\ddotted8{10}}
\COMMA
\be = \custpartn{1,...,7}{1,...,5}{\stline31\stline42\stline53\stline64\stline75\uarc12\ddotted6{10}\udotted8{10}}
\AND
\eta = \custpartn{1,...,7}{1,...,7}{\stline33\stline44\stline55\stline66\stline77\udotted8{10}\ddotted8{10}}
\]
from $\PB_\N$.
Then $\eta\in E\sub U$, and since $\al\be=1$ we have $\al\in R$.  We then have
\[
\al\eta = \custpartn{1,...,5}{1,...,7}{\stline13\stline24\stline35\stline46\stline57\udotted6{10}\ddotted8{10}}\in\al U,
\]
but we claim that $\al\eta\not\in U\al$, which then of course implies that $\al U\not\sub U\al$.  In fact, $\al\eta$ does not even belong to $\PB_\N\cdot\al$, as every element of~$\PB_\N\cdot\al$ contains at least one lower `hook' (inherited from $\al$), while $\al\eta$ does not.  

For more on monoids generated by idempotents and one-sided units, see \cite{JE2019a,JE2020a}.
\end{eg}

Example \ref{eg:FIM} involved a (strong) action pair $(U,S)$ for which $U$ and $S$ were both submonoids of $M$; in such cases, the monoid $US$ contains both $U$ and $S$.  The next example, which is one of our original sources of motivation, provides a contrast.

\begin{eg}\label{eg:Tn}
Let $n$ be a positive integer, and write $\bn=\{1,\ldots,n\}$.  In addition to the monoids $\PT_n$, $\T_n$, and so on, let $\E_n=\set{\id_A}{A\sub\bn}$ be the semilattice of partial identities, and write $\Sing(\E_n)=\E_n\sm\{\id_\bn\}$.  Then the following are all strong action pairs in $\PT_n$:
\[
(\E_n,\T_n) \COMMA (\E_n,\Sing(\T_n)) \COMMA (\Sing(\E_n),\T_n) \COMMA (\Sing(\E_n),\Sing(\T_n)).
\]
Indeed, conditions \ref{SA1} and \ref{SA2} follow from the laws
\begin{equation}\label{eq:Af}
\al\cdot\id_A=\id_{A\al^{-1}}\cdot \al \AND \id_A\cdot \al=\al\restr_A .
\end{equation}
One can check (though it is not entirely obvious, and in any case follows from results of later chapters; cf.~Remark \ref{rem:PTX}) that the subsemigroups $US$ of $\PT_n$ corresponding to the above pairs~$(U,S)$ are:
\bit
\item $\E_n\cdot\T_n = \PT_n$,
\item $\E_n\cdot\Sing(\T_n) = \Sing(\PT_n)$, and 
\item $\Sing(\E_n)\cdot\T_n = \Sing(\E_n)\cdot\Sing(\T_n) = \PTnTn$.
\eit
It is also worth noting that $\E_n\cdot\Sing(\T_n) = \Sing(\PT_n)$ contains $\Sing(\T_n)$ but not $\E_n$; it does however contain $\Sing(\E_n)=\E_n\sm\{\id_\bn\}$, a crucial property that features heavily in Chapter \ref{chap:pres}.  Presentations for the semigroups $\Sing(\T_n)$, $\Sing(\PT_n)$ and $\PTnTn$ can be found in \cite{JEtnsn,JEtnsn2,JEptnsn,JEptnsn2}.
Also, while $\Sing(\E_n)\cdot\T_n = \Sing(\E_n)\cdot\Sing(\T_n) = \PTnTn$ contains $\Sing(\E_n)$, it is disjoint from both $\T_n$ and $\Sing(\T_n)$.

We also observe that if $(U,S)=(\Sing(\E_n),\T_n)$ or $(\Sing(\E_n),\Sing(\T_n))$, then the condition $\al U^1\sub U^1\al$ in \ref{SA1} cannot be replaced by $\al U\sub U\al$ (cf.~Remark \ref{rem:AP4}).  Indeed, if $\al\in\T_n$ is any map with image $A\subsetneq\bn$, then $\al\cdot\id_A=\al \not= \id_B\cdot \al$ for any $B\subsetneq\bn$.
For the same reason, the action of $S$ on $U^1=\E_n$ does not restrict to an action of $S$ on $U=\Sing(\E_n)$, as ${}^\al\id_A=\id_\bn$ when $\im(\al)=A\subsetneq\bn$.

Similar considerations show that none of the above pairs satisfy the right-handed version $U^1\al\sub \al U^1$ of~\ref{SA1} when $n\geq2$.  Indeed, let $\al\in\T_n$ be an arbitrary constant map, and fix some $\es\subsetneq A\subsetneq\bn$.  Then for any $B\sub\bn$, $\al\cdot\id_B$ is equal to either $\al$ or else the empty map, and so is never equal to $\id_A\cdot \al$.  This all shows that $\id_A\cdot \al\in U\al\sm \al U^1$, and so $U^1\al\not\sub \al U^1$.  

We also have strong action pairs $(\E_n,\G_n)$ and $(\Sing(\E_n),\G_n)$ in $\PT_n$.  These lead to the (almost-)factorisable inverse semigroups
\[
\E_n\cdot\G_n = \I_n \AND \Sing(\E_n)\cdot\G_n = \Sing(\I_n),
\]
which are special cases of Example \ref{eg:FIM}.
\end{eg}

Generalisations of Example \ref{eg:Tn} will be discussed in Chapters \ref{chap:IA} and \ref{chap:wreath}.

\begin{eg}\label{eg:TX}
Extending the last example, consider again the full and partial transformation semigroups $\T_X$ and $\PT_X$ over an arbitrary set $X$.  Consider also the power set ${\P(X)=\set{A}{A\sub X}}$, as a semilattice under intersection. Note that $\T_X$ acts on $\P(X)$ via
\[
{}^\al\! A=A\al^{-1} \qquad\text{for $A\in\P(X)$ and $\al\in\T_X$.}
\]
Now let $S$ be a subsemigroup of $\T_X$, let $F$ be a subsemilattice of $\P(X)$ such that $F^1=F\cup\{X\}$ is closed under the action of $S$, and set $U=\set{\id_A}{A\in F}$.
Then $(U,S)$ is a strong action pair in $M=\PT_X$, as again follows from \eqref{eq:Af}.  So we have the subsemigroup
\[
US = \set{\id_A\cdot \al}{\al\in S,\ A\in F}= \set{\al\restr_A}{\al\in S,\ A\in F} \leq\PT_X,
\]
with composition given by
\[
\al\restr_A\cdot \be\restr_B = (\al \be) \restr_{A\cap B\al^{-1}} \qquad\text{for $\al,\be\in S$ and $A,B\in F$.}
\]
When $S=\T_X$ and $F=\P(X)$, we of course obtain $US=\PT_X$, but varying $S$ and $F$ gives rise to other interesting cases.
\ben
\item \label{it:TX1} Consider the case in which $X=\N=\{0,1,2,\ldots\}$.  Define transformations $\ga_1,\ga_2\in\T_\N$ by
\[
\ga_1 = \trans{0&1&2&3&\cdots\\1&2&3&4&\cdots}
\AND
\ga_2 = \trans{0&1&2&3&\cdots\\0&2&3&4&\cdots},
\]
and define the (infinite cyclic) subsemigroups $S_1=\la \ga_1\ra$ and $S_2=\la \ga_2\ra$ of $\PT_\N$.  Let $F$ be the set of all subsets of $\N$ of the form $x+\N=\{x,x+1,x+2,\ldots\}$ with $x\geq1$.  So $F$ is a subsemilattice of $\P(\N)$, and $F\cup\{\N\}$ is closed under the actions of both $S_1$ and $S_2$ ($F$ is closed under the action of $S_2$, but not of $S_1$).  So $(U,S_1)$ and $(U,S_2)$ are both strong action pairs in $\PT_\N$, and in fact we have $US_1=US_2$.  This shows that the `$S$~component' of a (strong) action pair $(U,S)$ is not uniquely determined by the semigroup $US$ itself.  Note also that $US_1=US_2$ contains no total maps, and also no idempotents.  So in general, a semigroup $US$ can be disjoint from both $U$ and $S$.
Clearly $US_1=US_2$ contains no isomorphic copy of $U$ (as $U$ consists entirely of idempotents), but any element of $US_1=US_2$ generates an infinite cyclic subsemigroup, which is isomorphic to both $S_1$ and $S_2$.

\item \label{it:TX2} Now let $X=\Z=\{0,\pm1,\pm2,\ldots\}$, and define 
\[
\ga_3 = \trans{\cdots&0&1&2&3&\cdots\\\cdots&1&2&3&4&\cdots} \in \T_\Z.
\]
Let $S_3=\la \ga_3\ra$, and let $F$ be as in \ref{it:TX1}, but now considered as a subsemilattice of $\P(\Z)$.  Note that $F\cup\{\Z\}$ is not closed under the action of $S_3$, and in fact $(U,S_3)$ is not a weak action pair at all, as even \ref{SA1} fails.  For example, with $\mathbb P=\{1,2,3,\ldots\}\in F$ and $\N={}^{\ga_3}\mathbb P=\{0,1,2,\ldots\}\not\in F\cup\{\Z\}$, we have $\ga_3\cdot\id_{\mathbb P}\in\ga_3 U^1$; on the other hand, $\dom(\ga_3\cdot\id_{\mathbb P})=\N$, yet every element $\de$ of $U^1\ga_3$ has $\dom(\de)=\Z$ or else $\dom(\de)\sub\mathbb P$.  However, we do have
\[
A\cap B\al^{-1}\in F \qquad\text{for all $A,B\in F$ and $\al\in S_3$,}
\]
and it follows from this and the rule $\al\restr_A\cdot \be\restr_B = (\al \be) \restr_{A\cap B\al^{-1}}$ that $US_3$ is still a semigroup.  In fact, by identifying $\PT_\N$ as a subset of $\PT_\Z$ (consisting of all partial maps with domain and range contained in $\N$), we have ${US_1=US_2=US_3}$, where $S_1$ and $S_2$ are as in \ref{it:TX1} above.  This time, note that $S_3$ consists of units of $\T_\Z$.

\item \label{it:TX3} Again consider $X=\N$, and let $S=\la \ga_1\ra$, with $\ga_1\in\T_\N$ as in \ref{it:TX1}.  This time we let $F$ consist of all finite subsets of $\N$, which is again a subsemilattice of $\P(\N)$ closed under the action of~$S$.  The resulting strong action pair $(U,S)$ leads to the subsemigroup~$US$ of $\PT_\N$, as usual.  This time $US$ contains a single idempotent, namely the empty map $\es$ (which is of course a zero of~$US$), and we note that $US$ is a \emph{nilsemigroup}, in the sense that every element of $US$ has a power equal to $\es$.  It follows that~$US$ does not contain any subsemigroup isomorphic to $U$ or to $S$.
\een
Among other things, the above examples show that if one is given a semigroup known to be of the form $US$ for some (strong) action pair $(U,S)$, it is not necessarily obvious what $U$ and/or~$S$ might be.
\end{eg}

It is also worth noting that the \emph{reflection monoids} of Everitt and Fountain \cite{EF2010,EF2013} are also special cases of the construction in Example \ref{eg:TX}, or indeed of Example \ref{eg:FIM}.

The next example features a \emph{weak} action pair $(U,S)$ in partial transformation monoids, and here we no longer assume that $U$ consists of partial identities.

\begin{eg}\label{eg:PTYZ}
Let $Y$ and $Z$ be disjoint non-empty sets, let $X=Y\cup Z$, and let $M=\PT_X$.  Let $S\leq\PT_Y$ and $V\leq\PT_Z$ be arbitrary subsemigroups.  We may also regard $S$ as a subsemigroup of $\PT_X$ (as any partial transformation of $Y$ is of course a partial transformation of $X$), and we additionally define $U = \set{\id_Y\cup\al}{\al\in V}$.  So $U$ is a subsemigroup of $\PT_X$, and clearly $U\cong V$.  It is easy to check that
\begin{equation}\label{eq:albe}
\al\be=\be\al=\al \qquad\text{for all $\al\in S$ and $\be\in U^1=U\cup\{\id_X\}$.}
\end{equation}
It immediately follows that condition \ref{SA1} holds for the pair $(U,S)$.  However, condition \ref{SA2} only holds if $|U|=1$, as $\al\be=\al'\be$ for all $\al,\al'\in U$ and $\be\in S$.  By contrast, the right-handed version of~\ref{SA2} holds: $\al\be=\al'\be'\implies\be=\be'$ for all $\al,\al'\in U$ and $\be,\be'\in S$.

It also follows from \eqref{eq:albe} that condition \ref{A1} holds with respect to the constant action: ${}^\al\be=\id_X$ for all $\al\in S$ and $\be\in U^1$.  Condition \ref{A2} is equivalent to \ref{SA2}, so still fails (unless $|U|=1$).  So the weak action pair $(U,S)$ is not an action pair.  In this case we have $US=S(=SU)$.
\end{eg}

In Example \ref{eg:FIM} we considered factorisable inverse monoids, which came from strong action pairs present in any inverse monoid.  As a foreshadowing of the next section, we now consider another action pair present in any inverse semigroup, though these are generally not strong.

\begin{eg}\label{eg:inv}
Let $S$ be an inverse semigroup with semilattice of idempotents $E=E(S)$.  Define the map $S\to E^1:s\mt s^+=ss^{-1}$.  It is easy to check that \ref{SA1} holds for the pair~$(E,S)$, and that the unary operation $s\mt s^+$ satisfies \ref{+1}--\ref{+4}.  For example, if $s,t\in S$, then
\[
st^+=stt^{-1} = ss^{-1}stt^{-1}=stt^{-1}s^{-1}s=(st)(st)^{-1}s=(st)^+s,
\]
giving \ref{+2}.  For \ref{+4}, suppose $es=ft$ for some $e,f\in E$ and $s,t\in S$.  Then certainly $(es)^+=(ft)^+$.  But
\[
(es)^+ = (es)(es)^{-1} = ess^{-1}e^{-1} = ess^{-1}e = eess^{-1} = es^+ \ANDSIM (ft)^+ = ft^+,
\]
so that $es^+=ft^+$, as required.  It follows from Proposition \ref{prop:+} that $(E,S)$ is an action pair in the monoid $S^1$.  By Lemma \ref{lem:SAPAP}\ref{SAPAP2}, this pair is strong if and only if $ss^{-1}=1$ for all $s\in S$, which occurs precisely when $S$ is a group.  

To understand the pair $(E,S)$ via the original Definition \ref{defn:AP}, we need to understand the action of $S$ on $E$.  For this, let $s\in S$ and $e\in E$.  Then $se=ss^{-1}se=ses^{-1}s=ses^{-1}\cdot s$.  It follows from the proof of Proposition \ref{prop:+} (see \eqref{eq:su}) that ${}^se=ses^{-1}\cdot s^+ = ses^{-1}\cdot ss^{-1} = ses^{-1}$.  This all shows that the action of $S$ on $E$ is given by conjugation, as in Example \ref{eg:act}\ref{act2}.  

As usual, the pair $(E,S)$ leads to the subsemigroup $ES$ of $S$, and of course we have $ES=S$.  Although this is not an interesting/new subsemigroup, it leads to a number of other interesting subsemigroups, as $(E,T)$ is an action pair for any $T\leq S$, by Lemma \ref{lem:US0}\ref{US01}.  We will explore this in more detail in the next section, where we extend our scope to left restriction monoids.
\end{eg}

\sect{Left restriction monoids}\label{sect:LRAP}

As just noted, another family of natural examples comes from the class of left restriction monoids, as defined in Section \ref{sect:lrs}.  For the duration of this section we fix some such monoid~$M$ with identity~$1$.  So $M$ has a (basic) unary operation $s\mt s^+$ satisfying 
\begin{enumerate}[label=\textup{(L\arabic*)}]
\bmc2
\item \label{L1'} $x^+x=x$,
\item \label{L2'} $x^+y^+=y^+x^+$,
\item \label{L3'} $(x^+y)^+=x^+y^+$,
\item \label{L4'} $xy^+=(xy)^+x$,
\item \label{L5'} $x^+x^+=x^+$,
\item \label{L6'} $(x^+)^+=x^+$.
\emc
\een
(Recall that the defining identities are \ref{L1'}--\ref{L4'}, while \ref{L5'} and \ref{L6'} are consequences.)
It follows quickly from \ref{L1'} that $1^+=1$.  As before, in the calculations to follow we write $=_1$ to indicate an application of \ref{L1'}, and similarly for $=_2$, and so on.  

The following subsets of $M$ will play an important role in all that follows:
\begin{align*}
E=E(M) &= \set{s\in M}{s=s^2}, & L=L(M)&=\set{s\in M}{1\in Ms},\\
P=P(M) &= \set{s^+}{s\in M}, & R=R(M)&=\set{s\in M}{1\in sM},\\
T=T(M) &= \set{s\in M}{s^+=1}, & G=G(M)&=L\cap R.
\end{align*}
The elements of $P$ are called \emph{projections}.  By \ref{L2'}, \ref{L3'} and \ref{L5'},~$P$ is a (monoid) semilattice.  By~\ref{L6'}, $P$ is fixed pointwise by the ${}^+$ operation, meaning that $u=u^+$ for all $u\in P$.  If~$M$ is inverse (with $s^+=ss^{-1}$), then $P=E$ is precisely the semilattice of (all) idempotents of~$M$.  However, this need not be the case for an arbitrary left restriction semigroup, as for example with~$\PT_X$ (with $\al^+=\id_{\dom(\al)}$ and $P=\set{\id_A}{A\sub X}$), although we always of course have~${P\sub E}$.  In fact, $E$ need not be a subsemigroup at all, let alone a semilattice; again consider~$\PT_X$.  The other subsets defined above \emph{are} submonoids, however.  The submonoids~$L$ and~$R$ consist of all left and right units of $M$, respectively, and~$G$ is the group of (two-sided) units.  We have~$L=R=G$ if~$M$ is finite (see for example \cite[Lemma 2.3]{JE2019a}), but this need not hold for infinite~$M$.  We will discuss~$T$ shortly, but our first goal is to show that $(P,M)$ is an action pair; see Proposition~\ref{prop:LR1}.  For this, and for later use, it is convenient to prove a simple fact, which we have already seen in Example~\ref{eg:inv} for the special case of inverse semigroups.

\begin{lemma}\label{lem:us+}
If $S$ is a left restriction semigroup with semilattice of projections $P=P(S)$, then $(us)^+=us^+$ for any $u\in P^1$ and~$s\in S$.
\end{lemma}

\pf
This is obvious for $u=1$; otherwise, $(us)^+=_6(u^+s)^+=_3u^+s^+=_6us^+$.
\epf

\begin{prop}\label{prop:LR1}
Let $M$ be a left restriction monoid, and write $P=P(M)$.  Then $(P,M)$ is an action pair.
\end{prop}

\pf
We use Proposition \ref{prop:+}.  For \ref{SA1}, let $s\in M$.  Then for any $u\in P$ we have
\begin{equation}\label{eq:su+s}
su=_6su^+=_4(su)^+s,
\end{equation}
which gives $sP\sub Ps$.

It remains to check that the map $M\to P:s\mt s^+$ (i.e., the ${}^+$ operation itself) satisfies conditions \ref{+1}--\ref{+4}.  Now, conditions \ref{+1} and \ref{+2} are simply \ref{L1'} and \ref{L4'}, respectively.  For \ref{+3}, let $s,t\in M$.  Then
\[
(st)^+ =_1 ((s^+s)t)^+ = (s^+(st))^+ =_3 s^+(st)^+ =_2 (st)^+s^+.
\]
Finally, for \ref{+4}, suppose $u,v\in P$ and $s,t\in M$ are such that $us=vt$.  Then certainly $(us)^+=(vt)^+$, and it follows from Lemma \ref{lem:us+} that $us^+=vt^+$.
\epf

\begin{rem}\label{rem:LR1}
If $S$ is a left restriction \emph{semigroup}, then $(P,S)$ is an action pair in the monoid~$S^1$, where again we write $P=P(S)=\set{s^+}{s\in S}$.  Indeed, this follows from the same proof as above.

Alternatively, when $S$ is a left restriction semigroup, $S^1$ becomes a left restriction \emph{monoid} in a unique way by putting $1^+=1$.  Proposition \ref{prop:LR1} says that $(P^1,S^1)$ is an action pair in $S^1$, and we then apply (both parts of) Lemma \ref{lem:US0}.
\end{rem}

\begin{rem}\label{rem:LR2}
Given any action pair $(U,S)$ we have the associated semigroup $US$.  In the case of the pair $(P,M)$, this subsemigroup is precisely $PM=M$ itself.  Thus, the pair $(P,M)$ does not lead to an `interesting' subsemigroup of $M$.  Similarly, the results of Chapter \ref{chap:pres} will be of no use in describing presentations for $M$ by means of the pair $(P,M)$, since they assume we already know a presentation for $M$ (and $P$).  Nevertheless, knowing that $(P,M)$ is an action pair is useful in describing many more such pairs below (using Lemma \ref{lem:US0}).  The pair $(P,M)$ will also be extremely important in Chapter \ref{chap:PAP}.
\end{rem}

The proof of Proposition \ref{prop:LR1} utilised the alternative formulation of action pairs from Proposition \ref{prop:+}, but it is convenient to understand the associated action of~$M$ on $P(=P^1)$.  To do so, fix some $u\in P$ and $s\in M$.  Then $su=(su)^+s$, by~\eqref{eq:su+s}.  By the proof of Proposition \ref{prop:+} (see \eqref{eq:su}), we therefore have ${}^su=(su)^+s^+=(su)^+$, where we used \ref{+3} in the last step.  To summarise, the action of $M$ on $P$ is given by
\begin{equation}\label{eq:MP}
{}^su = (su)^+ \qquad\text{for $s\in M$ and $u\in P$.}
\end{equation}
If $M$ is inverse (with $s^+=ss^{-1}$), then this action has the simple form ${{}^su=su(su)^{-1}=sus^{-1}}$; cf.~Examples \ref{eg:act}\ref{act2} and \ref{eg:inv}.  It is worth noting that the action \eqref{eq:MP} goes back to work of Fountain on the narrower class of \emph{adequate semigroups} \cite{Fountain1979}; see also \cite{GG2001} where it is used in the broader class of \emph{Ehresmann semigroups}.

Next we wish to prove a result concerning action pairs of the form $(Q,S)$, where $Q\leq P$ and $S\leq M$.  The set $T$ introduced above will play an important role in this, specifically in determining which of these pairs are strong.  

\begin{lemma}\label{lem:T}
Let $M$ be a left restriction monoid, and write $T=T(M)$ and $R=R(M)$.  Then
\ben
\item \label{T1} $T$ is a submonoid of~$M$,
\item \label{T2} $R\sub T$, with equality if $M$ is inverse,
\item \label{T3} $M\sm T=\set{s\in M}{s^+\not=1}$ is a right ideal (and hence a subsemigroup) of $M$.
\een
\end{lemma}

\pf
\firstpfitem{\ref{T1}}  Define the relation $\si=\bigset{(s,t)\in M\times M}{s^+=t^+}$.  Clearly $\si$ is an equivalence (it is the kernel of the map $s\mt s^+$), and it follows quickly from Lemma \ref{lem:+}\ref{lem+2} that $\si$ is a left congruence.  Lemma \ref{lem:RC} then tells us that $T=\set{s\in M}{s\mr\si1}$ is a submonoid.

\pfitem{\ref{T2}}  Lemma \ref{lem:group}\ref{group1} gives $R\sub T$.  If $M$ is inverse, and if $s\in T$, then $1=s^+=ss^{-1}$, so $s\in R$.  

\pfitem{\ref{T3}}  Let $s\in M\sm T$ and $t\in M$.  So $s^+\not=1$, and we must show that $(st)^+\not=1$.  But if $(st)^+=1$, then $s^+ = s^+\cdot1 = s^+(st)^+ =_3 (s^+\cdot st)^+ =_1 (st)^+ = 1$, a contradiction.
\epf

\begin{rem}
The equality $T=R$ need not hold if $M$ is not inverse.  For example, in $\PT_X$ we have
\[
T(\PT_X) = \T_X \AND R(\PT_X) = \set{\al\in\T_X}{\al \text{ is injective}}.
\]
The same example also shows that the right ideal $M\sm T$ (which is $\PT_X\sm\T_X$ in this case) is not always a left ideal.
\end{rem}

The next result follows immediately from Proposition \ref{prop:LR1}, and Lemmas~\ref{lem:SAPAP} and \ref{lem:US0}.

\begin{prop}\label{prop:LR2}
Let $M$ be a left restriction monoid, and write $P=P(M)$ and $T=T(M)$.  Then
\ben
\item \label{LR21} $(Q,S)$ is an action pair for any $Q\leq P$ and $S\leq M$ such that $Q^1$ is closed under the action of $S$ given in \eqref{eq:MP},
\item \label{LR22} such a pair $(Q,S)$ is strong if and only if $S\sub T$.  \epfres
\een
\end{prop}

\begin{rem}\label{rem:PS}
In particular, $(P,S)$ and $(\Pf,S)$ are both action pairs for any $S\leq M$, where we write $\Pf=P\sm\{1\}\leq P$.  Such pairs are strong if and only if $S\sub T$.  Thus, the following action pairs are all strong:
\[
(P,T) \COMMA
(P,R) \COMMA
(P,G) \COMMA
(\Pf,T) \COMMA
(\Pf,R) \AND
(\Pf,G) .
\]
As usual, these lead to the subsemigroups $PT$, $PR$ and so on.  
We first note that $PG(=GP)$ is in fact a factorisable inverse monoid.  Indeed, this follows from the easily checked facts that:
\bit
\item $E(PG)=P$, 
\item $g^{-1}u={}^{g^{-1}}ug^{-1}\in PG$ is an inverse of $ug$ for any $u\in P$ and $g\in G$.
\eit
The inverse semigroup $\Pf G(=G\Pf)=PG\sm G$ is \emph{almost-factorisable}.  
The submonoid $PT$ and the subsemigroup $\Pf T$ are \emph{left-factorisable} and \emph{almost-left-factorisable} in the terminology of \cite{Szendrei2013}.  Note that $PG=GP$ has a left-right duality, where the corresponding right action is given by ${e^g=g^{-1}eg}$, making it a restriction monoid; similar comments apply to the semigroup ${\Pf G=G\Pf}$.

(As we will see in Chapter \ref{chap:pres}, the very fact that $\Pf=P\sm\{1\}$ is a subsemigroup of $P$ is important in certain applications; see for example Assumption \ref{ass:US} and Remark \ref{rem:assUS}.)

The action pair $(P,L)$ is strong precisely when $L\sub T$.  When $M$ is inverse, this is equivalent to $L\sub R$ (cf.~Lemma \ref{lem:T}\ref{T2}), and hence to $L=R=G$.

Finally, the action pair $(P,M)$ itself is strong if and only if $M=T$, meaning that the left restriction structure of $M$ is trivial ($s^+=1$ for all $s\in M$).  When $M$ is inverse, this is equivalent to $M$ being a group ($ss^{-1}=1$ for all $s\in M$), as we observed in Example \ref{eg:inv}.  For the action pair $(\Pf,M)$, one can check that
\[
\Pf M = M\sm T = \set{s\in M}{s^+\not=1}.
\]
We showed in Lemma \ref{lem:T}\ref{T3} that this subsemigroup is in fact a right ideal of $M$.
\end{rem}

\sect{The first structure theorem}\label{sect:structure}

Having now discussed an extensive collection of examples, including the class of left restriction semigroups, we return to the general theory.  The main goal of this section is to describe the structure of a semigroup $US$ arising from an action pair~$(U,S)$, showing how such a semigroup can be built from a semidirect product $U\rtimes S$ and a special congruence.  We do this in Theorem~\ref{thm:FT} below, but we begin with a simple result that holds more generally for a \emph{weak} action pair $(U,S)$.  For such a pair, the action from \ref{A1} allows for the formation of the (external) semidirect product $U\rtimes S$, as in Definition \ref{defn:SD}.

\begin{prop}\label{prop:A1}
If $(U,S)$ is a weak action pair, then the map $\pi:U\rtimes S\to US$ given by $(u,s)\pi=us$ is a surmorphism with kernel
\[
\th = \bigset{((u,s),(v,t))}{u,v\in U,\ s,t\in S,\ us=vt}.
\]
Consequently, $\th$ is a congruence on $U\rtimes S$, and $US\cong(U\rtimes S)/\th$.
\end{prop}

\pf
If $u,v\in U$ and $s,t\in S$, then
\[
(u,s)\pi\cdot(v,t)\pi = us\cdot vt = u\cdot {}^sv\cdot st = (u\cdot {}^sv,st)\pi = ((u,s)\cdot(v,t))\pi,
\]
so that $\pi$ is a morphism.  Surjectivity is clear, as is the description of the kernel.  The final assertion follows from the First Isomorphism Theorem.
\epf

\begin{rem}\label{rem:A1}
Recall from Proposition \ref{prop:MM} that
\[
M = \bigset{(u,s)\in U\rtimes S}{u=us^+}
\]
is a subsemigroup of $U\rtimes S$, and that the map
\[
F:U\rtimes S\to M:(u,s)\mt(us^+,s)
\]
is a retraction.  (Note that to apply Proposition \ref{prop:MM} we need to keep Remark \ref{rem:AP2} in mind, which says that if $S$ is a monoid then it acts monoidally on $U^1$.)  Because of the identity $s=s^+s$ (cf.~Remark \ref{rem:AP2}), the following diagram commutes:
\begin{equation}\label{eq:USMUS}
\begin{tikzcd} 
 & M  \arrow{dr}{\pi\restr_M} & \\
U\rtimes S \arrow{ur}{F} \arrow[swap]{rr}{\pi} & & US
\end{tikzcd}
\end{equation}
In particular, $\pi\restr_M:M\to US$ is still a surmorphism, so we also have
\[
US\cong M/\vt \WHERE \vt = \ker(\pi\restr_M) = \th\restr_M.
\]
\end{rem}

Now consider an action pair $(U,S)$.  Of course Proposition \ref{prop:A1} holds for this pair, and we keep the meaning of the surmorphism $\pi:U\rtimes S\to US$, and the congruence $\th=\ker(\pi)$.  For any $u\in U^1$, we define the relation
\begin{equation}\label{eq:the}
\th_u = \bigset{(s,t)\in S\times S}{us=ut}.
\end{equation}
We also continue to write $s^+={}^s1$ for $s\in S$.

\begin{lemma}\label{lem:th}
Let $(U,S)$ be an action pair, and let $u,v\in U$ and $s,t\in S$.
\ben
\item \label{th1} We have $(u,s)\mr\th(v,t) \iff us^+=vt^+$ and $(s,t)\in\th_{us^+}$.
\item \label{th2} If $(U,S)$ is strong, then $(u,s)\mr\th(v,t) \iff u=v$ and $(s,t)\in\th_u$.
\een
\end{lemma}

\pf
The second part follows from the first, given Lemma \ref{lem:SAPAP}\ref{SAPAP2}.  For the first part, suppose $(u,s)\mr\th(v,t)$: i.e., $us=vt$.  Then $us^+=vt^+$ follows from \ref{+4}.  Together with \ref{+1}, it follows that $us^+s=us=vt=vt^+t=us^+t$, so $(s,t)\in\th_{us^+}$.  The converse is similar.
\epf

Proposition \ref{prop:A1} showed that, given an action pair $(U,S)$, the semigroup $US$ can be viewed abstractly as a homomorphic image of a semidirect product $U\rtimes S$.  Our next goal is to establish a converse of this fact; namely, we wish to describe \emph{all} homomorphic images of semidirect products that correspond to action pairs in this way.  To do so, we first establish some special properties of the above congruence $\th$, and the associated relations $\th_u$.  
Here and elsewhere, if $\si$ is a relation on a semigroup $T$, and if $x\in T$, we write $x\cdot\si=\bigset{(xs,xt)}{(s,t)\in\si}$.  The diagonal/equality relation on any set $X$ is denoted $\De_X=\bigset{(x,x)}{x\in X}$.

\begin{lemma}\label{lem:A2}
For an action pair $(U,S)$, and with the above notation, the following all hold:
\ben
\item \label{A21} For all $u\in U$ and $s\in S$, we have $(u,s)\mr\th(us^+,s)$.
\item \label{A22} For all $s,t\in S$, we have $(s^+,s)\mr\th(t^+,t)\implies s=t$.  
\item \label{A23} For all $u,v\in U$ and $s,t\in S$, we have $(u,s)\mr\th(v,t)\implies us^+ = vt^+$.
\item \label{A24} We have $\th_1=\De_S$.
\item \label{A25} For all $u\in U$, the relation $\th_u$ is a right congruence on $S$.
\item \label{A26} For all $u,w\in U$, we have $\th_u\sub\th_{wu}$.
\item \label{A27} For all $u\in U$ and $x\in S$, we have $x\cdot\th_u\sub\th_{{}^xu}$.
\item \label{A28} For all $u,w\in U$ and $(s,t)\in\th_u$, we have $u\cdot{}^sw=u\cdot{}^tw$ and $(s,t)\in\th_{u\cdot{}^sw}$.
\een
\end{lemma}

\pf
These are easily checked, using the definition of $\th$ and $\th_u$.  We additionally use~\ref{+1} for~\ref{A21} and~\ref{A22}, and Lemma \ref{lem:th} for \ref{A23}.  We just give the details for \ref{A27} and \ref{A28}.

%
%
%
%
%

\pfitem{\ref{A27}}  We have
\[
(s,t)\in\th_u \implies us=ut \implies xus=xut \implies {}^xu\cdot xs = {}^xu\cdot xt \implies (xs,xt)\in\th_{{}^xu}.
\]

\pfitem{\ref{A28}}  We begin by noting that
\[
(s,t)\in\th_u \implies us=ut \implies u\cdot{}^sw = u\cdot{}^tw,
\]
by Lemma \ref{lem:ap}.  Combining this with \ref{A1} and $us=ut$, it also follows that
\[
u\cdot{}^sw\cdot s = usw = utw = u\cdot{}^tw\cdot t = u\cdot{}^sw\cdot t,
\]
so that $(s,t)\in\th_{u\cdot{}^sw}$.
\epf

We now use these properties of $\th$ (and the $\th_u$) to construct a set of axioms to abstractly characterise the quotients of semidirect products giving rise to action pairs.  To keep the `concrete' pairs $(U,S)$ separate from the abstract, we will use the notation $V$ and $T$ for the latter.

\begin{defn}\label{defn:FT}
Suppose $V$ and $T$ are semigroups such that $T$ has a left action on $V^1$ (the monoid completion of $V$) by semigroup morphisms, written $(s,u)\mt{}^su$, for $s\in T$ and $u\in V^1$.  For $s\in T$, we write $s^+={}^s1\in V^1$.  Suppose $\si$ is a congruence on the semidirect product $V\rtimes T$, and for each $u\in V$ define the relation
\[
\si_u=\bigset{(s,t)\in T\times T}{(u,s)\mr\si(u,t)}.
\]
If $V$ is not a monoid (so $1\not\in V$), we additionally define $\si_1=\De_T$.
We say $\si$ is \emph{special} if the following all hold:  
\begin{enumerate}[label=\textup{(S\arabic*)}]
\item \label{S1} For all $u\in V$ and $s\in T$, we have $(u,s)\mr\si(us^+,s)$.
\item \label{S2} For all $s,t\in T$, we have $(s^+,s)\mr\si(t^+,t)\implies s=t$.
\item \label{S3} For all $u,v\in V$ and $s,t\in T$, we have $(u,s)\mr\si(v,t)\implies us^+ = vt^+$.
\item \label{S4} We have $\si_1=\De_T$.
\item \label{S5} For all $u\in V$, $\si_u$ is a right congruence on~$T$.
\item \label{S6} For all $u,w\in V$, we have $\si_u\sub\si_{wu}$.
\item \label{S7} For all $u\in V$ and $x\in T$, we have $x\cdot\si_u\sub\si_{{}^xu}$.
\item \label{S8} For all $u,w\in V$ and $(s,t)\in\si_u$, we have $u\cdot{}^sw=u\cdot{}^tw$ and $(s,t)\in\si_{u\cdot{}^sw}$.
\een
\end{defn}

If $T$ acts by monoid morphisms (i.e., if $s^+=1$ for all $s\in T$), then several simplifications arise in \ref{S1}--\ref{S8}.  For example, \ref{S1} is trivially true, and \ref{S3} says $(u,s)\mr\si(v,t)\implies u=v$.

Here is the main structural result of this chapter.

\begin{thm}\label{thm:FT}
Suppose $V$ and $T$ are semigroups such that $T$ has a left action on $V^1$ by semigroup morphisms.  Suppose also that $\si$ is a special congruence on the semidirect product $V\rtimes T$.  Then there is an action pair $(U,S)$ in some monoid $M$ such that
\[
U\cong V \COMMA S\cong T \AND US\cong(V\rtimes T)/\si.
\]
If $T$ acts on $V^1$ by monoid morphisms, then $(U,S)$ is strong.

Conversely, given any action pair $(U,S)$, we have $US\cong(U\rtimes S)/\th$ for some special congruence~$\th$ on the semidirect product $U\rtimes S$.
\end{thm}

\pf
The last assertion follows from Proposition \ref{prop:A1} and Lemma \ref{lem:A2}.  The remainder of the proof concerns the first two assertions, for which we fix $V$, $T$ and $\si$ with the specified properties.
Without loss of generality, we may assume that $T$ and $V^1$ are disjoint.  With a slight abuse of notation, we let $T^1=T\cup\{1\}$ be the monoid obtained from $T$ by adjoining $1$ (the identity of~$V^1$) as an identity, \emph{whether or not~$T$ was already a monoid}.  The reason for doing this is that $T$ might be a monoid acting non-monoidally on $V$, and we will soon construct an action of $T^1$ on $V^1$ that is required to be monoidal.  (But we stress that $V=V^1$ if $V$ is a monoid.)  

We denote the action of $T$ on $V^1$ by $(s,u)\mt{}^su$ ($s\in T$, $u\in V^1$), and we extend this to a monoidal action of $T^1$ on $V^1$ in the usual way, by further defining ${}^1u=u$ for all $u\in V^1$.
We then have the semidirect product $V^1\rtimes T^1$, which contains $V\rtimes T$ as a proper subsemigroup.  (Note that $V^1\rtimes T^1$ is not necessarily a monoid, although $(1,1)$ is always a left identity, as the action of~$T^1$ on~$V^1$ is monoidal; cf.~Lemma \ref{lem:USmon}.)  Throughout the proof, we write $s^+={}^s1$ for $s\in T^1$; we will make frequent use of the properties of these elements listed in Lemma \ref{lem:+}, often without explicit reference.  Since the action is monoidal, we have $1^+=1$.

The main effort in the proof goes into establishing the following:

\begin{lemma}\label{lem:Si}
There is a congruence $\Si$ on $V^1\rtimes T^1$ such that the following all hold:
\begin{enumerate}[label=\textup{($\Si$\arabic*)},leftmargin=12mm]
\item \label{Si1} $\si=\Si\restr_{V\rtimes T}$,
\item \label{Si2} For all $u,v\in V^1$, we have $(u,1) \mr\Si (v,1) \iff u=v$.
\item \label{Si3} For all $u\in V^1$ and $s\in T^1$, we have $(u,s) \mr\Si (us^+,s)$.
\item \label{Si4} For all $u,v\in V^1$ and $s,t\in T^1$, we have $(u,s) \mr\Si (v,t) \implies us^+=vt^+$.
\een
\end{lemma}

\pf
We define $\Si$ by specifying its equivalence classes.  Each $\si$-class in $V\rtimes T$ will determine a unique $\Si$-class, so consider some such $\si$-class~$K$.  By \ref{S2}, $K$ contains at most one element of the form $(s^+,s)$ with $s\in T$.  
\begin{enumerate}[label=\textup{(K\arabic*)}]
\item \label{K1} If $K$ does contain such an element $(s^+,s)$, then $K\cup\big\{(1,s)\big\}$ is a $\Si$-class.  (Note that this is just $K$ itself if $s^+=1\in V$.  If $s^+=1\not\in V$, then no $\si$-class contains $(s^+,s)$.)
\item \label{K2} Otherwise, $K$ is a $\Si$-class.  
\een
All elements of $V^1\rtimes T^1$ not yet assigned to $\Si$-classes are now assigned to singleton $\Si$-classes.  These elements are precisely:
\bit
\item $(u,1)$, for each $u\in V^1$, and
\item $(1,s)$, for each $s\in T$ for which $s^+=1\not\in V$ (if any).  (Note that such an element is still of the form $(1,s)=(s^+,s)$.)
\eit
It is clear from the construction that $\Si$ is an equivalence on $V^1\rtimes T^1$, and that \ref{Si1} and \ref{Si2} both hold.  So it remains to show that $\Si$ is compatible, and establish \ref{Si3} and \ref{Si4}.  We begin with these two conditions.

\pfitem{\ref{Si3}}  This follows from \ref{S1} when $u\in V$ and $s\in T$, and is trivial when $s=1$.  When $u=1$,~\ref{Si3} says $(1,s)\mr\Si(s^+,s)$; this is clear if $s^+=1$, and is true by definition of $\Si$ otherwise; cf.~\ref{K1}.

\pfitem{\ref{Si4}}  Suppose $(u,s) \mr\Si (v,t)$, where $u,v\in V^1$ and $s,t\in T^1$.  We must show that $us^+=vt^+$.  This is obvious if $(u,s)=(v,t)$, or follows from \ref{S3} if $(u,s) \mr\si (v,t)$, so we assume otherwise.  By construction (and by symmetry), it follows that $(u,s)=(1,s)$ and $(v,t)$ belongs to the $\si$-class of $(s^+,s)$, with $s^+\in V$; cf.~\ref{K1}.  But then $(v,t)\mr\si(s^+,s)$, and it follows from \ref{S3} and Lemma~\ref{lem:+}\ref{lem+1} that
\[
vt^+ = s^+s^+ = s^+ = 1\cdot s^+ = us^+,
\]
as required.

\medskip

As noted above, it remains (for the proof of the lemma) to establish the compatibility of $\Si$.  To do so, let $(\ba,\bb)\in\Si$ and~$\bc\in V^1\rtimes T^1$.  We must show that
\begin{equation}\label{eq:comp}
\bc\cdot(\ba,\bb)=(\bc\ba,\bc\bb)\in\Si \AND (\ba,\bb)\cdot\bc=(\ba\bc,\bb\bc)\in\Si.
\end{equation}
Clearly \eqref{eq:comp} holds if $\ba=\bb$, so we now assume that $\ba\not=\bb$.  Since $\ba\mr\Si\bb$, it follows that $\ba$ and~$\bb$ belong to a $\Si$-class of type \ref{K1} or \ref{K2}.  
By the form of these classes, we may write 
\[
\ba=(u,s)\COMMa \bb=(v,t) \ANd \bc=(w,x)\COMMA \text{where $u,v,w\in V^1$, $s,t\in T$ and $x\in T^1$,}
\]
and we have
\begin{equation}\label{eq:abc}
\bc\cdot(\ba,\bb) = ((w\cdot{}^xu,xs),(w\cdot{}^xv,xt)) \AND (\ba,\bb)\cdot\bc = ((u\cdot{}^sw,sx),(v\cdot{}^tw,tx)).
\end{equation}

\pfcase1  Suppose first that $(\ba,\bb)\in\si$; in particular, we have $u,v\in V$.  If $\bc\in V\rtimes T$, then \eqref{eq:comp} follows from the fact that $\si$ is compatible and $\si\sub\Si$.  So we are left to consider the cases in which $w=1$ and/or $x=1$.  We will soon consider these separately, but we first note that \ref{S3} gives $us^+=vt^+$.  Combining this with \ref{S1}, we have
\begin{equation}\label{eq:stus+}
(us^+,s) \mr\si (u,s) \mr\si (v,t) \mr\si (vt^+,t) = (us^+,t) \ \implies \ (s,t)\in\si_{us^+}.
\end{equation}

\pfcase{1.1}  Suppose $w=1$.  Examining \eqref{eq:abc}, we must show that
\bena\bmc2
\item \label{a} $({}^xu,xs) \mr\Si ({}^xv,xt)$, and
\item \label{b} $(us^+,sx) \mr\Si (vt^+,tx)$.
\emc\een
For \ref{b}, we note that \eqref{eq:stus+} and \ref{S5} give $(sx,tx)\in \si_{us^+}$, and so
\[
(us^+,sx) \mr\si (us^+,tx) = (vt^+,tx).
\]
For \ref{a}, we note that \eqref{eq:stus+} and \ref{S7} give $(xs,xt) \in \si_{{}^x(us^+)}$.  There is now a subtle point concerning the possibility that ${}^x(us^+)$ might equal $1$.
\bit
\item If ${}^x(us^+)=1$, then it follows from~\ref{S4} that $xs=xt$, and so $({}^x(us^+),xs) = ({}^x(us^+),xt)$.  
\item If ${}^x(us^+)\not=1$, then ${}^x(us^+)\in V$ and $({}^x(us^+),xs) \mr\si ({}^x(us^+),xt)$, by definition of $\si_{{}^x(us^+)}$, since $(xs,xt) \in \si_{{}^x(us^+)}$.  
\eit
Either way, we have $({}^x(us^+),xs) \mr\Si ({}^x(us^+),xt)$.  Now,
\[
{}^x(us^+) = {}^xu \cdot {}^x(s^+) = {}^xu(xs)^+ \ANDSIM {}^x(us^+) = {}^x(vt^+) = {}^xv(xt)^+.
\]
Combining all the above with \ref{Si3}, we conclude that
\[
({}^xu,xs) \mr\Si ({}^xu(xs)^+,xs) = ({}^x(us^+),xs) \mr\Si ({}^x(us^+),xt) = ({}^xv(xt)^+,xt) \mr\Si ({}^xv,xt),
\]
completing the proof of \ref{a}.

\pfcase{1.2}  Now suppose $x=1$.  Since we have already dealt with the $w=1$ case, we will also assume here that $w\not=1$.  Examining \eqref{eq:abc} again, this time we must show that
\bena\bmc2\addtocounter{enumi}{2}
\item \label{c} $(wu,s) \mr\Si (wv,t)$, and
\item \label{d} $(u\cdot{}^sw,s) \mr\Si (v\cdot{}^tw,t)$.
\emc\een
For \ref{c}, we note that \eqref{eq:stus+} and \ref{S6} give $(s,t)\in\si_{wus^+}$.  Since $wus^+\in V$ (as $u\in V$), it follows that $(wus^+,s)\mr\si(wus^+,t)$.  Combining this with \ref{S1} and $us^+=vt^+$, we have
\[
(wu,s) \mr\si (wus^+,s) \mr\si(wus^+,t) = (wvt^+,t) \mr\si (wv,t).
\]
For \ref{d}, and keeping $w\not=1$ in mind, we note that \eqref{eq:stus+} and \ref{S8} give
\begin{align*}
us^+\cdot{}^sw = us^+\cdot{}^tw &\AND (s,t)\in\si_{us^+\cdot{}^sw}.
\intertext{By Lemma \ref{lem:+}\ref{lem+4} we have $us^+\cdot{}^sw = u\cdot{}^sw$ and $us^+\cdot{}^tw=vt^+\cdot{}^tw=v\cdot{}^tw$.  So in fact,}
u\cdot{}^sw = v\cdot{}^tw &\AND (s,t)\in\si_{u\cdot{}^sw},
\end{align*}
which gives $(u\cdot{}^sw,s) \mr\si (u\cdot{}^sw,t) = (v\cdot{}^tw,t)$, completing the proof of \ref{d}.

\pfcase2  We now assume that $(\ba,\bb)\notin\si$ (but we still of course have $(\ba,\bb)\in\Si$).  By definition of $\Si$, and swapping $\ba$ and $\bb$ if necessary, we may assume that $\ba=(1,s)$, and that $\bb=(v,t)$ belongs to the $\si$-class of $\bd=(s^+,s)$, with $s^+\in V$.  Since $(\bb,\bd)\in\si$, it follows from Case 2 that $\bc\cdot(\bb,\bd)\in\Si$ and $(\bb,\bd)\cdot\bc\in\Si$.  So by transitivity, we can complete the proof of \eqref{eq:comp}, and hence of the lemma, by showing that
\begin{align*}
\bc\cdot(\ba,\bd)\in\Si &\AND (\ba,\bd)\cdot\bc\in\Si.
\intertext{For this, we first note that}
\bc\cdot(\ba,\bd)=((wx^+,xs),(w\cdot{}^x(s^+),xs) &\AND (\ba,\bd)\cdot\bc = (({}^sw,sx),(s^+\cdot{}^sw,sx)).
\end{align*}
By \ref{Si3}, we have
\[
\bc\ba = (wx^+,xs) \mr\Si (wx^+(xs)^+,xs) = (w(xs)^+,xs)=(w\cdot{}^x(s^+),xs)=\bc\bd.
\]
On the other hand, we have $\ba\bc=\bd\bc$, as ${}^sw=s^+\cdot{}^sw$.
We have now completed the proof of Lemma \ref{lem:Si}.
\epf

Returning now to the proof of the theorem, we fix a congruence $\Si$ on $V^1\rtimes T^1$ as in Lemma~\ref{lem:Si}, and we define $M=(V^1\rtimes T^1)/\Si$.  We denote the $\Si$-class of $(u,s)\in V^1\rtimes T^1$ by $[u,s]$, so that
\[
M = \bigset{[u,s]}{u\in V^1,\ s\in T^1} \qquad\text{with operation}\qquad [u,s]\cdot[v,t] = [u\cdot{}^sv,st].
\]
We already observed that $(1,1)$ is a left identity for $V^1\rtimes T^1$, so of course $[1,1]$ is a left identity for $M$.  It follows from \ref{Si3} that
\[
[u,s]\cdot[1,1] = [us^+,s] = [u,s] \qquad\text{for all $u\in V^1$ and $s\in T^1$,}
\]
so that $[1,1]$ is a right identity for $M$.  Consequently, $M$ is a monoid with identity $[1,1]$.

Now we define
\[
U = \bigset{[u,1]}{u\in V} \AND S = \bigset{[s^+,s]}{s\in T}.
\]
Clearly $U\leq M$, and it follows from \ref{Si2} that $U\cong V$.  We have $S\leq M$ because of the rule $[s^+,s]\cdot[t^+,t]=[(st)^+,st]$, which itself follows from $s^+\cdot{}^s(t^+) = {}^s(t^+)=(st)^+$.  We have $S\cong T$ because of \ref{S2}.  

Next we verify that $(U,S)$ is an action pair in $M$ by using Proposition \ref{prop:+}.  Property~\ref{SA1} follows from
\[
[s^+,s]\cdot[u,1] = [{}^su,s] = [{}^su,1]\cdot[s^+,s].
\]
We define the map $S\to U^1:[s^+,s]\mt[s^+,s]^+=[s^+,1]$.  When checking \ref{+1}--\ref{+3} we make extensive use of Lemma \ref{lem:+}.  For \ref{+4}, we must show that
\[
[u,1]\cdot[s^+,s] = [v,1]\cdot[t^+,t] \implies [u,1]\cdot[s^+,s]^+ = [v,1]\cdot[t^+,t]^+ \qquad\text{for all $u,v\in V^1$ and $s,t\in T$.}
\]
Expanding the products, this amounts to showing that
\[
[us^+,s] = [vt^+,t] \implies [us^+,1] = [vt^+,1]
\]
for all such $u,v,s,t$.  For this we use \ref{Si3} and \ref{Si4}:
\begin{align*}
[us^+,s] = [vt^+,t] 
\implies [u,s] = [v,t] 
\implies us^+ =vt^+ 
\implies [us^+,1] = [vt^+,1] .
\end{align*}
This completes the proof that $(U,S)$ is an action pair in $M$.

If $T$ acts on $V^1$ by monoid morphisms, then $s^+={}^s1=1$ for all $s\in T$, and it follows that $[s^+,s]^+=[1,1]$ for all $s$, so that $(U,S)$ is strong by Lemma~\ref{lem:SAPAP}\ref{SAPAP2}.

Next we note that
\[
US = \bigset{[u,1][s^+,s]}{u\in V,\ s\in T} = \bigset{[u,s]}{u\in V,\ s\in T},
\]
where in the last equality we use the fact that $[u,1][s^+,s] = [us^+,s] = [u,s]$ for $u\in V$ and $s\in T$, which itself follows from \ref{Si3}.  Finally, we define a map
\[
\Psi : V\rtimes T \to US \qquad\text{by}\qquad (u,s)\Psi = [u,s].
\]
It is clear that $\Psi$ is a surmorphism (it is the restriction to $V\rtimes T$ of the natural projection $V^1\rtimes T^1\to(V^1\rtimes T^1)/\Si=M$).  If $u,v\in V$ and $s,t\in T$, then by the definition of $\Psi$, and by~\ref{Si1}, we also have
\[
(u,s)\Psi = (v,t)\Psi \iff [u,s] = [v,t] \iff (u,s) \mr\Si (v,t) \iff (u,s) \mr\si (v,t) ,
\]
so that $\ker(\Psi)=\si$.  Thus, $US\cong(V\rtimes T)/\si$ by the First Isomorphism Theorem.
\epf

\begin{rem}\label{rem:A12}
Recall from Remark \ref{rem:A1} that the surmorphism $\pi:U\rtimes S\to US$ factors through the subsemigroup
\[
M = \bigset{(u,s)\in U\rtimes S}{u=us^+} \leq U\rtimes S.
\]
It is possible to use Theorem \ref{thm:FT} to classify the congruences $\si$ on semigroups of the form~$M$ such that $M/\si$ comes from an action pair.  We have chosen not to do this, however, because semidirect products are arguably the more natural/direct construction.  Nevertheless, the semigroup $M$ will play a very important role in Chapter \ref{chap:PAP}.
\end{rem}

\sect{Congruence conditions}\label{sect:cong}

We conclude Chapter \ref{chap:AP} with a sequence of results concerning generating sets for the congruence~$\th$ on $U\rtimes S$ defined in Proposition \ref{prop:A1}, where $(U,S)$ is an action pair in a monoid $M$.  For convenience, let us recall that
\[
\th = \bigset{((u,s),(v,t))}{u,v\in U,\ s,t\in S,\ us=vt}.
\]
The reason we are interested in describing such generating sets is that they feed into presentations for the semigroup $US\cong(U\rtimes S)/\th$ in Chapter \ref{chap:pres}.

The first result is the most general, and applies to any action pair $(U,S)$.  Roughly speaking, it shows that $\th$ is generated by pairs of two kinds: 
\bit
\item pairs with fixed $U$-coordinate, coming from the right congruences $\th_u$ in \eqref{eq:the}, and 
\item pairs with fixed $S$-coordinate, coming from the $s\mt s^+$ map in Proposition \ref{prop:+}.  
\eit

\begin{lemma}\label{lem:Om1}
If $(U,S)$ is an action pair, then $\th=\Om^\sharp$, where
\[
\Om = \bigset{((u,s),(u,t))}{u\in U,\ (s,t)\in\th_u} \cup \bigset{((u,s),(us^+,s))}{u\in U,\ s\in S}.
\]
\end{lemma}

\pf
Write $\th'=\Om^\sharp$.  It is easy to check that $\Om\sub\th$, so $\th'\sub\th$.  It remains to show that $\th\sub\th'$, so suppose $(u,s)\mr\th(v,t)$.  By Lemma \ref{lem:th}\ref{th1}, we have $us^+=vt^+$ and $(s,t)\in\th_{us^+}$.  But then
\[
(u,s) \mr\th' (us^+,s) \mr\th' (us^+,t) = (vt^+,t) \mr\th' (v,t). \qedhere
\]
\epf

If the pair $(U,S)$ is strong, then the previous result is essentially vacuous; cf.~Lemma \ref{lem:th}\ref{th2}.  The generating set $\Om$ from Lemma \ref{lem:Om1} can be reduced when $U$ and $S$ are both submonoids of the over-monoid~$M$, as we now show.  For the statement, recall from Remark \ref{rem:AP2} that $s=s^+s$, which tells us that $(1,s)\mr\th(s^+,s)$.

\begin{lemma}\label{lem:Om2}
If $(U,S)$ is an action pair, and if $U$ and $S$ are both submonoids, then $\th=\Om^\sharp$, where
\[
\Om = \bigset{((u,s),(u,t))}{u\in U,\ (s,t)\in\th_u} \cup \bigset{((1,s),(s^+,s))}{s\in S}.
\]
\end{lemma}

\pf
Write $\th'=\Om^\sharp$.  It again suffices to show that $\th\sub\th'$.  Given Lemma \ref{lem:Om1}, it suffices to show that $(u,s) \mr\th' (us^+,s)$ for all $u\in U$ and $s\in S$.  For this, we have
\[
(u,s) = (u,1)(1,s) \mr\th' (u,1)(s^+,s) =(us^+,s).  \qedhere
\]
\epf

Further simplifications are available when the pair $(U,S)$ is strong.  We begin by discussing certain conditions that might be satisfied by such a pair.

\begin{lemma}\label{lem:w}
Suppose $(U,S)$ is a strong action pair, and consider the following conditions:
\ben
\bmc2
\item \label{w1} $(\exists u\in U)\ (\forall s\in S)\ s=su$,
\item \label{w2} $(\forall s\in S)\ (\exists u\in U)\ s=su$,
\item \label{w3} $(\exists u\in U)\ (\forall s\in S)\ s=us$,
\item \label{w4} $(\forall s\in S)\ (\exists u\in U)\ s=us$,
\item \label{w5} $(\forall s\in S)\ (\exists u\in U)\ (\exists t\in S)\ s=ut$,
\item \label{w6} $(\exists s,t\in S)\ (\exists u\in U)\ s=ut$,
\item \label{w7} $(\exists v\in U)\ (\forall u\in U)\ (\forall s\in S)\ us=usv$,
\item \label{w8} $(\forall u\in U)\ (\forall s\in S)\ (\exists v\in U)\ us=usv$,
\item \label{w9} $(\exists v\in U)\ (\forall u\in U)\ (\forall s\in S)\ u=u\cdot{}^sv$,
\item \label{w10} $(\forall u\in U)\ (\forall s\in S)\ (\exists v\in U)\ u=u\cdot{}^sv$.
\emc
\een
The following implications hold:
\[
\begin{tikzpicture}
\nc\scl{1.6}
\foreach \x in {3,...,8} {\node () at (\x*\scl,0) {\ref{w\x}};}
\foreach \x/\y in {7/1,8/2} {\node () at (\x*\scl,\scl) {\ref{w\y}};}
\foreach \x/\y in {7/9,8/10} {\node () at (\x*\scl,-\scl) {\ref{w\y}};}
\foreach \x in {3,...,7} {\node () at (\x*\scl+.5*\scl,0) {$\Rightarrow$};}
\foreach \x in {7} {\node () at (\x*\scl+.5*\scl,\scl) {$\Rightarrow$};}
\foreach \x in {7} {\node () at (\x*\scl+.5*\scl,-\scl) {$\Rightarrow$};}
\foreach \x in {7,8} {\node () at (\x*\scl,.5*\scl) {$\Downarrow$};}
\foreach \x in {7,8} {\node () at (\x*\scl,-.5*\scl) {$\Updownarrow$};}
\end{tikzpicture}
\]
\end{lemma}

\pf
All implications other than \ref{w6}$\implies$\ref{w7}, \ref{w7}$\iff$\ref{w9} and \ref{w8}$\iff$\ref{w10} are obvious, and in fact hold for non-strong action pairs as well.  We treat the other implications now.

\pfitem{\ref{w6}$\implies$\ref{w7}}  If \ref{w6} holds, then we have $x=vy$ for some $x,y\in S$ and $v\in U$.  To show that \ref{w7} holds (with respect to this $v$), let $u\in U$ and $s\in S$ be arbitrary.  Then $u\cdot sx = usvy = u\cdot {}^sv\cdot sy$, so \ref{SA2} gives $u=u\cdot{}^sv$.  But then $us= u\cdot{}^sv\cdot s = usv$, as required.

\pfitem{\ref{w7}$\iff$\ref{w9} and \ref{w8}$\iff$\ref{w10}}  Using $usv = u\cdot {}^sv\cdot s$ and \ref{SA2}, we have ${us=usv \iff u=u\cdot{}^sv}$.
\epf

\begin{rem}\label{rem:w}
The conditions listed in Lemma \ref{lem:w} are quite natural.  Indeed, items \ref{w1}--\ref{w8} have the equivalent formulations:
\ben
\item \label{w1'} Some element of $U$ is a right identity for $S$.
\item \label{w2'} Every element of $S$ has a right identity from $U$.
\item \label{w3'} Some element of $U$ is a left identity for $S$.
\item \label{w4'} Every element of $S$ has a left identity from $U$.
\item \label{w5'} $S\sub US$.
\item \label{w6'} $S\cap US\not=\es$.
\item \label{w7'} Some element of $U$ is a right identity for $US$.
\item \label{w8'} Every element of $US$ has a right identity from $U$.
\een
While \ref{w9} and \ref{w10} may not appear to be quite as natural, we have shown that they are equivalent (when $(U,S)$ is strong) to \ref{w7} and \ref{w8}, and these alternative formulations will be useful in the next proof.

One important special case in which all of conditions \ref{w1}--\ref{w10} obviously hold is when $U$ is a submonoid of the over-monoid $M$; note that we only have to check \ref{w1} and \ref{w3}, which clearly hold in this case.
\end{rem}

\begin{rem}
It is also worth noting that the conditions listed in Lemma \ref{lem:w} are not all equivalent.  
For example, consider the partial transformation monoid $M=\PT_n$, where $n\geq2$, and the subsemigroups $U=\Sing(\E_n)$, $S_1=\T_n$ and $S_2=\Sing(\T_n)$.  As noted in Example \ref{eg:Tn}, the strong action pairs $(U,S_1)$ and $(U,S_2)$ give rise to the same subsemigroup:
\[
US_1=US_2=\Sing(\PT_n).
\]
One can check that:
\bit
\item $(U,S_1)$ satisfies only \ref{w8} and \ref{w10}, while
\item $(U,S_2)$ satisfies only \ref{w2}, \ref{w8} and \ref{w10}.
\eit

On the other hand, some of the conditions of Lemma \ref{lem:w} are equivalent under further assumptions.  For example, if every element of $U$ has a left identity (which happens for example if $U$ is regular), then \ref{w4}$\iff$\ref{w5}.  Indeed, suppose \ref{w5} holds, and let $s\in S$ be arbitrary.  Then~\ref{w5} gives $s=vt$ for some $v\in U$ and $t\in S$.  But then if $u\in U$ is a left identity for $v$, it follows that $s=vt=uvt=us$, showing that \ref{w4} holds.

A similar argument shows that if $U$ has a left identity, then \ref{w5}$\implies$\ref{w3}, so it follows that \ref{w3}$\iff$\ref{w4}$\iff$\ref{w5} in this case.
\end{rem}

\begin{rem}
Finally, we note that the implication \ref{w6}$\implies$\ref{w7} need not hold if $(U,S)$ is not strong.  For example, consider the symmetric inverse monoid $S=\I_X$ where $|X|\geq2$.  Then $(U,S)$ is an action pair, where $U=\Sing(\E_X)=\set{\id_A}{A\subsetneq X}$.  Indeed, this is easy to check directly, and it also follows from Remark \ref{rem:PS}, as $\I_X$ is inverse, and hence left restriction, and $U=\Pf$ in the notation of that remark.  Condition \ref{w6} holds trivially for this pair, as we can take $s=t=u=\id_A$ for any $A\subsetneq X$.  However, \ref{w7} does not hold since 
\[
US=\set{\al\in\I_X}{\dom(\al)\not=X}
\]
does not have a right identity, let alone a right identity from~$U$.
\end{rem}

For the rest of this section, we assume that for each $u\in U$, the right congruence $\th_u$ is generated (as a right congruence) by some set of pairs $\Om_u\sub S\times S$.  For any set $X$, and any set of pairs $\Si\sub X\times X$, we write
\[
\Si^{-1} = \bigset{(y,x)}{(x,y)\in\Si}.
\]

\begin{lemma}\label{lem:Om3}
If $(U,S)$ is a strong action pair, and if any of the conditions listed in Lemma~\ref{lem:w} hold (e.g., if $U$ is a submonoid), then $\th=\Om^\sharp$, where 
\[
\Om = \bigset{((u,s),(u,t))}{u\in U,\ (s,t)\in\Om_u}.
\]
\end{lemma}

\pf
By Lemma \ref{lem:w} we may assume that \ref{w10} holds.  As usual, it suffices to show that $\th\sub\th'$, where $\th'=\Om^\sharp$.  By Lemma \ref{lem:th}\ref{th2}, this amounts to showing that $(u,s) \mr\th' (u,t)$ for all $u\in U$ and $(s,t)\in\th_u$, so fix some such $u,s,t$.  Since $\th_u$ is generated (as a right congruence) by~$\Om_u$, there is a sequence
\[
s = s_1\to s_2\to\cdots\to s_k=t,
\]
where $s_1,\ldots,s_k\in S$, and such that for each $1\leq i<k$,
\[
s_i=a_ic_i \ANd s_{i+1}=b_ic_i \qquad\text{for some $(a_i,b_i)\in\Om_u\cup\Om_u^{-1}$ and $c_i\in S^1$.}
\]
Since $s=s_1$ and $t=s_k$, we can show that $(u,s)\mr\th'(u,t)$ by showing that $(u,s_i) \mr\th' (u,s_{i+1})$ for each $1\leq i<k$.  This is clear if $c_i=1$, since then
\[
((u,s_i),(u,s_{i+1})) = ((u,a_i),(u,b_i))\in\Om\cup\Om^{-1}\sub\th'.
\]
So now we assume that $c_i\not=1$.  By \ref{w10}, there exists $v\in U$ such that $u=u\cdot{}^{a_i}v$.  First note that
\[
u\cdot {}^{b_i}v\cdot b_i = u\cdot b_i\cdot v = u\cdot a_i\cdot v =u\cdot {}^{a_i}v\cdot a_i = u\cdot a_i,
\]
where we used $(a_i,b_i)\in\th_u$ (which implies $ua_i=ub_i$) in the second step; \ref{SA2} then gives $u\cdot {}^{b_i}v=u$.  It then follows that
\[
(u,s_i) = (u,a_i)\cdot(v,c_i) \mr\th' (u,b_i)\cdot(v,c_i) = (u,s_{i+1}),
\]
as required.
\epf

\begin{rem}
The above proof shows that $\th$ is in fact generated by $\Om$ as a \emph{right} congruence.
\end{rem}

The previous result showed that, under certain conditions, a generating set for the congruence~$\th$ can be constructed by incorporating generating sets for the right congruences $\th_u$.  We now aim to take this further, by additionally incorporating a generating set for the semigroup~$U$.  This of course necessitates stronger assumptions, and these include requiring that the $\th_u$ behave well with respect to joins.  Before stating the result, we first establish some notation, which will be of use throughout.  
For $u,v\in U$ we write
\begin{align}\label{eq:pre}
v\pre u \quad&\iff\quad u=wv &&\hspace{-1cm} \text{for some $w\in U^1$.}
\intertext{(So $\pre$ is Green's $\geq_{\L}$-preorder on~$U$.)  It is easy to see that}
\label{eq:uprev}
v\pre u \quad&\implies\quad \th_v\sub\th_u &&\hspace{-1cm}\text{for all $u,v\in U$.}
\intertext{It follows immediately that for any $V\sub U$,}
\label{eq:thv}
\bigvee_{v\in V, \atop v\pre u} \th_v &\sub \th_u  &&\hspace{-1cm}\text{for all $u\in U$.}
\intertext{Here and elsewhere, the join of equivalence relations is taken in the lattice of all equivalence relations (over the same set).  In particular, the join in \eqref{eq:thv} is the smallest equivalence relation on $S$ containing all of the $\th_v$ ($v\in V$, $v\pre u$); note that this join could involve infinitely many terms.  If $U$ is commutative, it also follows from~\eqref{eq:uprev} that}
\label{eq:thv1k}
\th_{u_1}\vee\cdots\vee\th_{u_k} &\sub \th_{u_1\cdots u_k} &&\hspace{-1cm}\text{for all $u_1,\ldots,u_k\in U$,}
\end{align}
as $u_i\pre u_1\cdots u_k$ for each $i$, by commutativity.

The next result concerns the cases in which we have equality in \eqref{eq:thv} or \eqref{eq:thv1k}.
The remaining results of this section assume $U$ and $S$ are submonoids of the over-monoid $M$; one could state them in greater generality, similarly to Lemma~\ref{lem:Om3}, but the combinations of conditions become unweildy, and we prefer the cleaner statements below.

\begin{lemma}\label{lem:Om4}
Suppose $(U,S)$ is a strong action pair in a monoid $M$, and suppose additionally that $U$ and $S$ are submonoids of $M$.  
\ben
\item \label{Om41} If there exists a subset $V\sub U$ such that 
\[
\th_u = \bigvee_{v\in V, \atop v\pre u} \th_v \qquad\text{for all $u\in U$,}
\]
then $\th=\Om^\sharp$, where
\[
\Om = \bigset{((v,s),(v,t))}{v\in V,\ (s,t)\in\Om_v}.
\]
\item \label{Om42} If $U$ is commutative, and if $\th_{uv}=\th_u\vee\th_v$ for all $u,v\in U$, then for any monoid generating set $V$ for $U$ we have $\th=\Om^\sharp$, where
\[
\Om = \bigset{((v,s),(v,t))}{v\in V,\ (s,t)\in\Om_v}.
\]
\een
\end{lemma}

\pf
\firstpfitem{\ref{Om41}}  Suppose $V\sub U$ satisfies the stated assumption.  As ever, it suffices to show that $\th\sub\th'$, where $\th'=\Om^\sharp$.  By Lemma \ref{lem:th}\ref{th2}, this amounts to showing that $(u,s) \mr\th' (u,t)$ for any $u\in U$ and $(s,t)\in\th_u$, so fix some such $u,s,t$.  By the assumption on $\th_u$, there is a sequence
\[
s=s_1\to s_2\to\cdots\to s_k=t
\]
such that each $(s_i,s_{i+1})\in\bigcup_{v\in V_u}(\Om_v\cup\Om_v^{-1})$, where $V_u=\set{v\in V}{v\pre u}$.  Since ${(u,s)=(u,s_1)}$ and $(u,t)=(u,s_k)$, it suffices to show that $(u,s_i) \mr\th' (u,s_{i+1})$ for each $1\leq i<k$, so fix some such $i$.  Then $(s_i,s_{i+1})\in\Om_v\cup\Om_v^{-1}$ for some $v\in V_u$, and so $((v,s_i),(v,s_{i+1}))\in\Om\cup\Om^{-1}$.  Since $v\pre u$ we have $u=wv$ for some $w\in U(=U^1)$, and then
\[
(u,s_i) = (w,1)\cdot(v,s_i) \mr\th' (w,1)\cdot(v,s_{i+1}) = (u,s_{i+1}).
\]

\pfitem{\ref{Om42}}  Suppose the assumptions of \ref{Om42} hold, and let $V$ be an arbitrary generating set for $U$.  It suffices to show that the assumption of \ref{Om41} holds with respect to this $V$.  By \eqref{eq:thv}, this amounts to showing that
\[
\bigvee_{v\in V, \atop v\pre u} \th_v \supseteq \th_u  \qquad\text{for all $u\in U$,}
\]
so fix some such $u$.  Since $U=\la V\ra$, we have $u=v_1\cdots v_k$ for some $v_1,\ldots,v_k\in V$, and so
\[
\th_u = \th_{v_1\cdots v_k} = \th_{v_1}\vee\cdots\vee\th_{v_k} \sub \bigvee_{v\in V, \atop v\pre u} \th_v,
\]
since each $v_i\pre u$ by commutativity of $U$.
\epf

The generating sets in Lemmas \ref{lem:Om3} and \ref{lem:Om4} can be further simplified in the (very) special case that $U$ and $S$ are submonoids of $M$, with $S$ a \emph{group}.  To deal with this case, we introduce another piece of notation.

For an arbitrary weak action pair $(U,S)$, we define the sets
\begin{equation}\label{eq:Ae}
\S_u = \set{s\in S}{us=u} \qquad\text{for each $u\in U$.}
\end{equation}
It is easy to see that each $\S_u$ is a (possibly empty) subsemigroup of $S$.  When $S$ is a submonoid of~$M$, each $\S_u$ is a (non-empty) submonoid of $S$; if $S$ is additionally a group, each $\S_u$ is a subgroup.  The next result concerns this group case, when $(U,S)$ is an action pair (not just a weak action pair).  For the statement, we fix a \emph{group} generating set $\Ga_u$ for~$\S_u$, for each $u\in U$; so each $\S_u$ is generated as a monoid by $\Ga_u\cup\Ga_u^{-1}$.  Note that in this case $S$ is contained in the group of units of $M$, so the action pair $(U,S)$ is automatically strong by Lemma \ref{lem:group}\ref{group2}.

\begin{lemma}\label{lem:Om5}
Suppose $(U,S)$ is a (strong) action pair in a monoid $M$, and suppose additionally that $U$ and $S$ are submonoids of $M$, with $S$ a group.  Then $\th=\Om^\sharp$, where
\[
\Om = \bigset{((u,1),(u,s))}{u\in U,\ s\in\Ga_u}.
\]
\end{lemma}

\pf
Once again, it suffices to show that $(u,s)\mr\th'(u,t)$ for all $u\in U$ and $(s,t)\in\th_u$, where $\th' = \Om^\sharp$.  First note that $(s,t)\in\th_u\implies us=ut\implies ust^{-1}=u$, so that $st^{-1}\in \S_u$.  It then follows that $st^{-1} = g_1\cdots g_k$ for some ${g_1,\ldots,g_k\in\Ga_u\cup\Ga_u^{-1}}$.  For each $1\leq i\leq k+1$ write $s_i=g_i\cdots g_kt$.  Since $s_1=s$ and $s_{k+1}=t$, it suffices to show that $(u,s_i) \mr\th' (u,s_{i+1})$ for each $1\leq i\leq k$.

So fix some such $1\leq i\leq k$.  If $g_i\in\Ga_u$, then clearly $(u,1)\mr\th'(u,g_i)$.  Otherwise, ${(u,1)\mr\th'(u,g_i^{-1})}$, so post-multiplying by $(1,g_i)$ gives $(u,g_i)\mr\th'(u,1)$.  Either way (and keeping $s_i=g_is_{i+1}$ in mind), it follows that
\[
(u,s_i) = (u,g_i)\cdot(1,s_{i+1}) \mr\th' (u,1)\cdot(1,s_{i+1}) = (u,s_{i+1}),
\]
as required.
\epf

As with Lemma \ref{lem:Om4}, the conclusion of Lemma \ref{lem:Om5} can be strengthened when the subgroups~$\S_u$ behave well with respect to joins.  (The join $\bigvee_iH_i$ of a collection of subgroups is the smallest subgroup containing the union $\bigcup_iH_i$.)
The statement again uses the relation $\pre$ on $U$ given in~\eqref{eq:pre}.

\begin{lemma}\label{lem:Om6}
Suppose $(U,S)$ is a (strong) action pair in a monoid $M$, and suppose additionally that $U$ and $S$ are submonoids of $M$, with $S$ a group.  
\ben
\item \label{Om61} If there exists a subset $V\sub U$ such that 
\[
\S_u = \bigvee_{v\in V, \atop v\pre u} \S_v \qquad\text{for all $u\in U$,}
\]
then $\th=\Om^\sharp$, where
\[
\Om = \bigset{((v,1),(v,s))}{v\in V,\ s\in\Ga_v}.
\]
\item \label{Om62} If $U$ is commutative, and if $\S_{uv}=\S_u\vee\S_v$ for all $u,v\in U$, then for any monoid generating set $V$ for $U$ we have $\th=\Om^\sharp$, where
\[
\Om = \bigset{((v,1),(v,s))}{v\in V,\ s\in\Ga_v}.
\]
\een
\end{lemma}

\pf
As in Lemma \ref{lem:Om4}, part \ref{Om62} follows from part \ref{Om61}.  

To prove \ref{Om61}, it suffices by Lemma~\ref{lem:Om5} to show that $(u,1) \mr\th' (u,s)$ for all $u\in U$ and $s\in\Ga_u$, where $\th'=\Om^\sharp$, so fix some such $u,s$.  Then $s \in \S_u = \bigvee_{v\in V_u} \S_v$, where again we write $V_u=\set{v\in V}{v\pre u}$, so $s=g_1\cdots g_k$ for some ${g_1,\ldots,g_k\in\bigcup_{v\in V_u}(\Ga_v\cup\Ga_v^{-1})}$.  (Recall that $\Ga_v$ is a group generating set for $\S_v$.)  Writing $s_i=g_i\cdots g_k$ for each $1\leq i\leq k+1$, we have $s_1=s$ and $s_{k+1}=1$, so it suffices to show that $(u,s_i)\mr\th'(u,s_{i+1})$ for each $1\leq i\leq k$.  But for any such $i$, we have $g_i\in\Ga_v\cup\Ga_v^{-1}$ for some $v\in V_u$, and $u=wv$ for some $w\in U$.  As in the proof of Lemma~\ref{lem:Om5}, we have $(v,1)\mr\th'(v,g_i)$, and (again remembering that $s_i=g_is_{i+1}$), it follows that
\[
(u,s_i) = (w,1)\cdot(v,g_i)\cdot(1,s_{i+1}) \mr\th' (w,1)\cdot(v,1)\cdot(1,s_{i+1}) = (u,s_{i+1}).  \qedhere
\]
\epf

\begin{rem}\label{rem:Om6}
As in \eqref{eq:uprev}, we have
\begin{align*}
v\pre u \implies \S_v&\sub\S_u &&\hspace{-2cm}\text{for all $u,v\in U$.}
\intertext{It follows that}
\bigvee_{v\in V, \atop v\pre u} \S_v &\sub \S_u &&\hspace{-2cm}\text{for all $u\in U$,}
\intertext{and if $U$ is commutative, that}
\S_{u_1}\vee\cdots\vee\S_{u_k} &\sub \S_{u_1\cdots u_k} &&\hspace{-2cm}\text{for all $u_1,\ldots,u_k\in U$.}
\end{align*}
Thus, in order to apply Lemma \ref{lem:Om6}, we only need to check the forwards inclusion in the relevant join equation.
\end{rem}

\chap{Proper action pairs}\label{chap:PAP}

Proposition \ref{prop:A1} showed that the semigroup $US$ arising from a weak action pair~$(U,S)$ is a homomorphic image of a semidirect product~$U\rtimes S$.  Theorem \ref{thm:FT} took this further in the case that $(U,S)$ is an action pair (not just a weak action pair), by additionally classifying the congruences on semidirect products that lead to action pairs in this way.
The purpose of the current chapter is to prove a number of additional structure theorems concerning semigroups arising from (weak) action pairs.  These results are inspired by important classical results on inverse semigroups, and wider classes such as left restriction semigroups, each of which involve the notion of a \emph{proper} semigroup from the relevant class.  

We discuss these classical results in Section \ref{sect:PAPback} (see Theorems \ref{thm:McA} and \ref{thm:OC}), as well as some of the backstory behind their development.  Section \ref{sect:AAP} shows that the mere existence of an action of $S$ on $U^1$ (for semigroups $U$ and $S$) allows one to construct an action pair $(\ul U,\ol S)$ in a suitable monoid $M$, with $\ul U\cong U$ and $\ol S\cong S$, and with the respective actions of $S$ and $\ol S$ on~$U^1$ and~$\ul U^1$ matching up appropriately.  Some aspects of this construction bear similarities with ideas in the proof of Theorem \ref{thm:FT}.  In Section \ref{sect:cover} we introduce the key notion of a \emph{proper} action pair, and show that a semigroup $US$ arising from an arbitrary weak action pair $(U,S)$ has a `proper cover'.  More specifically, Theorem \ref{thm:cover} shows that there exists a proper pair $(U',S')$ with $U'\cong U$ and~$S'\cong S$, and a natural surmorphism $U'S'\to US$.  As an application of this, we show how to deduce the classical Theorem \ref{thm:McA} as a corollary.
Sections \ref{sect:embed1}--\ref{sect:embed4} contain four embedding theorems for proper \emph{monoids}.  
\bit
\item Theorem \ref{thm:embed2} (in Section \ref{sect:embed2}) is the most general of these, and states that a monoid $US$ arising from a proper action pair $(U,S)$ can be embedded in a semidirect product $\U\rtimes(S/\si)$, where~$\U$ is a monoid containing $U$, and $\si$ is a certain natural congruence on $S$.  When $S$ is inverse,~$\si$ is the least group congruence.  The monoid $\U$ constructed in the proof is defined by means of a presentation that incorporates the structure of the monoids $U$ and $S/\si$, and the action of $S$ on~${U(=U^1)}$.  
\item Theorem~\ref{thm:embed1} (in Section~\ref{sect:embed1}) is actually a special case of Theorem \ref{thm:embed2}, stated and proved under the additional assumption that the elements $s^+={}^s1$ ($s\in S$) are central in~$U$.  This allows us to construct a very different monoid $\U$ during the proof, and the resulting semidirect product is in fact an (unrestricted) wreath product.  In the special case that $S$ is left restriction, this $\U$ turns out to be a semilattice (built in a natural way from the ideals of $U=P(S)$), and this allows us to deduce the classical Theorem \ref{thm:OC} as a corollary.  
\item Theorems \ref{thm:embed3} and \ref{thm:embed4} (in Sections \ref{sect:embed3} and \ref{sect:embed4}) are specialised versions of Theorem \ref{thm:embed2}.  They show that when $U$ is assumed to have certain additional structural properties (such as being a semilattice or a left-regular band), then $M=US$ embeds in $\U\rtimes(S/\si)$, as above, where $\U$ has the same structural property as $U$.
\eit

\sect{Background}\label{sect:PAPback}

In the early years of the development of semigroup theory, three directions of attack emerged to understand the structure and behaviour of an inverse semigroup $S$. These all had reference to the semilattice $E=E(S)=\set{x\in S}{x=x^2}$ of idempotents of $S$, and were  subsequently developed and applied to much broader classes. Specifically these directions are:
\bit
\item the use of Munn's \emph{fundamental inverse semigroup} $T_E$, built from order-isomorphisms of principal ideals of $E$ \cite{Munn1970}; 
\item the Ehresmann-Schein-Nambooripad approach, which uses the \emph{trace groupoid} of $S$ (which has identities $E$), together with the natural partial order on $S$, to build an inductive groupoid from which $S$ can be recovered; a comprehensive account of the genesis of this material is contained in Lawson's monograph \cite{Lawson1998};
\item McAlister's theory of \emph{proper} (a.k.a.~\emph{$E$-unitary}) inverse semigroups \cite{McAlister1974,McAlister1974b}. 
\eit
Our current inspiration stems from this final direction, and we now describe it in a little more detail, beginning with the key definition(s).  The following can be found, for example, on p55 of~\cite{CPbook2}.

\begin{defn}\label{defn:EU}
Let $S$ be a semigroup, and $A\sub S$ an arbitrary subset.  We say $S$ is \emph{$A$-unitary} if the following implications hold:
\begin{equation}\label{eq:EU}
ax\in A\implies x\in A \AND xa\in A\implies x\in A \qquad\text{for all $a\in A$ and $x\in S$.}
\end{equation}
\end{defn}

\begin{rem}\label{rem:EU}
Consider an \emph{inverse} semigroup $S$, with semilattice of idempotents $E=E(S)$.  Then~$S$ is $E$-unitary when the implications in \eqref{eq:EU} hold with respect to $A=E$.  It is easy to check (when $A=E$) that either of the implications in \eqref{eq:EU} follows from the other.  Thus, $S$ is $E$-unitary if and only if 
\begin{align}\label{eq:lr_proper}
\nonumber ex\in E &\implies x\in E &&\hspace{-2cm}\text{for all $e\in E$ and $x\in S$.}
\intertext{It is again easy to check that this condition is equivalent to}
e=ex &\implies x\in E &&\hspace{-2cm}\text{for all $e\in E$ and $x\in S$.}
\end{align}
\end{rem}

This all paves the way for the following:

\begin{defn}\label{defn:lr_proper}
An inverse semigroup $S$ is \emph{proper} if it is $E$-unitary, where $E=E(S)$: i.e., if~\eqref{eq:lr_proper} holds (cf.~Remark~\ref{rem:EU}).
\end{defn}

\begin{rem}
There are many other equivalent characterisations of proper ($E$-unitary) inverse semigroups; see for example \cite[Proposition 5.9.1]{Howie1995}.  One is for $\mathscr{R} \cap \sigma$ (or equivalently $\mathscr{L} \cap \sigma$) to be trivial, where~$\mathscr R$ and $\mathscr L$ are two of Green's relations, and where $\si$ is the least group congruence. Thus, one can coordinatise the elements of a proper inverse semigroup by an idempotent and a group element.  For broader classes, such as the left restriction semigroups we consider below, being $E$-unitary is a weaker condition than being proper; the definition of proper semigroups in these classes is given below.
\end{rem}

In two landmark papers, McAlister showed that:
\bit
\item any inverse semigroup has a \emph{proper cover} \cite{McAlister1974}: i.e., a preimage under an idempotent-separating morphism that is proper, and 
\item any proper inverse semigroup~$S$ is isomorphic to a so-called \emph{P-semigroup} \cite{McAlister1974b}. The latter is built from an action of the group $G=S/\sigma$ (where here $\sigma$ is the least group congruence on $S$) on a semilattice containing $E$.  Both the difficulty and the beauty of this approach is that the action of $G$ on $E$ is not in general closed.  It is also worth noting that while we know from the outset that $S$ acts on $E$ (cf.~Example \ref{eg:act}\ref{act2}), we also need $S/\si$ to act on (a semilattice containing) $E$.
\eit
Subsequently, O'Carroll \cite{OCarroll1976} showed that a proper inverse semigroup $S$ embeds into a semidirect product $F\rtimes G$, where $F$ is a semilattice containing $E$, and again $G=S/\si$. Immediately one sees a connection with (left) action pairs, since a proper inverse semigroup $S$ is coordinatised by two components, one from a semilattice and another from a group. Of course, $S$ might not be the internal product of these constituents.

There are many proofs of McAlister's `P-theorem', and of O'Carroll's subsequent embedding theorem.  In addition to the original strategies, we mention here the approaches of Munn \cite{Munn1976}, Billhardt \cite{Billhardt1997} and Steinberg \cite{Steinberg2003}, as well as the alternative strategy of Petrich and Reilly \cite{PR1979} using only partial actions (which yields a somewhat differerent formulation, but which nevertheless deserves mention here). A splendid account may be found in \cite{LM2007}, where \cite{PR1979}  is referred to as `the maverick alternative'.

The approach of McAlister and O'Carroll for  inverse semigroups prompted analagous work for  larger classes of semigroups; see for example \cite {FPW2004, Szendrei1982,Szendrei1980, Trotter1995,APW1992,MP1987,BGG2010,Lawson1986}.  The extensions are, roughly speaking, in two directions:
\bit
\item One is to keep the condition of $S$ being regular, but weaken the condition that $E(S)$ forms a semilattice: e.g., it could be a left-regular band (these satisfy the identities $x^2=x$ and~$xyx=xy$).
\item The other direction is to drop the condition that $S$ is regular, and here the main body of work  retains commutativity of (certain) idempotents. 
\eit
In particular, a theory analogous to that for  inverse semigroups has been developed for left restriction semigroups and  the special (and in some sense generating) classes of  left ample and weakly left ample semigroups; see for example \cite{Fountain1977, FG1993, GG2000, GG1999,BGG2010}.  The results of this chapter take this all much further, and apply to semigroups and monoids arising from (weak) action pairs.

Here is the appropriate extension to left restriction semigroups of the notion of a proper inverse semigroup \cite{FGG2009,GS2013,Lawson1986}.

\begin{defn}\label{defn:lr_proper2}
Let $S$ be a left restriction semigroup with semilattice of projections $P=P(S)$.  Define the congruence
\begin{align*}
\si = \si_S &= \bigset{(s,t)\in S\times S}{es=et\ (\exists e\in P)}\\
&= \bigset{(s,t)\in S\times S}{es=ft\ (\exists e,f\in P)}.
\end{align*}
We say $S$ is \emph{proper} if
\begin{equation}\label{eq:lr_proper2}
s=t \qquad\iff\qquad s^+=t^+ \ANd s\mr\si t \qquad\qquad\text{for all $s,t\in S$.}
\end{equation}
\end{defn}

\begin{rem}\label{rem:lr_proper2}
In fact, $\si=\nabla_P^\sharp$ is the least congruence on $S$ that identifies all projections.  As explained in \cite{Gould_notes}, $\si$ is also a unary semigroup congruence, meaning that ${(s,t)\in\si \implies (s^+,t^+)\in\si}$.  The latter also follows quickly from \ref{A2} and the definition of $\si$, as $(P,S)$ is an action pair; cf.~Proposition \ref{prop:LR1} and Remark \ref{rem:LR1}.

When $S$ is inverse (regarded as a left restriction semigroup with $s^+=ss^{-1}$), $\si$ is the least group congruence on~$S$.  Definitions \ref{defn:lr_proper} and \ref{defn:lr_proper2} are equivalent when $S$ is inverse; see for example~\cite{Lawson1998} or \cite[Chapter 5]{Howie1995}.  It is worth noting that $s^+=t^+$ is equivalent to $s\R t$ in an inverse semigroup, where $\R$ is one of Green's relations.

If $S$ is left ample (a left restriction semigroup that additionally satisfies the quasi-identity $xz=yz\Rightarrow xz^+=yz^+$), then $\sigma$  is the least right cancellative  monoid congruence on $S$ \cite{Fountain1977}.  If~$S$ is weakly left ample (a left restriction semigroup in which $P(S)=E(S)$), then $\sigma$ is the least unipotent monoid congruence on $S$ \cite{GG2000}.
\end{rem}

\begin{defn}
Let $S$ be a left restriction semigroup.  A \emph{cover} of~$S$ is a left restriction semigroup~$C$ together with a projection-separating surmorphism $\psi:C\to S$.  Such a $\psi$ is called a \emph{covering morphism}, and is required to be a unary semigroup morphism, in the sense that $(xy)\psi=(x\psi)(y\psi)$ and $x^+\psi=(x\psi)^+$ for all $x,y\in C$.   It follows quickly that $\psi$ restricts to an isomorphism $P(C)\to P(S)$.  It is worth noting that when $C$ and $S$ are inverse semigroups (regarded as left restriction semigroups with $x^+=xx^{-1}$), any semigroup morphism $C\to S$ automatically preserves the ${}^{-1}$ operation, and hence also the ${}^+$ operation.
\end{defn}

Here are the formal statements of the results mentioned above:

\begin{thm}[cf.~\cite{McAlister1974,FGG2009,BGG2010}]\label{thm:McA}
Every left restriction semigroup (and in particular every inverse semigroup) has a proper cover.
\end{thm}

\begin{thm}[cf.~\cite{OCarroll1976,Lawson1986,BGG2010}]\label{thm:OC}
Let $S$ be a proper left restriction semigroup (in particular an inverse semigroup) with semilattice of projections $P=P(S)$. Then $S$ embeds into a semidirect product~${Q\rtimes (S/\sigma)}$, where $Q$ is a semilattice containing $P$.
\end{thm}

As we have already stated, our main goal in the current chapter is to devise a suitable notion of a \emph{proper} (weak) action pair, and formulate and prove corresponding extensions of Theorems~\ref{thm:McA} and~\ref{thm:OC} in this context; see Theorems \ref{thm:cover}, \ref{thm:embed1}, \ref{thm:embed2}, \ref{thm:embed3} and \ref{thm:embed4}.
These more general results will also lead to new proofs of Theorems \ref{thm:McA} and \ref{thm:OC}.

\sect{From actions to action pairs}\label{sect:AAP}

(Weak) action pairs are defined in terms of an action of a semigroup $S$ on the monoid completion~$U^1$ of a semigroup $U$, where $U$ and $S$ are contained in a common over-monoid~$M$.  The purpose of the current section is to show that the mere existence of an action (by morphisms) leads to an action pair that mimicks the original action, in a technical sense made precise below.  This construction will be harnessed in the next section in order to show that any semigroup arising from a (weak) action pair has a `proper cover'.  

For convenience of reference, we will gather the standing assumptions for this section here, and also fix some notational conventions:

\begin{ass}\label{ass}
\ben
\item \label{a1} We assume that $U$ and $S$ are semigroups, and that $S$ has a left action on~$U^1$ via semigroup morphisms, denoted $(s,u)\mt {}^su$.
\een
\noindent For example, if $S$ has an action on $U$ itself, and if $U\not=U^1$, then this can be extended to an action on $U^1$ by additionally defining ${}^s1=1$ for all $s\in S$.  In this case, the action of $S$ on $U^1$ is by monoid morphisms.  More generally, however, we do not assume that the action in \ref{a1} is by monoid morphisms.
\ben\addtocounter{enumi}{1}
\item \label{a2} If $S$ happens to be a monoid, then we assume in addition that the action of $S$ on $U^1$ is monoidal: i.e., that the identity of $S$ acts identically.  
\een
Here, as usual, $U^1$ denotes $U$ if $U$ is a monoid or else $U^1=U\cup\{1\}$, where~$1$ is an adjoined identity.  
We use the same convention for $S^1$, and we assume without loss of generality (and for notational convenience) that the identities of $U^1$ and $S^1$ are the same.  But note that we do not assume that~$U$ and $S$ are contained in a common semigroup; in particular, we are not assuming at this point that $(U,S)$ is a (weak) action pair.  However, we will show that $(U,S)$ can indeed be identified with an action pair $(\ul U,\ol S)$ in a suitable monoid $M$.
\end{ass}

The requirement for $S$ to act monoidally (if it is a monoid) in Assumption \ref{ass}\ref{a2} is necessary for the results of this section.  In the next, we will apply these results in the case that $(U,S)$ is a weak action pair, and we noted in Remark \ref{rem:AP2} that whenever $S$ is a monoid, it acts monoidally on $U^1$ in this case.

The action of $S$ on $U^1$ allows for the formation of the (external) semidirect product $U\rtimes S$, as in Definition \ref{defn:SD}.  As usual, and if necessary, we can extend the action of $S$ on $U^1$ to a monoidal action of~$S^1$ on $U^1$ (it is already monoidal by assumption if $S=S^1$), and we then have the semidirect product $U^1\rtimes S^1$, which contains $U\rtimes S$ as a subsemigroup.  Although $U^1\rtimes S^1$ might not be a monoid (cf.~Lemma \ref{lem:USmon}), it follows from Proposition \ref{prop:MM} (and monoidality of the action of $S^1$ on $U^1$) that
\[
M = \bigset{(u,s)\in U^1\rtimes S^1}{u=us^+} \leq U^1\rtimes S^1
\]
is a monoid with identity $(1,1)$.  
Here, as usual, we write $s^+={}^s1$ for all $s\in S^1$, noting that~$1^+=1$.  

We begin by identifying two natural subsemigroups of $M$.  To this end, we define
\[
\ul u = (u,1) \AND \ol s = (s^+,s) \qquad\text{for $u\in U$ and $s\in S$,}
\]
and we set
\[
\ul U = \set{\ul u}{u\in U} \AND \ol S = \set{\ol s}{s\in S}.
\]
It is clear that $\ul U\sub M$, while $\ol S\sub M$ follows from Lemma \ref{lem:+}\ref{lem+1}.

\newpage

\begin{lemma}\label{lem:olul}
Given Assumption \ref{ass}, and with the above notation, we have 
\ben\bmc2
\item $\ul U\leq M$ and $\ul U\cong U$,
\item $\ol S\leq M$ and $\ol S\cong S$.  
\emc\een
\end{lemma}

\pf
The first part is clear.  The second follows quickly from the fact that for all $s,t\in S$,
\begin{equation}\label{eq:st^+}
\ol s\cdot\ol t = (s^+,s)\cdot(t^+,t) = (s^+\cdot{}^s(t^+),st) = ({}^s(t^+),st) = ((st)^+,st) = \ol{st},
\end{equation}
where we used parts \ref{lem+4} and \ref{lem+2} of Lemma \ref{lem:+}.
\epf

For convenience in what follows, we also write $\ul1=(1,1)$, even if $1\not\in U$.  On a number of occasions during the next proof (and later in the paper) we make use of the fact that ${\ul u\cdot\ol s = (us^+,s)}$ for all $u\in U$ and $s\in S$; in particular, $(u,s)=\ul u\cdot\ol s$ if $u=us^+$.

\begin{prop}\label{prop:olul}
Given Assumption \ref{ass}, and with the above notation, 
\ben
\item \label{olul1} $(\ul U,\ol S)$ is an action pair in $M = \bigset{(u,s)\in U^1\rtimes S^1}{u=us^+}$, 
\item \label{olul2} $(\ul U,\ol S)$ is strong if and only if the action of $S$ on $U^1$ is by monoid morphisms,
\item \label{olul3} $\ul U\cdot\ol S = \bigset{(u,s)\in U\rtimes S}{u=us^+}$.
\een
\end{prop}

\pf
\firstpfitem{\ref{olul1}}
By Lemma \ref{lem:olul}, $\ul U$ and $\ol S$ are subsemigroups of $M$.  Since $\ul U\cong U$ and $\ol S\cong S$, we have an action of $\ol S$ on $\ul U^1=\ul U\cup\big\{(1,1)\big\}$, given by
\[
{}^{\ol s}\ul u = \ul{{}^su} = ({}^su,1) \qquad\text{for $u\in U^1$ and $s\in S$.}
\]
Using Lemma \ref{lem:+}\ref{lem+4}, one can easily check that $\ol s\cdot\ul u = ({}^su,s) = {}^{\ol s}\ul u\cdot\ol s$ for all $u\in U$ and $s\in S$.  This shows that \ref{A1} holds.

For \ref{A2}, we must show that
\[
\ul u\cdot \ol s=\ul v\cdot \ol t \implies \ul u\cdot \ol s^+=\ul v\cdot \ol t^+ \qquad\text{for all $u,v\in U^1$ and $s,t\in S$.}
\]
Noting that
\begin{equation}\label{eq:ols+}
\ol s^+ = {}^{\ol s}\ul1 = \ul{{}^s1} = \ul{s^+} = (s^+,1),
\end{equation}
and similarly for $\ol t^+$, the required implication follows quickly from the fact that
\[
\ul u\cdot \ol s = (us^+,s) \COMMA \ul v\cdot \ol t = (vt^+,t) \COMMA \ul u\cdot \ol s^+ = (us^+,1) \AND \ul v\cdot \ol t^+= (vt^+,1).
\]

\pfitem{\ref{olul2}}
By Lemma \ref{lem:SAPAP}\ref{SAPAP2}, $(\ul U,\ol S)$ is strong if and only if $\ol s^+=\ul1$ for all $s\in S$.  But by \eqref{eq:ols+},
\[
\ol s^+=\ul1 \iff (s^+,1)=(1,1) \iff s^+=1 \iff {}^s1=1.
\]
The assertion quickly follows.

\pfitem{\ref{olul3}}
First suppose $u\in U$ and $s\in S$ are such that $u=us^+$.  Then
\[
(u,s)=(us^+,s)=\ul u\cdot\ol s\in\ul U\cdot\ol S.
\]
Conversely, suppose $\ba\in\ul U\cdot\ol S$, so that $\ba=\ul v\cdot\ol s$ for some $v\in U$ and $s\in S$.  Then with $u=vs^+\in U$, we have
\[
\ba = \ul v\cdot\ol s = (vs^+,s) = (u,s) \AND us^+ = vs^+s^+ = vs^+ = u,
\]
using Lemma \ref{lem:+}\ref{lem+1}.
\epf

\sect{Proper action pairs and a covering theorem}\label{sect:cover}

We are now almost ready to introduce the notion of a \emph{proper} action pair; see Definitions \ref{defn:si} and~\ref{defn:proper}.  The main result of the current section is Theorem \ref{thm:cover}, which shows that any semigroup arising from a (weak) action pair is a natural homomorphic image of a semigroup arising from a proper pair.  The proof of the theorem utilises the $(\ul U,\ol S)$ construction from the previous section.  Throughout the current section we comment on the special case of left restriction semigroups, and we end by showing how the classical Theorem \ref{thm:McA} follows from our Theorem \ref{thm:cover}.

Consider a weak action pair $(U,S)$ in a monoid $M$.  As ever, we denote the action of $S$ on~$U^1$ by $(s,u)\mt{}^su$.  We define
\[
S^+ = \set{s^+}{s\in S} = \set{{}^s1}{s\in S},
\]
noting that this subset of $U^1$ need not be a subsemigroup.  In any case, we also define
\[
P = P(U,S) = \la S^+\ra,
\]
which is the subsemigroup of $U^1$ generated by $S^+$.  Since $S^+$ consists of idempotents (cf.~Lemma~\ref{lem:+}\ref{lem+1}), it follows that $P$ is an idempotent-generated semigroup, and hence is a submonoid of $U^1$ if and only if $1\in S^+$; indeed, the latter claim follows from \cite[Lemma 2.1]{JE2020a}.  

\begin{lemma}\label{lem:P}
The subsemigroup $P$ of $U^1$ is closed under the action of $S$.
\end{lemma}

\pf
By Lemma \ref{lem:+}\ref{lem+2}, $S^+$ is closed under the action of $S$; so too therefore is $\la S^+\ra=P$.
\epf

\begin{rem}
By Lemma \ref{lem:P} (and consulting Definition \ref{defn:AP}), it follows immediately that $(P,S)$ is a weak action pair as well.  If $(U,S)$ is a (strong) action pair, then so too is $(P,S)$, by Lemma~\ref{lem:US0}\ref{US02}; when $(U,S)$ is strong, Lemma \ref{lem:SAPAP}\ref{SAPAP2} gives $P=\{1\}$.
\end{rem}

\begin{rem}\label{rem:P}
When $S$ is a left restriction semigroup, $S^+=\la S^+\ra$ is precisely the semilattice of projections $P(S)$.  This explains our choice of notation for $P=P(U,S)$.  Recall from Proposition~\ref{prop:LR1} (and Remark \ref{rem:LR1}) that $(P,S)$ is an action pair in $S$.
\end{rem}

In order to give the definition of a proper action pair, we need an appropriate congruence to play the role of $\si$ from Definition \ref{defn:lr_proper2}.

\begin{defn}\label{defn:si}
Let $(U,S)$ be a weak action pair in a monoid $M$, and write $P=\la S^+\ra$.  Define the relation
\begin{align*}
\ka = \ka(U,S) &= \bigset{(s,t)\in S\times S}{ps=qt\ (\exists p,q\in P)}\\
&= \bigset{(s,t)\in S\times S}{ps=qt\ (\exists p,q\in P^1)},
\end{align*}
and let $\si = \si(U,S) = \ka^\sharp$ be the congruence on $S$ generated by $\ka$.
\end{defn}

\begin{lemma}\label{lem:si}
The relation $\ka$ is reflexive, symmetric and compatible.  Consequently, $\si$ is the transitive closure of $\ka$.
\end{lemma}

\pf
We just prove the first statement, as the second immediately follows.  Reflexivity, symmetry and right-compatibility of $\si$ are clear.  For left-compatibility, let $(s,t)\in\ka$ and $x\in S$.  So $ps=qt$ for some $p,q\in P$.  
By \ref{A1} we have ${}^xp\cdot xs = x\cdot ps = x\cdot qt = {}^xq\cdot xt$.  Since ${}^xp,{}^xq\in P$, by Lemma~\ref{lem:P}, it follows that $(xs,xt)\in\ka$.
\epf

The next result concerns a special case in which we have $\si=\ka$.  For the statement, recall~\cite{Dubreil1941} that a semigroup $T$ is \emph{right-reversible} if any two left ideals (equivalently, any two principal left ideals) of $T$ have non-empty intersection: i.e., if for all $x,y\in T$, we have $ax=by$ for some $a,b\in T$.  Examples of right-reversible semigroups include semilattices (or in fact arbitrary commutative semigroups), or more generally left-regular bands; the latter satisfy the identity $xyx=xy$ (and $x=x^2$).  

\begin{lemma}\label{lem:si2}
If $P$ is right-reversible, then $\si=\ka$.
\end{lemma}

\pf
It suffices by Lemma \ref{lem:si} to show that $\ka$ is transitive (in this case).  So suppose $s \mr\ka t \mr\ka x$ for some $s,t,x\in S$.  Then $ps=qt$ and $p't=q'x$ for some $p,p',q,q'\in P$.  By right-reversibility we have $aq=bp'$ for some $a,b\in P$, and then $ap\cdot s=aqt=bp't=bq'\cdot x$, with $ap,bq'\in P$, so that $s\mr\ka x$, as required.
\epf

\begin{rem}\label{rem:si}
As a special case, when $S$ is a left restriction semigroup, $P=\la S^+\ra=P(S)$ is a semilattice, so $\si=\ka$ coincides with the relation $\si_S$ from Definition \ref{defn:lr_proper2}.
\end{rem}

We can now give the definition of proper pairs and semigroups.  

\begin{defn}\label{defn:proper}
We say a weak action pair $(U,S)$ is \emph{proper} if for all $u,v\in U^1$ and $s,t\in S$, 
\begin{equation}\label{eq:proper}
us=vt \qquad\iff\qquad us^+=vt^+ \ANd s\mr\si t.
\end{equation}
If $(U,S)$ is proper, the semigroup $US$ is then said to be \emph{$(U,S)$-proper}.
\end{defn}

The following simple fact will be used without explicit reference from now on:

\begin{lemma}
Any proper weak action pair is an action pair.
\end{lemma}

\pf
We need to show that \ref{A2} holds: i.e., that $us=vt \implies us^+=vt^+$ for all $u,v\in U^1$ and $s,t\in S$.  But this clearly follows from \eqref{eq:proper}.
\epf

\begin{rem}
Consider a \emph{strong} action pair $(U,S)$.  Since $s^+=1$ for all $s\in S$, by Lemma~\ref{lem:SAPAP}\ref{SAPAP2}, we have $P=S^+=\{1\}$, and it follows immediately that $\ka=\De_S$ (the equality relation on $S$), and hence that $\si=\De_S$.  Definition \ref{defn:proper} then says that $(U,S)$ is proper precisely when
\[
us=vt \qquad\iff\qquad u=v \ANd s= t \qquad\text{for all $u,v\in U^1$ and $s,t\in S$.}
\]
This then implies that the surmorphism $\pi:U\rtimes S\to US:(u,s)\mt us$ from Proposition~\ref{prop:A1} is injective, and hence an isomorphism.  

Figure \ref{fig:Venn} gives a Venn diagram displaying the various classes of pairs studied in the paper.
\end{rem}

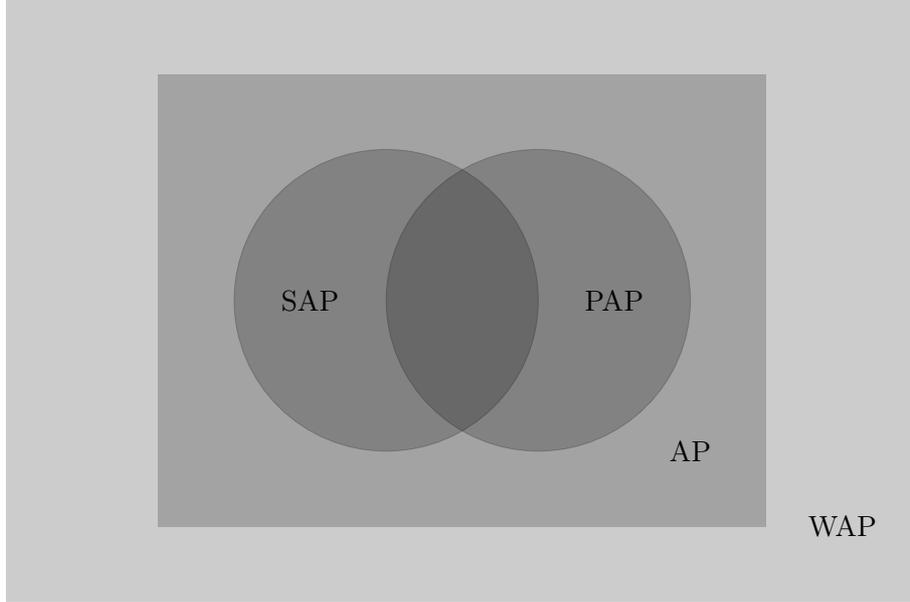
\begin{figure}[ht]
\begin{center}
\scalebox{1}{
\begin{tikzpicture}
  \tikzset{venn circle/.style={draw,circle,minimum width=4cm,fill=#1}}
  \node [venn circle = black,opacity=0.2] (C) at (-1,0) {$$};
  \node [venn circle = black,opacity=0.2] (A) at (1,0) {$$};
  \node at (-2,0) {SAP};
  \node at (2,0) {PAP};
  \node at (3,-2) {AP};
  \node at (5,-3) {WAP};
  \fill[opacity=0.2] (-4,-3)--(4,-3)--(4,3)--(-4,3);
  \fill[opacity=0.2] (-6,-4)--(6,-4)--(6,4)--(-6,4);
\end{tikzpicture} 
}
    \caption{Venn diagram indicating inclusions among the classes of strong action pairs (SAP), proper action pairs (PAP), action pairs (AP) and weak action pairs (WAP).}
    \label{fig:Venn}
   \end{center}
 \end{figure}

Now that we have defined proper action pairs, we begin by showing that these are indeed the `correct' generalisation of proper left restriction semigroups.

\begin{prop}\label{prop:lr}
If $S$ is a left restriction semigroup with semilattice of projections $P=P(S)$, then the following are equivalent:
\ben
\item \label{lr1} $S$ is proper, as in Definition \ref{defn:lr_proper2},
\item \label{lr2} $(P,S)$ is a proper action pair, as in Definition~\ref{defn:proper}.
\een
\end{prop}

\pf
We have already noted in Remark \ref{rem:P} that $P=S^+=\la S^+\ra$ in this case, and in Remark~\ref{rem:si} that the definitions of $\si$ from Definitions \ref{defn:lr_proper2} and \ref{defn:si} coincide.

\pfitem{\ref{lr1}$\implies$\ref{lr2}}  Suppose first that $S$ is proper, and let $u,v\in P^1$ and $s,t\in S$.  We must show that 
\begin{align}
\label{eq:p1} us=vt \qquad&\iff\qquad us^+=vt^+ \ANd s\mr\si t .
\intertext{Since $S$ is proper, and since $us,vt\in S$, it follows from \eqref{eq:lr_proper2} that}
\label{eq:p2} us=vt \qquad&\iff\qquad (us)^+=(vt)^+ \ANd us\mr\si vt .
\end{align}
By Lemma \ref{lem:us+}, we have ${us^+=vt^+ \iff (us)^+=(vt)^+}$.  It also follows quickly from the definition of $\si$ (in Definition \ref{defn:lr_proper2}) that $s\mr\si t \iff us\mr\si vt$, keeping in mind the fact that $P$ is a semilattice.  Equivalence of \eqref{eq:p1} and \eqref{eq:p2} is now clear.

\pfitem{\ref{lr2}$\implies$\ref{lr1}}  Conversely, suppose $(P,S)$ is proper, and let $s,t\in S$.  We must show that 
\[
s=t \qquad\iff\qquad s^+=t^+ \ANd s\mr\si t .
\]
As the forwards implication is clear, suppose $s^+=t^+$ and $s\mr\si t$.  Then with $u=v=s^+=t^+\in P$, we have $us^+=s^+=t^+=vt^+$.  Since $(P,S)$ is proper, and since also $s\mr\si t$, it follows from~\eqref{eq:proper} that $us=vt$: i.e.,~$s=t$.
\epf

The next result shows that when the action pair $(U,S)$ is proper, and when $P=\la S^+\ra$ satisfies a certain kind of semigroup identity, the congruence $\si$ can be described equationally, rather than by asserting the existence of suitable elements of $P$.

\begin{lemma}\label{lem:si3}
Let $(U,S)$ be a proper action pair, and suppose $P$ satisfies a semigroup identity of the form $f(x,y)x=g(x,y)y$, where $f(x,y),g(x,y)\in\{x,y\}^*$.  Then
\[
\si = \ka = \bigset{(s,t)\in S\times S}{f(s^+,t^+)s=g(s^+,t^+)t}.
\]
\end{lemma}

\pf
The identity implies that $P$ is right-reversible, so Lemma \ref{lem:si2} gives $\si=\ka$.  For the rest of the proof we write
\[
\si'=\bigset{(s,t)\in S\times S}{f(s^+,t^+)s=g(s^+,t^+)t}.
\]
Since $f(s^+,t^+),g(s^+,t^+)\in P^1$ for any $s,t\in S$, we clearly have $\si'\sub\si$.  

Conversely, suppose $s\mr\si t$.  Since $s^+,t^+\in P$, we have $f(s^+,t^+)\cdot s^+=g(s^+,t^+)\cdot t^+$.  Combining this with $s\mr\si t$, it follows from \eqref{eq:proper} that $f(s^+,t^+)\cdot s=g(s^+,t^+)\cdot t$, whence $s\mr\si't$.
\epf

\begin{rem}\label{rem:si3}
Since semilattices satisfy the identity $xy=yx$, it follows from Lemma \ref{lem:si3} that
\begin{align*}
\si = \ka &= \bigset{(s,t)\in S\times S}{t^+s = s^+t}
\intertext{when $P$ is a semilattice.  Similarly,}
\si = \ka &= \bigset{(s,t)\in S\times S}{s^+t^+s = s^+t}
\end{align*}
if $P$ is a left-regular band.

If $P$ satisfies an identity of the form $f(x,y,z_1,\ldots,z_k)x=g(x,y,z_1,\ldots,z_k)y$, then one can arbitrarily substitute $x$ and/or $y$ in place of the $z_i$ to obtain an identity of the form ${f'(x,y)x=g'(x,y)y}$, and Lemma \ref{lem:si3} then applies.
\end{rem}

In Sections \ref{sect:embed1}--\ref{sect:embed4} we will be concerned with proper action pairs $(U,S)$, when $U$ and $S$ are both submonoids of the over-monoid $M$.  The $\si$-class of the identity of $S$ will play an important role on many occasions during these sections.  The next lemma characterises the elements of this class, and shows that they act on $U^1$ in a particularly simple way.  In the statement, we do not assume that~$U$ is a submonoid of $M$, and we use the notation $\wh1$ for consistency with later use.

\begin{lemma}\label{lem:wh1}
Let $(U,S)$ be a proper action pair, with $S$ a submonoid of the over-monoid~$M$, and let $\wh1$ be the $\si$-class of $1\in S$.
\ben
\item \label{wh11}  We have $\wh1 = S\cap S^+ = \set{s\in S}{s=s^+}$.
\item \label{wh12}  For any $s\in \wh1$, and for any $u\in U^1$, we have ${}^su = su = s^+u$.
\een
\end{lemma}

\pf
\firstpfitem{\ref{wh11}}  We prove this part by showing that for all $s\in S$, the following are equivalent:
\bena
\bmc3
\item \label{wh11a} $s\mr\si1$,
\item \label{wh11c} $s\in S^+$,
\item \label{wh11b} $s=s^+$.
\emc
\een
\firstpfitem{\ref{wh11a}$\implies$\ref{wh11b}}  Taking $u=t=1$ and $v=s^+$, we have $us^+=vt^+$ and $s\mr\si t$.  Since $(U,S)$ is proper, it follows from \eqref{eq:proper} that $us=vt$: i.e., $s=s^+$.

\pfitem{\ref{wh11b}$\implies$\ref{wh11c}}  This is obvious.

\pfitem{\ref{wh11c}$\implies$\ref{wh11a}}  Suppose $s=t^+$ for some $t\in S$.  Then since $1\cdot s=t^+\cdot1$, with $1,t^+\in S^+\sub P$, it follows that $(s,1)\in\ka\sub\si$; cf.~Definition \ref{defn:si}.  

\pfitem{\ref{wh12}}  We have $s=s^+$ by part \ref{wh11}, so of course $su=s^+u$.  Combining $s=s^+$ with~\ref{A1} and Lemma \ref{lem:+}\ref{lem+4}, we have $su = {}^su\cdot s = {}^su\cdot s^+ = {}^su$.
\epf

The next lemma provides a handy criterion for an action pair to be proper.  For the statement, we say that a subsemigroup $V$ of a semigroup $T$ is \emph{left-dense} in $T$ if the following condition holds:
\begin{align}
\label{eq:WLA} (\forall t\in T)\ (\exists a\in V)& \quad at\in V.
\intertext{Note, for example, that any right ideal is left-dense, but the converse need not hold.  (The term `left-dense' is taken from \cite[p.~98]{BP1985}, where the meaning is slightly different.  There `$(\exists a\in V)$' is replaced by `$(\exists a\in T)$', and this results in a weaker condition.  For example, take~$T$ to be an arbitrary group, $V$ a proper subgroup, and $t\in T\sm V$; any $a\in T$ such that $at\in V$ necessarily belongs to $T\sm V$.)  Before stating the lemma, we claim that left-density is equivalent to the ostensibly stronger condition:}
\label{eq:WLA'} (\forall s,t\in T^1)\ (\exists a\in V)& \quad as,at\in V.
\end{align}
Obviously \eqref{eq:WLA'} implies \eqref{eq:WLA}.  For the converse, suppose \eqref{eq:WLA} holds, and let $s,t\in T^1$.  There exists $c\in V$ such that $ct\in V$ (this is obvious if $t=1$, or follows from \eqref{eq:WLA} otherwise).  Since certainly $cs\in T$, there also exists $b\in V$ such that $b\cdot cs\in V$.  We then have $as,at\in V$ for $a=bc\in V$.
(The weaker version of left-density from \cite{BP1985} does not imply the corresponding version of \eqref{eq:WLA'}, as seen again with $T$ a group, $V$ a proper subgroup, $s\in V$ and $t\in T\sm V$.)

\begin{lemma}\label{lem:De}
Suppose $(U,S)$ is an action pair.
\ben
\item \label{De1} If $P$ is left-dense in $U$, then 
\[
us=vt \qquad\Rightarrow\qquad us^+=vt^+ \ANd s\mr\si t \qquad\qquad\text{for all $u,v\in U^1$ and $s,t\in S$.}
\]
\item \label{De2} If $\si=\De_S$, then
\[
us=vt \qquad\Leftarrow\qquad us^+=vt^+ \ANd s\mr\si t \qquad\qquad\text{for all $u,v\in U^1$ and $s,t\in S$.}
\]
\item \label{De3} If $P$ is left-dense in $U$, and if $\si=\De_S$, then $(U,S)$ is proper.
\een
\end{lemma}

\pf
It clearly suffices to prove \ref{De1} and \ref{De2}.  For these, let $u,v\in U^1$ and $s,t\in S$.  

\pfitem{\ref{De1}}  Suppose $u s=v t$.  It follows from \ref{A2} that $u s^+=v t^+$.  By \eqref{eq:WLA'} there exists $p\in P$ such that $pu,pv\in P$, and it then follows from $pu s=pv t$ that $(s,t)\in\ka\sub\si$.

\pfitem{\ref{De2}}  Suppose $u s^+=v t^+$ and $ s\mr\si  t$.  Since $\si = \De_{ S}$, we have $ s= t$.  Combining this with \ref{+1}, it follows that $u s = u s^+\cdot  s = v t^+\cdot  t = v t$.  \qedhere
\epf

Consider again semigroups $S$ and $U$ for which $S$ has a left action on $U^1$ by semigroup morphisms, and assume that $S$ acts monoidally if it happens to be monoid; cf.~Assumption \ref{ass}.  As in Section \ref{sect:AAP}, we have the semidirect products $U\rtimes S$ and $U^1\rtimes S^1$, and the subsemigroups
\[
\ul U = \set{\ul u}{u\in U} \AND \ol S = \set{\ol s}{s\in S}
\]
of $U^1\rtimes S^1$, where $\ul u=(u,1)$ and $\ol s=(s^+,s)$, for $u\in U$ and $s\in S$.  We showed in Proposition~\ref{prop:olul} that $(\ul U,\ol S)$ is an action pair in the monoid
\[
M=\bigset{(u,s)\in U^1\rtimes S^1}{u=us^+},
\]
and the next result takes this further.  For the proof, we use the fact that $\ol s^+=\ul{s^+}$ for $s\in S$; cf.~\eqref{eq:ols+}.  Together with the fact that $S\to\ol S:s\mt\ol s$ is an isomorphism, it follows that
\[
\ol s_1^+\cdots \ol s_k^+=\ul{s_1^+\cdots s_k^+} \qquad\text{for all $s_1,\ldots,s_k\in S$.}
\]

\begin{prop}\label{prop:proper}
Given Assumption \ref{ass}, and with the above notation, we have $\si(\ul U,\ol S) = \De_{\ol S}$.  Moreover, the action pair~$(\ul U,\ol S)$ is proper, so the semigroup
\[
\ul U\cdot\ol S = \bigset{(u,s)\in U\rtimes S}{u=us^+}
\]
is $(\ul U,\ol S)$-proper.
\end{prop}

\pf
Throughout the proof, we write $\ka=\ka(\ul U,\ol S)$ and $\si=\si(\ul U,\ol S)$, and also
\[
P = P(\ul U,\ol S) = \la \ol S^+\ra \WHERE \ol S^+ = \set{\ol s^+}{s\in S} = \set{\ul{s^+}}{s\in S}.
\]
To show that $\si=\De_{\ol S}$ we must show that $\ka=\De_{\ol S}$, so suppose $\ol s \mr\ka \ol t$ for some $s,t\in S$.  By definition, this means that $\bp\cdot \ol s= \bq\cdot \ol t$ for some $\bp,\bq\in P$.  Now, $\bp=\ol s_1^+\cdots \ol s_k^+=\ul{s_1^+\cdots s_k^+}$ for some $s_1,\ldots,s_k\in S$, so that $\bp=\ul u$, where $u=s_1^+\cdots s_k^+\in U^1$.  Similarly, $\bq=\ul v$ for some $v\in U^1$.  But then
\[
(u s^+,s) = \bp\cdot\ol s = \bq\cdot\ol t = (v t^+,t).
\]
It follows from this that $s=t$, whence $\ol s=\ol t$, and this completes the proof that $\ka=\De_{\ol S}$.

Given Proposition \ref{prop:olul}, it remains to show that $(\ul U,\ol S)$ is proper: i.e., that
\[
\ul u\cdot\ol s=\ul v\cdot\ol t \qquad\iff\qquad \ul u\cdot\ol s^+=\ul v\cdot\ol t^+ \ANd \ol s\mr\si \ol t \qquad\text{for all $u,v\in U^1$ and $s,t\in S$.}
\]
The backwards implication follows from Lemma \ref{lem:De}\ref{De2}, as $\si(\ul U,\ol S) = \De_{\ol S}$.  For the forwards implication, suppose $\ul u\cdot\ol s=\ul v\cdot\ol t$.  Since $(\ul U,\ol S)$ is an action pair, it follows from~\ref{A2} that $\ul u\cdot\ol s^+=\ul v\cdot\ol t^+$.  Also, $(us^+,s) = \ul u\cdot\ol s=\ul v\cdot\ol t = (vt^+,t)$ gives $s=t$, so certainly~$\ol s\mr\si \ol t$.
\epf

We can now give the following covering theorem, which is the main result of this section.  In essence, it says that any semigroup $US$ arising from a weak action pair $(U,S)$ is a homomorphic image of a semigroup $U'S'$ arising from a \emph{proper} action pair $(U',S')$, with $U'\cong U$ and $S'\cong S$.

\begin{thm}\label{thm:cover}
Suppose $(U,S)$ is a weak action pair in a monoid $M$.  
\ben
\item \label{cover1} There exists a proper action pair $(U',S')$ in a monoid $M'$, with $U'\cong U$ and $S'\cong S$, and a surmorphism $\psi:U'S'\to US$.  
\item \label{cover2} If $S$ is a submonoid of $M$, then such a pair and surmorphism exist for which:
\bit
\item $S'$ is a submonoid of $M'$, 
\item $U'$ is a subsemigroup of $U'S'$, and 
\item the restriction $\psi\restr_{U'}$ is an isomorphism $U'\to U$.
\eit
\item \label{cover3} If $s^+\in U$ for all $s\in S$ (including the case that $U$ is a submonoid of $M$), then such a pair and surmorphism exist for which:
\bit
\item $S'$ is a subsemigroup of $U'S'$, and 
\item the restriction $\psi\restr_{S'}$ is an isomorphism $S'\to S$.
\eit
\een
\end{thm}

\pf
\firstpfitem{\ref{cover1}}  Since $(U,S)$ is a weak action pair, $S$ (and hence also $S^1$) acts on $U^1$, and $S^1$ acts monoidally (cf.~Remark \ref{rem:AP2}).  We then take
\[
M'=\bigset{(u,s)\in U^1\rtimes S^1}{u=us^+} \COMMa U' = \ul U =\set{\ul u}{u\in U} \ANd S' = \ol S = \set{\ol s}{s\in S},
\]
as above.  By Lemma \ref{lem:olul}, we have $U'\cong U$ and $S'\cong S$.  By Proposition \ref{prop:proper}, $(U',S')$ is a proper action pair in $M'$, and 
\[
U'S' = \bigset{(u,s)\in U\rtimes S}{u=us^+}.
\]

Define the map $\psi:U'S'\to US:(u,s)\mt  us$.  Since this is the restriction to $U'S'$ of the morphism $\pi:U\rtimes S\to US$ from Proposition \ref{prop:A1}, it follows that $\psi$ is a morphism.
To demonstrate surjectivity, let $x\in US$, so that $x=us$ for some $u\in U$ and $s\in S$.  As in Remark~\ref{rem:AP2}, we have $s=s^+s$, and so $(\ul u\cdot\ol s)\psi = (us^+,s)\psi = us^+\cdot s=us = x$, with $\ul u\cdot\ol s\in U'S'$.

\pfitem{\ref{cover2}}  Suppose $S$ is a submonoid of $M$.  Clearly then $S'=\ol S$ is a submonoid of~$M'$, so that ${U'=U'\cdot\ol1\leq U'S'}$.  Since $\psi\restr_{U'}$ maps each $\ul u=(u,1)\mt u\cdot1=u$, it is clearly an isomorphism~$U'\to U$.

\pfitem{\ref{cover3}}  The stated assumption ensures that $U'S'$ contains each $\ol s = (s^+,s)$, with $s\in S$: i.e., that $S'\leq U'S'$.  Since $\psi$ maps each $\ol s=(s^+,s)\mt s^+s=s$, the final assertion follows.
\epf

\begin{rem}
We have already noted that a semigroup $US$ arising from a (weak) action pair $(U,S)$ need not contain $U$ itself, or even an isomorphic copy of $U$; see Example \ref{eg:TX}\ref{it:TX3}.  Moreover,~$U'$ is not contained in $U'S'$ (in the above construction) if $S$ is not a submonoid of $M$.  Thus,~$S$ being a submonoid is a necessary assumption in part \ref{cover2} of Theorem \ref{thm:cover}, at least for the particular choice of $U'=\ul U$ and $S'=\ol S$ used in the proof.

In part \ref{cover3}, note that the stronger assumption of $U$ being a submonoid of $M$ implies that~$U'$ is a submonoid of $M'$.
\end{rem}

As an application of Theorem \ref{thm:cover}, we show how to deduce Theorem \ref{thm:McA}.  We begin with the case of left restriction \emph{monoids}:

\begin{cor}[cf.~{\cite[Theorem 6.4]{BGG2010}}]\label{cor:cover}
Every left restriction monoid has a proper cover.
\end{cor}

\pf
Consider a left restriction monoid $M$, and the action pair $(P,M)$ in $M$, where $P$ is the semilattice of projections (cf.~Proposition \ref{prop:LR1}).  Let $(P',M')=(\ul P,\ol M)$ be the proper action pair constructed in the proof of Theorem \ref{thm:cover}, and let
\[
C = P'M' = \bigset{(u,s)\in P\rtimes M}{u=us^+}.
\]
So $\psi:C\to PM=M:(u,s)\mt us$ is a (monoid) surmorphism, and $\psi\restr_{P'}$ is an isomorphism~${P'\to P}$.

It is routine to check that $C$ is left restriction under the operation $(u,s)^+=(u,1)$.  Thus, $P(C) = \bigset{(u,1)}{u\in P} = \set{\ul u}{u\in P} = P'$.  Consequently, $\psi$ is projection-separating, since $\psi\restr_{P'}$ is an isomorphism.  

It remains to check that $\psi$ respects the unary operations of $S$ and $C$.  But for any $(u,s)\in C$ we use Lemma \ref{lem:us+} to calculate
\[
(u,s)^+\psi = (u,1)\psi = u = us^+ = (us)^+ = ((u,s)\psi)^+,
\]
as required.  
\epf

We now show how the general case follows.

\pf[\bf Proof of Theorem \ref{thm:McA}]
Let $S$ be a left restriction semigroup.  If $S$ is a monoid, then Corollary~\ref{cor:cover} applies, so we assume otherwise.  Let $M=S^1=S\sqcup\{1\}$, so that $M$ is also left restriction (where we additionally define $1^+=1$).  Let $C=P'M'$ be the proper cover of $M$ from the proof of Corollary \ref{cor:cover}, with covering morphism $\psi:C\to M:(u,s)\mt us$.  Since $M\sm\{1\}=S\leq M$, is follows that
\begin{equation}\label{eq:c}
\bc\psi=1 \iff \bc=\bone \qquad\text{for all $\bc\in C$,}
\end{equation}
where we write $\bone=(1,1)$ for the identity of $C$.  

Now we let $D=\set{\bc\in C}{\bc^+\not=\bone}$.  By Lemma \ref{lem:T}\ref{T3}, $D$ is a subsemigroup of $C$.  Moreover,~$D$ is closed under ${}^+$, as 
\[
\bc\in D \implies (\bc^+)^+=\bc^+\not=\bone \implies \bc^+\in D,
\]
by \ref{L6}.  It follows that $D$ is a proper left restriction semigroup.  By \eqref{eq:c}, the restriction $\Psi=\psi\restr_D$ maps $D$ into $S=M\sm\{1\}$, as
\[
\bc\psi=1 \implies 1=1^+=(\bc\psi)^+=\bc^+\psi \implies \bc^+=\bone \implies \bc\not\in D.
\]
In fact, $\Psi:D\to S$ is surjective, since for any $s\in S$ we have $s=\bc\psi$ for some $\bc\in C$, and we must in fact have $\bc\in D$; indeed, if $\bc^+=\bone$, then we would have $1=\bone\psi=\bc^+\psi=(\bc\psi)^+=s^+$, contradicting the fact that $S=M\sm\{1\}$ is closed under ${}^+$.  Thus, $\Psi$ is a surmorphism $D\to S$, and it is projection-separating since $\psi$ is.
\epf

\begin{rem}\label{rem:invcover}
Consider an \emph{inverse} monoid $M$, regarded as a left restriction monoid under $s^+=ss^{-1}$.  So $P=P(M)=E(M)$ is the semilattice of (all) idempotents of $M$; as usual, we denote this semilattice by $E$.  Let $C$ and $\psi$ be as in the proof of Corollary \ref{cor:cover}, so that
\[
C = \bigset{(e,s)\in E\rtimes M}{e=es^+} = \bigset{(e,s)\in E\rtimes M}{e\leq s^+}
\]
is a proper left restriction monoid under $(e,s)^+=(e,1)$, and $\psi:C\to M:(e,s)\mt es$ is a covering morphism.  (In the above, $\leq$ is the \emph{natural partial order} on $M$, defined by $s\leq t\iff s\in Et$.  Note that for $e,f\in E$ we have $e\leq f\iff e=ef$.)  It is routine to show the following:
\bit
\item $C$ is inverse, with inversion given by $(e,s)^{-1} = ({}^{s^{-1}}e,s^{-1})$,
\item $E(C) = \bigset{(e,s)}{e,s\in E,\ e=es} = \bigset{(e,s)}{e,s\in E,\ e\leq s}$.
\eit
(To show that $C$ is inverse, one shows that it is regular, and $E(C)$ is commutative.  To verify the above claim about $(e,s)^{-1}$, one needs to note that 
\[
(e,s)\in C \implies e=es^+ \implies {}^s({}^{s^{-1}}e) = {}^{ss^{-1}}e = {}^{s^+}e = s^+e(s^+)^{-1} = s^+es^+ = es^+s^+ = e.)
\]
We note, however, that $C$ is not necessarily a proper inverse monoid, even though it is a proper left restriction monoid.  In fact, one can easily show that $C$ is a proper inverse monoid if and only if~$S$ is a proper inverse monoid.  

The key point here is that $P(C) = \bigset{(e,1)}{e\in E}$ is generally not equal to $E(C)$, so that $(E(C),C)$ might not be a proper action pair, even though $(P(C),C)$ always is; cf.~Proposition~\ref{prop:LR1}.  (In fact, $E(C)=P(C) \iff E=\{1\}$, as $(e,e)\in E(C)\sm P(C)$ for any $e\in E\sm\{1\}$.)  Similarly, $\psi$ is generally not idempotent-separating, even though it is projection-separating, as we have $(e,s)\psi=e$ for any idempotent $(e,s)$ of $C$.  (So again, $\psi$ is idempotent-separating if and only if $E=\{1\}$.  As we have noted on a number of occasions, $E=\{1\}$ is equivalent to $M$ being a group.)
\end{rem}

\sect{Embedding theorems for proper monoids I}\label{sect:embed1}

In this section, and the subsequent three sections, we turn to the task of finding an appropriate generalisation of Theorem \ref{thm:OC} in the context of proper action pairs.  That is, we consider a semigroup $US$ arising from a proper pair $(U,S)$, and seek to embed $US$ in a semidirect product $\U\rtimes(S/\si)$, where~$\U$ is some semigroup containing $U$, and where $\si$ is the congruence on $S$ given in Definition~\ref{defn:si}.  For technical reasons, we will need to assume that $U$ and $S$ are both submonoids of the over-monoid~$M$.  However, this does not end up hampering us, in the sense that we can still derive (the full semigroup version of) Theorem \ref{thm:OC} as a corollary.

The main result of the current section is Theorem \ref{thm:embed1}.  This states that a $(U,S)$-proper monoid $M=US$ embeds in a suitable semidirect product $\U\rtimes(S/\si)$, under the assumption that the submonoid $P=\la S^+\ra$ of $U$ (see Definition \ref{defn:si}) is central in $U$.  This centrality condition is quite natural.  For example, when $M$ is a proper left restriction monoid, arising from the proper action pair $(P,M)$, we have already observed that $P=U$ is the semilattice of projections of~$M$, and is therefore commutative, and hence central in itself.  After proving Theorem \ref{thm:embed1}, we deduce the monoid version of Theorem \ref{thm:OC} in Corollary \ref{cor:embed1}, and then show that the general case follows.

The centrality assumption just discussed is in fact not necessary to obtain an embedding of the desired kind.  Indeed, in the next section we prove Theorem \ref{thm:embed2}, which is essentially the same as Theorem \ref{thm:embed1} but with the centrality assumption removed.  We do, however, have strong motivation for stating and proving the two separate embedding results.  First, the construction of the monoid~$\U$ is somewhat more abstract in the general case.  By contrast, the monoid $\U$ used in the current section is more transparent, and the corresponding semidirect product $\U\rtimes(S/\si)$ is in fact an (unrestricted) wreath product $I\operatorname{wr} {} (S/\si)$, in the sense of \cite{RRT2003,Meldrum1995}, for a suitable monoid~$I$, as we explain in a little more detail below.  Another advantage of the current construction is that when $M$ is a proper left restriction monoid, $\U$ is easily seen to be a semilattice, and as noted above this allows us to deduce Theorem \ref{thm:OC} as a corollary.

In this section, and the next three, an important role will be played by certain special factorisations of elements of a semigroup~$US$ arising from a (weak) action pair $(U,S)$:

\begin{defn}\label{defn:nf}
Let $(U,S)$ be a weak action pair.  A \emph{natural factorisation} for an element $a\in US$ is a pair $(u,s)\in U\times S$ such that $a=us$ and $u=us^+$.
\end{defn}

The next result lists some special properties of these factorisations.  Note that we can think of a natural factorisation of $a\in US$ as an element of the semigroup
\[
\bigset{(u,s)\in U\rtimes S}{u=us^+}
\]
that maps to $a$ under the surmorphism $\pi:U\rtimes S\to US$ from Proposition \ref{prop:A1}; cf.~Remark \ref{rem:A1}.  The first two parts of the next result can be proved with reference to this viewpoint, but we give simple direct proofs for convenience.

\begin{lemma}\label{lem:nf}
Let $(U,S)$ be a weak action pair.
\ben
\item \label{nf1} Every element of $US$ has a natural factorisation.
\item \label{nf2} If $(u,s)$ and $(v,t)$ are natural factorisations of elements $a$ and $b$ of $US$, respectively, then $(u\cdot{}^sv,st)$ is a natural factorisation of $ab$.
\item \label{nf3} If $(U,S)$ is an action pair, and if $(u,s)$ and $(v,t)$ are natural factorisations for the same element of $US$, then $u=v$.
\item \label{nf4} If $(U,S)$ is a proper action pair, and if $(u,s)$ and $(v,t)$ are natural factorisations for the same element of $US$, then $s\mr\si t$.
\een
\end{lemma}

\pf
\firstpfitem{\ref{nf1}}  Let $a\in US$.  We certainly have $a=vs$ for some $v\in U$ and $s\in S$, and we take $u=vs^+$.  Then since $s=s^+s$ (cf.~Remark \ref{rem:AP2}) we have $us=vs^+s=vs=a$, while $us^+=vs^+s^+=vs^+=u$, by Lemma~\ref{lem:+}\ref{lem+1}.

\pfitem{\ref{nf2}}  We obtain $ab = u\cdot{}^sv\cdot st$ right from \ref{A1}.  Combining Lemma \ref{lem:+}\ref{lem+2} with $v=vt^+$, we have
\[
u\cdot{}^sv\cdot(st)^+ = u\cdot{}^sv\cdot{}^s(t^+) = u\cdot{}^s(vt^+) = u\cdot{}^sv.
\]

\pfitem{\ref{nf3}}  Suppose $us=vt$, where $(u,s)$ and $(v,t)$ are natural.  By \ref{A2} we have $us^+=vt^+$, and since $u=us^+$ and $v=vt^+$ (by naturality), it follows that $u=v$.

\pfitem{\ref{nf4}}  This follows immediately from \eqref{eq:proper}.
\epf

\begin{rem}\label{rem:nf}
On many occasions, Lemma \ref{lem:nf} will allow us to define a morphism $\psi$ whose domain $US$ arises from a proper action pair $(U,S)$.  Given $a\in US$, parts \ref{nf1}, \ref{nf3} and \ref{nf4} allow us to unambiguously define $a\psi$, by explaining how to construct it from $u$ and the $\si$-class of~$s$ for \emph{any} natural factorisation $(u,s)$ of $a$.  Part \ref{nf2} will be useful in showing that such a map $\psi$ is a morphism.
\end{rem}

We are now almost ready to state our first embedding theorem.  The theorem concerns a proper action pair $(U,S)$ in a monoid $M$, where $U$ and $S$ are both submonoids.  Without loss of generality, we may assume that $M=US$, so that $M$ is $(U,S)$-proper.  We continue to use the notation~$S^+$, $P=\la S^+\ra$ and $\si=\ka^\sharp$ from Section \ref{sect:cover}.  So $P$ is an idempotent-generated submonoid of $U$, and $\si$ is a congruence on $S$.  
For $s\in S$ we write $\wh s$ for the $\si$-class of $s$.

As we have already noted, the following theorem concerns the case in which $P=\la S^+\ra$ is central in~$U$, meaning that $pu=up$ for all $p\in P$ and $u\in U$.  This is of course equivalent to the generating set $S^+$ being central in $U$.  

\newpage

\begin{thm}\label{thm:embed1}
Let $M=US$ be a $(U,S)$-proper monoid, and suppose $P=\la S^+\ra$ is central in~$U$.  Then there exists a semidirect product $\M = \U\rtimes (S/\si)$ and an embedding $\psi:M\to\M$ such that:
\ben
\item $\U$ contains a subsemigroup $U'$ isomorphic to $U$,
\item $\psi\restr_U:U\to U'\rtimes\{\wh1\}$ is an isomorphism,
\item $\U$ is commutative if $U$ is commutative.
\een
\end{thm}

We build towards the proof of the theorem with a series of preliminary results.  For the remainder of this section, we assume that $M=US$ is a $(U,S)$-proper monoid (with~$U$ and~$S$ submonoids), and that $P=\la S^+\ra$ is central in $U$.

We begin with a result that gives a convenient equational formulation of $\si$.

\begin{lemma}\label{lem:Psi}
The monoid $P$ is a semilattice, and consequently $\si = \bigset{(s,t)\in S\times S}{t^+s = s^+t}$.
\end{lemma}

\pf
Commutativity of $P$ follows from centrality.  Since every element of $S^+$ is an idempotent (cf.~Lemma \ref{lem:+}\ref{lem+1}), so too therefore is every element of $P$.  The assertion regarding $\si$ follows from Lemma \ref{lem:si3}; cf.~Remark \ref{rem:si3}.
\epf

In what follows, it will be convenient to write
\[
\cS=S/\si=\set{\wh s}{s\in S}.
\]
So we must construct a monoid $\U$, a semidirect product $\M=\U\rtimes\cS$, and a suitable embedding~${\psi:M\to\M}$.  The monoid $\U$ will in fact be of the form $I^\cS$ for some monoid $I$, and then~$\U\rtimes\S$ will be an \emph{unrestricted wreath product} $I\operatorname{wr}\cS = I^\cS\rtimes\cS$ in the sense of \cite{RRT2003,Meldrum1995}.  Since we will be studying a different kind of (transformational) wreath product in Chapter \ref{chap:wreath}, and since we wish to emphasise the monoid $\U$ itself in our construction (rather than $I$), we will not speak explicitly of (any kind of) wreath products in the current chapter.  However, in order to prepare the reader for our construction we briefly recall the definition of the unrestricted wreath product~$V\operatorname{wr} T$, for semigroups $V$ and $T$.  First, $V^T$ denotes as usual the semigroup of all functions~$T\to V$ under componentwise product, which we denote by $\star$.  Given $f,g\in V^T$, the product $f\star g\in V^T$ is defined by~$x(f\star g) = xf\cdot xg$ for $x\in T$ (with the product $xf\cdot xg$ taken in $V$).  Next, $T$ has a left action on~$V^T$ by semigroup morphisms, given by $(t,f)\mt {}^t\!f$, where the latter is defined by $x({}^t\!f) = (xt)f$ for $x\in T$.  This allows for the formation of the semidirect product $V^T\rtimes T$, as in Definition \ref{defn:SD}, and this is taken as the definition of $V\operatorname{wr} T$.  As noted above, we will no longer refer explicitly to this kind of wreath product, even though the semidirect product $\M=\U\rtimes\S$ will be of this form.

The construction of $\U$ begins with the \emph{power semigroup} of $U$:
\begin{equation}\label{eq:PU}
\P(U)=\set{V}{V\sub U} \qquad\text{with set product}\qquad V\cdot W = \set{vw}{v\in V,\ w\in W} \qquad\text{for $V,W\sub U$.}
\end{equation}
We then take $I$ to be the principal ideal of $\P(U)$ generated by $P$.  Since $\P(U)$ is a monoid with identity $\{1\}$, the ideal $I = \P(U)\cdot P\cdot\P(U)$ consists of all subsets of $U$ of the form $V\cdot P\cdot W$, with $V,W\sub U$.  In fact, since $P$ is a central submonoid of $U$, we have
\begin{align*}
I &= \P(U) \cdot P = P \cdot \P(U) = \set{V\cdot P=P\cdot V}{V\sub U},
\intertext{and moreover}
I &= P\cdot\P(U)\cdot P = \set{V\in\P(U)}{V=V\cdot P=P\cdot V}.
\end{align*}
Thus, $I$ is a (local) monoid with identity $P$ (and zero~$\es$).
We note, however, that $I$ is generally not a \emph{submonoid} of $\P(U)$, as the respective identities of these monoids are $P$ and $\{1\}$.

We now let
\[
\U = I^\cS
\]
be the set of all functions $\cS\to I$.  We denote the componentwise product in $\U$ by $\star$.  So for $f,g\in\U$, we have $f\star g:\cS\to I$, defined by
\[
\wh x (f\star g) = \wh xf\cdot\wh xg \qquad\text{for each $x\in S$.}
\]
Note here that $\wh xf$ and $\wh xg$ are both elements of $I$, and in particular of $\P(U)$, so $\wh xf\cdot\wh xg$ denotes the set product in~$I$.  It is easy to see that $\star$ is associative, so that $\U$ is a semigroup.  In fact, $\U$ is a monoid, whose identity~$\bone$ is the constant map with image $\{P\}$: i.e., $\wh x\bone=P$ for all $x\in S$.

The first part of the next result verifies the third item in Theorem \ref{thm:embed1}.  The second part concerns the special case in which $U=P$; since $P$ is central in $U$, certainly $U=P$ is commutative in this case.

\begin{lemma}\label{lem:UUcom}
\ben
\item \label{UUcom1} If $U$ is commutative, then so too is $\U$.
\item \label{UUcom2} If $U=P$ (is commutative), then $U$ and $\U$ are both semilattices.
\een
\end{lemma}

\pf
\firstpfitem{\ref{UUcom1}}  If $U$ is commutative, then obviously $\P(U)$ is as well.  But then so too is its subsemigroup~$I$, and hence also $\U=I^\cS$.

\pfitem{\ref{UUcom2}}  We already observed in Lemma \ref{lem:Psi} that $P$ is a semilattice.  Since $\U=I^\cS$, it is enough to show that $I$ is a semilattice.  An arbitrary element of $I$ is a set of the form $P\cdot V\cdot P$ for some $V\sub U=P$.  These sets are precisely the ideals of $P$, and since the set product of two ideals of a semilattice is their intersection, we see that $I$ is the $\cap$-semilattice of ideals of $P$.
\epf

To define the semidirect product $\M=\U\rtimes\cS$, we need an action of $\cS$ on $\U$, and this is also componentwise.  For $s\in S$ and $f\in\U$, we define ${}^{\wh s}\!f\in\U$ by
\[
\wh x({}^{\wh s}\!f) = (\wh x\wh s)f \qquad\text{for all $x\in S$.}
\]
Clearly ${}^{\wh s}\!f$ is well defined: i.e., the above definition does not depend on the choice of representative~$s$ of the $\si$-class $\wh s$.  It is also a routine matter to check that we have a monoidal action by monoid morphisms: i.e., that
\[
{}^{\wh s}({}^{\wh t}\!f) = {}^{\wh s\wh t}\!f \COMMA {}^{\wh s}(f\star g) = {}^{\wh s}\!f\star{}^{\wh s}\!g \COMMA {}^{\wh1}\!f=f \AND {}^{\wh s}\bone=\bone \qquad\text{for all $s,t\in S$ and $f,g\in \U$.}
\]
We may therefore form the semidirect product
\[
\M = \U\rtimes\cS = \bigset{(f,\wh s)}{f\in\U,\ s\in S} \qquad\text{with operation}\qquad (f,\wh s)\cdot(g,\wh t) = (f\star{}^{\wh s}\!g,\wh s\wh t).
\]
By Lemma \ref{lem:USmon}\ref{USmon3}, $\M$ is a monoid with identity $(\bone,\wh 1)$.  (As noted above, our monoid $\M$ is an unrestricted wreath product $I\wr\cS=I^{\cS}\rtimes\cS$.)

Now that we have defined the monoid $\M=\U\rtimes\S$, we wish to construct an embedding~${\psi:M\to\M}$.  As in Remark \ref{rem:nf}, $\psi$ will be defined in terms of natural factorisations.  Specifically, given a natural factorisation $(u,s)$ for an element $a\in M$, we will define $a\psi = (f_u,\wh s)$, for a certain $f_u\in\U$.  We now turn to the definition of these elements $f_u$, and establish the important properties of the mapping $u\mt f_u$.

To this end, let $u\in U$.  For $s\in S$, we define the subset
\begin{equation}\label{eq:Vsu}
V_{\wh s,u} = \set{{}^tu}{t\in\wh s} \sub U.
\end{equation}
We then define the function $f_u\in\U$ by
\[
\wh sf_u = V_{\wh s,u}\cdot P \qquad\text{for all $s\in S$.}
\]
Finally, we define the map
\begin{equation}\label{eq:ufu}
\phi:U\to\U:u\mt f_u.
\end{equation}

\begin{lemma}\label{lem:phi_mono}
The map $\phi:U\to\U$ is a monomorphism.
\end{lemma}

\pf
We begin by showing that $\phi$ is injective.  To do so, suppose $f_u=f_v$ for some $u,v\in U$.  In particular, $\wh1f_u=\wh1f_v$, so that
\[
u = {}^1u\cdot1 \in V_{\wh1,u}\cdot P = \wh1f_u=\wh1f_v = V_{\wh1,v}\cdot P.
\]
It follows that $u = {}^tv\cdot p$ for some $t\in\wh1$ and $p\in P$.  Combining this with Lemma \ref{lem:wh1}\ref{wh12}, and centrality of $P$, it follows that $u = t^+v\cdot p = t^+p\cdot v\in Pv$, as $t^+p\in P$.  By symmetry, $v \in Pu$.  Since $P$ is a semilattice (cf.~Lemma \ref{lem:Psi}), it follows from Lemma \ref{lem:sl} that $u=v$, as required.

To show that $\phi$ is a morphism, suppose $u,v\in U$.  We must show that $f_{uv}=f_u\star f_v$: i.e., that
\begin{align}
\nonumber \wh sf_{uv} &= \wh sf_u \cdot \wh sf_v &&\hspace{-2cm}\text{for all $s\in S$.}
\intertext{By definition, and using centrality of $P$, this reduces to showing that}
\label{eq:Vsuv}
V_{\wh s,uv} \cdot P &= V_{\wh s,u}\cdot V_{\wh s,v}\cdot P &&\hspace{-2cm}\text{for all $s\in S$.}
\end{align}
For the fowards inclusion, let $x\in V_{\wh s,uv} \cdot P$.  Then $x={}^t(uv)\cdot p$ for some $t\in\wh s$ and $p\in P$, and so 
\[
x = {}^tu\cdot{}^tv\cdot p \in V_{\wh s,u}\cdot V_{\wh s,v}\cdot P.
\]
Conversely, let $x\in V_{\wh s,u}\cdot V_{\wh s,v}\cdot P$, so that $x={}^tu\cdot{}^zv\cdot p$ for some $t,z\in \wh s$ and $p\in P$.  Since $t\mr\si z$ (as $t,z\in \wh s$), it follows from Lemma \ref{lem:Psi} that $z^+\cdot t=t^+\cdot z$.  Lemma \ref{lem:ap} then gives 
\begin{equation}\label{eq:z+tv}
z^+\cdot {}^tv=t^+\cdot {}^zv.
\end{equation}
But then
\begin{align*}
x = {}^tu\cdot{}^zv\cdot p &= {}^tu\cdot t^+\cdot{}^zv\cdot p &&\text{by Lemma \ref{lem:+}\ref{lem+4}}\\
&= {}^tu\cdot z^+\cdot{}^tv\cdot p &&\text{by \eqref{eq:z+tv}}\\
&= {}^tu\cdot {}^tv\cdot z^+p &&\text{by centrality of $P$}\\
&= {}^t(uv)\cdot z^+p \in V_{\wh s,uv}\cdot P.
\end{align*}
This completes the proof of \eqref{eq:Vsuv}, and hence of the lemma.
\epf

It follows that $\phi$ is an isomorphism from $U$ to its image, which is the subsemigroup
\[
U' = \set{f_u}{u\in U} \leq \U.
\]
Note that while $U'$ is a monoid (as $U$ is), it is generally not a submonoid of $\U$, as the identities of $U'$ and $\U$ are $f_1$ and $\bone$, respectively.

We can now tie together the loose ends.  The next proof uses the natural factorisations from Definition \ref{defn:nf}.

\pf[\bf Proof of Theorem \ref{thm:embed1}]
We define the map
\[
\psi:M=US \to \M=\U\rtimes\cS
\]
by
\[
a\psi = (f_u,\wh s) \qquad\text{for any natural factorisation $(u,s)$ of $a\in M=US$.}
\]
This is well defined by Lemma \ref{lem:nf}; cf.~Remark \ref{rem:nf}.

We begin by showing that $\psi$ is injective.  To do so, suppose $a\psi=b\psi$ for some $a,b\in M$, and fix natural factorisations $(u,s)$ and $(v,t)$ for $a$ and $b$, respectively.  Then $(f_u,\wh s)=a\psi=b\psi=(f_v,\wh t)$, so that $f_u=f_v$ and $\wh s=\wh t$.  The former implies $u=v$, by Lemma \ref{lem:phi_mono}.  Combining this with naturality, we obtain $us^+=u=v=vt^+$.  Since also $s\mr\si t$ (as $\wh s=\wh t$), and since $(U,S)$ is proper, it follows from \eqref{eq:proper} that $us=vt$: i.e., $a=b$.

Clearly ${u\psi=(f_u,\wh1)}$ for all $u\in U$, so $\psi\restr_U$ maps $U$ bijectively onto $U'\rtimes\{\wh1\}$, and we have already observed that $U'\cong U$.  We also observed in Lemma \ref{lem:UUcom} that $\U$ is commutative if $U$ is.  

It remains to show that $\psi$ is a morphism.  So fix some $a,b\in M$, with natural factorisations $(u,s)$ and $(v,t)$.  We must show that $(ab)\psi = a\psi\cdot b\psi$.  We first observe that
\[
a\psi\cdot b\psi = (f_u,\wh s)\cdot(f_v,\wh t) = (f_u\star{}^{\wh s}\!f_v,\wh s\wh t).
\]
By Lemma \ref{lem:nf}\ref{nf2}, $ab$ has natural factorisation $(u\cdot{}^sv,st)$, and so
\[
(ab)\psi = (f_{u\cdot{}^sv},\wh s\wh t).
\]
Thus, the proof will be complete if we can show that $f_u\star{}^{\wh s}\!f_v = f_{u\cdot{}^sv}$.  To do so, let $x\in S$ be arbitrary.  Then
\[
\wh x f_{u\cdot{}^sv} = V_{\wh x,u\cdot{}^sv} \cdot P \AND \wh x(f_u\star{}^{\wh s}\!f_v) = \wh xf_u\cdot (\wh x\wh s)f_v = V_{\wh x,u} \cdot V_{\wh x\wh s,v} \cdot P,
\]
again using centrality of $P$.  So we must show that
\begin{equation}\label{eq:VP}
V_{\wh x,u\cdot{}^sv} \cdot P = V_{\wh x,u} \cdot V_{\wh x\wh s,v} \cdot P.
\end{equation}

For the forwards inclusion, let $c\in V_{\wh x,u\cdot{}^sv} \cdot P$, so that $c={}^y(u\cdot{}^sv)\cdot p$ for some $y\in\wh x$ and $p\in P$.  Then $c={}^yu\cdot{}^{ys}v\cdot p$, with ${}^yu\in V_{\wh x,u}$ and ${}^{ys}v\in V_{\wh x\wh s,v}$ (as $y\mr\si x\implies ys\mr\si xs$).

Conversely, let $d\in V_{\wh x,u} \cdot V_{\wh x\wh s,v} \cdot P$, so that $d={}^yu\cdot{}^zv\cdot p$ for some $y\in\wh x$, $z\in\wh x\wh s$ and $p\in P$.  Since $z\mr\si xs \mr\si ys$, Lemma \ref{lem:Psi} gives $(ys)^+\cdot z = z^+\cdot ys$.  It then follows from Lemma \ref{lem:ap} that
\begin{equation}\label{eq:ys+v}
(ys)^+\cdot {}^zv = z^+\cdot {}^{ys}v.
\end{equation}
But then
\begin{align*}
d = {}^yu\cdot{}^zv\cdot p &= {}^y(us^+)\cdot{}^zv\cdot p &&\text{as $u=us^+$, by naturality}\\
&= {}^yu\cdot {}^y(s^+)\cdot{}^zv\cdot p \\
&= {}^yu\cdot (ys)^+\cdot{}^zv\cdot p &&\text{by Lemma \ref{lem:+}\ref{lem+2}}\\
&= {}^yu\cdot z^+\cdot {}^{ys}v\cdot p &&\text{by \eqref{eq:ys+v}}\\
&= {}^yu\cdot {}^{ys}v\cdot z^+p &&\text{by centrality of $P$}\\
&= {}^y(u\cdot {}^sv)\cdot z^+p .
\end{align*}
This completes the proof of \eqref{eq:VP}, and hence of the lemma, since $ {}^y(u\cdot {}^sv)\in V_{\wh x,u\cdot{}^sv}$ and~${z^+p\in P}$.
\epf

This time we may deduce Theorem \ref{thm:OC} as a corollary.  Again we begin with the monoid case:

\begin{cor}\label{cor:embed1}
Let $S$ be a proper left restriction monoid with semilattice of projections ${P=P(S)}$.  Then $S$ can be embedded in a semidirect product $Q\rtimes(S/\si)$, where $Q$ is a semilattice containing~$P$.
\end{cor}

\pf
We aim to apply Theorem \ref{thm:embed1} to the pair $(U,S)$ where $U=P$.  We first observe that~$(P,S)$ is indeed a proper action pair, by Proposition \ref{prop:lr}.  Since $P$ is a semilattice, it is certainly central in $U=P$.  Theorem \ref{thm:embed1} then does indeed apply, and it tells us that $S=PS$ embeds in a semidirect product $\U\rtimes\cS$, where $\U$ is a monoid containing (an isomorphic copy) of $U=P$, and where $\cS=S/\si$.  By Lemma \ref{lem:UUcom}\ref{UUcom2}, $\U$ is a semilattice.
\epf

The general case now follows very quickly:

\pf[\bf Proof of Theorem \ref{thm:OC}]
Starting with an arbitrary proper left restriction \emph{semigroup} $S$, we embed this into $S^1$, and then apply Corollary \ref{cor:embed1} to $S^1$.
\epf

\begin{rem}\label{rem:embed1}
Theorem \ref{thm:embed1} concerned the special case that the submonoid $P=\la S^+\ra$ is central in $U$.  In the very special case that $U$ is commutative, the monoid $\U$ we constructed is also commutative.  In the very very special case that $U$ is a semilattice, unfortunately $\U$ is not necessarily a semilattice; it \emph{is}, however, in the very very very special case that $U=P$ is a semilattice, as we showed in Lemma \ref{lem:UUcom}\ref{UUcom2}.  

In the next section we prove Theorem \ref{thm:embed2}.  In a sense, this is a more general version of Theorem \ref{thm:embed1}, as we drop the assumption that $P$ is central in $U$.  In the proof we construct a very different monoid $\U$ to deal with this greater level of generality, and in fact $\U=\Mpres XR$ is defined by a monoid presentation.  This construction does not allow us to deduce that~$\U$ is commutative when~$U$ is commutative; thus, Theorem \ref{thm:embed2} does not completely subsume Theorem~\ref{thm:embed1}.  However, in Section~\ref{sect:embed3} we prove Theorem \ref{thm:embed3}, which is a specialisation of Theorem~\ref{thm:embed2} in the case that $U$ is commutative.  The monoid $\U$ constructed in the proof of Theorem~\ref{thm:embed3} is of the form $\U=\Mpres X{R\cup C}$, where $C$ is an additional set of relations forcing $\U$ to be commutative.  As we will see, this does have the additional benefit that $\U$ is a semilattice when~$U$ is.  Yet another specialisation is given in Theorem \ref{thm:embed4}, which treats the case that $U$ is a left-regular band.  Although a semilattice is a special case of a left-regular band, the differing technicalities in each case lead us to consider the two separately.
\end{rem}

\sect{Embedding theorems for proper monoids II}\label{sect:embed2}

Consider again a proper action pair $(U,S)$ in a monoid $M$, where $U$ and $S$ are submonoids, and assume again that $M=US$.  
As discussed above, our aim now is to prove a stronger version of Theorem \ref{thm:embed1}, in which we drop the assumption that the submonoid ${P=\la S^+\ra}$ is central in $U$.
In fact, the monoid $P$ will play no role at all in the current section, apart from its involvement in the construction of the congruence $\si=\ka^\sharp$ from Definition~\ref{defn:si}.  In the following statement, we continue to write $\wh s$ for the $\si$-class of $s\in S$.

\begin{thm}\label{thm:embed2}
Let $M=US$ be a $(U,S)$-proper monoid.  Then there exists a semidirect product $\M = \U\rtimes (S/\si)$ and an embedding $\psi:M\to\M$ such that:
\ben
\item \label{embed21} $\U$ contains a subsemigroup $U'$ isomorphic to $U$,
\item \label{embed22} $\psi\restr_U:U\to U'\rtimes\{\wh1\}$ is an isomorphism.
\een
\end{thm}

The proof of the theorem occupies the remainder of the section.  Because we no longer assume that $P$ is central in $U$, we are forced to adopt a completely different strategy.  Accordingly, the construction of $\U$ (and $\M$) is somewhat more involved than in the previous section.
As before, we write
\[
\cS = S/\si = \set{\wh s}{s\in S}.
\]

The monoid $\U$ will in fact be defined by means of a presentation $\Mpres XR$, where $X$ and~$R$ are defined as follows.  First,
\[
X = \set{x_{\wh s,u}}{s\in S,\ u\in U}
\]
is an alphabet in one-one correspondence with the cartesian product $\cS \times U$.  We then define $R$ to be the set consisting of all the following relations over $X$, displayed here as equations for clarity:
\begin{equation}\label{eq:R}
x_{\wh t,u\cdot{}^sv} = x_{\wh t,u}x_{\wh t\wh s,v} \qquad\text{for $s,t\in S$ and $u,v\in U$ with $u=us^+$.}
\end{equation}
These relations may appear mysterious at this stage, but the reader may notice a resemblance to \eqref{eq:VP}.  In any case, we will see that the relations contain (just) enough information about products in $U$ and $\cS$, and about the action of $S$ on $U(=U^1)$, to make everything work.

For the duration of this section we write ${\sim}=R^\sharp$ for the congruence on the free monoid $X^*$ generated by~$R$, and we denote the $\sim$-class of $w\in X^*$ by $\ol w$.  (There should be no confusion with our earlier use of the over-line notation in Sections \ref{sect:AAP} and \ref{sect:cover}, as we will never refer to the~$(\ul U,\ol S)$ construction here.)
We then define the monoid
\[
\U = \Mpres XR = X^*/{\sim} = \set{\ol w}{w\in X^*}.
\]
It is worth noting at this point that 
\begin{equation}\label{eq:suv}
x_{\wh t,uv} \sim x_{\wh t,u}x_{\wh t,v} \qquad\text{for all $t\in S$ and $u,v\in U$.}
\end{equation}
Indeed, taking $s=1$ in \eqref{eq:R}, we see that $R$ contains the relation $x_{\wh t,uv} = x_{\wh t,u}x_{\wh t,v}$.

In order to define a semidirect product $\U\rtimes\cS$, we need an action of $\cS$ on $\U$.  We begin by defining an action of $\cS$ on $X^*$.  For this, let $s\in S$ and $w\in X^*$.  Then $w=x_{\wh t_1,u_1}\cdots x_{\wh t_k,u_k}$ for some $t_1,\ldots,t_k\in S$ and $u_1,\ldots,u_k\in U$, and we define
\begin{equation}\label{eq:SactX*}
{}^{\wh s}w = x_{\wh s\wh t_1,u_1}\cdots x_{\wh s\wh t_k,u_k}.
\end{equation}
This is trivially a monoidal action by monoid morphisms: i.e.,
\[
{}^{\wh s}({}^{\wh t}w)={}^{\wh s\wh t}w \COMMa {}^{\wh s}(ww')={}^{\wh s}w\cdot{}^{\wh s}w' \COMMa {}^{\wh1}w=w \ANd {}^{\wh s}\ew=\ew \qquad\text{for all $s,t\in S$ and $w,w'\in X^*$.}
\]
Here $\ew$ denotes the empty word, which is the identity of $X^*$.  Given any relation $(w,w')$ from~$R$, as in \eqref{eq:R}, and any $s\in S$, it is easy to see that $({}^{\wh s}w,{}^{\wh s}w')$ is again a relation from $R$.  It immediately follows that we have an induced action of $\cS$ on $\U$ given by
\[
{}^{\wh s} \ol w = \ol{{}^{\wh s}w} \qquad\text{for $w\in X^*$ and $s\in S$.}
\]
Of course this is still monoidal, and by monoid morphisms.  Consequently, we may form the semidirect product
\[
\M=\U\rtimes\S.
\]
By Corollary \ref{cor:USmon}, $\M$ is a monoid with identity $(\ol\ew,\wh1)$.

We now define a function
\[
\phi:U\to\U : u \mt \ol x_{\wh1,u}.
\]
Taking $t=1$ in \eqref{eq:suv}, it follows quickly that $\phi$ is a semigroup morphism.  (But note that $\phi$ is not a \emph{monoid} morphism, as $1\phi = \ol x_{\wh1,1}$, while the identity of $\U$ is $\ol\ew$, the $\sim$-class of the empty word $\ew$.  Since $R$ consists entirely of relations with non-empty words on both sides (see \eqref{eq:R}), it follows immediately that $\ol\ew=\{\ew\}$, and so $\ol x_{\wh1,1}\not=\ol\ew$.)

The next result demonstrates a crucial relationship between the respective actions of~$S$ and~$\cS$ on~$U$ and~$\U$.

\begin{lemma}\label{lem:phi}
If $u,v\in U$ and $s\in S$ are such that $u=us^+$, then
\[
(u\cdot{}^sv)\phi = (u\phi)\cdot{}^{\wh s}(v\phi).
\]
\end{lemma}

\pf
Taking $t=1$ in \eqref{eq:R}, we see that 
\[
\ol x_{\wh1,u\cdot{}^sv} = \ol x_{\wh 1,u}\cdot \ol x_{\wh s,v}.
\]
But of course $\ol x_{\wh1,u\cdot{}^sv} = (u\cdot{}^sv)\phi$ and $\ol x_{\wh 1,u} = u\phi$, while $\ol x_{\wh s,v} = \ol x_{\wh s\wh 1,v} = {}^{\wh s}\ol x_{\wh1,v} = {}^{\wh s}(v\phi)$.
\epf

The next result shows that $\U$ and $\phi$ are, in a sense, `universal' with respect to the property just established.

\begin{prop}\label{prop:xi1}
Let $M=US$ be a $(U,S)$-proper monoid.
\ben
\item \label{xi11} Suppose there is a semigroup $\V$ on which $\cS=S/\si$ acts monoidally by semigroup morphisms, and a morphism $\xi:U\to \V$ such that
\begin{equation}\label{eq:xi}
(u\cdot{}^sv)\xi = (u\xi)\cdot{}^{\wh s}(v\xi) \qquad\text{for all $u,v\in U$ and $s\in S$ with $u=us^+$.}
\end{equation}
Then there exists a morphism $\Xi:\U\to \V$ such that the following diagram commutes:
\[
\begin{tikzcd}[row sep=2.5em]
 & \U \arrow{dr}{\Xi} \\
U \arrow{ur}{\phi} \arrow{rr}{\xi} && \V
\end{tikzcd}
\]
\item \label{xi12} The morphism $\phi:U\to\U$ is an embedding if and only if some $\V$ and $\xi$ exist as in part \ref{xi11}, with $\xi$ an embedding.
\een
\end{prop}

\pf
\firstpfitem{\ref{xi11}}  Suppose $\V$ and $\xi:U\to \V$ exist with the stated properties.  We begin by defining a morphism
\[
\ze:X^*\to \V:x_{\wh s,u} \mt {}^{\wh s}(u\xi).
\]
Using \eqref{eq:xi}, it is easy to check that $\ze$ preserves the relations from $R$: i.e., that
\[
(x_{\wh t,u\cdot{}^sv})\ze = (x_{\wh t,u}x_{\wh t\wh s,v})\ze \qquad\text{for all $s,t\in S$ and $u,v\in U$ with $u=us^+$.}
\]
It follows that $R^\sharp\sub\ker(\ze)$, so $\ze$ induces a morphism $\Xi : \U=X^*/R^\sharp\to \V:\ol w\mt w\ze$.  This $\Xi$ is defined on generators of $\U$ by
\[
\Xi : \U\to \V : \ol x_{\wh s,u} \mt {}^{\wh s}(u\xi).
\]
But then for any $u\in U$ we have
\[
u(\phi\circ\Xi) = \ol x_{\wh1,u} \Xi = {}^{\wh1}(u\xi) = u\xi,
\]
since $\cS$ acts monoidally on $\V$.

\pfitem{\ref{xi12}}  If $\phi$ is an embedding, then we just take $\V=\U$ and $\xi=\phi$; cf.~Lemma \ref{lem:phi}.  Conversely, if $\xi=\phi\circ\Xi$ is an embedding, then of course $\phi$ is as well.
\epf

Recall that we wish to show that the $(U,S)$-proper monoid $M=US$ embeds in the semidirect product $\M=\U\rtimes\cS$.  We will soon see that a certain very natural kind of embedding exists, and it is convenient to explicitly formulate what we mean by `natural' here.  The following definition refers to the natural factorisations from Definition \ref{defn:nf}.

\begin{defn}\label{defn:nm}
Let $M=US$ be a $(U,S)$-proper monoid, and let $\V$ be a semigroup on which $\cS=S/\si$ acts monoidally via semigroup morphisms.  A morphism
\[
\psi:M\to \V\rtimes \cS
\]
is \emph{natural} if there exists a morphism $\xi:U\to \V$ such that
\[
a\psi = (u\xi,\wh s) \qquad\text{for any natural factorisation $(u,s)$ of $a\in M=US$.}
\]
In this case we say $\psi$ is \emph{natural along $\xi$}.
\end{defn}

Strictly speaking, the next two results are not essential for our main purpose, which is to construct an embedding of $M=US$ into $\M=\U\rtimes\cS$.  However, we include them as they demonstrate that $\U$ is somehow `canonical' with respect to the existence of a \emph{natural} such embedding.

\begin{prop}\label{prop:xi2}
Let $M=US$ be a $(U,S)$-proper monoid, and let $\V$ be a semigroup on which $\cS=S/\si$ acts monoidally via semigroup morphisms.  For a morphism $\xi:U\to \V$, define the map
\[
\psi:M\to \V\rtimes\cS \qquad\text{by}\qquad a\psi = (u\xi,\wh s) \qquad\text{for any natural factorisation $(u,s)$ of $a\in M=US$.}
\]
Then
\ben
\item \label{xi21} $\psi$ is a morphism (and is hence natural along $\xi$) if and only if \eqref{eq:xi} holds,
\item \label{xi22} $\psi$ is injective if and only if $\xi$ is injective.
\een
\end{prop}

\pf
\firstpfitem{\ref{xi21}}  Let $a$ and $b$ be elements of $US$, with natural factorisations $(u,s)$ and $(v,t)$, respectively.  By Lemma \ref{lem:nf}\ref{nf2}, $ab$ has natural factorisation $(u\cdot{}^sv,st)$, so we have
\[
(ab)\psi = ((u\cdot{}^sv)\xi,\wh s\wh t) \AND a\psi\cdot b\psi = (u\xi,\wh s)\cdot(v\xi,\wh t) = ((u\xi)\cdot{}^{\wh s}(v\xi),\wh s\wh t).
\]
The claim follows quickly.

\pfitem{\ref{xi22}}  Suppose first that $\psi$ is injective, and suppose $u,v\in U$ are such that $u\xi=v\xi$.  Now, $(u,1)$ and $(v,1)$ are clearly natural factorisations of $u,v\in U\sub M$, so
\[
u\psi = (u\xi,\wh1) = (v\xi,\wh1) = v\psi .
\]
It follows from injectivity of $\psi$ that $u=v$.

Conversely, suppose $\xi$ is injective, and suppose $a,b\in M$ are such that $a\psi=b\psi$.  Fix natural factorisations $(u,s)$ and $(v,t)$ of $a$ and $b$, respectively.  Then
\[
(u\xi,\wh s) = a\psi = b\psi = (v\xi,\wh t),
\]
so it follows that $u\xi=v\xi$ and $s\mr\si t$.  From the former, and injectivity of $\xi$, we deduce $u=v$.  By naturality, we have $us^+=u=v=vt^+$.  Since $s\mr\si t$, and since $(U,S)$ is proper, it follows from~\eqref{eq:proper} that $us=vt$: i.e., $a=b$.
\epf

The next result summarises our progress so far:

\begin{cor}\label{cor:phi}
Let $M=US$ be a $(U,S)$-proper monoid.  Then the following are equivalent:
\ben
\item \label{phi1} there is a natural embedding of $M$ in a semidirect product $\V\rtimes\cS$ for some semigroup $\V$,
\item \label{phi2} $\phi:U\to\U$ is an embedding,
\item \label{phi3} there is a natural embedding of $M$ in the semidirect product $\U\rtimes\cS$ along $\phi$.
\een
\end{cor}

\pf
\firstpfitem{\ref{phi1}$\implies$\ref{phi2}}  If such a natural embedding exists, say along $\xi:U\to\V$, then by Proposition~\ref{prop:xi2}, $\xi$ is an embedding and satisfies \eqref{eq:xi}.  It then follows from Proposition \ref{prop:xi1}\ref{xi12} that $\phi$ is an embedding.

\pfitem{\ref{phi2}$\implies$\ref{phi3}}  This follows immediately from Proposition \ref{prop:xi2}, taking $\V=\U$ and $\xi=\phi$, and keeping Lemma \ref{lem:phi} in mind.

\pfitem{\ref{phi3}$\implies$\ref{phi1}}  This is obvious.
\epf

Thus, the main step remaining in the proof of Theorem \ref{thm:embed2} is to show that the morphism~${\phi:U\to\U}$ is injective: i.e., that
\[
\ol x_{\wh1,u}=\ol x_{\wh1,v} \implies u=v \qquad\text{for all $u,v\in U$.}
\]
This turns out to be quite involved, and is finally accomplished in Proposition \ref{prop:phi} below.  The proof involves a certain invariant associated to the $\sim$-class of (some) words over $X$, including those of the form $x_{\wh1,u}$.  This invariant comes from what we will call a \emph{trace}.  Working towards its definition, we begin by defining a map $\lam : X^+ \to \cS$ as follows.  (Here, as usual, $X^+=X^*\sm\{\ew\}$ is the free \emph{semigroup} over $X$.)  For ${w=x_{\wh s_1,u_1}\cdots x_{\wh s_k,u_k}\in X^+}$, we define $w\lam=\wh s_1$.  So $w\lam$ is the `$\cS$-coordinate' of the first letter from $w$.  

\begin{lemma}\label{lem:lam}
If $w,w'\in X^+$ and $w\sim w'$, then $w\lam=w'\lam$.
\end{lemma}

\pf
This is true for every relation $(w,w')\in R$, as in \eqref{eq:R}, and the result quickly follows.
\epf

We now define the set
\[
W = \set{w\in X^+}{w\lam=\wh1}
\]
of all (non-empty) words over $X$ whose first letter has `$\cS$-coordinate' $\wh1$.  It follows from Lemma~\ref{lem:lam} that $W$ is closed under $\sim$, in the sense that two $\sim$-equivalent words from $X^*$ are either both contained in $W$ or both not.  
We now introduce the idea of a \emph{trace} of (some) words from $W$.  

\begin{defn}\label{defn:trace}
Consider a word $w\in W$, so that
\[
w=x_{\wh 1,u_1}x_{\wh s_2,u_2}\cdots x_{\wh s_k,u_k}  \qquad\text{for some $k\geq1$ and some $u_1,\ldots,u_k\in U$ and $s_2,\ldots,s_k\in S$.}
\]
Consider also a tuple $\bt=(t_2,\ldots,t_k)\in \wh s_2\times\cdots\times\wh s_k$, so that $t_i\in S$ and $t_i\mr\si s_i$ for all $2\leq i\leq k$.  For each $1\leq i\leq k$, we define
\[
p_i = u_1\cdot {}^{t_2}u_2 \cdots {}^{t_i}u_i \in U.
\]
(As usual, we interpret $p_1=u_1$.)  We say the tuple $\bp=(p_1,\ldots,p_k)\in U^k$ is a \emph{trace} of $w$ if
\[
p_{i-1} = p_{i-1}t_i^+ \qquad\text{for all $2\leq i\leq k$.}
\]
In this case, we say the trace $\bp=(p_1,\ldots,p_k)$ of $w$ is \emph{witnessed} by the tuple $\bt=(t_2,\ldots,t_k)$.
\end{defn}

\begin{rem}\label{rem:trace}
It is not always obvious whether a given word from $W$ has any trace at all.  However, it is easy to see that a word $w=x_{\wh1,u}$ from $W$ of length $1$ has trace $\bp=(u)$, as witnessed (vacuously) by the empty tuple $\bt=\es$.  

On the other hand, consider a word $w=x_{\wh1,u}x_{\wh s,v}\in W$ of length $2$.  The possible tuples $\bt$ in Definition \ref{defn:trace} have the form $\bt=(t)$, where $t\in \wh s$.  And then $\bp=(u,u\cdot{}^tv)$ is a trace for~$w$ (witnessed by $\bt$) if and only if $p_1=p_1t_2^+$: i.e., $u=ut^+$.  To summarise, $w=x_{\wh1,u}x_{\wh s,v}$ has a trace if and only if $u=ut^+$ for some $t\in\wh s$.
\end{rem}

Nevertheless, the next result shows that if a word has a trace, then this is unique.

\begin{lemma}\label{lem:trace1}
A word $w\in W$ has at most one trace.
\end{lemma}

\pf
Consider a word $w =x_{\wh 1,u_1}x_{\wh s_2,u_2}\cdots x_{\wh s_k,u_k}\in W$, and suppose $w$ has traces ${\bp=(p_1,\ldots,p_k)}$ and $\bq=(q_1,\ldots,q_k)$, witnessed by tuples $\bt=(t_2,\ldots,t_k)$ and $\bz=(z_2,\ldots,z_k)$, respectively.  So by definition, $\bt,\bz\in\wh s_2\times\cdots\times\wh s_k$, and we have
\begin{align*}
p_i = u_1\cdot {}^{t_2}u_2 \cdots {}^{t_i}u_i &\AND q_i = u_1\cdot {}^{z_2}u_2 \cdots {}^{z_i}u_i &&\text{for all $1\leq i\leq k$,}
\intertext{and also}
p_{i-1} = p_{i-1}t_i^+ &\AND q_{i-1} = q_{i-1}z_i^+ &&\text{for all $2\leq i\leq k$.}
\end{align*}
We must show that $\bp=\bq$, and for this we show by induction that $p_i=q_i$ for all $1\leq i\leq k$.  Certainly $p_1=u_1=q_1$.  We now assume that $i\geq2$, and that $p_{i-1}=q_{i-1}$.  Then
\[
p_{i-1}t_i^+ = p_{i-1} = q_{i-1} = q_{i-1}z_i^+.
\]
Since also $t_i\mr\si z_i$ (as $t_i,z_i\in\wh s_i$), and since $(U,S)$ is proper, it follows from \eqref{eq:proper} that~${p_{i-1}t_i = q_{i-1}z_i}$.  Lemma \ref{lem:ap} then gives
\[
p_{i-1}\cdot{}^{t_i}u_i = q_{i-1}\cdot{}^{z_i}u_i.
\]
But by definition, we have $p_{i-1}\cdot{}^{t_i}u_i = p_i$ and $q_{i-1}\cdot{}^{z_i}u_i = q_i$, so this completes the inductive step.
\epf

Here is the main technical result we need:

\begin{lemma}\label{lem:trace2}
Suppose $w,w'\in W$ are such that $w\sim w'$.
\ben
\item \label{trace21} If either of $w$ or $w'$ has a trace, then both do.
\item \label{trace22} If $w$ and $w'$ have traces $\bp=(p_1,\ldots,p_k)$ and $\bq=(q_1,\ldots,q_l)$, then $p_k=q_l$.
\een
\end{lemma}

\pf
We prove both parts at the same time, and it suffices to assume that $w$ and $w'$ differ by a single application of a relation from $R$.  That is, we assume that
\[
w = w_1 \cdot x_{\wh t,u\cdot{}^sv} \cdot w_2 \AND w' = w_1\cdot x_{\wh t,u}x_{\wh t\wh s,v}\cdot w_2,
\]
for some $w_1,w_2\in X^*$, and some $s,t\in S$ and $u,v\in U$ with $u=us^+$.  We also write
\[
w_1 = x_{\wh s_1,u_1}\cdots x_{\wh s_m,u_m} \AND w_2 = x_{\wh s_{m+2},u_{m+2}}\cdots x_{\wh s_k,u_k} \qquad\text{where each $s_i\in S$ and $u_i\in U$,}
\]
noting that either (or both) of $w_1$ or $w_2$ could be empty.  It is also convenient to define
\[
s_{m+1} = t \AND u_{m+1} = u\cdot{}^sv.
\]
So then
\[
w = x_{\wh s_1,u_1}\cdots x_{\wh s_k,u_k} \AND w' = x_{\wh s_1,u_1}\cdots x_{\wh s_m,u_m} \cdot x_{\wh t,u}x_{\wh t\wh s,v} \cdot x_{\wh s_{m+2},u_{m+2}}\cdots x_{\wh s_k,u_k}.
\]
Since $w\in W$, we must have $\wh s_1=\wh1$; in the case that $w_1=\ew$ (i.e., $m=0$), this says that $\wh t=\wh1$.  We also define
\begin{equation}\label{eq:v}
(v_1,\ldots,v_{k+1}) = (u_1,u_2,\ldots,u_m,u,v,u_{m+2},\ldots,u_k).
\end{equation}
Note that $v_i$ is the `$U$-coordinate' of the $i$th letter of $w'$.

We begin with the forwards implication.

\pfitem{($\Rightarrow$)}  Suppose first that $w$ has trace $\bp=(p_1,\ldots,p_k)$, witnessed by $\bt=(t_2,\ldots,t_k)\in\wh s_2\times\cdots\times\wh s_k$, meaning that
\bit
\item $p_i = u_1\cdot{}^{t_2}u_2\cdots {}^{t_i}u_i$ for all $1\leq i\leq k$, and
\item $p_{i-1}=p_{i-1}t_i^+$ for all $2\leq i\leq k$.
\eit
We claim that $w'$ has trace $\bq$, witnessed by $\bz$, for
\[
\bz = (t_2,\ldots,t_m,t_{m+1},t_{m+1}\cdot s,t_{m+2},\ldots,t_k)
\AND
\bq = (p_1,\ldots,p_m,p_m\cdot{}^{t_{m+1}}u,p_{m+1},\ldots,p_k).
\]
Note that $\bz$ and $\bq$ are obtained from $\bt$ and $\bp$, respectively, by inserting single entries.  Since the final entries of $\bp$ and $\bq$ are equal, this will complete the proof of the forwards implication.  To assist with the following calculations, it is convenient to line up the entries of $\bz$ and $\bq$, along with the letters of $w'$:
\begin{center}
\begin{tabular}{c|ccccccccc}
&$1$&$2$&$\cdots$&$m$&$m+1$&$m+2$&$m+3$&$\cdots$&$k+1$\\
\hline
$\bz$&&$t_2$&$\cdots$&$t_m$&$t_{m+1}$&$t_{m+1}\cdot s$&$t_{m+2}$&$\cdots$&$t_k$\\
$\bq$&$p_1$&$p_2$&$\cdots$&$p_m$&$p_m\cdot{}^{t_{m+1}}u$&$p_{m+1}$&$p_{m+2}$&$\cdots$&$p_k$\\
$w'$&$x_{\wh s_1,u_1}$&$x_{\wh s_2,u_2}$&$\cdots$&$ x_{\wh s_m,u_m} $&$ x_{\wh t,u}$&$x_{\wh t\wh s,v} $&$ x_{\wh s_{m+2},u_{m+2}}$&$\cdots $&$x_{\wh s_k,u_k}$
\end{tabular}
\end{center}
To prove the claim, we follow Definition \ref{defn:trace} carefully, and we proceed in three steps.

\pfstep1  We first check that $\bz\in\wh s_2\times\cdots\times\wh s_m \times \wh t \times \wh t\wh s \times \wh s_{m+2}\times\cdots\times\wh s_k$.  
\bit
\item For $2\leq i\leq m$ we have $z_i=t_i \mr\si s_i$.
\item For $i=m+1$ we have $z_{m+1} = t_{m+1} \mr\si s_{m+1} = t$.
\item For $i=m+2$ we have $z_{m+2} = t_{m+1}\cdot s \mr\si s_{m+1}\cdot s = ts$.
\item For $m+3\leq i\leq k+1$ we have $z_i=t_{i-1} \mr\si s_{i-1}$.
\eit

\pfstep2  Next we check that $q_i = v_1\cdot {}^{z_2}v_2\cdots{}^{z_i}v_i$ for all $1\leq i\leq k+1$, where the $v_j$ are as in~\eqref{eq:v}.
\bit
\item For $1\leq i\leq m$ we have $q_i = p_i = u_1\cdot{}^{t_2}u_2\cdots {}^{t_i}u_i = v_1\cdot{}^{z_2}v_2\cdots {}^{z_i}v_i$.
\item For $i=m+1$ we have $q_{m+1} = p_m\cdot{}^{t_{m+1}}u = u_1\cdot{}^{t_2}u_2\cdots {}^{t_m}u_m \cdot {}^{t_{m+1}}u = v_1\cdot{}^{z_2}v_2\cdots {}^{z_m}v_m{}^{z_{m+1}}v_{m+1}$.
\item For $m+2\leq i\leq k+1$ we have 
\begin{align*}
q_i = p_{i-1} &= (u_1\cdot{}^{t_2}u_2\cdots {}^{t_m}u_m)\cdot {}^{t_{m+1}}u_{m+1} \cdot ({}^{t_{m+2}}u_{m+2}\cdots {}^{t_{i-1}}u_{i-1})\\
&= (u_1\cdot{}^{t_2}u_2\cdots {}^{t_m}u_m)\cdot {}^{t_{m+1}}(u\cdot{}^sv) \cdot ({}^{t_{m+2}}u_{m+2}\cdots {}^{t_{i-1}}u_{i-1})\\
&= (u_1\cdot{}^{t_2}u_2\cdots {}^{t_m}u_m)\cdot {}^{t_{m+1}}u\cdot {}^{t_{m+1}\cdot s}v \cdot ({}^{t_{m+2}}u_{m+2}\cdots {}^{t_{i-1}}u_{i-1})\\
&= (v_1\cdot{}^{z_2}v_2\cdots {}^{z_m}v_m)\cdot {}^{z_{m+1}}v_{m+1}\cdot {}^{z_{m+2}}v_{m+2} \cdot ({}^{z_{m+3}}v_{m+3}\cdots {}^{z_i}v_i).
\end{align*}
\eit

\pfstep3
Finally, we check that $q_{i-1}=q_{i-1}z_i^+$ for all $2\leq i\leq k+1$.
\bit
\item For $2\leq i\leq m+1$ we have $q_{i-1} = p_{i-1} = p_{i-1}t_i^+ = q_{i-1}z_i^+$.
\item For $i=m+2$ we have
\begin{align*}
q_{m+1}\cdot z_{m+2}^+ &= p_m\cdot{}^{t_{m+1}}u\cdot(t_{m+1}\cdot s)^+\\
&= p_m\cdot{}^{t_{m+1}}u\cdot{}^{t_{m+1}}(s^+) &&\text{by Lemma \ref{lem:+}\ref{lem+2}}\\
&= p_m\cdot{}^{t_{m+1}}(us^+) \\
&= p_m\cdot{}^{t_{m+1}}u &&\text{as $u=us^+$}\\
&= q_{m+1}.
\end{align*}
\item For $m+3\leq i\leq k+1$ we have $q_{i-1} = p_{i-2} = p_{i-2}t_{i-1}^+ = q_{i-1}z_i^+$.
\eit

\pfitem{($\Leftarrow$)}  Conversely, we now assume $w'$ has trace $\bq=(q_1,\ldots,q_{k+1})$, witnessed by
\begin{equation}\label{eq:z}
\bz=(z_2,\ldots,z_{k+1})\in\wh s_2\times\cdots\times\wh s_m \times \wh t \times \wh t\wh s \times \wh s_{m+2}\times\cdots\times\wh s_k,
\end{equation}
meaning that
\bit
\item $q_i = v_1\cdot{}^{z_2}v_2\cdots {}^{z_i}v_i$ for all $1\leq i\leq k+1$, where again the $v_j$ are as in \eqref{eq:v}, and
\item $q_{i-1}=q_{i-1}z_i^+$ for all $2\leq i\leq k+1$.
\eit
We claim that $w$ has trace $\bp$, witnessed by $\bt$, where
\[
\bt = (z_2,\ldots,z_{m+1},z_{m+3},\ldots,z_{k+1}) \AND \bp = (q_1,\ldots,q_m,q_{m+2},\ldots,q_{k+1}).
\]
This time $\bt$ and $\bp$ are obtained by deleting an appropriate entry of $\bz$ and $\bq$.
Again, since the final entries of $\bp$ and $\bq$ are equal, this will complete the proof of the backwards implication, and hence of the lemma.  
(The reader might worry that $\bp$ and $\bq$ could perhaps not have the same final entry if $m+1=k+1$: i.e., if $\bp$ is obtained by deleting the \emph{final} entry of $\bq$.  However, this case does not arise, as by the form of $w=x_{\wh s_1,u_1}\cdots x_{\wh s_m,u_m} \cdot x_{\wh t,u\cdot{}^sv} \cdot x_{\wh s_{m+2},u_{m+2}}\cdots x_{\wh s_k,u_k}$, we have $k\geq m+1$.)
This time the entries of $\bt$ and $\bp$, and the letters from $w$, are as follows:
\begin{center}
\begin{tabular}{c|ccccccccc}
&$1$&$2$&$\cdots$&$m$&$m+1$&$m+2$&$\cdots$&$k$\\
\hline
$\bt$&&$z_2$&$\cdots$&$z_m$&$z_{m+1}$&$z_{m+3}$&$\cdots$&$z_{k+1}$\\
$\bp$&$q_1$&$q_2$&$\cdots$&$q_m$&$q_{m+2}$&$q_{m+3}$&$\cdots$&$q_{k+1}$\\
$w$&$x_{\wh s_1,u_1}$&$x_{\wh s_2,u_2}$&$\cdots$&$ x_{\wh s_m,u_m} $&$ x_{\wh t,u\cdot{}^sv}$&$ x_{\wh s_{m+2},u_{m+2}}$&$\cdots $&$x_{\wh s_k,u_k}$\\
\multicolumn{5}{c}{}&$\rotatebox[origin=c]{90}{=}$\\
\multicolumn{5}{c}{}&$x_{\wh s_{m+1},u_{m+1}}$
\end{tabular}
\end{center}
Again we prove the claim in three steps.

\pfstep1
We first note that $\bt\in\wh s_2\times\cdots\times\wh s_m\times\wh t\times\wh s_{m+2}\times\cdots\times\wh s_k$ follows immediately from~\eqref{eq:z}.

\pfstep2
Next we check that $p_i = u_1\cdot {}^{z_2}u_2\cdots{}^{z_i}u_i$ for all $1\leq i\leq k$.
\bit
\item For $1\leq i\leq m$ we have $p_i=q_i=v_1\cdot{}^{z_2}v_2\cdots{}^{z_i}v_i=u_1\cdot{}^{t_2}u_2\cdots{}^{t_i}u_i$.
\item The $m+1\leq i\leq k$ case is rather more involved.  For this we first claim that 
\begin{equation}\label{eq:qzqzs}
q_{m+1}z_{m+2} = q_{m+1}z_{m+1}s.
\end{equation}
To prove this, it suffices by \eqref{eq:proper} to show that
\begin{equation}\label{eq:qzqzs1}
q_{m+1}z_{m+2}^+ = q_{m+1}(z_{m+1}s)^+ \AND z_{m+2} \mr\si z_{m+1}s.
\end{equation}
The latter follows from $z_{m+2} \mr\si ts \mr\si z_{m+1}s$; cf.~\eqref{eq:z}.  For the former, we first note that ${q_{m+1}z_{m+2}^+=q_{m+1}}$.  On the other hand, since $q_{m+1} = q_m \cdot {}^{z_{m+1}}v_{m+1}$, we have
\begin{align}
\nonumber q_{m+1}(z_{m+1}s)^+ &= q_m \cdot {}^{z_{m+1}}v_{m+1} \cdot {}^{z_{m+1}}(s^+) &&\text{by Lemma \ref{lem:+}\ref{lem+2}}\\
\nonumber &= q_m \cdot {}^{z_{m+1}}(v_{m+1}s^+) \\
\nonumber &= q_m \cdot {}^{z_{m+1}}v_{m+1} &&\text{as $v_{m+1}s^+=us^+=u=v_{m+1}$}\\
\nonumber &= q_{m+1}.
\intertext{This completes the proof of \eqref{eq:qzqzs1}, and hence of \eqref{eq:qzqzs}.   We then have}
\nonumber q_{m+1}\cdot {}^{z_{m+2}}v_{m+2} &= q_{m+1}\cdot{}^{z_{m+1}s}v_{m+2} &&\text{by \eqref{eq:qzqzs} and Lemma \ref{lem:ap}}\\
\nonumber &= q_m \cdot {}^{z_{m+1}}v_{m+1}\cdot{}^{z_{m+1}s}v_{m+2} \\
\nonumber &= q_m \cdot {}^{z_{m+1}}(v_{m+1}\cdot{}^sv_{m+2}) \\
\nonumber &= q_m \cdot {}^{z_{m+1}}(u\cdot{}^sv) \\
\label{eq:qzqzs2} &= q_m \cdot {}^{t_{m+1}}u_{m+1} .
\end{align}
But then for any $m+1\leq i\leq k$ we have
\begin{align*}
p_i = q_{i+1} &= q_{m+1} \cdot {}^{z_{m+2}}v_{m+2} \cdot {}^{z_{m+3}}v_{m+3}\cdots {}^{z_{i+1}}v_{i+1} \\
&= q_m \cdot {}^{t_{m+1}}u_{m+1} \cdot {}^{z_{m+3}}v_{m+3}\cdots {}^{z_{i+1}}v_{i+1} &&\text{by \eqref{eq:qzqzs2}}\\
&= v_1\cdot{}^{z_2}v_2\cdots{}^{z_m}v_m \cdot {}^{t_{m+1}}u_{m+1} \cdot {}^{z_{m+3}}v_{m+3}\cdots {}^{z_{i+1}}v_{i+1} \\
&= u_1\cdot{}^{t_2}u_2\cdots{}^{t_m}u_m \cdot {}^{t_{m+1}}u_{m+1} \cdot {}^{t_{m+2}}u_{m+2}\cdots {}^{t_i}u_i ,
\end{align*}
as required.
\eit

\pfstep3
Finally, we check that $p_{i-1}=p_{i-1}t_i^+$ for all $2\leq i\leq k$.
\bit
\item For $2\leq i\leq m+1$ we have $p_{i-1}=q_{i-1}=q_{i-1}z_i^+=p_{i-1}t_i^+$.
\item For $m+2\leq i\leq k$ we have $p_{i-1}=q_i=q_iz_{i+1}^+=p_{i-1}t_i^+$.  
\eit
We have now finally completed the proof.
\epf

It is now relatively straightforward to show the following.

\begin{prop}\label{prop:phi}
The morphism $\phi:U\to\U:u\mt\ol x_{\wh1,u}$ is an embedding.
\end{prop}

\pf
Suppose $u,v\in U$ are such that $\ol x_{\wh1,u} = \ol x_{\wh1,v}$.  This means that $w \sim w'$, where $w=x_{\wh1,u}$ and $w'=x_{\wh1,v}$.  As in Remark \ref{rem:trace}, $w$ and $w'$ have traces $(u)$ and $(v)$, respectively, so it follows from Lemma \ref{lem:trace2} that $u=v$.
\epf

And we can now finally tie together the loose ends.

\pf[\bf Proof of Theorem \ref{thm:embed2}]
We take $\U=\Mpres XR$, as constructed above.  By Corollary \ref{cor:phi} and Proposition \ref{prop:phi}, there is a natural embedding $\psi$ of $M$ in $\M=\U\rtimes\cS$ along $\phi$.  By definition, this $\psi$ is defined by
\[
a\psi = (u\phi,\wh s) \qquad\text{for any natural factorisation $(u,s)$ of $a\in M=US$.}
\]
It remains to check items \ref{embed21} and \ref{embed22} in the statement of the theorem.  But
\[
u\psi = (u\phi,\wh1) \qquad\text{for all $u\in U$,}
\]
as $(u,1)$ is a natural factorisation for any $u\in U$.  Thus, we can take $U'=\im(\phi) \cong U$.
\epf

\begin{rem}\label{rem:embed2}
Theorem \ref{thm:embed1} (minus the clause concerning commutativity of $U$ and $\U$) is of course a special case of Theorem \ref{thm:embed2}, but the proofs we have given for these two results are completely different.  In particular, the monoids $\U$ constructed in the proofs are not at all alike.  One advantage of the proof of Theorem \ref{thm:embed1} is that the monoid~$\U$ turns out to be a semilattice when $M$ is a proper left restriction monoid; see the proof of Corollary~\ref{cor:embed1}.  
This is not the case, however, for the monoid $\U=\Mpres XR$ constructed in the current section, though it is worth noting that if $U$ is a band, then $\U$ is idempotent-generated.  Indeed, taking $v=u(=u^2)$ in \eqref{eq:suv}, we have $x_{\wh t,u} \sim x_{\wh t,u}\cdot x_{\wh t,u}$ for all $t\in S$ and $u\in U$.

Nevertheless, some modifications may be made in certain special cases, which lead to specialised versions of Theorem \ref{thm:embed2}.  Specifically, if we assume that $U$ is commutative, or a semilattice, or a left-regular band, then we can modify the construction to ensure that $\U$ has the same property (commutative, semilattice, left-regular band, respectively).  This will be the subject of the next two sections.
\end{rem}

\sect{Embedding theorems for proper monoids III}\label{sect:embed3}

The main result of the previous section was Theorem \ref{thm:embed2}, which shows that any $(U,S)$-proper monoid $M=US$ embeds naturally in a suitable semidirect product $\U\rtimes(S/\si)$.  In this section and the next we specialise this theorem in two important cases: when $U$ is commutative (but see Remark \ref{rem:embed3}), and when $U$ is a left-regular band, respectively.

\newpage

\begin{thm}\label{thm:embed3}
Let $M=US$ be a $(U,S)$-proper monoid, with $U$ commutative (respectively, a semilattice).  Then there exists a semidirect product $\M = \UC\rtimes (S/\si)$ and an embedding $\psi:M\to\M$ such that:
\ben
\item \label{embed30} $\UC$ is commutative (respectively, a semilattice), 
\item \label{embed31} $\UC$ contains a subsemigroup $U'$ isomorphic to $U$,
\item \label{embed32} $\psi\restr_U:U\to U'\rtimes\{\wh1\}$ is an isomorphism.
\een
\end{thm}

\begin{rem}\label{rem:embed3}
As explained in Remark \ref{rem:embed1}, if we removed the semilattice clauses from Theorem~\ref{thm:embed3} then the result would follow from (or is indeed contained in) Theorem \ref{thm:embed1}.  However, as explained in the same remark, the semilattice version of Theorem \ref{thm:embed3} does not follow from the proof of Theorem \ref{thm:embed1} given above.  The proof of Theorem \ref{thm:embed3} that we give below assumes only that $U$ is commutative, as it does not simplify at all if we assume that it is a semilattice.

There are two other reasons for including the coming proof of Theorem \ref{thm:embed3}.  The first is to show how the method of the previous section can be adapted to the commutative case; we add additional relations to the presentation $\U=\Mpres XR$ from the previous proof to obtain a commutative monoid $\UC=\Mpres X{R\cup C}$, which will be a semilattice if $U$ is.  The second additional reason is to pave the way for the next section, in which we carry out the same task in the more complex case that $U$ is a left-regular band; there we construct a suitable left-regular band $\UL=\Mpres X{R\cup L}$.
\end{rem}

For the rest of this section, we assume that $M=US$ is a $(U,S)$-proper monoid, with $U$ commutative.  We will not assume that $U$ is a semilattice, but will occasionally make a comment about the special case in which it is.  It is possible to prove results analogous to Proposition~\ref{prop:xi1} and Corollary~\ref{cor:phi}, demonstrating the `universality' of our construction of $\UC$.  These are omitted, however, as the statements and proofs are almost identical to those of the results just quoted.

The next result refers to the submonoid $P=\la S^+\ra\leq U$.  The proof is essentially the same as for Lemma \ref{lem:Psi}.  (Note that $P$ is commutative because~$U$ is.)

\begin{lemma}\label{lem:PsiC}
The monoid $P$ is a semilattice, and consequently $\si = \bigset{(s,t)\in S\times S}{t^+s = s^+t}$.  \epfres
\end{lemma}

As ever, we write
\[
\cS = S/\si = \set{\wh s}{s\in S}.
\]
The monoid $\UC$ will again be defined by means of a presentation $\Mpres X{R\cup C}$.  The alphabet~$X$ and the relations $R$ will have the same meaning as in the previous section, but for convenience we repeat the definitions here.  Specifically, we have
\[
X = \set{x_{\wh s,u}}{s\in S,\ u\in U},
\]
while $R$ consists of all the relations 
\begin{align}
\label{eq:RR} x_{\wh t,u\cdot{}^sv} &= x_{\wh t,u}x_{\wh t\wh s,v} &&\text{for $s,t\in S$ and $u,v\in U$ with $u=us^+$.}
\intertext{We also define $C$ to be the set of all relations}
\label{eq:CC} x_{\wh s,u}x_{\wh t,v} &= x_{\wh t,v}x_{\wh s,u} &&\text{for $s,t\in S$ and $u,v\in U$.}
\end{align}
For the duration of this section we write ${\sim}=(R\cup C)^\sharp$ for the congruence on the free monoid~$X^*$ generated by~$R\cup C$, and we denote the $\sim$-class of $w\in X^*$ by $\ol w$.  
We then define the monoid
\[
\UC = \Mpres X{R\cup C} = X^*/{\sim} = \set{\ol w}{w\in X^*}.
\]

\begin{lemma}\label{lem:UCcom}
\ben
\item The monoid $\UC$ is commutative.
\item If $U$ is a semilattice, then so too is $\UC$.
\een
\end{lemma}

\pf
This follows from Lemma \ref{lem:pres}.  (As in Remark \ref{rem:embed2}, when $U$ is a semilattice, $R$ contains the relations $x_{\wh t,u}= x_{\wh t,u}\cdot x_{\wh t,u}$ for all $t\in S$ and $u\in U$.)
\epf

As before, $\cS$ acts monoidally on $X^*$ by monoid morphisms via
\[
{}^{\wh s}w = x_{\wh s\wh t_1,u_1}\cdots x_{\wh s\wh t_k,u_k} \qquad\text{for $w=x_{\wh t_1,u_1}\cdots x_{\wh t_k,u_k}\in X^*$ and $s\in S$.}
\]
We noted previously that for any relation $(w,w')$ from~$R$, and any $s\in S$, $({}^{\wh s}w,{}^{\wh s}w')$ is again a relation from $R$.  This is also of course true if instead $(w,w')$ is a relation from $C$.  So we have an induced action of $\cS$ on $\UC$ given by
\[
{}^{\wh s} \ol w = \ol{{}^{\wh s}w} \qquad\text{for $w\in X^*$ and $s\in S$,}
\]
which is still monoidal, and by monoid morphisms.  We therefore have the (monoid) semidirect product
\[
\M = \UC\rtimes\cS.
\]

We now define a function
\[
\phiC:U\to\UC : u \mt \ol x_{\wh1,u}.
\]
Taking $s=t=1$ in \eqref{eq:RR}, we have
\[
x_{\wh1,uv} \sim x_{\wh1,u}x_{\wh1,v} \qquad\text{for all $u,v\in U$,}
\]
and it follows that $\phiC$ is a semigroup morphism.  We also have the following, proved in the same way as Lemma \ref{lem:phi}.

\begin{lemma}\label{lem:phiC}
If $u,v\in U$ and $s\in S$ are such that $u=us^+$, then
\[
(u\cdot{}^sv)\phiC = (u\phiC)\cdot{}^{\wh s}(v\phiC).  \epfreseq
\]
\end{lemma}

We have now defined the monoid $\UC$ and the semidirect product $\M=\UC\rtimes\cS$, and we have verified part \ref{embed30} of Theorem \ref{thm:embed3}.  The next result defines the map $\psi$ from the theorem.

\begin{prop}\label{prop:phiC1}
\ben
\item \label{phiC1} We have a morphism $\psi:M\to\M$ given by
\[
a\psi = (u\phiC,\wh s) \qquad\text{for any natural factorisation $(u,s)$ of $a\in M=US$.}
\]
\item \label{phiC2} $\psi:M\to\M$ is injective if and only if $\phiC:U\to\UC$ is injective.
\een
\end{prop}

\pf
\firstpfitem{\ref{phiC1}}  The map $\psi$ is well defined because of Lemma \ref{lem:nf} (cf.~Remark \ref{rem:nf}).  The proof of Proposition \ref{prop:xi2}\ref{xi21} is easily adapted to show that $\psi$ is a morphism.

\pfitem{\ref{phiC2}}  This is essentially the same as the proof of Proposition \ref{prop:xi2}\ref{xi22}.
\epf

Thus, to complete the proof of Theorem \ref{thm:embed3}, it remains to show that $\phiC$ is injective: i.e., that
\[
\ol x_{\wh1,u}=\ol x_{\wh1,v} \implies u=v \qquad\text{for all $u,v\in U$.}
\]
Indeed, it will then follow that the restriction
\[
\psi\restr_U:U\to\M:u\mt(u\phiC,\wh1)=(\ol x_{\wh1,u},\wh1)
\]
maps $U$ isomorphically onto its image, which is precisely $U'\rtimes\{\wh1\}$ for $U'=\set{\ol x_{\wh1,u}}{u\in U}\leq\UC$.

The proof in the previous section that $\phi: U\to\Mpres XR$ was injective was quite involved, relying on the concept of the \emph{trace} of certain words over $X$; see Definition \ref{defn:trace}.  Because of the commutativity of $\UC$, we can replace traces with the somewhat simpler notion of a \emph{shadow}:

\begin{defn}\label{defn:shadow}
For a word $w=x_{\wh s_1,u_1}\cdots x_{\wh s_k,u_k}\in X^*$, we define the function
\[
\vs_w : \wh s_1\times\cdots\times\wh s_k \to U \qquad\text{by}\qquad \vs_w(\bt) = {}^{t_1}u_1\cdots{}^{t_k}u_k \qquad\text{for $\bt=(t_1,\ldots,t_k)\in \wh s_1\times\cdots\times\wh s_k$.}
\]
(Note that it is convenient to write $\vs_w$ to the left of its argument, hence the notation $\vs_w(\bt)$ instead of $\bt\vs_w$.)
We call $\vs_w(\bt)$ a \emph{shadow} of $w$, and we say it is \emph{witnessed} by $\bt$.  We write
\[
\Si_w = \im(\vs_w) = \set{\vs_w(\bt)}{\bt\in\wh s_1\times\cdots\times\wh s_k}
\]
for the set of all shadows of $w$.
\end{defn}

\begin{rem}\label{rem:shadow}
When $k=0$, so that $w=\ew$ is empty, the domain of $\vs_\ew$ is $\wh s_1\times\cdots\times\wh s_k = \{\es\}$, by the standard set-theoretic interpretation of a cartesian product of an empty family.  By the usual convention regarding empty products we have $\vs_\ew(\es) = 1$.  Thus, $\Si_\ew=\{1\}$.

It is also worth noting a relationship between shadows and traces.  Recall from Definition~\ref{defn:trace} that the latter are defined for (some) words of the form $w=x_{\wh1,u_1}x_{\wh s_2,u_2}\cdots x_{\wh s_k,u_k}$.  A trace of such a word~$w$ is a certain special kind of tuple $\bp=(p_1,\ldots,p_k)\in U^k$, and we claim that if such a trace exists, then $p_k$ is a shadow of $w$.  Indeed, among the defining properties of the trace $\bp$, we have
\[
p_k=u_1\cdot{}^{t_2}u_2\cdots{}^{t_k}u_k \qquad\text{for some $(t_2,\ldots,t_k)\in\wh s_2\times\cdots\times\wh s_k$.}
\]
But then $p_k = {}^1u_1\cdot{}^{t_2}u_2\cdots{}^{t_k}u_k = \vs_w(\bt)$ for $\bt=(1,t_2,\ldots,t_k)\in\wh1\times\wh s_2\times\cdots\times\wh s_k$.

Finally, we also note that $P=\la S^+\ra$ is precisely the set of all shadows of all words of the form $w=x_{\wh s_1,1}\cdots x_{\wh s_k,1}$ ($s_1,\ldots,s_k\in S$).  Indeed, given any tuple ${\bt=(t_1,\ldots,t_k)\in \wh s_1\times\cdots\times\wh s_k}$, we have
\[
\vs_w(\bt) = {}^{t_1}1\cdots{}^{t_k}1 = t_1^+\cdots t_k^+\in P.
\]
Conversely, any element $z_1^+\cdots z_l^+$ of $P$ is a shadow of the word $x_{\wh z_1,1}\cdots x_{\wh z_l,1}$.
\end{rem}

Since $\Si_w\sub U$ for all $w\in X^*$, we can think of $\Si_w$ as an element of the power semigroup~$\P(U)$, as defined in \eqref{eq:PU}.  

\begin{lemma}\label{lem:sm}
The map $\xi_1:X^*\to\P(U):w\mt\Si_w$ is a monoid morphism.
\end{lemma}

\pf
We observed in Remark \ref{rem:shadow} that $\Si_\ew=\{1\}$, which is the identity of $\P(U)$.  So it remains to show that
\[
\Si_{ww'} = \Si_w\cdot\Si_{w'} \qquad\text{for all $w,w'\in X^*$.}
\]
But this follows quickly upon writing $w = x_{\wh s_1,u_1}\cdots x_{\wh s_k,u_k} $ and $ w' = x_{\wh s_1',u_1'}\cdots x_{\wh s_l',u_l'}$, and noting that:
\bit
\item any shadow of $w$ has the form ${}^{t_1}u_1\cdots{}^{t_k}u_k$ for some $(t_1,\ldots,t_k)\in \wh s_1\times\cdots\times\wh s_k$,
\item any shadow of $w'$ has the form ${}^{t_1'}u_1'\cdots{}^{t_l'}u_l'$ for some $(t_1',\ldots,t_l')\in \wh s_1'\times\cdots\times\wh s_l'$,
\item any shadow of $ww'$ has the form ${}^{t_1}u_1\cdots{}^{t_k}u_k\cdot {}^{t_1'}u_1'\cdots{}^{t_l'}u_l'$ for some $(t_1,\ldots,t_k)\in \wh s_1\times\cdots\times\wh s_k$ and $(t_1',\ldots,t_l')\in \wh s_1'\times\cdots\times\wh s_l'$.  \qedhere
\eit
\epf

Since $U$ is commutative, and since $P$ is a submonoid, we have a morphism
\[
\xi_2:\P(U)\to\P(U):V\mt P\cdot V.
\]
The image of this morphism is contained in $P\cdot \P(U) = P\cdot\P(U)\cdot P$, the local monoid of $\P(U)$ with identity $P$. (This local monoid played an important role in Section \ref{sect:embed1} as well, where it was denoted $I$.)  Of particular importance to us is the composite $\xi_1\circ\xi_2$, where $\xi_1$ is the morphism from Lemma \ref{lem:sm}.  We denote this composite by
\[
\xi = \xi_1\circ\xi_2 :X^*\to\P(U): w\mt P\cdot\Si_w.
\]

\begin{lemma}\label{lem:shadow3}
We have $R\cup C\sub\ker(\xi)$.
\end{lemma}

\pf
Consider some $(w,w')\in R\cup C$.  We must show that $w\xi=w'\xi$: i.e., that $P\cdot\Si_w=P\cdot\Si_{w'}$.  Since $P$ is a submonoid of $U$, it is enough to show that
\begin{equation}\label{eq:Siww'}
\Si_w\sub P\cdot\Si_{w'} \AND \Si_{w'}\sub P\cdot\Si_w.
\end{equation}
We consider separate cases according to whether $(w,w')$ belongs to $R$ or $C$, starting with the simpler of the two.

\pfcase1  Suppose first that $(w,w')\in C$, so as in \eqref{eq:CC} we have
\[
w = x_{\wh s,u}x_{\wh t,v} \AND w' = x_{\wh t,v}x_{\wh s,u} \qquad\text{for some $s,t\in S$ and $u,v\in U$.}
\]
By symmetry, it is enough to establish the first inclusion in \eqref{eq:Siww'}.  Since $1\in P$, we can do this by showing that
\[
\Si_w\sub\Si_{w'}.
\]
To do so, let $a\in\Si_w$ be an arbitrary shadow of $w$; we must show that $a\in\Si_{w'}$.  Since $w$ has length $2$, the shadow $a$ of $w$ is witnessed by a tuple of the form $\bz=(z_1,z_2)\in\wh s\times\wh t$, and we have $a = \vs_{w}(\bz) = {}^{z_1}u\cdot{}^{z_2}v$.  Since $U$ is commutative it follows that $a = {}^{z_2}v\cdot{}^{z_1}u = \vs_{w'}(\bz')$, where $\bz'=(z_2,z_1)\in\wh t\times\wh s$.  This shows that $a\in\Si_{w'}$, as required.

\pfcase2  Next, suppose $(w,w')\in R$, so as in \eqref{eq:RR} we have
\[
w = x_{\wh t,u\cdot{}^sv} \AND w' = x_{\wh t,u}x_{\wh t\wh s,v} \qquad\text{for some $s,t\in S$ and $u,v\in U$ with $u=us^+$.}
\]
This time we do not have symmetry, so we must demonstrate both inclusions in \eqref{eq:Siww'}.

To establish the first of these inclusions, we will again show that
\[
\Si_w\sub\Si_{w'}.
\]
To prove this, let $a\in\Si_w$.  This time $a$ is witnessed by a tuple of the form $\bz=(z)$, meaning that $z\in\wh t$ and $a=\vs_{w}(\bz) = {}^z(u\cdot{}^sv)$.  Since $\si$ is a congruence, we have
\[
z\in\wh t \implies z\mr\si t \implies zs \mr\si ts \implies zs\in\wh t\wh s,
\]
so it follows that $\bz' = (z,zs) \in \wh t\times \wh t\wh s$.  But then $w'$ has shadow
\[
\vs_{w'}(\bz') = {}^zu \cdot {}^{zs}v = {}^zu \cdot{}^z({}^sv) = {}^z(u\cdot{}^sv) = a,
\]
showing that $a\in\Si_w'$.

To demonstrate the second inclusion in \eqref{eq:Siww'}, let $b\in\Si_{w'}$ be an arbitrary shadow of $w'$.  We must show that
\begin{equation}\label{eq:apb}
b=pc \qquad\text{for some $p\in P$ and some shadow $c\in\Si_{w}$.}
\end{equation}
The shadow $b$ of $w'$ is witnessed by some tuple $\bz=(z_1,z_2)$.  This means that $z_1\in\wh t$, $z_2\in\wh t\wh s$ and 
\begin{align}
\label{eq:az1uz2v}
b=\vs_{w'}(\bz)&={}^{z_1}u\cdot{}^{z_2}v.
\intertext{Now we put $\bz'=(z_1)$.  Since $z_1\in\wh t$, it follows that }
\label{eq:vswbz}
\vs_{w}(\bz') &= {}^{z_1}(u\cdot{}^sv)
\intertext{is a shadow of $w$.  From $z_1\in\wh t$ and $z_2\in\wh t\wh s$ we deduce that $z_2 \mr\si ts \mr\si z_1s$, so it follows from Lemma \ref{lem:PsiC} that}
\nonumber z_2^+\cdot z_1s &= (z_1s)^+\cdot z_2.
\intertext{Since $z_2^+,(z_1s)^+\in U$ and $z_1s,z_2\in S$ (and $v\in U$), Lemma \ref{lem:ap} then gives}
\label{eq:z1s+}
z_2^+\cdot {}^{z_1s}v &= (z_1s)^+\cdot {}^{z_2}v.
\end{align}
Putting everything together, and using properties of the action, we then have
\begin{align*}
z_2^+\cdot \vs_{w}(\bz') &= z_2^+\cdot {}^{z_1}(u\cdot{}^sv) &&\text{by \eqref{eq:vswbz}}\\
&= z_2^+\cdot {}^{z_1}u\cdot{}^{z_1s}v \\
&= {}^{z_1}u\cdot z_2^+\cdot {}^{z_1s}v &&\text{by commutativity of $U$}\\
&= {}^{z_1}u\cdot (z_1s)^+\cdot {}^{z_2}v &&\text{by \eqref{eq:z1s+}}\\
&= {}^{z_1}u\cdot {}^{z_1}(s^+)\cdot {}^{z_2}v &&\text{by Lemma \ref{lem:+}\ref{lem+2}}\\
&= {}^{z_1}(u s^+)\cdot {}^{z_2}v \\
&= {}^{z_1}u\cdot {}^{z_2}v &&\text{as $u=us^+$}\\
&= b &&\text{by \eqref{eq:az1uz2v}.}
\end{align*}
So \eqref{eq:apb} holds with $p=z_2^+\in P$ and $c=\vs_{w}(\bz')\in\Si_{w}$.
\epf

Because of Lemma \ref{lem:shadow3}, it follows that we have a well-defined morphism
\[
\Xi : \UC = \Mpres{X}{R\cup C} \to \P(U) \qquad\text{given by} \qquad \ol w\Xi = w\xi = P\cdot\Si_w \qquad\text{for $w\in X^*$.}
\]
We can now prove the following:

\begin{prop}\label{prop:phiC}
The morphism $\phiC:U\to\UC:u\mt\ol x_{\wh1,u}$ is an embedding.
\end{prop}

\pf
We prove the proposition by showing that the composite
\[
\phiC\circ\Xi:U\to\P(U):u\mt P\cdot\Si_{x_{\wh1,u}}
\]
is injective.  To do so, fix some ${(u,v)\in\ker(\phiC\circ\Xi)}$, meaning that $u,v\in U$ and
\[
P\cdot\Si_{x_{\wh1,u}} = P\cdot\Si_{x_{\wh1,v}}.
\]
Now, $u={}^1u$ is a shadow of $x_{\wh1,u}$ (witnessed by $(1)$), and so $u=1\cdot u\in P\cdot\Si_{x_{\wh1,u}} = P\cdot\Si_{x_{\wh1,v}}$.  It follows that $u=pa$ for some $p\in P$ and some shadow $a$ of $x_{\wh1,v}$.  By definition, this shadow has the form $a={}^tv$ for some $t\in\wh 1$.  By Lemma \ref{lem:wh1}\ref{wh12} we have $a={}^tv = t^+v$, and so $u = pa = pt^+v \in Pv$, as $pt^+\in P$.  By symmetry we also have $v\in Pu$.  It then follows from Lemma~\ref{lem:sl} that $u=v$, and the proof is complete.
\epf

As we have already observed, this completes the proof of Theorem \ref{thm:embed3}.

\begin{rem}
We noted above that the power semigroup $\P(U)$ played an important role in Section \ref{sect:embed1} as well, as did its local monoid $P\cdot\P(U)\cdot P$.  Comparing the above proof of Proposition~\ref{prop:phiC} with that of Lemma \ref{lem:phi_mono}, there is a sense in which we have come full circle.  Indeed, for $s\in S$ and $u\in U$, the shadow set
\[
\Si_{x_{\wh s,u}} = \set{{}^tu}{t\in\wh s}
\]
is precisely the set $V_{\wh s,u}$ constructed in Section \ref{sect:embed1}; see \eqref{eq:Vsu}.  In the proof of Lemma \ref{lem:phi_mono}, the injectivity of the map $u\mt f_u$ defined in \eqref{eq:ufu} boiled down to the implication
\[
\wh1f_u = \wh1f_v\implies u=v \qquad\text{for $u,v\in U$.}
\]
This is of course equivalent to the injectivity of the map
\begin{equation}\label{eq:UPU}
U \to \P(U) : u\mt \wh1f_u = V_{\wh1,u}\cdot P.
\end{equation}
Since
\[
V_{\wh1,u}\cdot P = P\cdot V_{\wh1,u} = P\cdot\Si_{x_{\wh1,u}} = u (\phiC\circ\Xi),
\]
it follows that the map in \eqref{eq:UPU} is precisely $\phiC\circ\Xi$ from the proof of Proposition \ref{prop:phiC}.  (But recall that we did not assume $U$ itself was commutative in Section \ref{sect:embed1}, only that $P$ was central in $U$.)

On the other hand, it does not seem possible to prove Theorem \ref{thm:embed4} (the main result of the next section, which concerns left-regular bands) using techniques similar to those of Sections \ref{sect:embed1} and \ref{sect:embed3}.
\end{rem}

\sect{Embedding theorems for proper monoids IV}\label{sect:embed4}

As already discussed, our main goal in this final section of the chapter is to prove the following:

\begin{thm}\label{thm:embed4}
Let $M=US$ be a $(U,S)$-proper monoid, with $U$ a left-regular band.  Then there exists a semidirect product $\M = \UL\rtimes (S/\si)$ and an embedding $\psi:M\to\M$ such that:
\ben
\item \label{embed40} $\UL$ is a left-regular band,  
\item \label{embed41} $\UL$ contains a subsemigroup $U'$ isomorphic to $U$,
\item \label{embed42} $\psi\restr_U:U\to U'\rtimes\{\wh1\}$ is an isomorphism.
\een
\end{thm}

Working towards the proof, for the rest of the section we fix a $(U,S)$-proper monoid $M=US$, and we assume that $U$ is a left-regular band.  As usual we write $\cS=S/\si=\set{\wh s}{s\in S}$.  Since~$U$ is left-regular, so too is its submonoid $P=\la S^+\ra$.  It then follows from Lemma~\ref{lem:si3} (cf.~Remark~\ref{rem:si3}) that
\begin{equation}\label{eq:siLRB}
\si = \bigset{(s,t)\in S\times S}{s^+t^+s=s^+t}.
\end{equation}

This time we define $\UL = \Mpres X{R\cup L}$, where again $X = \set{x_{\wh s,u}}{s\in S,\ u\in U}$, $R$ consists of all the relations 
\begin{align}
\label{eq:RRR} 
x_{\wh t,u\cdot{}^sv} &= x_{\wh t,u}x_{\wh t\wh s,v} &&\text{for $s,t\in S$ and $u,v\in U$ with $u=us^+$,}
\intertext{and additionally $L$ is the set of all relations}
\label{eq:LLL} 
x_{\wh s,u}x_{\wh t,v}x_{\wh s,u} &= x_{\wh s,u}x_{\wh t,v} &&\text{for $s,t\in S$ and $u,v\in U$.}
\end{align}
By Lemma \ref{lem:pres}\ref{pres3}, and remembering that $R$ contains the relations $x_{\wh t,u}=x_{\wh t,u}\cdot x_{\wh t,u}$ for each $t\in S$ and $u\in U$ (cf.~Remark~\ref{rem:embed2}), $\UL$ is a left-regular band.  We write ${\sim}=(R\cup L)^\sharp$, and denote the $\sim$-class of $w\in X^*$ by $\ol w$, so that
\[
\UL = \Mpres X{R\cup L} = X^*/{\sim} = \set{\ol w}{w\in X^*}.
\]
As ever, the action of $\cS$ on $X^*$ given in \eqref{eq:SactX*} induces a monoidal action of $\cS$ on $\UL$ by monoid morphisms; indeed, we just need to observe that $(w,w')\in L \implies ({}^{\wh s}w,{}^{\wh s}w')\in L$ for all $s\in S$.  This then allows us to define the semidirect product
\[
\M=\UL\rtimes\cS.
\]

We define the map
\[
\phiL:U\to\UL:u\mt\ol x_{\wh 1,u}.
\]
As in Lemma \ref{lem:phiC}, we have
\[
(u\cdot{}^sv)\phiL = (u\phiL)\cdot{}^{\wh s}(v\phiL) \qquad\text{for $u,v\in U$ and $s\in S$ with $u=us^+$.}
\]
As in Proposition \ref{prop:phiC1}, this allows us to define a morphism
\[
\psi:M\to\M \qquad\text{by}\qquad a\psi = (u\phiL,\wh s) \qquad\text{for any natural factorisation $(u,s)$ of $a\in M=US$,}
\]
which is injective if and only if $\phiL$ is.

Thus, we can complete the proof of Theorem \ref{thm:embed4} by showing that $\phiL$ is injective.  The simpler method involving shadows from the previous section does not apply here, as several arguments there relied crucially on commutativity of $U$ (see especially the proof of Lemma~\ref{lem:shadow3}).  Thus, we will once again use the \emph{traces} from Definition \ref{defn:trace}.  Recall that these are defined for (some) words from the set $W$, which consists of all non-empty words over $X$ whose first letter has the form $x_{\wh 1,u}$ for some $u\in U$.  Examining \eqref{eq:RRR} and \eqref{eq:LLL}, note that $W$ is still closed under ${\sim}=(R\cup L)^\sharp$, in the sense that a pair of $\sim$-equivalent words either both belong to $W$ or else both do not (cf.~Lemma~\ref{lem:lam}).  By Lemma \ref{lem:trace1}, any word from $W$ has at most one trace.  The main technical result we need is the following version of Lemma \ref{lem:trace2}:

\begin{lemma}\label{lem:trace3}
Suppose $w,w'\in W$ are such that $w\sim w'$.
\ben
\item \label{trace31} If either of $w$ or $w'$ has a trace, then both do.
\item \label{trace32} If $w$ and $w'$ have traces $\bp=(p_1,\ldots,p_k)$ and $\bq=(q_1,\ldots,q_l)$, then $p_k=q_l$.
\een
\end{lemma}

\pf
As with Lemma \ref{lem:trace2}, it suffices to assume that $w$ and $w'$ differ by a single application of a relation from $R\cup C$.  Lemma \ref{lem:trace2} has already taken care of the case that the relation is from~$R$, so for the rest of the proof we assume that the relation is from $C$.  Thus, up to symmetry we have
\[
w = w_1 \cdot x_{\wh s,u}x_{\wh t,v} \cdot w_2 \AND w' = w_1\cdot x_{\wh s,u}x_{\wh t,v}x_{\wh s,u} \cdot w_2,
\]
for some $w_1,w_2\in X^*$, $s,t\in S$ and $u,v\in U$.  We also write
\[
w_1 = x_{\wh s_1,u_1}\cdots x_{\wh s_m,u_m} \AND w_2 = x_{\wh s_{m+3},u_{m+3}}\cdots x_{\wh s_k,u_k} \qquad\text{where each $s_i\in S$ and $u_i\in U$,}
\]
noting that either (or both) of $w_1$ or $w_2$ could be empty.  It is also convenient to define
\[
s_{m+1} = s \COMMA u_{m+1} = u \COMMA s_{m+2}=t \AND u_{m+2}=v.
\]
So then
\begin{align*}
w &= x_{\wh s_1,u_1}\cdots x_{\wh s_k,u_k} \\
&= x_{\wh s_1,u_1}\cdots x_{\wh s_m,u_m} \cdot x_{\wh s_{m+1},u_{m+1}}x_{\wh s_{m+2},u_{m+2}} \cdot x_{\wh s_{m+3},u_{m+3}}\cdots x_{\wh s_k,u_k},
\intertext{and}
w' &= x_{\wh s_1,u_1}\cdots x_{\wh s_m,u_m} \cdot x_{\wh s,u}x_{\wh t,v}x_{\wh s,u} \cdot x_{\wh s_{m+3},u_{m+3}}\cdots x_{\wh s_k,u_k} \\
&= x_{\wh s_1,u_1}\cdots x_{\wh s_m,u_m} \cdot x_{\wh s_{m+1},u_{m+1}}x_{\wh s_{m+2},u_{m+2}}x_{\wh s_{m+1},u_{m+1}} \cdot x_{\wh s_{m+3},u_{m+3}}\cdots x_{\wh s_k,u_k}.
\end{align*}
Since $w\in W$, we must have $\wh s_1=\wh1$; in the case that $w_1=\ew$ (i.e., $m=0$), this says that $\wh s=\wh1$.  We also define
\begin{align}
\nonumber (v_1,\ldots,v_{k+1}) &= (u_1,u_2,\ldots,u_m,u,v,u,u_{m+3},\ldots,u_k) \\
\label{eq:vv} &= (u_1,u_2,\ldots,u_m,u_{m+1},u_{m+2},u_{m+1},u_{m+3},\ldots,u_k).
\end{align}
Note that $v_i$ is the `$U$-coordinate' of the $i$th letter of $w'$.

We begin with the forwards implication.

\pfitem{($\Rightarrow$)}  Suppose first that $w$ has trace $\bp=(p_1,\ldots,p_k)$, witnessed by $\bt=(t_2,\ldots,t_k)\in\wh s_2\times\cdots\times\wh s_k$, meaning that
\bit
\item $p_i = u_1\cdot{}^{t_2}u_2\cdots {}^{t_i}u_i$ for all $1\leq i\leq k$, and
\item $p_{i-1}=p_{i-1}t_i^+$ for all $2\leq i\leq k$.
\eit
We claim that $w'$ has trace $\bq$, witnessed by $\bz$, for
\begin{align*}
\bz &= (t_2,\ldots,t_m,t_{m+1},t_{m+2},t_{m+1},t_{m+3},\ldots,t_k)\\
\text{and}\qquad
\bq &= (p_1,\ldots,p_m,p_{m+1},p_{m+2},p_{m+2},p_{m+3},\ldots,p_k).
\end{align*}
Note that $\bz$ and $\bq$ are obtained from $\bt$ and $\bp$, respectively, by inserting single entries.  Since the final entries of $\bp$ and $\bq$ are equal, this will complete the proof of the forwards implication.  To assist with the following calculations, it is convenient to line up the entries of $\bz$ and $\bq$, along with the letters of $w'$:
\begin{center}\setlength{\tabcolsep}{4.5pt}
\begin{tabular}{c|cccccccccc}
&$1$&$2$&$\cdots$&$m$&$m+1$&$m+2$&$m+3$&$m+4$&$\cdots$&$k+1$\\
\hline
$\bz$&&$t_2$&$\cdots$&$t_m$&$t_{m+1}$&$t_{m+2}$&$t_{m+1}$&$t_{m+3}$&$\cdots$&$t_k$\\
$\bq$&$p_1$&$p_2$&$\cdots$&$p_m$&$p_{m+1}$&$p_{m+2}$&$p_{m+2}$&$p_{m+3}$&$\cdots$&$p_k$\\
$w'$&$x_{\wh s_1,u_1}$&$x_{\wh s_2,u_2}$&$\cdots$&$ x_{\wh s_m,u_m} $&$ x_{\wh s,u}$&$x_{\wh t,v} $&$ x_{\wh s,u}$&$ x_{\wh s_{m+3},u_{m+3}}$&$\cdots $&$x_{\wh s_k,u_k}$\\
\multicolumn{5}{c}{}&$\rotatebox[origin=c]{90}{=}$&$\rotatebox[origin=c]{90}{=}$&$\rotatebox[origin=c]{90}{=}$\\
\multicolumn{5}{c}{}&$x_{\wh s_{m+1},u_{m+1}}$&$x_{\wh s_{m+2},u_{m+2}}$&$x_{\wh s_{m+1},u_{m+1}}$
\end{tabular}
\end{center}
To prove the claim, we follow Definition \ref{defn:trace} carefully, and we proceed in three steps.

\pfstep1  First, it is clear that $\bz\in\wh s_2\times\cdots\times\wh s_m \times \wh s\times \wh t\times \wh s \times \wh s_{m+3}\times\cdots\times\wh s_k$.  For this, note that $t_{m+1}\in\wh s_{m+1}=\wh s$ and $t_{m+2}\in\wh s_{m+2}=\wh t$.

\pfstep2  Next we check that $q_i = v_1\cdot {}^{z_2}v_2\cdots{}^{z_i}v_i$ for all $1\leq i\leq k+1$, where the $v_j$ are as in~\eqref{eq:vv}.
\bit
\item For $1\leq i\leq m+2$ we have $q_i = p_i = u_1\cdot{}^{t_2}u_2\cdots {}^{t_i}u_i = v_1\cdot{}^{z_2}v_2\cdots {}^{z_i}v_i$.
\item For $m+3\leq i\leq k+1$, and remembering that $U$ is a left-regular band, we have 
\begin{align*}
q_i = p_{i-1} &= (u_1\cdot{}^{t_2}u_2\cdots {}^{t_m}u_m)\cdot {}^{t_{m+1}}u_{m+1}{}^{t_{m+2}}u_{m+2} \cdot ({}^{t_{m+3}}u_{m+3}\cdots {}^{t_{i-1}}u_{i-1})\\
&= (u_1\cdot{}^{t_2}u_2\cdots {}^{t_m}u_m)\cdot {}^{t_{m+1}}u_{m+1}{}^{t_{m+2}}u_{m+2} {}^{t_{m+1}}u_{m+1} \cdot ({}^{t_{m+3}}u_{m+3}\cdots {}^{t_{i-1}}u_{i-1})\\
&= (v_1\cdot{}^{z_2}v_2\cdots {}^{z_m}v_m)\cdot {}^{z_{m+1}}v_{m+1}\cdot {}^{z_{m+2}}v_{m+2}{}^{z_{m+3}}v_{m+3} \cdot ({}^{z_{m+4}}v_{m+4}\cdots {}^{z_i}v_i).
\end{align*}
\eit

\pfstep3
Finally, we check that $q_{i-1}=q_{i-1}z_i^+$ for all $2\leq i\leq k+1$.
\bit
\item For $2\leq i\leq m+2$ we have $q_{i-1} = p_{i-1} = p_{i-1}t_i^+ = q_{i-1}z_i^+$.
\item For $i=m+3$ we have
\begin{align*}
q_{m+2} = p_{m+2} &= p_m\cdot {}^{t_{m+1}}u_{m+1}{}^{t_{m+2}}u_{m+2} \\
&= p_m\cdot {}^{t_{m+1}}u_{m+1}{}^{t_{m+2}}u_{m+2}{}^{t_{m+1}}u_{m+1} &&\text{by left-regularity}\\
&= p_m\cdot {}^{t_{m+1}}u_{m+1}{}^{t_{m+2}}u_{m+2}{}^{t_{m+1}}u_{m+1}\cdot t_{m+1}^+ &&\text{by Lemma \ref{lem:+}\ref{lem+4}}\\
&= p_m\cdot {}^{t_{m+1}}u_{m+1}{}^{t_{m+2}}u_{m+2}\cdot t_{m+1}^+ &&\text{by left-regularity}\\
&= p_{m+2}t_{m+1}^+ = q_{m+2}z_{m+3}^+.
\end{align*}
\item For $m+4\leq i\leq k+1$ we have $q_{i-1} = p_{i-2} = p_{i-2}t_{i-1}^+ = q_{i-1}z_i^+$.
\eit

\pfitem{($\Leftarrow$)}  Conversely, we now assume $w'$ has trace $\bq=(q_1,\ldots,q_{k+1})$, witnessed by
\begin{equation}\label{eq:zz}
\bz=(z_2,\ldots,z_{k+1})\in\wh s_2\times\cdots\times\wh s_m \times \wh s \times \wh t\times \wh s \times \wh s_{m+3}\times\cdots\times\wh s_k,
\end{equation}
meaning that
\bit
\item $q_i = v_1\cdot{}^{z_2}v_2\cdots {}^{z_i}v_i$ for all $1\leq i\leq k+1$, where again the $v_j$ are as in \eqref{eq:vv}, and
\item $q_{i-1}=q_{i-1}z_i^+$ for all $2\leq i\leq k+1$.
\eit
We will shortly construct a suitable trace for $w$.  However, in order to deal with a subtle point that did not arise in the proof of Lemma \ref{lem:trace2}, we delay the construction, and instead begin by showing that
\begin{equation}\label{eq:qm2=qm3}
q_{m+2}=q_{m+3}.
\end{equation}
Working towards this, and defining $z_1=1$ for convenience, we first claim that 
\begin{equation}\label{eq:qizj+}
q_i = q_i z_j^+ \qquad\text{for all $2\leq i\leq k$ and $1\leq j\leq i+1$.}
\end{equation}
Indeed, this is is obvious for $j=1$, and is true for $j=i+1$ by the properties of the tuples $\bq$ and $\bz$ (stated above).  For $2\leq j\leq i$ we have
\begin{align*}
q_i &= q_{j-1}\cdot {}^{z_j}v_j \cdot {}^{z_{j+1}}{v_{j+1}}\cdots {}^{z_i}v_i \\
&= q_{j-1}\cdot {}^{z_j}v_j\cdot z_j^+ \cdot {}^{z_{j+1}}{v_{j+1}}\cdots {}^{z_i}v_i &&\text{by Lemma \ref{lem:+}\ref{lem+4}}\\
&= q_{j-1}\cdot {}^{z_j}v_j\cdot z_j^+ \cdot {}^{z_{j+1}}{v_{j+1}}\cdots {}^{z_i}v_i \cdot z_j^+&&\text{by left-regularity}\\
&= q_{j-1}\cdot {}^{z_j}v_j \cdot {}^{z_{j+1}}{v_{j+1}}\cdots {}^{z_i}v_i \cdot z_j^+&&\text{by Lemma \ref{lem:+}\ref{lem+4} again}\\
&= q_iz_j^+.
\end{align*}
Next we note that
\begin{equation}\label{eq:zm1sizm3}
z_{m+1}\mr\si z_{m+3}.
\end{equation}
Indeed, for $m\geq1$, this follows from \eqref{eq:zz}, as then $z_{m+1},z_{m+3}\in\wh s$.  When $m=0$ (which occurs when $w_1=\ew$), \eqref{eq:zm1sizm3} says that $1=z_1\mr\si z_3$.  But we observed just before \eqref{eq:vv} that $\wh s=\wh 1$ when $m=0$, so \eqref{eq:zm1sizm3} still holds in this case as $z_3=z_{m+3}\in\wh s=\wh1=\wh z_1$.  Now that we have completed the proof of \eqref{eq:zm1sizm3}, it follows from \eqref{eq:siLRB} that
\[
z_{m+1}^+z_{m+3}^+z_{m+1} = z_{m+1}^+z_{m+3}.
\]
Lemma \ref{lem:ap} then gives
\begin{equation}\label{eq:zm1zm3}
z_{m+1}^+z_{m+3}^+\cdot {}^{z_{m+1}}u = z_{m+1}^+\cdot{}^{z_{m+3}}u.
\end{equation}
Keeping $(v_{m+1},v_{m+2},v_{m+3}) = (u,v,u)$ in mind (and remembering that $z_1=1$, which is relevant for the $m=0$ case), we then calculate
\begin{align*}
q_{m+3} = q_{m+2}\cdot{}^{z_{m+3}}u &= q_{m+2}\cdot z_{m+1}^+\cdot{}^{z_{m+3}}u &&\text{by \eqref{eq:qizj+}}\\
&= q_{m+2}\cdot z_{m+1}^+z_{m+3}^+\cdot {}^{z_{m+1}}u &&\text{by \eqref{eq:zm1zm3}}\\
&= q_{m+2}\cdot {}^{z_{m+1}}u &&\text{by \eqref{eq:qizj+}}\\
&= q_m\cdot {}^{z_{m+1}}u\cdot{}^{z_{m+2}}v\cdot {}^{z_{m+1}}u \\
&= q_m\cdot {}^{z_{m+1}}u\cdot{}^{z_{m+2}}v &&\text{by left-regularity}\\
&= q_{m+2},
\end{align*}
completing the proof of \eqref{eq:qm2=qm3}.

We now claim that $w$ has trace $\bp$, witnessed by $\bt$, where
\[
\bt = (z_2,\ldots,z_{m+2},z_{m+4},\ldots,z_{k+1}) \AND \bp = (q_1,\ldots,q_{m+2},q_{m+4},\ldots,q_{k+1}).
\]
This time $\bt$ and $\bp$ are obtained by deleting an appropriate entry of $\bz$ and $\bq$.
However, unlike the situation in Lemma \ref{lem:trace2}, it is possible that $\bp$ is obtained by deleting the \emph{final} entry of~$\bq$.  This occurs precisely when $m+3=k+1$ (i.e., when $w_2=\ew$).  In this case, we have ${\bq=(q_1,\ldots,q_{m+2},q_{m+3})}$ and $\bp=(q_1,\ldots,q_{m+2})$, and these still have the same final entries because of \eqref{eq:qm2=qm3}.

Thus, we can complete the proof of the lemma by showing that $w$ does indeed have trace $\bp$, witnessed by $\bt$.  
This time the entries of $\bt$ and $\bp$, and the letters from $w$, are as follows:
\begin{center}
\begin{tabular}{c|ccccccccc}
&$1$&$2$&$\cdots$&$m$&$m+1$&$m+2$&$m+3$&$\cdots$&$k$\\
\hline
$\bt$&&$z_2$&$\cdots$&$z_m$&$z_{m+1}$&$z_{m+2}$&$z_{m+4}$&$\cdots$&$z_{k+1}$\\
$\bp$&$q_1$&$q_2$&$\cdots$&$q_m$&$q_{m+1}$&$q_{m+2}$&$q_{m+4}$&$\cdots$&$q_{k+1}$\\
$w$&$x_{\wh s_1,u_1}$&$x_{\wh s_2,u_2}$&$\cdots$&$ x_{\wh s_m,u_m} $&$ x_{\wh s,u}$&$ x_{\wh t,v}$&$ x_{\wh s_{m+3},u_{m+3}}$&$\cdots $&$x_{\wh s_k,u_k}$\\
\multicolumn{5}{c}{}&$\rotatebox[origin=c]{90}{=}$&$\rotatebox[origin=c]{90}{=}$\\
\multicolumn{5}{c}{}&$x_{\wh s_{m+1},u_{m+1}}$&$x_{\wh s_{m+2},u_{m+2}}$
\end{tabular}
\end{center}
As ever, we proceed in three steps.

\pfstep1
We first note that $\bt\in\wh s_2\times\cdots\times\wh s_m\times\wh s\times\wh t\times\wh s_{m+2}\times\cdots\times\wh s_k$ follows immediately from~\eqref{eq:zz}.

\pfstep2
Next we check that $p_i = u_1\cdot {}^{z_2}u_2\cdots{}^{z_i}u_i$ for all $1\leq i\leq k$.
\bit
\item For $1\leq i\leq m+2$ we have $p_i=q_i=v_1\cdot{}^{z_2}v_2\cdots{}^{z_i}v_i=u_1\cdot{}^{t_2}u_2\cdots{}^{t_i}u_i$.
\item For $m+3\leq i\leq k$ we have
\begin{align*}
p_i = q_{i+1} &= q_{m+3} \cdot {}^{z_{m+4}}v_{m+4}\cdots {}^{z_{i+1}}v_{i+1} \\
&= q_{m+2} \cdot {}^{z_{m+4}}v_{m+4}\cdots {}^{z_{i+1}}v_{i+1} &&\text{by \eqref{eq:qm2=qm3}}\\
&= v_1\cdot{}^{z_2}v_2\cdots{}^{z_{m+2}}v_{m+2} \cdot {}^{z_{m+4}}v_{m+4}\cdots {}^{z_{i+1}}v_{i+1} \\
&= u_1\cdot{}^{t_2}u_2\cdots{}^{t_{m+2}}u_{m+2} \cdot {}^{t_{m+3}}u_{m+3}\cdots {}^{t_i}u_i ,
\end{align*}
as required.
\eit

\pfstep3
Finally, we check that $p_{i-1}=p_{i-1}t_i^+$ for all $2\leq i\leq k$.
\bit
\item For $2\leq i\leq m+2$ we have $p_{i-1}=q_{i-1}=q_{i-1}z_i^+=p_{i-1}t_i^+$.
\item For $i=m+3$, we use \eqref{eq:qm2=qm3} to calculate
\[
p_{m+2} = q_{m+2} = q_{m+3} = q_{m+3}z_{m+4}^+ = p_{m+2}t_{m+3}^+.
\]
\item For $m+4\leq i\leq k$ we have $p_{i-1}=q_i=q_iz_{i+1}^+=p_{i-1}t_i^+$.  
\eit
We have now finally completed the proof.
\epf

As in Section \ref{sect:embed3} (cf.~Proposition \ref{prop:phi}) the injectivity of $\phiL:U\to\UL$ follows quickly from Lemma \ref{lem:trace3}.  As observed above, this completes the proof of Theorem \ref{thm:embed4}.

\begin{rem}
Theorem \ref{thm:embed3} (proved in the previous section) concerns two special cases:
\bit
\item $U$ is commutative, meaning that it satisfies the identity $xy=yx$,
\item $U$ is a semilattice, meaning that it satisfies the identities $xy=yx$ and $x^2=x$.
\eit
Similarly, Theorem \ref{thm:embed4} concerns the special case in which:
\bit
\item $U$ is a left-regular band, meaning that it satisfies the identities $xyx=xy$ and $x^2=x$.  
\eit
The reader might wonder why we did not give a statement for the case that $U$ is assumed only to satisfy $xyx=xy$.  The simple reason for this is that $U$ being a monoid means that the identity $xyx=xy$ implies $x^2=x$ upon substituting $y=1$.  On the other hand, the identities $xy=yx$ and $x^2=x$ are independent of each other.

We leave it as an open problem to classify the (collections of) identities on $U$ that lead to analogous embedding results.  For example, if $M=US$ is a  $(U,S)$-proper monoid, with $U$ a band, does there exist an embedding $M\to\U_B\rtimes\cS$, with $\U_B$ a band containing a copy of $U$?

Even in the case of a \emph{right}-regular band (satisfying $xyx=yx$ and $x^2=x$), it is not clear whether such an embedding will exist.  If it does, then it is likely that a very different proof strategy would be required.  Indeed, the method of Sections \ref{sect:embed1}, \ref{sect:embed3} and \ref{sect:embed4} relied crucially on~$P$ satisfying an identity of the form $f(x,y)x = g(x,y)y$, so that Lemma \ref{lem:si3} applied and gave an equational formulation of the congruence $\si$.  But neither of the defining identities for right-regular bands are of this form.
\end{rem}

\newpage

\part{Presentations}\label{part:II}

This part of the paper explicitly focusses on presentations by generators and relations.  
Since so many naturally occuring semigroups arise from action pairs, our purpose here is to develop a body of general results concerning presentations for such semigroups.  We then apply these to a number of important examples.

Chapter \ref{chap:pres} contains the main theoretical results of this part of the paper, each of which gives a presentation for a semigroup arising from an action pair, modulo various natural assumptions.  Roughly speaking, stricter assumptions lead to smaller sets of relations.  In particular, the results and/or proofs tend to be simpler when the semigroup in question is in fact a monoid.  \mbox{Chapters~\ref{chap:LR}--\ref{chap:wreath}} then apply the general machinery developed in Chapter \ref{chap:pres} to several examples.  Specifically, we treat:
\bit
\item free left restriction monoids in Chapter \ref{chap:LR}, 
\item several monoids and semigroups of (partial) endomorphisms of an independence algebra in Chapter \ref{chap:IA}, and 
\item a number of transformational wreath products in Chapter \ref{chap:wreath}.  
\eit
The introductions to these chapters contain full summaries of the results they contain, so we now turn to the general results promised above.

\chap{Presentations}\label{chap:pres}

As just discussed, the purpose of this chapter is prove a number of general results on presentations.  Each such result concerns a semigroup of the form $US$ arising from an action pair $(U,S)$ in a monoid~$M$.  (In fact, some results hold more generally for \emph{weak} action pairs.)  In an ideal world, one would hope that a presentation for $US$ could be `pieced together' from presentations for $U$ and~$S$, in the following sense:
\bit
\item We start with presentations $\pres{X_U}{R_U}$ and $\pres{X_S}{R_S}$ for~$U$ and $S$.
\item We then hope to obtain a presentation for $US$ of the form $\pres{X_U\cup X_S}{R_U\cup R_S\cup R}$ for some additional set of relations $R$ over $X_U\cup X_S$.
\eit
This is indeed possible in some cases, a number of which we treat in the current chapter.  
\bit
\item Sections~\ref{sect:2M1}--\ref{sect:2M3} deal with the case in which $U$ and $S$ are both submonoids of $M$, the main results here being Theorems~\ref{thm:ESmon0}, \ref{thm:ESmon}, \ref{thm:Msimp}, \ref{thm:M} and \ref{thm:ESmon00}.  In Remark \ref{rem:ESmon0} we apply the results of Section~\ref{sect:cong} to give especially neat presentations when certain natural conditions hold, and in particular we are able to deduce the main results of \cite{EEF2005} as a (very) special case.  
\item The case in which only $U$ is assumed to be a submonoid (and certain other conditions hold) is covered in Section~\ref{sect:1M1}; see Theorems~\ref{thm:ES} and \ref{thm:simp}.  It turns out that this case, and the previous one, include a great many important and natural examples, several of which will be covered in subsequent chapters.  
\item Section~\ref{sect:1M2} considers the case in which only $S$ is assumed to be a submonoid, and we will see that this is rather more complicated; see Theorem \ref{thm:US_sd}.
\eit
Unfortunately, the hope described above cannot be realised in full generality.  Indeed, Example~\ref{eg:TX}\ref{it:TX3} shows that~$US$ might not even contain isomorphic copies of $U$ or $S$.  Moreover, even when $US$ does contain copies of both $U$ and $S$, the hope is still unrealistic, even for extremely simple semigroups $U$ and $S$.  For example, if $U$ and~$S$ are both free monogenic semigroups (isomorphic to the positive integers under addition), then the direct product $U\times S$ (the simplest possible case of the $(U\rtimes S)/\th$ construction) is not even finitely generated \cite{HR1994,RRW1998,East2021}.

Given the above considerations, the next result seems to be the most general one could hope to give concerning presentations for $US$.  It is a special case of (the semigroup version of) Lemma~\ref{lem:ST}.  The statement utilises the surmorphism $\pi:U\rtimes S\to US$ from Proposition \ref{prop:A1}, with kernel $\th$.

\begin{prop}\label{prop:ES_sd}
Suppose $(U,S)$ is a weak action pair, and suppose the semidirect product $U\rtimes S$ has presentation $\Spres XR$ via $\phi:X^+\to U\rtimes S$.  Let $N:U\rtimes S\to X^+$ be a normal form function, and suppose $\th=\Om^\sharp$.  Then $US$ has presentation $\Spres X{R\cup R_\Om}$ via $\phi\pi:X^+\to US$, where $R_\Om=\bigset{(N(\ba),N(\bb))}{(\ba,\bb)\in\Om}$.  \epfres
\end{prop}

\begin{rem}\label{rem:SD}
The previous result uses a presentation for the semidirect product $U\rtimes S$ as a `black box'.  Presentations for $U\rtimes S$ are not always readily constructable from presentations for~$U$ and $S$, as already discussed (even for the simplest case of \emph{direct} products).  However, when~$U\rtimes S$ is a monoid, neat constructions are given by Lavers \cite{Lavers1998}; we will utilise these in Section \ref{sect:2M1}.

The case in which $U$ is a monoid, and the semigroup $S$ acts on $U(=U^1)$ by monoid morphisms, is considered in \cite{FADEG2019}, where a presentation for $U\rtimes S$ is constructed from a presentation for $S$ and the entire multiplication table of $U$.  As we will see in Section \ref{sect:1M2}, similar presentations exist in the case that $S$ is a monoid with a monoidal action on $U^1$; see Theorem \ref{thm:US_sd}.
\end{rem}

\sect{Two submonoids I}\label{sect:2M1}

In this section we consider the simplest special case, where $U\rtimes S$ is a monoid.  By Corollary~\ref{cor:USmon}, this occurs precisely when $U$ and $S$ are both monoids, with $S$ acting on monoidally on $U(=U^1)$ by monoid morphisms.  

For the duration of this section, we fix a weak action pair $(U,S)$ in a monoid $M$, with $U$ and~$S$ submonoids, and we assume the action in \ref{A1} is by monoid morphisms.  In other words, we assume that $U$ and $S$ are submonoids of $M$, and that \ref{SA1'} holds.
We will later additionally assume that $(U,S)$ is an action pair, not just a weak action pair.  In this case, $(U,S)$ is strong, by Lemma~\ref{lem:SAPAP}\ref{SAPAP2}.

We also fix monoid presentations for $U$ and $S$:
\[
\text{$\Mpres{X_U}{R_U}$ \ \ via \ \  $\phi_U:X_U^*\to U$ \AND $\Mpres{X_S}{R_S}$ \ \ via \ \ $\phi_S:X_S^*\to S$.}
\]
For convenience, and without loss of generality, we assume that:
\bit
\item $X_U$ and $X_S$ are disjoint, 
\item no letter of $X_U$ or $X_S$ maps to~$1$, and
\item $\phi_U\restr_{X_U}$ and $\phi_S\restr_{X_S}$ are injective.
\eit
We write $\ol w=w\phi_U$ and $\ol v=v\phi_S$ for all $w\in X_U^*$ and $v\in X_S^*$.  There is no chance of confusion here, as the only word common to~$X_U^*$ and $X_S^*$ is the empty word $\ew$, and we have $\ew\phi_U=\ew\phi_S=1$.  We also fix normal form functions
\[
N_U:U\to X_U^* \AND N_S:S\to X_S^*.
\]
We may assume that $N_U(1)=N_S(1)=\ew$, and also that $N_U(\ol x)=x$ and $N_S(\ol y)=y$ for all $x\in X_U$ and $y\in X_S$.   

Since $U$, $S$ and $U\rtimes S$ are all monoids, a presentation for the semidirect product $U\rtimes S$ can be constructed using a result of Lavers \cite{Lavers1998}, which we state below.  Lavers' presentation is built out of the presentations $\Mpres{X_U}{R_U}$ and $\Mpres{X_S}{R_S}$, and an additional set~$R_1$ of relations over $X_U\cup X_S$ that capture the action of $S$ on~$U(=U^1)$ from \ref{A1}.  To motivate these additional relations, consider letters $x\in X_S$ and~$y\in X_U$.  Since $\ol x\in S$ and $\ol y\in U$, within $M$ we have~${\ol x\cdot \ol y = {}^{\ol x}\ol y\cdot \ol x}$, with ${}^{\ol x}\ol y\in U(=U^1)$.  The set $R_1$ represents each such equation as a relation, and is defined by
\begin{equation}\label{eq:R1}
R_1 = \bigset{(xy,{}^xy\cdot x)}{x\in X_S,\ y\in X_U},
\end{equation}
where for convenience we write ${}^xy=N_U({}^{\ol x}\ol y)\in X_U^*$.
The next result is \cite[Corollary~2]{Lavers1998}; it is also a special case of Theorem \ref{thm:M} below.

\begin{thm}\label{thm:Lavers}
If $U$ and $S$ are monoids, and if $S$ acts monoidally on $U(=U^1)$ by monoid morphisms, then with the above notation, the monoid $U\rtimes S$ has presentation
\[
\Mpres{X_U\cup X_S}{R_U\cup R_S\cup R_1}
\]
via
\[
\phi:(X_U\cup X_S)^*\to U\rtimes S:x\mt \begin{cases}
(x\phi_U,1) &\text{if $x\in X_U$,}\\
(1,x\phi_S) &\text{if $x\in X_S$.}  
\end{cases}
\epfreseq
\]
\end{thm}

By Proposition \ref{prop:A1}, $US\cong(U\rtimes S)/\th$ is a homomorphic image of $U\rtimes S$.  Thus, we can use Lemma \ref{lem:ST} to extend the above presentation for $U\rtimes S$ to a presentation for~$US$ via the surmorphism
\[
\Phi = \phi\pi :(X_U\cup X_S)^*\to US:x\mt \begin{cases}
x\phi_U &\text{if $x\in X_U$,}\\
x\phi_S &\text{if $x\in X_S$.}  
\end{cases}
\]  
In the next statement we use the canonical normal form function $N:U\rtimes S\to(X_U\cup X_S)^*$ defined by $N(u,s)=N_U(u)N_S(s)$.

\begin{thm}\label{thm:ESmon0}
Suppose $(U,S)$ is a weak action pair, with $U$ and $S$ both submonoids, and with~$S$ acting on $U(=U^1)$ by monoid morphisms.  Suppose also that $\th=\Om^\sharp$.  Then with the above notation, the monoid~$US$ has presentation
\[
\Mpres{X_U\cup X_S}{R_U\cup R_S\cup R_1\cup R_\Om}
\]
via $\Phi$, where $R_\Om= \bigset{(N(\ba),N(\bb))}{(\ba,\bb)\in\Om}$.  \epfres
\end{thm}

\begin{rem}\label{rem:ESmon0}
In the special case that $(U,S)$ is an action pair (and not just a weak action pair), we have already observed that $(U,S)$ is in fact strong.  In this case, we can take $\Om$ to be the generating set for $\th$ from Lemma \ref{lem:Om3}, and then
\begin{align}
\label{eq:ROm1} R_\Om &= \bigset{(N_U(u)N_S(s),N_U(u)N_S(t))}{u\in U,\ (s,t)\in\Om_u},
\longintertext{where the $\Om_u$ generate the right congruences $\th_u$ from \eqref{eq:the}.  \endgraf
Lemmas \ref{lem:Om4}, \ref{lem:Om5} and~\ref{lem:Om6} allow us to replace the above $R_\Om$ with a smaller set of relations if the pair $(U,S)$ satisfies any of the additional conditions listed in these lemmas.  When Lemma~\ref{lem:Om4}\ref{Om41} applies, we can take}
\label{eq:ROm2} R_\Om &= \bigset{(N_U(v)N_S(s),N_U(v)N_S(t))}{v\in V,\ (s,t)\in\Om_v},
\intertext{for a suitable subset $V\sub U$.  If Lemma \ref{lem:Om4}\ref{Om42} applies, then since $U=\la\ol X_U\ra$, we can take}
\label{eq:ROm3} R_\Om &= \bigset{(xN_S(s),xN_S(t))}{x\in X_U,\ (s,t)\in\Om_{\ol x}}.
\intertext{In the (very) special case that $S$ is a group, $(U,S)$ is automatically strong by Lemma \ref{lem:group}\ref{group2}.  Lemma \ref{lem:Om5} then applies, and we can take}
\label{eq:ROm4} R_\Om &= \bigset{(N_U(u),N_U(u)N_S(s))}{u\in U,\ s\in\Ga_u},
\intertext{where the $\Ga_u$ generate the subgroups $\S_u$ from \eqref{eq:Ae}.  Again, this can be further simplified if either of the conditions of Lemma \ref{lem:Om6} hold.  When Lemma \ref{lem:Om6}\ref{Om61} applies, we can take}
\label{eq:ROm5} R_\Om &= \bigset{(N_U(v),N_U(v)N_S(s))}{v\in V,\ s\in\Ga_v}
\intertext{for a suitable subset $V\sub U$.  When Lemma \ref{lem:Om6}\ref{Om62} applies, we can take}
\label{eq:ROm6} R_\Om &= \bigset{(x,xN_S(s))}{x\in X_U,\ s\in\Ga_{\ol x}}
\end{align}
In particular, these versions of Theorem \ref{thm:ESmon0} include the (very very) special case of factorisable inverse monoids (cf.~Example \ref{eg:FIM}), where $U$ and $S$ are respectively the semilattice of idempotents and group of units of $US$, and we therefore obtain the main results of \cite{EEF2005} as corollaries.  
\end{rem}

\sect{Two submonoids II}\label{sect:2M2}

We now consider the more general case in which $U$ and $S$ are still submonoids of the over-monoid~$M$, but we no longer assume $U\rtimes S$ is a monoid.  Even though $US$ is still a submonoid of~$M$ in this case, we cannot hope to obtain a presentation for $US$ via Lavers' Theorem \ref{thm:Lavers}, as this obviously does not apply when $U\rtimes S$ is not a monoid.  We will see in Theorem \ref{thm:ESmon}, however, that we can still obtain a rather neat presentation for $US$ that bears some resemblence to Theorem \ref{thm:ESmon0}, and the simplifications discussed in Remark \ref{rem:ESmon0}.  We explore the reason behind this resemblance in Section \ref{sect:2M3}.

Throughout this section we fix an action pair $(U,S)$ in a monoid $M$, and we assume that $U$ and $S$ are submonoids.  Again we fix presentations for $U$ and $S$:
\[
\text{$\Mpres{X_U}{R_U}$ \ \ via \ \  $\phi_U:X_U^*\to U$ \AND $\Mpres{X_S}{R_S}$ \ \ via \ \ $\phi_S:X_S^*\to S$.}
\]
We assume again that:
\bit
\item $X_U$ and $X_S$ are disjoint, 
\item no letter of $X_U$ or $X_S$ maps to~$1$, and
\item $\phi_U\restr_{X_U}$ and $\phi_S\restr_{X_S}$ are injective.
\eit
We again write $\ol w=w\phi_U$ and $\ol v=v\phi_S$ for $w\in X_U^*$ and $v\in X_S^*$, and fix normal form functions
\[
N_U:U\to X_U^* \AND N_S:S\to X_S^*,
\]
assuming that $N_U(1)=N_S(1)=\ew$, and that $N_U(\ol x)=x$ and $N_S(\ol y)=y$ for $x\in X_U$ and $y\in X_S$.  Since $US=\la\ol X_U\cup\ol X_S\ra$, as $U$ and $S$ are submonoids, we still have a surmorphism
\[
\Phi :(X_U\cup X_S)^*\to US:x\mt \begin{cases}
x\phi_U &\text{if $x\in X_U$,}\\
x\phi_S &\text{if $x\in X_S$.}  
\end{cases}
\]  
As in the previous section, we define
\begin{align*}
R_1 &= \bigset{(xy,{}^xy\cdot x)}{x\in X_S,\ y\in X_U},
\intertext{where again we write ${}^xy=N_U({}^{\ol x}\ol y)\in X_U^*$.  We also define}
R_2 &= \bigset{(N_U(u)N_S(s),N_U(u)N_S(t))}{u\in U,\ (s,t)\in\Om_u},
\intertext{where each $\Om_u$ generates $\th_u$ ($u\in U$) as a right congruence.  This time, since $(U,S)$ might not be strong, we also need the relations}
R_3 &= \bigset{(x,x^+x)}{x\in X_S,\ \ol x^+\not=1},
\end{align*}
where here we write $x^+=N_U(\ol x^+)$ for each~$x\in X_S$.  Note that $x^+=\ew$ if $\ol x^+=1$, so that we do not need to include a relation $(x,x^+x)$ for such an $x$.  In particular, $R_3=\es$ if $(U,S)$ is strong; cf.~Lemma \ref{lem:SAPAP}\ref{SAPAP2}.  It is also worth noting that $\ol{x^+}=\ol x^+$, by definition of $N_U$.

\begin{thm}\label{thm:ESmon}
If $(U,S)$ is an action pair, and if $U$ and $S$ are both submonoids, then with the above notation, the monoid~$US$ has presentation
\[
\Mpres{X_U\cup X_S}{R_U\cup R_S\cup R_1\cup R_2\cup R_3}
\]
via $\Phi$.
\end{thm}

\pf
We have already noted that $\Phi$ is a surmorphism, and it is routine to check that
\[
R_U\cup R_S\cup R_1\cup R_2\cup R_3 \sub \ker(\Phi).
\]
It therefore remains to show that $\ker(\Phi)\sub(R_U\cup R_S\cup R_1\cup R_2\cup R_3)^\sharp$, and for this we require some technical lemmas.  But first we fix notation.

We write $\ol w=w\Phi$ for $w\in(X_U\cup X_S)^*$.  We also use the abbreviations 
\[
{}^wv = N_U({}^{\ol w}\ol v) \AND w^+ = N_U(\ol w^+) \qquad\text{for all $w\in X_S^*$ and $v\in X_U^*$.}
\]
In particular, ${}^\ew v=N_U({}^1\ol v)=N_U(\ol v)$ for all $v\in X_U^*$.  By definition of $N_U$, we have
\begin{equation}\label{eq:olwv}
\ol{{}^wv}={}^{\ol w}\ol v \AND \ol{w^+}=\ol w^+ \qquad\text{for all $w\in X_S^*$ and $v\in X_U^*$.}
\end{equation}
We also note that
\begin{equation}\label{eq:we}
w^+ = N_U(\ol w^+) = N_U({}^{\ol w}1) = {}^w\ew \qquad\text{for all $w\in X_S^*$.}
\end{equation}
We write ${\sim} = (R_U\cup R_S\cup R_1\cup R_2\cup R_3)^\sharp$ for the congruence on $(X_U\cup X_S)^*$ generated by the relations.  We also write ${\sim_U}=R_U^\sharp$, and similarly for $\sim_S$, $\sim_1$ and so on.  This will allow us to indicate which of the various sets of relations are used in the manipulations of words to follow.

\begin{lemma}\label{lem:xv}
For any $x\in X_S$ and $v\in X_U^*$, we have $xv\sim{}^xv\cdot x$.
\end{lemma}

\pf
We prove the result by induction on $k=\ell(v)$, the length (number of letters) of $v$.  If $k=0$, then $v=\ew$, and we have $xv=x \sim_3 x^+x={}^x\ew\cdot x={}^xv\cdot x$, where we used \eqref{eq:we} in the third step.  (In the second step, note that $x=x^+x$ if $\ol x^+=1$.)

We now assume that $k\geq1$, so that $v=wy$ for some $w\in X_U^*$ and $y\in X_U$ with $\ell(w)=k-1$.  Then by induction, $xv = xwy \sim {}^xw\cdot xy \sim_1 {}^xw\cdot{}^xy\cdot x$, so it remains to show that ${}^xw\cdot{}^xy \sim {}^xv$.  But this follows quickly from properties of the action and of normal form functions (cf.~\eqref{eq:N}), as
\[
{}^xw\cdot{}^xy = N_U({}^{\ol x}\ol w)\cdot N_U({}^{\ol x}\ol y) \sim_U N_U({}^{\ol x}\ol w\cdot {}^{\ol x}\ol y) = N_U({}^{\ol x}(\ol w\cdot \ol y)) = N_U({}^{\ol x}\ol v) = {}^xv.  \qedhere
\]
\epf

\begin{lemma}\label{lem:wv}
For any $w\in X_S^*$ and $v\in X_U^*$, we have $wv\sim{}^wv\cdot w$.
\end{lemma}

\pf
We use induction on $k=\ell(w)$.  The $k=0$ case is trivial, as ${{}^\ew v=N_U(\ol v)\sim_U v}$.  For $k\geq1$, we write $w=w'x$, where $w'\in X_S^*$ and $x\in X_S$ with $\ell(w')=k-1$.  Then by induction and Lemma \ref{lem:xv}, we have $wv = w'xv \sim w'\cdot{}^xv \cdot x \sim {}^{w'}\!({}^xv)\cdot w'x =  {}^{w'}\!({}^xv)\cdot w$, so it remains to show that ${}^{w'}\!({}^xv) \sim {}^wv$.  For this we have
\[
{}^{w'}\!({}^xv) = N_U \big({}^{\ol w'}\! (\ol{{}^x v}) \big) = N_U \big({}^{\ol w'}\! ({}^{\ol x}\ol v)\big) = N_U({}^{\ol w'\ol x}\ol v) = N_U({}^{\ol w}\ol v) = {}^wv,
\]
where we used \eqref{eq:olwv} in the second step.
\epf

\begin{lemma}\label{lem:w+}
For any $w\in X_S^*$, we have $w\sim w^+w$.
\end{lemma}

\pf
By Lemma \ref{lem:wv} and \eqref{eq:we}, we have $w=w\cdot\ew\sim{}^w\ew\cdot w=w^+w$.
\epf

\begin{rem}\label{rem:no2}
For later reference, it is important to note that the proof so far has not used any of the relations from $R_2$.
\end{rem}

Returning now to the proof of the theorem, we recall that the last remaining step is to show that $\ker(\Phi)\sub{\sim}$.  To do this, let $(w,v)\in\ker(\Phi)$: i.e., $w,v\in(X_U\cup X_S)^*$ and $\ol w=\ol v$.  We must show that $w\sim v$.  

We first observe that
\begin{equation}\label{eq:uv}
w \sim_1 w_1w_2 \AND v \sim_1 v_1v_2 \qquad\text{for some $w_1,v_1\in X_U^*$ and $w_2,v_2\in X_S^*$.}
\end{equation}
It follows from Lemma \ref{lem:w+} that
\begin{equation}\label{eq:uv+}
w \sim w_1w_2^+w_2 \AND v \sim v_1v_2^+v_2.
\end{equation}
From \eqref{eq:uv}, we have
\[
\ol w_1\ol w_2 = \ol w = \ol v = \ol v_1 \ol v_2 \qquad\text{with $\ol w_1,\ol v_1\in U$ and $\ol w_2,\ol v_2\in S$.}
\]
It follows from \ref{A2} that $\ol w_1\ol w_2^+=\ol v_1\ol v_2^+$, and we denote this element of $U$ by $u$.  Now,
\[
\ol{w_1w_2^+} = \ol w_1 \ol w_2^+ = u \ANDSIM \ol{v_1v_2^+} = u,
\]
so it follows that $w_1w_2^+\sim N_U(u) \sim v_1v_2^+$.  Combining this with \eqref{eq:uv+}, we have
\begin{equation}\label{eq:uv2}
w \sim N_U(u)\cdot w_2 \AND v \sim N_U(u)\cdot v_2,
\end{equation}
so it remains to show that
\begin{equation}\label{eq:uv3}
N_U(u)\cdot w_2 \sim N_U(u)\cdot v_2.
\end{equation}
Applying $\Phi$ to \eqref{eq:uv2}, we have $u\cdot\ol w_2=\ol w=\ol v=u\cdot\ol v_2$, so that $(\ol w_2,\ol v_2)\in\th_u$.  Since $\th_u$ is generated as a right congruence by $\Om_u$, it follows that there is a sequence
\[
\ol w_2 = s_1 \to s_2\to\cdots\to s_k=\ol v_2
\]
where $s_1,\ldots,s_k\in S$, and such that for each $1\leq i<k$,
\[
s_i=a_ic_i \ANd s_{i+1}=b_ic_i \qquad\text{for some $(a_i,b_i)\in\Om_u\cup\Om_u^{-1}$ and $c_i\in S(=S^1)$.}
\]
Since $N_U(u)\cdot w_2\sim_S N_U(u)\cdot N_S(s_1)$ and $N_U(u)\cdot v_2\sim_S N_U(u)\cdot N_S(s_k)$, we can complete the proof of \eqref{eq:uv3}, and hence of the theorem, by showing that
\[
N_U(u)\cdot N_S(s_i) \sim N_U(u)\cdot N_S(s_{i+1}) \qquad\text{for each $1\leq i<k$.}
\]
But for any such $i$ we have
\[
N_U(u)\cdot N_S(s_i) \sim_S N_U(u)\cdot N_S(a_i)\cdot N_S(c_i) \sim_2 N_U(u)\cdot N_S(b_i)\cdot N_S(c_i) \sim_S N_U(u)\cdot N_S(s_{i+1}),
\]
as required.  
\epf

\begin{rem}\label{rem:ESmon1}
When the pair $(U,S)$ is strong, Theorem \ref{thm:ESmon} reduces to Theorem \ref{thm:ESmon0}, or more specifically to the version of Theorem \ref{thm:ESmon0} where $R_\Om$ is the set of relations in \eqref{eq:ROm1}.
\end{rem}

\begin{rem}\label{rem:pres12}
The presentation 
\begin{align}\label{eq:pres1}
&\Mpres{X_U\cup X_S}{R_U\cup R_S\cup R_1\cup R_2\cup R_3}
\intertext{for $US$ from Theorem \ref{thm:ESmon} contains/extends the presentation}
\label{eq:pres2}
&\Mpres{X_U\cup X_S}{R_U\cup R_S\cup R_1}.
\intertext{Lavers' Theorem \ref{thm:Lavers} says that \eqref{eq:pres2} is a presentation for the semidirect product $U\rtimes S$ in the case that the latter is a monoid.  This is obviously not true when ${U\rtimes S}$ is not a monoid, so the reader may wonder `why' our presentation \eqref{eq:pres1} for $US$ still extends the presentation \eqref{eq:pres2}.  We will explore this question in Section \ref{sect:2M3}, where we will in fact see (in Theorem \ref{thm:M}) that the intermediate}
\nonumber &\Mpres{X_U\cup X_S}{R_U\cup R_S\cup R_1\cup R_3}
\end{align}
is a presentation for the (local) monoid $\bigset{(u,s)\in U\rtimes S}{u=us^+}\leq U\rtimes S$.
\end{rem}

In the presentation $\Mpres{X_U\cup X_S}{R_U\cup R_S\cup R_1\cup R_2\cup R_3}$ from Theorem \ref{thm:ESmon}, the sets of relations $R_1$ and $R_3$ are quantified over letters from $X_U\cup X_S$: i.e., generators for $U$ and~$S$.  However, $R_2$ is quantified over all of $U$.  Ideally, one would like to replace $R_2$ with a set of relations quantified over $X_U$.  This is sometimes possible, but not always.  The second part of the next result describes one such situation, and the first part describes a more general situation in which some degree of simplification of $R_2$ is still possible.  The conditions are similar to those discussed in Section \ref{sect:cong}, but note that we do not need to assume the pair $(U,S)$ is strong here.  In the statement, $\pre$ is the relation on $U$ given in \eqref{eq:pre}.

\begin{thm}\label{thm:Msimp}
Suppose $(U,S)$ is an action pair, with $U$ and $S$ submonoids.
\ben
\item \label{Msimp1}  If there exists a subset $V\sub U$ such that
\[
\th_u = \bigvee_{v\in V, \atop v\pre u} \th_v \qquad\text{for all $u\in U$,}
\]
then $US$ has presentation $\Mpres{X_U\cup X_S}{R_U\cup R_S\cup R_1\cup R_2'\cup R_3}$ via $\Phi$, where
\[
R_2' = \bigset{(N_U(v)N_S(s),N_U(v)N_S(t))}{v\in V,\ (s,t)\in\Om_v}.
\]
\item \label{Msimp2}  If $U$ is commutative, and if $\th_{uv}=\th_u\vee\th_v$ for all $u,v\in U$, then $US$ has presentation
\[
\Mpres{X_U\cup X_S}{R_U\cup R_S\cup R_1\cup R_2''\cup R_3}
\]
via $\Phi$, where $R_2'' = \bigset{(xN_S(s),xN_S(t))}{x\in X_U,\ (s,t)\in\Om_{\ol x}}$.
\een
\end{thm}

\pf
As with Lemma \ref{lem:Om4}, it suffices to prove \ref{Msimp1}.  For this, we proceed in similar fashion to the proof of Lemma \ref{lem:Om4}\ref{Om41}.

By Theorem \ref{thm:ESmon}, and since $R_2'\sub R_2$, it suffices to show that $R_2\sub{\approx}$, where this time we write ${\approx}=(R_U\cup R_S\cup R_1\cup R_2'\cup R_3)^\sharp$.  For this, let $u\in U$ and $(s,t)\in\Om_u$.  We must show that $N_U(u)N_S(s)\approx N_U(u)N_S(t)$.  
Again we write $V_u=\set{v\in V}{v\pre u}$.  Since $(s,t)\in\th_u=\bigvee_{v\in V_u}\th_v$, there is a sequence
\[
s = s_1 \to s_2 \to \cdots \to s_k = t,
\]
such that each $(s_i,s_{i+1})\in\bigcup_{v\in V_u}(\Om_v\cup\Om_v^{-1})$, and it is enough to show that
\[
N_U(u)N_S(s_i) \approx N_U(u)N_S(s_{i+1}) \qquad\text{for all $1\leq i<k$.}
\]
Fix some such $i$, so that $(s_i,s_{i+1})\in\Om_v\cup\Om_v^{-1}$ for some $v\in V_u$.  Then $u=wv$ for some $w\in U(=U^1)$, and 
\[
N_U(u)N_S(s_i)\approx N_U(w)\cdot N_U(v)N_S(s_i) \approx N_U(w)\cdot N_U(v)N_S(s_{i+1}) \approx N_U(u)N_S(s_{i+1}).  \qedhere
\]
\epf

\begin{rem}
In the very special case that $S$ is additionally assumed to be a group, the pair $(U,S)$ is then strong.  Simplifications arising in this case have already been discussed in Remark~\ref{rem:ESmon0}.
\end{rem}

\sect{Two submonoids III}\label{sect:2M3}

In this section we briefly discuss an alternative way to prove Theorem \ref{thm:ESmon}, and also give a generalisation of Theorem \ref{thm:ESmon0}; see Theorem \ref{thm:ESmon00}.  This involves proving an additional general result, Theorem \ref{thm:M}, which is of independent interest, and which will in fact be directly applicable when we consider free left restriction monoids in Chapter \ref{chap:LR}.

So again we consider a weak action pair $(U,S)$, where $U$ and $S$ are both submonoids of the over-monoid (which shall not be named in this section).  By Proposition~\ref{prop:A1}, the monoid $US$ is a homomorphic image of $U\rtimes S$, though the latter need not be a monoid; cf.~Lemma \ref{lem:USmon}.  In the case that $U\rtimes S$ \emph{is} a monoid, Theorem~\ref{thm:ESmon0} built a presentation for $US$, starting with Lavers' (monoid) presentation for $U\rtimes S$ from Theorem~\ref{thm:Lavers}.  As we noted in the previous section, however, this approach is obviously not available in the case that $U\rtimes S$ is not a monoid.  Despite this, the presentation for $US$ from Theorem \ref{thm:ESmon} (which applies when $(U,S)$ is an action pair, not just a weak action pair) certainly `looks like' it has been built in this way (cf.~Remark \ref{rem:pres12}), and we explore the reason for this in the current section.  

Recall from Remark \ref{rem:A1} that the natural surmorphism $\pi:U\rtimes S\to US:(u,s)\mt us$ factors through the local monoid
\[
M = (1,1)\cdot(U\rtimes S)\cdot(1,1) = \bigset{(u,s)\in U\rtimes S}{u=us^+} \leq U\rtimes S,
\]
so that $US\cong M/\vt$, where $\vt=\ker(\pi\restr_M)$.  It follows that one could obtain a presentation for~$US$ by first finding a (monoid) presentation for $M$, then finding a generating set for the congruence~$\vt$, and then finally applying Lemma \ref{lem:ST}.  We show here how to do the first of these tasks; see Theorem \ref{thm:M}.

The next result concerns a monoid $S$ acting monoidally on a monoid $U$ by semigroup morphisms.  In particular, this is the case whenever $(U,S)$ is a weak action pair; cf.~Remark \ref{rem:AP2}.
For the statement and proof we keep the notation of Section \ref{sect:2M2}, in particular:
\bit
\item the presentations and normal form functions for $U$ and $S$, 
\item the over-line notation for words over $X_U$ and $X_S$, and 
\item the additional sets of relations $R_1$ and $R_3$ (the relations $R_2$ are not required).
\eit
We also define the morphism
\[
\varphi:(X_U\cup X_S)^*\to M:x\mt \begin{cases}
(\ol x,1) &\text{if $x\in X_U$,}\\
(\ol x^+,\ol x) &\text{if $x\in X_S$.}
\end{cases}
\]

\begin{thm}\label{thm:M}
If the monoid $S$ acts monoidally on the monoid $U$ by semigroup morphisms, then the (local) monoid $M = \bigset{(u,s)\in U\rtimes S}{u=us^+}$ has presentation
\[
\Mpres{X_U\cup X_S}{R_U\cup R_S\cup R_1\cup R_3}
\]
via $\varphi$.
\end{thm}

\pf
First, for any $(u,s)\in M$ we have
\[
(u,s) = (us^+,s) = (u,1)\cdot(s^+,s) = N_U(u)\varphi \cdot N_S(s)\varphi = (N_U(u)\cdot N_S(s))\varphi,
\]
and this shows that $\varphi$ is surjective.  (In the above calculation, it is clear that $N_U(u)\varphi = (u,1)$.  For $N_S(s)\varphi = (s^+,s)$, we apply \eqref{eq:st^+}.)

Next, it is easy to check that $\varphi$ preserves the relations, so that ${\approx}\sub\ker(\varphi)$, where we write ${{\approx}=(R_U\cup R_S\cup R_1\cup R_3)^\sharp}$.  Here we just show that $R_1\sub\ker(\varphi)$.  For this, let $x\in X_S$ and $y\in X_U$.  Then, using Lemma \ref{lem:+}\ref{lem+4}, we have
\[
(xy)\varphi = (\ol x^+,\ol x)\cdot(\ol y,1) = (\ol x^+\cdot{}^{\ol x}\ol y,\ol x) = ({}^{\ol x}\ol y,\ol x) = ({}^{\ol x}\ol y\cdot\ol x^+,\ol x) = ({}^{\ol x}\ol y,1)\cdot(\ol x^+,\ol x) = ({}^xy\cdot x)\varphi.
\]

To complete the proof we must show that $\ker(\varphi)\sub{\approx}$, so suppose $(w,v)\in\ker(\varphi)$.  As in Lemmas \ref{lem:xv}--\ref{lem:w+} (cf.~Remark~\ref{rem:no2}), we have
\[
w \approx w_1w_2 \approx w_1w_2^+w_2 \ANd v \approx v_1v_2 \approx v_1v_2^+v_2 \qquad\text{for some $w_1,v_1\in X_U^*$ and $w_2,v_2\in X_S^*$.}
\]
But then
\[
w\varphi = w_1\varphi\cdot w_2\varphi = (\ol w_1,1)\cdot(\ol w_2^+,\ol w_2) = (\ol w_1\ol w_2^+,\ol w_2) \ANDSIM v\varphi = (\ol v_1\ol v_2^+,\ol v_2).
\]
Since $w\varphi=v\varphi$, it follows that $\ol w_1\ol w_2^+ = \ol v_1\ol v_2^+$ and $\ol w_2 = \ol v_2$, so that $w_1w_2^+\sim_U v_1v_2^+$ and $w_2\sim_S v_2$.  But then
\[
w \approx w_1w_2^+\cdot w_2 \approx v_1v_2^+\cdot v_2 \approx v,
\]
and the proof is complete.
\epf

\begin{rem}
When the action of $S$ on $U$ is additionally by monoid morphisms (i.e., when each $s^+=1$), the semidirect product $U\rtimes S$ is a monoid, and $M=U\rtimes S$; cf.~Corollary \ref{cor:USmon} and Proposition \ref{prop:MM}.  In this case we have $R_3=\es$, and Theorem \ref{thm:M} reduces to Lavers' Theorem~\ref{thm:Lavers}.
\end{rem}

An application of Lemma \ref{lem:ST} yields the following more general version of Theorem \ref{thm:ESmon0}.  In the statement we use the surmorphism
\[
\Phi = \varphi\circ\pi\restr_M:(X_U\cup X_S)^*\to US:x\mt \begin{cases}
x\phi_U &\text{if $x\in X_U$,}\\
x\phi_S &\text{if $x\in X_S$,}  
\end{cases}
\]
where $\varphi:(X_U\cup X_S)^*\to M$ is as in Theorem \ref{thm:M}.  We also refer to the congruence $\vt=\ker(\pi\restr_M)$, and the canonical normal form function $N:M\to(X_U\cup X_S)^*$ defined by $N(u,s)=N_U(u)N_S(s)$.

\begin{thm}\label{thm:ESmon00}
Suppose $(U,S)$ is a weak action pair, with $U$ and $S$ both submonoids, and suppose also that $\vt=\Om^\sharp$.  Then with the above notation, the monoid~$US$ has presentation
\[
\Mpres{X_U\cup X_S}{R_U\cup R_S\cup R_1\cup R_3\cup R_\Om}
\]
via $\Phi$, where $R_\Om= \bigset{(N(\ba),N(\bb))}{(\ba,\bb)\in\Om}$.  \epfres
\end{thm}

\sect{One submonoid I}\label{sect:1M1}

Theorem \ref{thm:ESmon} gave a presentation for the monoid $US$ arising from an action pair $(U,S)$, when~$U$ and $S$ are both submonoids of the over-monoid.  However, in some of our intended applications,~$S$ and $US$ are not a submonoids, even though $U$ is a (commutative) submonoid.  This section deals with such cases, where we must make some additional assumptions.  It is convenient to list these up front:

\newpage

\begin{ass}\label{ass:US}
We assume that $(U,S)$ is an action pair in a monoid $M$ such that:
\ben
\item \label{ass1} $U$ is a submonoid of $M$ (which implies that $S\sub US$),
\item \label{ass2} $U\sm\{1\}$ is a subsemigroup of $U$, and $U\sm\{1\}\sub US$,
\item \label{ass3} $(U,S^1)$ is an action pair extending $(U,S)$, in the sense described before Lemma \ref{lem:US1}.
\een
\end{ass}

\begin{rem}\label{rem:assUS}
Item \ref{ass2} of Assumption \ref{ass:US} is a fairly strong condition, but holds in all of our motivating examples in later chapters.  For example, it is known that $U\sm\{1\}\leq U$ if $U$ is an idempotent-generated monoid \cite[Lemma 4.9]{Sandwiches1}.  In particular, this holds for monoid bands, including monoid semilattices.  An important special case therefore occurs when $U=P(M)$ is the semilattice of projections of a left restriction monoid $M$; cf.~Propositions \ref{prop:LR1} and \ref{prop:LR2} (and Remark \ref{rem:PS}).  

By Lemma \ref{lem:US1}, item \ref{ass3} of Assumption \ref{ass:US} holds if and only if:
\bit
\item $us=v \implies us^+=v$ for all $u,v\in U^1$ and $s\in S$.
\eit
By the same lemma, this holds automatically if $(U,S)$ is strong, and then $(U,S^1)$ is also strong.

It is also worth noting that (under Assumption \ref{ass:US}) $US$ is a submonoid of $M$ if and only if~$S$ is a submonoid, as follows from Lemma \ref{lem:submon}.  Of course if $S$ (and hence $US$) happens to be a submonoid, then we could just apply Theorem \ref{thm:ESmon} to obtain a (monoid) presentation for $US$.
\end{rem}

For the duration of this section, we fix a pair $(U,S)$ satisfying Assumption \ref{ass:US}.
We have already noted in Remark \ref{rem:SD} that presentations for $U\rtimes S$ exist in the case that $S$ acts on $U(=U^1)$ by monoid morphisms \cite{FADEG2019}: i.e., when the pair $(U,S)$ is strong.  However, since the presentation from \cite{FADEG2019} utilises the entire multiplication table for~$U$, attempting to use this (and Proposition~\ref{prop:ES_sd}) to obtain presentations for $US\cong(U\rtimes S)/\th$ in this case is undesirable.  It also does not apply to the more general situation in which $(U,S)$ is not strong.  
In any case, we will see in Theorem \ref{thm:ES} that it is still possible to construct a presentation for~$US$ (under Assumption~\ref{ass:US}) very much like that from Theorem \ref{thm:ESmon}.

For the rest of this section, we fix semigroup presentations for $U\sm\{1\}$ and $S$:
\[
\text{$\Spres{X_U}{R_U}$ \ \ via \ \  $\phi_U:X_U^+\to U\sm\{1\}$ \AND $\Spres{X_S}{R_S}$ \ \ via \ \ $\phi_S:X_S^+\to S$.}
\]
(The assumption that $U\sm\{1\}\leq U$ means that $U$ itself has presentation $\Mpres{X_U}{R_U}$.)
Again we assume without loss of generality that $X_U$ and $X_S$ are disjoint, that $\phi_U\restr_{X_U}$ and $\phi_S\restr_{X_S}$ are both injective, and we write $\ol w=w\phi_U$ and $\ol v=v\phi_S$ for $w\in X_U^+$ and $v\in X_S^+$.  We also fix normal form functions
\[
N_U:U\sm\{1\}\to X_U^+ \AND N_S:S\to X_S^+,
\]
assuming that $N_U(\ol x)=x$ and $N_S(\ol y)=y$ for $x\in X_U$ and $y\in X_S$.  It is important to note the following convenitions regarding the identity $1$.
\bit
\item Even though $1$ obviously does not belong to $U\sm\{1\}$, it will still be convenient to write $N_U(1)=\ew$, which as usual denotes the empty word.  In the calculations to follow, $N_U(1)=\ew$ will only ever appear as a sub-word of a non-empty word.
\item It is possible that $1\in S$, in which case $N_S(1)$ is (by definition) a fixed non-empty word over~$X_S$ mapping to $1$.  
\item On the other hand, if $1\not\in S$, it is convenient to write $N_S(1)=\ew$ as well.  In this case, $N_S(1)=\ew$ will again only ever appear below as a sub-word of a non-empty word.
\eit

Since $US$ contains both $U\sm\{1\}$ and $S$, there is a well-defined morphism
\[
\Phi :(X_U\cup X_S)^+\to US : x\mt \begin{cases}
x\phi_U &\text{if $x\in X_U$,}\\
x\phi_S &\text{if $x\in X_S$,}
\end{cases}
\]
which is easily seen to be surjective.

As in Section~\ref{sect:2M2}, we define
\begin{align*}
R_1 = \bigset{(xy,{}^xy\cdot x)}{x\in X_S,\ y\in X_U} 
\AND R_3 = \bigset{(x,x^+x)}{x\in X_S,\ \ol x^+\not=1},
\end{align*}
where again we write 
\[
{}^xy=N_U({}^{\ol x}\ol y)\in X_U^* \AND x^+ = N_U(\ol x^+) \in X_U^* \qquad\text{for $x\in X_S$ and $y\in X_U$.}
\]
Note that it is possible to have ${}^{\ol x}\ol y=1$ for some $x\in X_S$ and $y\in X_U$, in which case ${{}^xy=N_U(1)=\ew}$, and the corresponding relation from $R_1$ is simply $(xy,x)$.
Similarly, it is possible to have $x^+=\ew$ for some $x\in X_S$ (when $\ol x^+=1\not\in U$), but in this case $R_3$ does \emph{not} contain the relation $(x,x^+x)$.

To define the final set of relations, we first require another piece of notation.  For each $u\in U$ we extend the relation $\th_u$ from \eqref{eq:the} to
\begin{equation}\label{eq:The}
\Th_u = \bigset{(s,t)\in S^1\times S^1}{us=ut},
\end{equation}
which is again a right congruence on $S^1$.  The reason we need to work with $\Th_u$ in place of~$\th_u$ is that it is possible to have $u=us$ for $u\in U\sm\{1\}$ and $s\in S$, even if $S$ is not itself a monoid.  In fact, such an equation holds for \emph{every} $u\in U$.  Indeed, due to the assumption $U\sm\{1\}\sub US$, each $u\in U\sm\{1\}$ satisfies $u=vs$ for some $v\in U$ and $s\in S$.  Since $(U,S^1)$ is an action pair by Assumption \ref{ass:US}\ref{ass3}, it follows from Lemma \ref{lem:US1}\ref{US12} that $u=vs^+$, and so $u=vs=vs^+s=us$.  Of course if~$S$ is also a submonoid, then $\Th_u=\th_u$ for all $u$, and $u=us\iff(1,s)\in\th_u$; but in this case we already have the (monoid) presentation for $US$ in Theorem \ref{thm:ESmon}.

For each $u\in U\sm\{1\}$ we fix a set of pairs $\Om_u\sub S^1\times S^1$ that generates $\Th_u$ as a right congruence, and we define
\[
R_2 = \bigset{(N_U(u)N_S(s),N_U(u)N_S(t))}{u\in U\sm\{1\},\ (s,t)\in\Om_u}.
\]
Note that it is possible that $s=1$ (or $t=1$) for some $(s,t)\in\Om_u$, and we might have $1\not\in S$.  In this case, we interpret $N_S(s)=\ew$ as above, and we note that the word $N_U(u)N_S(s)$ is still non-empty.

\begin{thm}\label{thm:ES}
Suppose $(U,S)$ is an action pair satisfying Assumption~\ref{ass:US}.  Then with the above notation, the semigroup $US$ has presentation
\[
\Spres{X_U\cup X_S}{R_U\cup R_S\cup R_1\cup R_2\cup R_3}
\]
via $\Phi$.
\end{thm}

\pf
The proof is mostly the same as for Theorem \ref{thm:ESmon}, but with a little extra care required in some places since we do not assume $S$ is a submonoid.

Since $\Phi$ is a surmorphism and preserves the relations, it remains to show that ${\ker(\Phi)\sub{\sim}}$, where ${\sim} = (R_U\cup R_S\cup R_1\cup R_2\cup R_3)^\sharp$.
We again write ${\sim_U}=R_U^\sharp$, and similarly for $\sim_S$, $\sim_1$, and so on.  By convention, we also allow ourselves to write $\ew\sim\ew$.

We write $\ol w=w\Phi$ for $w\in(X_U\cup X_S)^+$.  It will also be convenient to define $\ol\ew=1$, even though $\ew$ does not belong to $(X_U\cup X_S)^+$, and $1$ might not belong to $US=\im(\Phi)$.  We also use the abbreviations 
\[
{}^wv = N_U({}^{\ol w}\ol v) \AND w^+ = N_U(\ol w^+) \qquad\text{for $w\in X_S^*$ and $v\in X_U^*$,}
\]
keeping in mind the above conventions regarding adjoined identities and empty words.  In particular, we again have
\[
\ol{{}^wv}={}^{\ol w}\ol v \COMMA 
\ol{w^+}=\ol w^+ \COMMA
{}^\ew v=N_U(\ol v) \AND
w^+ = {}^w\ew
\qquad\text{for all $w\in X_S^*$ and $v\in X_U^*$.}
\]
We then have the following three lemmas, whose proofs are exactly as for Lemmas \ref{lem:xv}--\ref{lem:w+}:

\begin{lemma}\label{lem:xv'}
For any $x\in X_S$ and $v\in X_U^*$, we have $xv\sim{}^xv\cdot x$.  \epfres
\end{lemma}

\begin{lemma}\label{lem:wv'}
For any $w\in X_S^*$ and $v\in X_U^*$, we have $wv\sim{}^wv\cdot w$.  \epfres
\end{lemma}

\begin{lemma}\label{lem:w+'}
For any $w\in X_S^*$, we have $w\sim w^+w$.  \epfres
\end{lemma}

Returning now to the proof of the theorem, let $(w,v)\in\ker(\Phi)$: i.e., $w,v\in(X_U\cup X_S)^+$ and $\ol w=\ol v$.  We must show that $w\sim v$.  

We first observe that
\begin{equation}\label{eq:uv'}
w \sim_1 w_1w_2 \AND v \sim_1 v_1v_2 \qquad\text{for some $w_1,v_1\in X_U^*$ and $w_2,v_2\in X_S^*$.}
\end{equation}
Note that one of $w_1$ or $w_2$ might be empty, but not both (and similarly for $v_1$ and $v_2$).  It follows from Lemma \ref{lem:w+'} that
\begin{equation}\label{eq:uv+'}
w \sim w_1w_2^+w_2 \AND v \sim v_1v_2^+v_2.
\end{equation}
From \eqref{eq:uv'}, and keeping $\ol\ew=1$ in mind, we have
\[
\ol w_1\ol w_2 = \ol w = \ol v = \ol v_1 \ol v_2 \qquad\text{with $\ol w_1,\ol v_1\in U(=U^1)$ and $\ol w_2,\ol v_2\in S^1$.}
\]
Since $(U,S^1)$ is an action pair (cf.~Assumption \ref{ass:US}), it follows that $\ol w_1\ol w_2^+=\ol v_1\ol v_2^+$.  Keeping in mind that $U\sm\{1\}$ is a semigroup, this implies that the words $w_1w_2^+$ and $v_1v_2^+$ (both from $X_U^*$) are either both empty, or else both non-empty.  We consider these possibilities separately.

\pfcase1 Suppose first that $w_1w_2^+=\ew=v_1v_2^+$.  Then from \eqref{eq:uv+'} we have $w\sim w_2$ and $v\sim v_2$, with ${w_2,v_2\in X_S^+}$.  But then $\ol w_2=\ol w=\ol v=\ol v_2$, so that $w_2\sim_S v_2$, and so $w\sim v$.

\pfcase2  Now suppose $w_1w_2^+\not=\ew\not=v_1v_2^+$.  Let $u=\ol w_1\ol w_2^+=\ol v_1\ol v_2^+\in U\sm\{1\}$, so that ${w_1w_2^+\sim_U N_U(u)\sim_U v_1v_2^+}$.  Combined with \eqref{eq:uv+'} it follows that
\begin{equation}\label{eq:uv2'}
w \sim N_U(u)\cdot w_2 \AND v \sim N_U(u)\cdot v_2.
\end{equation}
Applying $\Phi$ to \eqref{eq:uv2'}, we have $u\cdot\ol w_2=\ol w=\ol v=u\cdot\ol v_2$, so that $(\ol w_2,\ol v_2)\in\Th_u$.  Since $\Th_u$ is generated as a right congruence by $\Om_u$, it follows that there is a sequence
\[
\ol w_2 = s_1 \to s_2\to\cdots\to s_k=\ol v_2
\]
where $s_1,\ldots,s_k\in S^1$, and such that for each $1\leq i<k$,
\[
s_i=a_ic_i \ANd s_{i+1}=b_ic_i \qquad\text{for some $(a_i,b_i)\in\Om_u\cup\Om_u^{-1}$ and $c_i\in S^1$.}
\]
Since $w\sim N_U(u)\cdot w_2\sim N_U(u)\cdot N_S(s_1)$ and similarly $v\sim N_U(u)\cdot N_S(s_k)$, it suffices to show that
\[
N_U(u)\cdot N_S(s_i) \sim N_U(u)\cdot N_S(s_{i+1}) \qquad\text{for each $1\leq i<k$.}
\]
But for any such $i$ we have
\[
N_U(u)\cdot N_S(s_i) \sim_S N_U(u)\cdot N_S(a_i)\cdot N_S(c_i) \sim_2 N_U(u)\cdot N_S(b_i)\cdot N_S(c_i) \sim_S N_U(u)\cdot N_S(s_{i+1}),
\]
as required.  
\epf

\begin{rem}
Despite the obvious similarities, Theorem \ref{thm:ESmon} does not follow from Theorem~\ref{thm:ES} (in the special case that $S$ is a submonoid), as the former does not assume that $U\sm\{1\}$ is a subsemigroup of $U$.
\end{rem}

As with Theorem \ref{thm:Msimp}, and with an essentially identical proof, the presentation from Theorem~\ref{thm:ES} can be simplified in certain special cases:

\begin{thm}\label{thm:simp}
Suppose $(U,S)$ is an action pair satisfying Assumption \ref{ass:US}.  
\ben
\item \label{simp1}  If there exists a subset $V\sub U\sm\{1\}$ such that
\[
\Th_u = \bigvee_{v\in V, \atop v\pre u} \Th_v \qquad\text{for all $u\in U\sm\{1\}$,}
\]
then $US$ has presentation $\Spres{X_U\cup X_S}{R_U\cup R_S\cup R_1\cup R_2'\cup R_3}$ via $\Phi$, where
\[
R_2' = \bigset{(N_U(v)N_S(s),N_U(v)N_S(t))}{v\in V,\ (s,t)\in\Om_v}.
\]
\item \label{simp2}  If $U$ is commutative, and if $\Th_{uv}=\Th_u\vee\Th_v$ for all $u,v\in U$, then $US$ has presentation
\[
\Spres{X_U\cup X_S}{R_U\cup R_S\cup R_1\cup R_2''\cup R_3}
\]
via $\Phi$, where $R_2'' = \bigset{(xN_S(s),xN_S(t))}{x\in X_U,\ (s,t)\in\Om_{\ol x}}$.  \epfres
\een
\end{thm}

\sect{One submonoid II}\label{sect:1M2}

In the previous section we gave presentations for the semigroup $US$ arising from an action pair~$(U,S)$ in which only $U$ was assumed to be a submonoid of the over-monoid $M$, modulo other conditions listed in Assumption~\ref{ass:US}.  In this section we consider the case in which only $S$ is assumed to be a submonoid.  As we will see, this is a much more complicated matter, and the `hope' described at the beginning of Chapter \ref{chap:pres} is far from realisable.  

An example of a semigroup of this form is $\Sing(\I_n)=\I_n\sm\G_n$, the singular part of a finite symmetric inverse monoid; cf.~Example \ref{eg:Tn}.  Here $\Sing(\I_n)=US$ arises from the (strong) action pair ${(U,S)=(\Sing(\E_n),\G_n)}$.  Two different presentations for $\Sing(\I_n)$ may be found in \cite{JEinsn,JEinsn2}, and it is important to note that these presentations are not `built' from presentations for~$U$ and~$S$.  Neither \emph{could} any such presentation be built in this way, as $\Sing(\I_n)$ does not contain any copy of~${S=\G_n}$.  On the other hand, $\Sing(\I_n)$ does of course contain $U=\Sing(\E_n)$, and also (isomorphic copies of)~$\G_{n-1}$, which played an important role in both papers.  Another example is~${\PT_n\sm\T_n}$, the semigroup of all strictly partial transformations~\cite{JEptnsn}, which arises from the (strong) action pair $(\Sing(\E_n),\T_n)$.  These two examples generalise to \emph{almost-(left-)factorisable} inverse and left restriction semigroups; cf.~Examples \ref{eg:FIM} and~\ref{eg:Tn} and Remark~\ref{rem:PS}.

Although it is not possible to give general results for building presentations for $US$ out of presentations for $U$ and $S$ when we assume only $S$ is a submonoid, it is however possible to adapt ideas from \cite{FADEG2019} in order to give a general presentation for $U\rtimes S$ in this case; see Theorem~\ref{thm:US_sd} below.  (Note that \cite[Theorem 3.1]{FADEG2019} applies to semidirect products $U\rtimes S$ where $U$ is a monoid, which is the opposite of our current focus.)  One can then apply Proposition \ref{prop:ES_sd} to extend this to a presentation for $US\cong(U\rtimes S)/\th$, modulo a generating set for the congruence $\th$.  

We assume throughout this section that the monoid $S$ acts monoidally on $U^1$ for some semigroup $U$, allowing for the formation of the semidirect product $U\rtimes S$ as in Definition \ref{defn:SD}.  (This monoidal assumption is indeed satisfied when $(U,S)$ is a weak action pair; cf.~Remark~\ref{rem:AP2}.)  We write $1$ for the identity of both $S$ and $U^1$.  So ${}^1u=u$ for all $u\in U^1$.

We also assume that~$U$ has presentation $\Spres XR$ via $\phi:X^+\to U$, assuming as usual that $\phi\restr_X$ is injective.  (We will not be referring to a presentation for $S$.)  For $w\in X^+$, we write $\ol w=w\phi\in U$, and we write ${\sim}=R^\sharp$.  We also fix a normal form function $N:U\to X^+$, assuming that $N(\ol x)=x$ for all $x\in X$.  If $1\in U$, then by definition $N(1)$ is some non-empty word over~$X$ mapping to $1$.  If $1\not\in U$, then by convention we define $N(1)=\ew$, even though this does not belong to $X^+$.

We now introduce an alphabet
\[
Y = \set{x_s}{x\in X,\ s\in S}
\]
in one-one correspondence with the cartesian product $X\times S$, and we define the morphism
\[
\Phi:Y^+\to U\rtimes S:x_s \mt (\ol x,s).
\]
We identify each $x\in X$ with $x_1\in Y$, and in this way identify $X$ with the subset $\set{x_1}{x\in X}$ of $Y$.  Note then that $x\Phi= x_1\Phi=(\ol x,1)$ for all $x\in X$.

For any (non-empty) word $w=x_1\cdots x_k\in X^+$, and for any $s\in S$, we write
\[
w_s = x_1\cdots x_{k-1}(x_k)_s\in Y^+.
\]
Note that $(wv)_s=w\cdot v_s$ for any $w\in X^*$, $v\in X^+$ and $s\in S$.  We define two sets of relations as follows:
\[
R_1 = \bigset{(w_s,v_s)}{(w,v)\in R,\ s\in S}
\AND
R_2 = \bigset{(x_sy_t,(xN({}^s\ol y))_{st})}{x,y\in X,\ s,t\in S}.
\]
By identifying $X$ with a subset of $Y$, as above, we also have $R\sub R_1$, as $(w,v)\equiv(w_1,v_1)\in R_1$ for $(w,v)\in R$.
Note also that in $R_2$, we have
\[
(xN({}^s\ol y))_{st} = \begin{cases}
x_{st} &\text{if ${}^s\ol y=1\not\in U$, as then $N({}^s\ol y)=\ew$,}\\
x\cdot N({}^s\ol y)_{st} &\text{otherwise.}
\end{cases}
\]
In the first case, the corresponding relation from $R_2$ is just $(x_sy_t,x_{st})$.  

\begin{lemma}\label{lem:ws}
For any $w\in X^+$ and $s\in S$, we have $w_s\Phi=(\ol w,s)$.  Consequently, $\Phi$ is a surmorphism.
\end{lemma}

\pf
We prove the first assertion by induction on $k=\ell(w)$, the length of $w$.  If $k=1$, then $w\in X$ and $w_s\in Y$, so $w_s\Phi=(\ol w,s)$ by definition.  We now assume that $k\geq2$, so that $w=xv$ for some $x\in X$ and $v\in X^+$ with $\ell(v)=k-1$, and we note that $w_s=x\cdot v_s=x_1\cdot v_s$.  By induction, and monoidality of the action, it follows that
\[
w_s\Phi = (x_1\Phi)\cdot(v_s\Phi) = (\ol x,1)\cdot(\ol v,s) = (\ol x\cdot\ol v,s) = (\ol w,s).
\]
The second assertion now follows, because $(u,s)=(\ol{N(u)},s)=N(u)_s\Phi$ for all $u\in U$ and~${s\in S}$.
\epf

\begin{lemma}\label{lem:R1R2}
We have $R_1\cup R_2\sub\ker(\Phi)$.
\end{lemma}

\pf
Consider first a relation $(w_s,v_s)$ from $R_1$.  So $s\in S$ and $(w,v)\in R$.  The latter gives $\ol w=\ol v$, and it then follows from Lemma \ref{lem:ws} that $w_s\Phi = (\ol w,s) = (\ol v,s) = v_s\Phi$.

Now consider a relation $(x_sy_t,(xN({}^s\ol y))_{st})$ from $R_2$.  Then again using Lemma \ref{lem:ws} we have
\[
(xN({}^s\ol y))_{st}\Phi = (\ol{xN({}^s\ol y)},st) = (\ol x \cdot \ol{N({}^s\ol y)},st) = (\ol x \cdot {}^s\ol y,st) = (\ol x,s)\cdot(\ol y,t) = (x_sy_t)\Phi.  \qedhere
\]
\epf

We now write ${\approx}=(R_1\cup R_2)^\sharp$ for the congruence on $Y^+$ generated by the relations $R_1\cup R_2$.  We also write ${\approx_1}=R_1^\sharp$ and ${\approx_2}=R_2^\sharp$.  Recall that ${\sim}=R^\sharp$ is the congruence on $X^+$ generated by $R$.

\begin{lemma}\label{lem:wsvs}
If $w,v\in X^+$ and $w\sim v$, then $w_s\approx v_s$ for any $s\in S$.
\end{lemma}

\pf
It suffices to assume that $w$ and $v$ differ by a single application of a relation from $R$.  So by symmetry we have $w=aw'b$ and $v=av'b$ for some $a,b\in X^*$ and $(w',v')\in R$.  Now, $R_1$ contains the relations $(w'_s,v'_s)$ and $(w',v')\equiv(w'_1,v'_1)$.  So in the cases $b=\ew$ and $b\not=\ew$, we have
\[
w_s = a\cdot w'_s \approx_1 a\cdot v'_s = v_s
\AND
w_s = a\cdot w' \cdot b_s \approx_1 a\cdot v'\cdot b_s = v_s,
\]
respectively.
\epf

\begin{lemma}\label{lem:wvt}
If $w\in Y^+$, then $w\approx v_t$ for some $v\in X^+$ and $t\in S$.
\end{lemma}

\pf
We prove the lemma by induction on $k=\ell(w)$, the length of $w$.  If $k=1$, then $w\in Y$, so $w=x_s$ for some $x\in X$ and $s\in S$, so we just take $v=x$ and $t=s$.  We now assume that $k\geq2$, so that $w=w'y_s$ for some $w'\in Y^+$, $y\in X$ and $s\in S$.  By induction, since $\ell(w')=k-1$, we have $w'\approx v'_{t'}$ for some $v'\in X^+$ and $t'\in S$.  By definition, $v'_{t'} = v'' x_{t'}$ for some $v''\in X^*$ and $x\in X$.  But then
\[
w = w' y_s \approx v'_{t'} y_s = v''x_{t'} y_s \approx_2 v''(xN({}^{t'}\ol y))_{t's} = (v'' xN({}^{t'}\ol y))_{t's},
\]
and we take $v=v'' xN({}^{t'}\ol y)$ and $t=t's$.
\epf

\begin{lemma}\label{lem:wus}
If $w\in Y^+$ and $w\Phi=(u,s)$, then $w\approx N(u)_s$.
\end{lemma}

\pf
By Lemma \ref{lem:wvt}, we have $w\approx v_t$ for some $v\in X^+$ and $t\in S$.  It then follows from Lemmas~\ref{lem:ws} and \ref{lem:R1R2} that
\[
(u,s) = w\Phi = v_t\Phi = (\ol v,t),
\]
so that $t=s$ and $\ol v=u=\ol{N(u)}$.  It follows from the latter that $v\sim N(u)$.  Combining all of this with Lemma \ref{lem:wsvs} we obtain
\[
w\approx v_t \approx N(u)_t = N(u)_s.  \qedhere
\]
\epf

We can now prove the main result of this section:

\begin{thm}\label{thm:US_sd}
Suppose $U$ is a semigroup, and $S$ a monoid with a monoidal action on $U^1$.  Then with the above notation, the semidirect product $U\rtimes S$ has presentation $\Spres Y{R_1\cup R_2}$ via $\Phi$.
\end{thm}

\pf
By Lemmas \ref{lem:ws} and \ref{lem:R1R2}, it remains only to show that $\ker(\Phi)\sub{\approx}$.  So suppose ${(w,v)\in\ker(\Phi)}$, and write $(u,s)=w\Phi=v\Phi$.  Then by Lemma \ref{lem:wus} we have $w\approx N(u)_s\approx v$.
\epf

\begin{rem}
The presentation $\Spres Y{R_1\cup R_2}$ from Theorem \ref{thm:US_sd} is not as `compact' as that of Lavers' Theorem \ref{thm:Lavers} (which only applies in much more special circumstances):
\bit
\item The generating set $Y$ is essentially $|S|$ copies of $X$, and $R_1$ is $|S|$ copies of $R$.
\item Moreover, $R_2$ contains $|X|^2$ copies of the entire multiplication table of $S$, in the sense that for every $x,y\in X$ and $s,t\in S$, we have a relation of the form $(x_sy_t,w_{st})$.
\eit
In general, however, one cannot hope to obtain a presentation for $U\rtimes S$ in terms of a smaller generating set.  For example, consider the \emph{direct} product $U\times S$, where $U=X^+$ is a free semigroup over a non-empty set $X$, and $S$ is an arbitrary monoid.  For any $x\in X$ and $s\in S$, it is clear that an expression $(x,s) = (u_1,t_1)\cdots(u_k,t_k)$, with each $u_i\in U$ and $t_i\in S$ can only exist with $k=1$ (and $(u_1,t_1)=(x,s)$).  This shows that any generating set for $U\times S$ must contain $X\times S$.  
\end{rem}

\begin{rem}
We have not included any applications of Theorem \ref{thm:US_sd} in the current paper.  While we believe it would be very worthwhile to investigate such applications, especially to almost-(left-)factorisable inverse or left restriction semigroups (cf.~Example \ref{eg:FIM} and Remark~\ref{rem:PS}), it is beyond the scope of the current work.
For example, the papers \cite{JEinsn,JEinsn2,JEptnsn}, which concern the semigroups $\Sing(\I_n)$ and $\PT_n\sm\T_n$ discussed above, total more than~$60$ pages.
As another relevant comparison, the 41-page paper \cite{FADEG2019} is essentially devoted to a single semigroup of the form $U\rtimes S$, with $U$ a monoid and $S$ a non-monoid semigroup (the opposite configuration to that covered by Theorem \ref{thm:US_sd}).
A number of wreath products including $M\wr\Sing(\I_n)$ are treated in \cite{CE2022}, using entirely different methods.
We leave it as an open problem to investigate the above-mentioned applications of Theorem \ref{thm:US_sd}.
\end{rem}

We conclude the chapter with the following observation.

\begin{cor}\label{cor:US_sd}
Suppose $U$ is a semigroup, and $S$ a monoid with a monoidal action on $U^1$.  If~$U$ is finitely presented, and if $S$ is finite, then $U\rtimes S$ is finitely presented.
\end{cor}

\pf
Examining the definitions, we see that the assumptions ensure that $Y$, $R_1$ and $R_2$ are all finite, so the result follows immediately from Theorem \ref{thm:US_sd}.
\epf

\begin{rem}
Corollary \ref{cor:US_sd} gave a sufficient condition for a semidirect product $U\rtimes S$ (with~$S$ a monoid, acting monoidally on $U^1$) to be finitely presented.  We leave it as an open problem to find necessary and sufficient conditions; for \emph{direct} products $U\times S$, see \cite{RRW1998}.
\end{rem}

\chap{Free left restriction monoids}\label{chap:LR}

We now come to the first of our applications of the general results of Chapter \ref{chap:pres}.  Our goal in this chapter is to obtain a monoid presentation for the free left restriction monoid $\LR_X$ over an arbitrary set $X$.  (Of course $\LR_X$ has trivial presentation $\pres X\es$ in the signature of left restriction monoids.  Monoid presentations will inevitably require more relations, yet presentations in this simpler signature are valuable.)

We recall the definition of $\LR_X$ in Section \ref{sect:LRprelim}, where we also discuss its structure and decomposition in terms of action pairs.  We then give the presentation in Section \ref{sect:LRpres}; see Theorem~\ref{thm:LX}.  As we note in Remark \ref{rem:LX}, this presentation utilises the unique minimum (monoid) generating set for $\LR_X$.  Theorem \ref{thm:PX} gives a presentation for the semilattice of projections of $\LR_X$.

\sect{Preliminaries}\label{sect:LRprelim}

\emph{Left restriction semigroups} were defined in Section \ref{sect:lrs} as a variety of unary semigroups, and consequently free objects exist.  It transpires that the free left restriction semigroup on a set~$X$ coincides with the free left ample semigroup on~$X$; see for example \cite{GG2000}.  The latter form a quasi-variety, and a sub-class of the class of left restriction semigroups.  Note that left ample semigroups were originally called \emph{left type $A$ semigroups} in the literature.  The original description of free \emph{right} type $A$ semigroups was given by Fountain in \cite{Fountain1991b}, and these are in fact subsemigroups of the free inverse semigroup \cite{Munn1974,Scheiblich1973} over the same base set.  See \cite{Gould1996,GG2000,FGG2009,Kambites2011} for connections with free objects in other varieties and quasi-varieties.  An equivalent form of Fountain's construction, but taking the left-right dual, is given below.  The free left restriction \emph{monoid} is simply the monoid completion of the corresponding semigroup, so we treat monoids here for convenience.  All of the results in this section can be translated to results for free left restriction \emph{semigroups} by replacing monoid presentations with semigroup presentations (with the same generators and relations).  

Let $X$ be an arbitrary set, and as usual let $X^*$ be the free monoid over $X$, with empty word denoted $\ew$.  For $w\in X^*$ we write $w^\da$ for the set of all prefixes of $w$ (including $\ew$ and $w$), and for $A\sub X^*$ we write $A^\da=\bigcup_{w\in A}w^\da$ for the set of all prefixes of all words from $A$.  Let
\[
\P_X = \set{A\sub X^*}{A^\da=A,\ 0<|A|<\aleph_0}
\]
be the set of all non-empty, finite, prefix-closed subsets of $X^*$.  Then $\P_X$ is a monoid semilattice under~$\cup$, with identity $\{\ew\}$.  The \emph{free left restriction monoid over $X$} is the set
\[
\LR_X = \bigset{(A,w)}{A\in \P_X,\ w\in A} ,
\]
with:
\bit
\item product $(A,w)\cdot(B,v)=(A\cup wB,wv)$, where $wB=\set{wv}{v\in B}$, and
\item unary operation $(A,w)^+=(A,\ew)$.  
\eit
The monoid $\LR_X$ contains isomorphic copies of both $\P_X$ and $X^*$, which we identify with the submonoids
\[
\P_X \equiv P(\LR_X) = \bigset{(A,\ew)}{A\in \P_X} \AND X^* \equiv \bigset{(w^\da,w)}{w\in X^*}.
\]
By identifying $A\in\P_X$ and $w\in X^*$ with $A\equiv(A,\ew)$ and $w\equiv(w^\da,w)$, we have 
\begin{align}
\label{eq:Aw} A\cdot w & \equiv(A\cup w^\da,w) &&\hspace{-2cm}\text{for all $A\in\P_X$ and $w\in X^*$.}
\intertext{In particular,}
\nonumber (A,w) &\equiv A\cdot w &&\hspace{-2cm}\text{for all $(A,w)\in\LR_X$.}
\end{align}
It follows that $\LR_X=\P_X\cdot X^*$.

It of course follows from Proposition \ref{prop:LR1} (and the identification $\P_X\equiv P(\LR_X)$) that $(\P_X,\LR_X)$ is an action pair in $\LR_X$, though this does not help us to find a presentation; cf.~Remark \ref{rem:LR2}.  However, Proposition \ref{prop:LR2}\ref{LR21} immediately gives us the following action pair, which will prove to be much more useful:

\begin{prop}\label{prop:LRX}
For any set $X$, $(\P_X,X^*)$ is an action pair in $\LR_X$, and $\LR_X=\P_X\cdot X^*$.  \epfres
\end{prop}

Before turning to presentations, we conclude this section with a series of remarks considering how the monoid~$\LR_X$ fits into the context of the various results and constructions from Part \ref{part:I} of the paper.

\begin{rem}\label{rem:LRX1}
It is instructive to consider a direct proof of Proposition \ref{prop:LRX}, starting from Definition \ref{defn:AP}:
\begin{enumerate}[label=\textup{(A\arabic*)}]
\item The required action of $X^*$ on $\P_X$ is given by
\begin{equation}\label{eq:uA}
{}^w\!A = (wA)^\da = w^\da\cup wA \qquad\text{for $A\in \P_X$ and $w\in X^*$.}
\end{equation}
It is then routine to show that
\[
{}^{wv}\!A={}^w({}^v\!A) \COMMA {}^w(A\cup B)={}^w\!A\cup {}^w\!B \AND w\cdot A={}^w\!A\cdot w ,
\]
for all $A,B\in \P_X$ and $w,v\in X^*$.  (We also of course have ${}^\ew\!A=A$ for all $A\in\P_X$.)
\item As noted above, the identity of $\P_X$ is $\{\ew\}$, and we have
\begin{equation}\label{eq:wda}
w^+ = {}^w\{\ew\}=w^\da \qquad\text{for any $w\in X^*$.}
\end{equation}
We then quickly deduce the required implication
\[
A\cdot w = B\cdot v \implies A\cdot w^+ = B\cdot v^+ \qquad\text{for all $A,B\in \P_X$ and $w,v\in X^*$}
\]
from the expansions $A\cdot w=(A\cup w^\da,w)$, $A\cdot w^+=(A\cup w^\da,\ew)$, and so on.
\een
It also follows from \eqref{eq:wda} that the associated map ${X^*\to\P_X:w\mt w^+}$ from Proposition \ref{prop:+} is given by ${w^+=w^\da\equiv(w^\da,\ew)}$.  This is of course just the restriction to $X^*$ of the ${}^+$ operation on~$\LR_X$.
\end{rem}

\begin{rem}\label{rem:LRX2}
It follows from Lemma \ref{lem:SAPAP}\ref{SAPAP2} that the pair $(\P_X,X^*)$ is not strong, as $w^\da\not=\{\ew\}$ when $w\in X^*$ is non-empty.
Hence, the left-uniqueness property \ref{SA2} does not hold.  In contrast to this, it follows immediately from \eqref{eq:Aw} that the pair $(\P_X,X^*)$ has the \emph{right-uniqueness} property:
\[
A\cdot w = B\cdot v \implies w=v \qquad\text{for all $A,B\in \P_X$ and $w,v\in X^*$.}
\]
One important consequence of this is that the congruence $\si=\si(\P_X,X^*)$ on $X^*$ from Definition~\ref{defn:si} is trivial: i.e., $\si=\De_{X^*}$.  Together with Lemma \ref{lem:De} (and $P=U=\P_X$), it follows immediately from this that the pair $(\P_X,X^*)$ is proper.  Consequently, $\LR_X=\P_X\cdot X^*$ is $(\P_X,X^*)$-proper.
\end{rem}

\begin{rem}\label{rem:LRX3}
It does not follow from Remark \ref{rem:LRX2} that $\LR_X$ is a proper left restriction monoid in the sense of Definition \ref{defn:lr_proper2}, although this does turn out to be the case.  (To the best of our knowledge, this fact has not been explicitly stated in the literature, but it follows quickly from known results; see for example \cite[Proposition 3.3]{GG2000}.)  To demonstrate this we could show that the pair $(\P_X,\LR_X)$ is proper, and then apply Proposition~\ref{prop:lr}.  But this is no easier than directly using Definition \ref{defn:lr_proper2}.  To use this, we need to understand the congruence $\si=\si_{\LR_X} (=\si(\P_X,\LR_X))$ on $\LR_X$, and it is easy to see that for any elements $(A,w)$ and $(B,v)$ of $\LR_X$,
\begin{equation}\label{eq:siLX}
(A,w)\mr\si(B,v) \iff w=v .
\end{equation}
Indeed, the forwards implication is clear, while if $w=v$, then $(B,\ew)\cdot(A,w)=(A,\ew)\cdot(B,v)$.  We then have
\begin{align*}
(A,w) = (B,v) \quad\iff\quad A=B \text{ and } w=v 
\quad\iff\quad (A,w)^+=(B,v)^+ \text{ and } (A,w)\mr\si(B,v).
\end{align*}
\end{rem}

\begin{rem}\label{rem:LRX4}
Consider again the pair $(\P_X,X^*)$.  It follows from Proposition \ref{prop:A1} that $\LR_X$ is a homomorphic image of the semidirect product
\[
\P_X\rtimes X^* = \bigset{(A,w)}{A\in \P_X,\ w\in X^*} \qquad\text{with operation}\qquad (A,w)\cdot(B,v)=(A\cup {}^w\!B,wv),
\]
where the action is given in \eqref{eq:uA}.  Explicitly, the surmorphism $\pi:\P_X\rtimes X^*\to\LR_X$ from Proposition \ref{prop:A1} is given by
\[
(A,w)\pi = A\cdot w=(A\cup w^\da,w).
\]
Note that $\P_X\rtimes X^*$ is not a monoid, even though $\P_X$ and $X^*$ both are.  This is because the action of $X^*$ on $\P_X$ is not by monoid morphisms (cf.~\eqref{eq:wda}), even though it is of course monoidal; cf.~Lemma \ref{lem:USmon}.
On the other hand, it follows immediately from its definition that $\LR_X$ is also the \emph{subsemigroup} of $\P_X\rtimes X^*$ consisting of all pairs $(A,w)$ with $w\in A$.  (For this, note that when $w\in A$, we have~${A\cup{}^w\!B=A\cup wB}$.)  

This is of course reminiscent of the embedding theorems from Chapter \ref{chap:PAP}.  Indeed, Theorem~\ref{thm:embed2} tells us that $\LR_X=\P_X\cdot X^*$ embeds in a semidirect product $\U\rtimes X^*$ for some semigroup~$\U$ containing $\P_X$.  (For this we also need to remember that $\si(\P_X,X^*)=\De_{X^*}$; cf Remark \ref{rem:LRX2}.)  As noted in the previous paragraph, we can in fact take $\U$ to be $\P_X$ itself.

Alternatively, we could apply Theorem \ref{thm:embed2} to the pair $(\P_X,\LR_X)$, or even apply Theorem~\ref{thm:OC} to $\LR_X$ itself.  In either case, the relevant theorem tells us that $\LR_X(=\P_X\cdot\LR_X)$ embeds in a semidirect product $\U\rtimes(\LR_X/\si)$ for some semigroup $\U$ containing $\P_X$, where $\si=\si_{\LR_X}=\si(\P_X,\LR_X)$.  It is easy to see from \eqref{eq:siLX} that $\LR_X/\si \cong X^*$, as has also been shown in \cite[Theorem~5.1(iii)]{BGG2011}.  
\end{rem}

\begin{rem}\label{rem:LRX5}
Finally, we consider the covering Theorem \ref{thm:cover} in the context of~$\LR_X$.  Since~$\LR_X$ is proper, it is trivially covered by itself, and it turns out that this is precisely what we recover from our proof of Theorem \ref{thm:cover}.  Indeed, the $(\ul U,\ol S)$ construction applied to the pair ${(U,S)=(\P_X,X^*)}$ leads to the (monoid) subsemigroups
\[
\ul{\P_X} = \bigset{(A,\ew)}{A\in\P_X} \equiv \P_X \AND \ol{X^*} = \bigset{(w^\da,w)}{w\in X^*} \equiv X^*
\]
of the semidirect product $U\rtimes S(=U^1\rtimes S^1)=\P_X\rtimes X^*$.  The proper monoid
\[
\ul U\cdot\ol S = \bigset{(u,s)\in U\rtimes S}{u=us^+} = \bigset{(A,w)\in\P_X\rtimes X^*}{A=A\cup w^+}
\]
is then precisely $\LR_X$ itself, since
\[
A = A\cup w^+ \iff A = A\cup w^\da \iff w\in A \qquad\text{for all $A\in\P_X$ and $w\in X^*$.}
\]
\end{rem}

\sect{Presentations}\label{sect:LRpres}

Our goal now is to give a monoid presentation for the free left restriction monoid~$\LR_X$ over an arbitrary set $X$.  We keep the notation from the previous section.  In particular, $\LR_X = \P_X\cdot X^*$ arises from the action pair $(\P_X,X^*)$, where $\P_X$ and $X^*$ are identified with submonoids of $\LR_X$, as explained above.  We \emph{could} apply Theorem \ref{thm:ESmon} to this pair, but it actually follows from Remark~\ref{rem:LRX5} that Theorem \ref{thm:M} applies, and this will lead us to a slightly quicker derivation.  In any case, we need presentations for $\P_X$ and $X^*$; with the latter of course being trivial, we turn to the former.

The next result gives a monoid generating set for $\P_X$.  Recall that an \emph{atom} of a semigroup~$S$ is an element $a$ such that $a=xy \implies a\in\{x,y\}$ for all $x,y\in S$.  It is easy to see that any generating set for a semigroup must contain all atoms.

\begin{prop}\label{prop:GaX}
For any set $X$, $\Ga_X=\set{w^\da}{w\in X^+}$ is the (unique) minimum monoid generating set for~$\P_X$.
\end{prop}

\pf
If $A\in\P_X$, then
\[
A = A^\da = \bigcup_{w\in A\sm\{\ew\}}w^\da \in \la\Ga_X\ra,
\]
keeping $|A|<\aleph_0$ in mind.  This shows that $\P_X=\la\Ga_X\ra$.

To complete the proof, it suffices to show that every element of $\Ga_X$ is an atom.  So let $w\in X^+$, and suppose $w^\da=A\cup B$ for some $A,B\in\P_X$.  Then $A,B\sub w^\da$.  
Since $w\in w^\da=A\cup B$, we may assume by symmetry that $w\in A$, and then $w^\da\sub A^\da=A$, so that $A=w^\da$.
\epf

We now wish to give a presentation for $\P_X$ in terms of the generating set $\Ga_X$.  With this in mind, we define an alphabet
\[
A_X = \set{a_w}{w\in X^+}
\]
in one-one correspondence with $\Ga_X$.  By Proposition \ref{prop:GaX}, we have a surmorphism
\[
\phi_X : A_X^* \to \P_X : a_w \mt w^\da.
\]
We then define $R_X=R_X^1\cup R_X^2$, where
\[
R_X^1 = \bigset{(a_wa_v,a_va_w)}{w,v\in X^+} \AND R_X^2 = \bigset{(a_{wv}a_w,a_{wv})}{w\in X^+,\ v\in X^*}.
\]
Note that we have the alternative description
\[
R_X^2 = \bigset{(a_wa_v,a_w)}{w,v\in X^+,\ v\text{ is a prefix of } w}.
\]

\begin{thm}\label{thm:PX}
For any set $X$, $\P_X$ has presentation $\Mpres{A_X}{R_X}$ via $\phi_X$.
\end{thm}

\pf
We have already noted that $\phi_X$ is a surmorphism, and it is clear that ${R_X\sub\ker(\phi_X)}$.  Indeed, we have $R_X^1\sub\ker(\phi_X)$ because $\P_X$ is commutative, and $R_X^2\sub\ker(\phi_X)$ follows from the observation that
\[
(wv)^\da \cup w^\da = (wv)^\da \qquad\text{for all $w,v\in X^*$.}
\]
It therefore remains to show that $\ker(\phi_X)\sub{\sim}$, where ${\sim}=R_X^\sharp$.  So suppose $(p,q)\in\ker(\phi_X)$, and write
\[
p = a_{w_1}\cdots a_{w_k} \AND q = a_{v_1}\cdots a_{v_l} \WHERE w_1,\ldots,w_k,v_1,\ldots,v_l\in X^+.
\]
For any $1\leq i\leq l$, we have
\[
v_i \in v_i^\da \sub v_1^\da\cup\cdots\cup v_l^\da = q\phi_X = p\phi_X = w_1^\da\cup\cdots\cup w_k^\da.
\]
It follows from this that $v_i\in w_j^\da$ for some $1\leq j\leq k$, which means that $v_i$ is a prefix of $w_j$.  But then $R_X^2$ contains the relation $(a_{w_j}a_{v_i},a_{w_j})$.  Combining this with $R_X^1$, it follows that $p\sim pa_{v_i}$.  Since this is true for all $1\leq i\leq l$, it follows that $p\sim pa_{v_1}\cdots a_{v_l} = pq$.  By symmetry we also have $q\sim qp$.  Since $pq\sim qp$ by $R_X^1$, it follows that $p\sim q$, as required.
\epf

We now apply Theorem \ref{thm:M} to obtain a (monoid) presentation for $\LR_X$.

\begin{thm}\label{thm:LX}
For any set $X$, the free left restriction monoid $\LR_X$ has presentation
\[
\Mpres{A_X\cup X}{R_X\cup R_1\cup R_2}
\]
via
\[
\Phi_X : (A_X\cup X)^* \to \LR_X : \begin{cases}
a_w \mt w^\da \equiv(w^\da,\ew) &\text{for $w\in X^+$,}\\
x \mt x \equiv (x^\da,x) &\text{for $x\in X$,}
\end{cases}
\]
where
\[
R_1 = \bigset{(xa_w,a_{xw}x)}{x\in X,\ w\in X^+} \AND R_2 = \bigset{(x,a_xx)}{x\in X}.
\]
\end{thm}

\pf
We begin with the presentation $\Mpres{A_X}{R_X}$ for $\P_X$ from Theorem \ref{thm:PX}, and the canonical presentation $\Mpres X\es$ for $X^*$ (via the identity map $X^*\to X^*$).  Theorem \ref{thm:M} then tells us that $\LR_X=\bigset{(A,w)\in\P_X\rtimes X^*}{A=A\cup w^+}$ has presentation 
\[
\Mpres{A_X\cup X}{R_X\cup R_1\cup R_2}
\]
via $\Phi_X$, where
\[
R_1 =  \bigset{(xa_w,{}^x(a_w)\cdot x)}{x\in X,\ w\in X^+} \AND R_2 = \bigset{(x,x^+x)}{x\in X}.
\]
(The second set of relations was denoted $R_3$ in Theorem \ref{thm:M}.)
In $R_1$, ${}^x(a_w)$ is a word mapping to ${}^x(w^\da) = (xw)^\da$, so of course we can take ${}^x(a_w)=a_{xw}$.  Similarly, we can take $x^+=a_x$ in~$R_2$.
\epf

\begin{rem}\label{rem:LX}
The presentation in Theorem \ref{thm:LX} involves the generating set $\Ga_X\cup X$ for $\LR_X$.  It is easy to see (using Proposition \ref{prop:GaX}) that this is the unique minimum generating set.
\end{rem}

\chap{Independence algebras}\label{chap:IA}

In this chapter we apply the theory developed so far to a number of monoids and semigroups of partial endomorphisms of independence algebras.  In Sections \ref{sect:Aprelim} and \ref{sect:EndA} we recall the necessary definitions, and prove some preliminary results on independence algebras, and the partial endomorphism monoid~$\PEnd(A)$ of such an algebra $A$.  A key point is that $\PEnd(A)$ is a left restriction monoid.  Section \ref{sect:MaxA} contains a number of results concerning maximal subalgebras, and Section \ref{sect:fca} treats subalgebras of arbitrary finite codimension; the main result here is Theorem \ref{thm:SubA}, which gives a presentation for the $\cap$-semilattice $\FSub(A)$ of all finite-codimensional subalgebras of a strong independence algebra $A$.  Section \ref{sect:AutA} provides natural generating sets for the automorphism groups $\Aut(A)$ in the case that~$A$ is finite-dimensional, and also for the group $\FAut(A)$ of finitary automorphisms of an arbitrary strong independence algebra.
In Section~\ref{sect:AES} we identify six strong action pairs $(E,S)$ in $M=\PEnd(A)$, each leading in the usual way to a naturally occurring subsemigroup~$ES$ of $\PEnd(A)$, the elements of which we also characterise.  We then look, in Section~\ref{sect:Acong}, at the various `congruence conditions' from Section~\ref{sect:cong} and Chapter~\ref{chap:pres}, as they apply to each of our pairs.  Finally, we discuss presentations for the semigroups~$ES$ in Section \ref{sect:EA}.  Presentations for $S$ are not known in general (even for relatively simple classes of the algebra~$A$), so it is not possible to give explicit presentations for $ES$ in full generality.  However, a presentation for~$E$ can be deduced from Theorem \ref{thm:SubA} in the case that~$A$ is finite-dimensional and strong, as $E$ is then isomorphic to $\Sub(A)=\FSub(A)$.  We conclude by explaining how the presentation for $E$ can be extended to a presentation for $ES$ modulo an appropriate presentation for $S$.  For background on universal algebra, see \cite{BS1981}.

\sect{Preliminaries}\label{sect:Aprelim}

Let $A$ be a (universal) algebra.  We typically identify $A$ with its underlying set, and as usual we identify any nullary operation on $A$ with a constant: i.e., the unique element in the image of the operation.  
If $X\subseteq A$ then $\langle X\rangle$ denotes the subalgebra generated by $X$.  The least subalgebra of~$A$ will be denoted
\[
\CA=\la\es\ra.
\]
Note that $\CA=\es$ if and only if $A$ has no constants.  
We write $B\leq A$ ($B<A$) to indicate that $B$ is a (proper) subalgebra of~$A$, and we write $\Sub(A)=\set{B}{B\leq A}$ for the lattice of all subalgebras.  The meet and join of $B,C\in\Sub(A)$ are $B\cap C$ and $B\vee C=\la B\cup C\ra$, respectively.

A subset $X\sub A$ is \emph{independent} if for any $x\in X$ we have $x\notin \langle X\setminus x\rangle$.  (Here and elsewhere in this chapter, we avoid clutter by identifying one-element subsets of $A$ with their unique element, so $X\sm x$ stands for $X\sm\{x\}$, and so on.)  A \emph{basis} is an independent generating set for $A$.
Note that a singleton subset $\{x\}$ is independent if and only if $x\not\in\CA$.

\begin{defn}\label{defn:IA}
Let $A$ be a (universal) algebra. Then $A$ is an \emph{independence algebra} if it satisfies the \emph{exchange property} and the \emph{free basis property}:
\begin{enumerate}[label=\textup{(EP)},leftmargin=13mm]
\item \label{EP} For all $X\sub A$ and $x,y\in A$, if $x\in \langle X\cup  y\rangle$ and $x\notin \langle X\rangle$, then $y\in \langle X\cup  x\rangle$.
\een
\begin{enumerate}[label=\textup{(FBP)},leftmargin=13mm]
\item \label{FBP} If $X$ is basis for $A$, then any map $X\rightarrow A$ can be extended (uniquely) to an endomorphism of $A$.
\een
\end{defn}

Archetypal examples of independence algebras include:
\bit
\item vector spaces over division rings, where we have a binary operation (addition), a unary operation (negation), a nullary operation (the zero), and an additional unary operation for each element of the ground ring (scalar multiplication), 
\item sets, regarded as algebras with no operations, and 
\item free group-acts, including the special case of sets when the group is trivial; see Example~\ref{eg:FXM}.
\eit
As noted in \cite{Gould1995}, independence algebras were first studied under the name of \emph{$v^*$-algebras}.  They appear first in an article of Narkiewicz \cite{Narkiewicz1961}, and were inspired by Marczewski's study of notions of independence, initiated in \cite{Marczewski1959}.  A classification into five families, together with one finite exception, was given by Urbanik in \cite{Urbanik1966}.

The next (well known) result is a simple consequence of \ref{EP}.

\begin{lemma}\label{lem:Xy}
If $X$ is an independent subset of an independence algebra $A$, then $X\cup y$ is independent for any $y\in A\sm\la X\ra$.  \epfres
\end{lemma}

The next result follows from \cite[Lemma 3.7]{Gould1995}, and will typically be used without explicit reference.

\begin{lemma}\label{lem:auto}
If $X$ and $Y$ are bases of an independence algebra $A$, then any endomorphism of $A$ extending a bijection $X\to Y$ is an automorphism.  \epfres
\end{lemma}

The following facts concerning an arbitrary independence algebra $A$ are well known; see for example \cite[Section 3]{Gould1995}.
\bit
\item Every independent subset of $A$ is contained in a maximal independent subset.  Since $\es$ is independent, maximal independent subsets exist.  Any maximal independent subset generates $A$, and is hence a basis.  So any independent set can be extended to a basis.
\item All bases for $A$ have the same cardinality, which we call the \emph{dimension} of $A$, and denote~$\dim(A)$.  This coincides with the usual meaning of dimension when $A$ is a vector space; when $A$ is a set (with no operations), $\dim(A)=|A|$.  Note that the dimension of~$A$ is also called the \emph{rank} of $A$ in the literature, and denoted $\rank(A)$.
\item Any subalgebra $B\leq A$ is also an independence algebra, and any basis of $B$ may be extended to a basis for $A$. Thus, $\dim(B)$ exists and $\dim(B)\leq \dim(A)$.  
\item If $\dim(A)=n<\aleph_0$, then $A$ is the unique subalgebra of dimension $n$, and the maximal (proper) subalgebras are precisely those of dimension $n-1$.  (Recall that $B\leq A$ is maximal if $B<C\leq A\implies C=A$.)  When $\dim(A)\geq\aleph_0$ we require the notion of \emph{codimension} to describe the maximal subalgebras; see Definition \ref{defn:codim} and Proposition \ref{prop:Max}.
\item If $X\sub A$ is any independent set, then any map $X\to A$ can be extended to an endomorphism of $A$.  (First extend $X$ to a basis $X\sqcup Y$ of $A$, and extend the given map $X\to A$ arbitrarily to $X\sqcup Y\to A$; by \ref{FBP}, the latter extends to an endomorphism.)
\eit
Here and elsewhere, $\sqcup$ denotes disjoint union.

The notion of \emph{codimension} of subalgebras was used in \cite{FL1993}, and goes back to ideas from Jones \cite{Jones1989}.  The following definition is a special case of that given in \cite{Jones1989}, which involved more general \emph{closure operators} on arbitrary sets.  

\begin{defn}\label{defn:Bind}
Let $B$ be a subalgebra of an independence algebra $A$.  A subset $X$ of $A$ is \emph{$B$-independent} if $x\not\in\big\la B\cup (X\sm x)\big\ra$ for all $x\in X$.  The set $X$ is then called a \emph{$B$-basis} for~$\la B\cup X\ra$.  
\end{defn}

Clearly any $B$-independent set is disjoint from $B$.  Note that independent sets are precisely the $\CA$-independent sets in the sense of Definition \ref{defn:Bind}, where as usual $\CA=\la\es\ra$.  The next result is a special case of the main theorem of \cite{Jones1989}:

\begin{lemma}
For any subalgebra $B$ of an independence algebra $A$, all $B$-bases for $A$ have the same cardinality.  \epfres
\end{lemma}

We may therefore make the following definition:

\begin{defn}\label{defn:codim}
Given a subalgebra $B$ of an independence algebra $A$, the \emph{codimension of $B$ in~$A$}, denoted $\codim(B)$, is defined to be the cardinality of any $B$-basis of $A$.
\end{defn}

Strictly speaking, one should perhaps use notation such as $\codim_A(B)$, to explicitly indicate the over-algebra $A$.  But since the context should always be clear, we simply write $\codim(B)$ to avoid clutter.

We will also need the following basic facts.

\begin{lemma}\label{lem:XY}
If $B$ is a subalgebra of an independence algebra $A$, then $X\cup Y$ is a basis of $A$ for any basis $X$ of $B$, and any $B$-basis $Y$ of $A$.
\end{lemma}

\pf
Since $A=\la B\cup Y\ra=\la\la X\ra\cup Y\ra=\la X\cup Y\ra$, it remains to show that $X\cup Y$ is independent.  For this, and aiming for a contradiction, suppose $z\in\big\la(X\cup Y)\sm z\big\ra$ for some $z\in X\cup Y$.  

\pfcase1  If $z\in Y$, then $z\in\big\la X\cup (Y\sm z)\big\ra=\big\la B\cup (Y\sm z)\big\ra$, contradicting $B$-independence of $Y$.

\pfcase2  Now suppose $z\in X$, so that $z\in\big\la(X\sm z)\cup Y\big\ra$.  
Choose finite subsets $X'\sub X\sm z$ and $Y'\sub Y$ such that $z\in\la X'\cup Y'\ra$, and suppose $|X'|+|Y'|$ is minimal among all such $X',Y'$.  Since~$X$ is a basis of $B$, it follows from independence that $z\not\in\la X'\ra$, and so $Y'\not=\es$.  Then for any $y\in Y'$, and with $Z=(X'\cup Y')\sm y$, we have
\[
z \in \la X'\cup Y'\ra = \la Z\cup y\ra \AND z\not\in \la Z\ra,
\]
using minimality of $|X'|+|Y'|$ for the latter.  It then follows from \ref{EP} that
\[
y\in\la Z\cup z\ra\sub\big\la(X\cup Y)\sm y\big\ra.
\]
Since $y\in Y'\sub Y$, we have therefore reduced to Case 1.
\epf

It follows that $\dim(B)+\codim(B)=\dim(A)$ for all $B\leq A$.  In particular, if $\dim(A)=n<\aleph_0$, then $\codim(B)=n-\dim(B)$ for all $B\leq A$.
We will also use the following converse of Lemma~\ref{lem:XY}, often without reference.

\begin{lemma}\label{lem:BAXY}
If $B$ is a subalgebra of an independence algebra $A$, and if $X\sqcup Y$ is independent for some basis $X$ of $B$, then $Y$ is $B$-independent.
\end{lemma}

\pf
Aiming to prove the contrapositive, suppose $Y$ is not $B$-independent, so ${y\in \big\la B\cup(Y\sm y)\big\ra}$ for some $y\in Y$.  But $B=\la X\ra$, so $y\in \big\la B\cup(Y\sm y)\big\ra=\big\la X\cup(Y\sm y)\big\ra=\big\la (X\cup Y)\sm y\big\ra$, where we used the fact that $X$ and $Y$ are disjoint in the last equality.  This shows that $X\cup Y$ is not independent.
\epf

The first part of the next result is \cite[Lemma 1.3]{FL1993}.

\begin{lemma}\label{lem:BCA}
Let $A$ be an independence algebra.
\ben
\item \label{BCA1} If $B\leq C\leq A$, then $\codim(B)\geq\codim(C)$.
\item \label{BCA2} If $B\leq C\leq A$ and $\codim(B)=\codim(C)<\aleph_0$, then $B=C$.
\een
\end{lemma}

\pf
Let $X$ be a basis of $B$, and extend this to bases $X\sqcup Y$ and $X\sqcup Y\sqcup Z$ of $C$ and $A$, respectively.  
By Lemma \ref{lem:BAXY}, $Y\sqcup Z$ is a $B$-basis of $A$, and $Z$ is a $C$-basis of $A$, so
\[
\codim(B)=|Y\sqcup Z| \AND \codim(C)=|Z|.
\]
Part \ref{BCA1} immediately follows.  For \ref{BCA2}, we have
\[
\codim(B)=\codim(C)<\aleph_0 \implies |Y\sqcup Z|=|Z|<\aleph_0 \implies Y=\es \implies B=\la X\ra=\la X\sqcup Y\ra=C.  \qedhere
\]
\epf

Some of the results in Chapter \ref{chap:IA} apply only to so-called \emph{strong} independence algebras:

\begin{defn}\label{defn:strong}
An independence algebra $A$ is \emph{strong} if it satisfies the \emph{strong property}:
\begin{enumerate}[label=\textup{(SP)},leftmargin=13mm]
\item \label{SP} For independent subsets $X,Y\sub A$, if $\la X\ra\cap\la Y\ra=\CA$, then $X\cup Y$ is independent.
\een
\end{defn}

\begin{rem}\label{rem:NS1}
Sets, vector spaces and free group-acts are all strong independence algebras.  It is also easy to show that every independence algebra of dimension $\leq2$ is strong.

A specific example of a three-dimensional non-strong independence algebra is given in \cite[Example~1.5]{FL1993}, and it will be convenient to recall the construction here.  Let $A=\{a_1,a_2,a_3,a_4\}$, with a ternary operation $f$ defined by
\begin{align*}
f(a_i,a_i,a_i) &= a_i &&\text{for all $i$}\\
f(a_i,a_i,a_j) = f(a_i,a_j,a_i) = f(a_j,a_i,a_i) &= a_j &&\text{if $i,j$ are distinct}\\
f(a_i,a_j,a_k) &= a_l &&\text{if $i,j,k,l$ are distinct.}
\end{align*}
All proper subsets of $A$ are independent, but $A$ is not itself independent, and~\ref{SP} fails (for example) for $X=\{a_1,a_2\}$ and $Y=\{a_3,a_4\}$.
\end{rem}

Of crucial importance to us is the following result, which is \cite[Lemma~1.6]{FL1993}.  In fact, the property in the lemma can be taken as an equivalent definition of strong independence algebras.

\begin{lemma}\label{lem:strong}
Let $A$ be a strong  independence algebra, and let $B,C\leq A$.  Let $X$ be a basis for $B\cap C$, and extend $X$ to bases $X\sqcup Y$ and $X\sqcup Z$ for $B$ and $C$, respectively.  Then $X\sqcup Y \sqcup Z$ is a  basis for $B\vee C$.  \epfres
\end{lemma}

We will also need the following simple result, which is a generalisation of the inclusion-exclusion formula for sets.

\begin{lemma}\label{lem:BC1}
If $B$ and $C$ are subalgebras of a strong independence algebra $A$, then 
\ben
\item \label{BC11} $\dim(B\vee C) + \dim(B\cap C) = \dim(B)+\dim(C)$,
\item \label{BC12} $\codim(B\vee C) + \codim(B\cap C) = \codim(B)+\codim(C)$.
\een
\end{lemma}

\pf
Let $X$ be a basis for $B\cap C$, and extend this to bases $X\sqcup Y$ and $X\sqcup Z$ for $B$ and $C$.  By Lemma \ref{lem:strong}, $X\sqcup Y\sqcup Z$ is a basis for $B\vee C$, and we extend this to a basis $X\sqcup Y\sqcup Z\sqcup W$ for $A$.  Both parts follow quickly upon writing down expressions for the various (co)dimensions in terms of the cardinalities of $X$, $Y$, $Z$ and $W$. 
\epf

\begin{rem}\label{rem:NS2}
The assumption that $A$ is strong cannot be removed from Lemma \ref{lem:BC1}.  For example, let $A$ be the algebra from \cite{FL1993} discussed in Remark~\ref{rem:NS1}.  Then for the subalgebras $B=\{a_1,a_2\}$ and $C=\{a_3,a_4\}$, we have
\[
\dim(B)=\dim(C)=2 \COMMA \dim(B\cap C)=0 \AND \dim(B\vee C)=\dim(A)=3,
\]
which witnesses the failure of the first identity from Lemma \ref{lem:BC1} in $A$.  (The same (sub)algebras work to demonstrate the failure of the second identity.)
\end{rem}

\sect{(Partial) endomorphisms}\label{sect:EndA}

A \emph{partial endomorphism} of an algebra $A$ is a morphism $B\to A$ for some subalgebra $B\leq A$; if this morphism is injective, it is a \emph{partial automorphism}, and it is then an isomorphism from its domain onto its image.  We also have \emph{endomorphisms} and \emph{automorphisms}, with their usual meaning.  We write
\begin{equation}\label{eq:PEndA}
\PEnd(A) \COMMA \End(A) \COMMA \PAut(A) \AND \Aut(A)
\end{equation}
for the sets of all partial endomorphisms, endomorphisms, partial automorphisms and automorphisms of $A$, respectively.  These are all monoids under composition, with identity element $\id_A$.  Note that for $\al,\be\in\PEnd(A)$, we have
\[
\dom(\al\be) = (\im(\al)\cap\dom(\be))\al^{-1} \AND \im(\al\be) = (\im(\al)\cap\dom(\be))\be,
\]
both of which are subalgebras of $A$ by standard facts from universal algebra.

The inclusions among the monoids listed in \eqref{eq:PEndA} are indicated as follows:
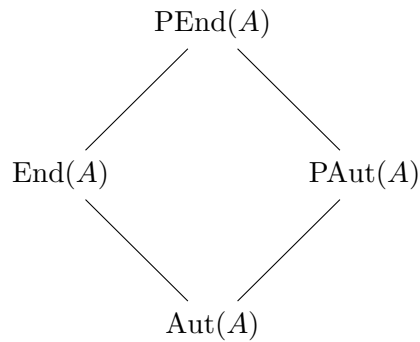
\[
\begin{tikzpicture}[scale=2]
\node (PE) at (0,2) {$\PEnd(A)$};
\node (E) at (-1,1) {$\End(A)$};
\node (PA) at (1,1) {$\PAut(A)$};
\node (A) at (0,0) {$\Aut(A)$};
\draw (PE)--(E)--(A)--(PA)--(PE);
\end{tikzpicture}
\]
Moreover, $\PAut(A)$ is an inverse monoid, and $\Aut(A)$ is the group of units of all of the above monoids.  
If $A$ is an independence algebra, then we also have $\End(A)\vee\PAut(A)=\PEnd(A)$, as follows from Theorem \ref{thm:ESA2}\ref{ESA21} below.  On the other hand, $\End(A)\cap\PAut(A)$ consists of all injective endomorphisms of $A$; this is equal to $\Aut(A)$ if and only if $\dim(A)<\aleph_0$.
When $A$ is simply a set (an algebra with no operations), the monoids listed in \eqref{eq:PEndA} are $\PT_A$, $\T_A$, $\I_A$ and~$\G_A$, respectively.  These monoids of (partial) transformations were introduced in Section \ref{sect:ts}.

The next result is of fundamental importance.  It is well known, and follows immediately from the fact that left restriction monoids form a variety.  Specifically, $\PT_A$ is a left restriction monoid under $\al^+=\id_{\dom(\al)}$, and the submonoid ${\PEnd(A)\leq\PT_A}$ is closed under ${}^+$.

\begin{prop}\label{prop:LRIA}
For any algebra $A$, the partial endomorphism monoid $\PEnd(A)$ is a left restriction monoid under the unary operation $\al^+=\id_{\dom(\al)}$.  \epfres
\end{prop}

Recall from Section \ref{sect:LRAP} that for any left restriction monoid $M$ we have a number of important submonoids, such as $P(M)=\set{s^+}{s\in M}$ and $T(M)=\set{s\in M}{s^+=1}$.
For any algebra $A$, we denote the semilattice of projections of $\PEnd(A)$ by
\begin{equation}\label{eq:EA}
\E_A = P(\PEnd(A)) = \set{\al^+}{\al\in\PEnd(A)} = \set{\id_B}{B\leq A}.
\end{equation}
So $\E_A$ is precisely the semilattice of \emph{partial identities} of $A$: i.e., all identity maps on subalgebras of $A$.  Since
\[
\id_B\cdot\id_C = \id_{B\cap C} \AND \id_B=\id_C\iff B=C \qquad\text{for all $B,C\leq A$,}
\]
it follows that $\E_A$ is in fact isomorphic to $\Sub(A)=\set{B}{B\leq A}$, considered as an $\cap$-semilattice.

It is also worth noting that for any algebra $A$,
\begin{equation}\label{eq:TA}
T(\PEnd(A)) = \set{\al\in\PEnd(A)}{\al^+=\id_A} = \set{\al\in\PEnd(A)}{\dom(\al)=A} = \End(A).
\end{equation}

For an independence algebra $A$, we also define the subsets
\begin{align*}
\Sing(\PEnd(A)) &= \PEnd(A)\sm\Aut(A) ,\\
\Sing(\End(A)) &= \End(A)\sm\Aut(A) ,\\
\Sing(\PAut(A)) &= \PAut(A)\sm\Aut(A) ,
\intertext{and we refer to the elements of these sets as \emph{singular}.  These are subsemigroups (indeed, two-sided ideals) of $\PEnd(A)$, $\End(A)$ and $\PAut(A)$, respectively, if and only if $\dim(A)<\aleph_0$.  On the other hand,}
\PEnd(A)\sm\End(A) &= \set{\al\in\PEnd(A)}{\dom(\al)<A}
\intertext{is a subsemigroup (indeed, right ideal) of $\PEnd(A)$, consisting of all \emph{strictly} partial endomorphisms, for any algebra $A$.  However,}
\PEnd(A)\sm\PAut(A) &= \set{\al\in\PEnd(A)}{\al\text{ is not injective}}
\end{align*}
is not a subsemigroup if $\dim(A)\geq3$.  For example, let $\{x,y,z\}\sub A$ be independent, let $B=\la x,y\ra$ and $C=\la y,z\ra$, and consider $\al:B\to A$ and $\be:C\to A$ defined by $x\al=y\al=x$ and $y\be=z\be=y$; then $\al,\be\in\PEnd(A)\sm\PAut(A)$, yet $\al\be=\id_{\CA}\in\PAut(A)$.

Now suppose $A$ is an independence algebra.  If $\al\in\PEnd(A)$, then the dimension of the subalgebra $\im(\al)\leq A$ is taken as the definition of the \emph{rank} of $\al$:
\[
\rank(\al) = \dim(\im(\al)).
\]
Note that $0\leq\rank(\al)\leq\dim(A)$.  We also have
\begin{equation}\label{eq:rank}
\rank(\al\be)\leq\min(\rank(\al),\rank(\be)) \qquad\text{for all $\al,\be\in\PEnd(A)$.}
\end{equation}
Indeed, this was proved for $\al,\be\in\End(A)$ in \cite[Lemma 4.1]{Gould1995}, and the argument there works unmodified for partial endomorphisms.  At times it will also be convenient to speak of the \emph{corank} of $\al\in\PEnd(A)$, defined by
\[
\corank(\al) = \codim(\im(\al)).
\]

When $\dim(A)=n<\aleph_0$, we have
\begin{align*}
\Aut(A)&=\set{\al\in\PEnd(A)}{\rank(\al)=n}\\
&=\set{\al\in\End(A)}{\rank(\al)=n}\\
&=\set{\al\in\PAut(A)}{\rank(\al)=n},
\end{align*}
so that $\Sing(\PEnd(A)) = \set{\al\in\PEnd(A)}{\rank(\al)<n}$, with similar statements for $\Sing(\End(A))$ and $\Sing(\PAut(A))$.

At this point it is worth recording some observations concerning small ranks and dimensions. 

\begin{rem}\label{rem:small}
For the duration of this remark, we write $C=\CA=\la\es\ra$ for the unique subalgebra of~$A$ of dimension $0$.  We have already noted that $C$ is non-empty if and only if $A$ has constants (nullary operations).
\ben
\item \label{small1} Since every subalgebra of $A$ contains $C$, and since every (partial) endomorphism fixes each constant, every (partial) endomorphism fixes $C$ pointwise.  
\item \label{small2} A (partial) endomorphism has rank $0$ if and only if its image is $C$.  It follows that if $A$ has no constants, then $A$ has no rank-$0$ endomorphisms (unless $A$ is itself empty).  On the other hand, $\id_C$ is always a rank-$0$ \emph{partial} endomorphism.
\item \label{small3} If $A$ has constants, then any function from a basis of $A$ into $C$ can be extended to an endomorphism of rank $0$.  
\item \label{small4} By \eqref{eq:rank}, the sets
\[
I_0 = \set{\al\in\End(A)}{\rank(\al)=0} \AND J_0 = \set{\al\in\PEnd(A)}{\rank(\al)=0}
\]
are ideals of $\End(A)$ and $\PEnd(A)$.  As above $I_0=\es$ if $A$ is non-empty and has no constants; in this case we also have $J_0=\{\id_C=\id_\es=\es\}$.  

\item \label{small5}  If $\al\in J_0$ and $\be\in\PEnd(A)$, then since $\be$ fixes $C=\im(\al)$ pointwise, $\al\be=\al$.  That is, each element of $J_0$ is a left zero for $\PEnd(A)$.  Consequently, the ideals $I_0$ and $J_0$ are left zero semigroups.

\item \label{small6} If $\dim(A)=0$, then $A=C$ and $\PEnd(A)=\End(A)=\PAut(A)=\Aut(A)=\{\id_A\}$.

\item \label{small7} If $\dim(A)=1$, then 
\begin{align*}
\End(A)&=\Aut(A)\cup I_0 ,\\
\PEnd(A)&= \Aut(A)\cup J_0=\End(A)\cup\{\id_C\} ,\\
\PAut(A) &= \Aut(A)\cup\{\id_C\} .
\end{align*}
In particular, if $\dim(A)=1$ and $A$ has no constants, then $C=\es$, $I_0=\es$ and
\[
\End(A)=\Aut(A) \AND \PEnd(A)=\PAut(A)=\Aut(A)\cup\{\es\}.
\]
\een
\end{rem}

It will also be convenient to prove the following:

\begin{lemma}\label{lem:Sing_es}
For any independence algebra $A$, the following are equivalent:
\ben
\item \label{Ses1} $\End(A)=\Aut(A)$,
\item \label{Ses2} $\Sing(\End(A))=\es$,
\item \label{Ses3} $\dim(A)=0$ or $[\dim(A)=1$ and $A$ has no constants$]$.
\een
\end{lemma}

\pf
Clearly \ref{Ses1} and \ref{Ses2} are equivalent.  

\pfitem{\ref{Ses3}$\implies$\ref{Ses1}}  This implication was covered in items \ref{small6} and \ref{small7} of Remark \ref{rem:small}.

\pfitem{\ref{Ses1}$\implies$\ref{Ses3}}  Aiming to prove the contrapositive, suppose \ref{Ses3} does not hold, so that either $\dim(A)\geq2$ or else $\dim(A)=1$ and $A$ has constants.  In either case, there exists $x\in A$ such that $\la x\ra<A$.  It follows that for any basis $X$ of $A$, the endomorphism $\al$ extending the constant map $X\to\{x\}$ is not an automorphism.  
\epf

\sect{Maximal subalgebras}\label{sect:MaxA}

In this section we prove a number of results concerning maximal (proper) subalgebras of an independence algebra $A$.   These will be used in a number of places in the remainder of Chapter~\ref{chap:IA}.  We denote the set of all maximal subalgebras of $A$ by $\Max(A)$.

\newpage

\begin{prop}\label{prop:Max}
For a proper subalgebra $B$ of an independence algebra $A$, the following are equivalent:
\ben
\item \label{Max1} $A=\la B\sqcup x\ra$ for some $x\in A\sm B$,
\item \label{Max2} $A=\la B\sqcup x\ra$ for all $x\in A\sm B$ (i.e., $B$ is maximal),
\item \label{Max3} $\codim(B)=1$.
\een
Consequently, $\Max(A)=\set{B\leq A}{\codim(A)=1}$.
\end{prop}

\pf
\firstpfitem{\ref{Max1}$\implies$\ref{Max2}}  Suppose $A=\la B\sqcup x\ra$ for some $x\in A\sm B$, and let $y\in A\sm B$ be arbitrary; we must show that $A=\la B\sqcup y\ra$.  But $y\in A\sm B=\la B\sqcup x\ra\sm\la B\ra$, so \ref{EP} gives $x\in\la B\sqcup y\ra$.  It then follows that $A=\la B\sqcup x\ra\sub\la B\sqcup y\ra$.

\pfitem{\ref{Max2}$\implies$\ref{Max3}}  Since $B<A$, we have $\codim(B)\geq1$.  It follows directly from \ref{Max2} that $\codim(B)\leq1$.

\pfitem{\ref{Max3}$\implies$\ref{Max1}}  This is clear.
\epf

In particular, if $\dim(A)=n<\aleph_0$, then $\Max(A)=\set{B\leq A}{\dim(B)=n-1}$, as observed in Section \ref{sect:Aprelim}.  Our next result shows that every finite-codimensional subalgebra of $A$ is a finite intersection of maximal subalgebras.  
In the statement, the expression $C_1\cap\cdots\cap C_k$ is understood to represent $A$ if $k=0$: i.e., if the list $C_1,\ldots,C_k$ is empty.

\begin{lemma}\label{lem:BCC}
Let $A$ be an independence algebra, and let $B\leq A$ with $\codim(B)=k<\aleph_0$.  Then $B=C_1\cap\cdots\cap C_k$ for some $C_1,\ldots,C_k\in\Max(A)$.
\end{lemma}

\pf
Let~$X$ be a basis of $B$, and $\{y_1,\ldots,y_k\}$ a $B$-basis of $A$, so that $Y=X\sqcup\{y_1,\ldots,y_k\}$ is a basis of $A$ by Lemma \ref{lem:XY}.  For each $1\leq i\leq k$, let $C_i=\la Y\sm y_i\ra$; it follows from Proposition~\ref{prop:Max} that each $C_i$ is maximal.  We will show by induction on $k$ that $B=C_1\cap\cdots\cap C_k$.  This is vacuously true for $k=0$, so we assume that $k\geq1$.  

Let $X'=X\sqcup y_k$, and put $B'=\la X'\ra$, so that $\codim(B')=k-1$, and $Y=X'\sqcup\{y_1,\ldots,y_{k-1}\}$.  It follows inductively that $B'=C_1\cap\cdots\cap C_{k-1}$, so it remains to show that $B=B'\cap C_k$.  
From $X\sub X'$ and $X\sub Y\sm y_k$, we have $B\sub B'$ and $B\sub C_k$, so that $B\leq B'\cap C_k\leq B'$.  Since $B'=\la B\sqcup y_k\ra$, it follows from Proposition \ref{prop:Max} that $B$ is maximal in $B'$, so we can complete the proof that $B=B'\cap C_k$ by showing that $B'\cap C_k<B'$.  But this is clear:  we obviously have~$y_k\in B'$, yet $y_k\not\in\la Y\sm y_k\ra=C_k$, as $Y$ is independent.
\epf

\begin{rem}\label{rem:BCC}
In the case that $\dim(A)<\aleph_0$, it follows from Lemma \ref{lem:BCC} that $\Sub(A)$ is generated (as an $\cap$-monoid) by the set $\Max(A)$.  We will see in Proposition \ref{prop:MaxA} that $\Max(A)$ is in fact the unique minimum generating set.  Actually, Proposition \ref{prop:MaxA} also treats the case that~$A$ is infinite-dimensional and strong, but with $\Sub(A)$ replaced by the lattice of finite-codimensional subalgebras.
\end{rem}

For subalgebras of arbitrary codimension we have the following:

\begin{lemma}\label{lem:BCC2}
If $A$ is an independence algebra, and if $B\leq A$, then
\[
B=\bigcap_{C\in \Si}C \WHERE \Si=\set{C\in\Max(A)}{B\leq C}.
\]
\end{lemma}

\pf
Write $D=\bigcap_{C\in \Si}C$.  Since $B\leq C$ for each $C\in \Si$, certainly $B\leq D$.

To show that $D\leq B$, it suffices to show that $A\sm B\sub A\sm D$, so fix some $y\in A\sm B$.  We will show that $y\not\in D$ by finding $C\in \Si$ such that $y\not\in C$.  Fix a basis $X$ of $B$, so that $X\sqcup y$ is independent, by Lemma \ref{lem:Xy}.  Extend $X\sqcup y$ to a basis $X\sqcup y\sqcup Z$ of $A$, and define $C=\la X\sqcup Z\ra$.  By Proposition \ref{prop:Max}, $C$ is maximal, and certainly $B=\la X\ra\sub C$, so that $C\in \Si$.  Since $X\sqcup y\sqcup Z$ is independent, $y\not\in\la X\sqcup Z\ra=C$.
\epf

The next result concerns the situation in which $A$ has subalgebras missing only one element, and shows that this is an extremely restrictive property.  Since every subalgebra contains $\CA$, we could only possibly have $A\sm x\leq A$ if $x\not\in\CA$.  Such an element $x$ could only exist if $A\not=\CA$: i.e., if $\dim(A)\not=0$.

\begin{prop}\label{prop:sub1}
Let $A$ be an independence algebra with $\dim(A)\not=0$, and write $C=\CA$ and $X=A\sm C$.  Then the following are equivalent:
\ben
\item \label{sub11} $A\sm x\leq A$ for some $x\in X$,
\item \label{sub12} $A\sm x\leq A$ for all $x\in X$,
\item \label{sub13} $X$ is independent,
\item \label{sub14} $X$ is the unique basis of $A$,
\item \label{sub15} $\Sub(A)=\set{B}{ C\sub B\sub A}$,
\item \label{sub16} $\Max(A)=\set{A\sm x}{x\in X}$.
\een
\end{prop}

\pf
\firstpfitem{\ref{sub11}$\implies$\ref{sub12}}  Suppose $A\sm x\leq A$ for some $x\in X$, and let $y\in X$ be arbitrary; we must show that $A\sm y\leq A$.  Since $x,y\not\in C$, we may extend $x$ and $y$ to bases $U\sqcup x$ and $V\sqcup y$ for $A$, noting that $|U|=|V|$.  We then let $\al$ be the endomorphism of~$A$ extending any bijection $U\sqcup x\to V\sqcup y$ that maps $x\mt y$.  Then $\al$ is an automorphism (as it maps a basis bijectively to a basis), and since $A\sm x\leq A$, it follows that $A\sm y=(A\sm x)\al\leq A$.  

\pfitem{\ref{sub12}$\implies$\ref{sub13}}  If we had $x\in\la X\sm x\ra$ for some $x\in X$, then since $A\sm x\leq A$, it would follow that $x\in\la X\sm x\ra\sub\la A\sm x\ra=A\sm x$, a contradiction.

\pfitem{\ref{sub13}$\implies$\ref{sub14}}  If $X$ is independent, then it is clearly the maximum independent subset of $A$.

\pfitem{\ref{sub14}$\implies$\ref{sub15}}  Since every subalgebra contains $C$, it suffices to show that $B\leq A$ for any $C\sub B\sub A$, so fix some such $B$.  To show that $B\leq A$, it suffices to show that no element of $A\sm B$ belongs to~$\la B\ra$.  But for any $x\in A\sm B$, we have $(B\sm C)\sqcup x\sub A\sm C=X$.  Since $X$ is independent, so too is $(B\sm C)\sqcup x$.  Thus, $x\not\in\la B\sm C\ra=\la B\ra$.

\pfitem{\ref{sub15}$\implies$\ref{sub16} and \ref{sub16}$\implies$\ref{sub11}}  These are both clear.
\epf

Consider an independence algebra $A$, and again write $C=\CA$ and $X=A\sm C$.  Any automorphism of $A$ fixes $C$ pointwise, and is therefore of the form $\id_C\cup\al$ for some permutation $\al\in\G_X$.  It follows that $\Aut(A)$ is a subgroup of $\set{\id_C\cup\al}{\al\in\G_X}$, and the latter is of course isomorphic to $\G_X$.  
It turns out that the conditions of Proposition \ref{prop:sub1} lead to the greatest `freedom' in automorphisms, in the sense that $\Aut(A)$ is precisely this copy of $\G_X$.  Indeed, if $X$ is a basis for~$A$, then any permutation $\al\in\G_X$ extends to an automorphism (cf.~Lemma \ref{lem:auto}), which must be $\id_C\cup\al$.  
However, it turns out that $X$ being independent is not (quite) necessary for this `maximum freedom' to occur, and the next result characterises the algebras for which it does.

\begin{prop}\label{prop:sub2}
Let $A$ be an independence algebra with $\dim(A)\not=0$, and write $C=\CA$ and $X=A\sm C$.  Then the following are equivalent:
\ben
\item \label{sub21} $X\sm x$ is independent for some $x\in X$,
\item \label{sub22} $X\sm x$ is independent for all $x\in X$,
\item \label{sub23} $\Aut(A) = \set{\id_C\cup\al}{\al\in\G_X}$.
\een
When the above conditions hold, we have $\Aut(A)\cong\G_X$.
\end{prop}

\newpage 

\pf
If $X$ is independent, then all three conditions hold, so we assume otherwise for the rest of the proof.  In particular, it follows from Proposition~\ref{prop:sub1} that $A\sm x\not\leq A$ for all $x\in A$.  

\pfitem{\ref{sub21}$\implies$\ref{sub22}}  Suppose $X\sm x$ is independent for some $x\in X$, and let $y\in X$ be arbitrary; we must show that $X\sm y$ is independent.  This is clear if $x=y$, so suppose otherwise, and write $U=X\sm\{x,y\}$.  Since $U\sqcup y=X\sm x$ is independent we have $y\not\in\la U\ra$.  

Next we claim that $x\not\in\la U\ra$.  To prove this, suppose to the contrary that $x\in\la U\ra$.  Then $\la U\ra = \la U\sqcup x\ra = \la C\sqcup U\sqcup x\ra = \la A\sm y\ra = A$, since $A\sm y\not\leq A$.  It follows that $y\in A=\la U\ra$, a contradiction.  
Now that the claim has been proved, it follows from Lemma \ref{lem:Xy} (and the fact that $U$ is independent) that $U\sqcup x$ is independent.  But $U\sqcup x=X\sm y$, so we are done.

\pfitem{\ref{sub22}$\implies$\ref{sub23}}  It suffices to show that $\id_C\cup\al$ is an automorphism for any $\al\in\G_X$.  To do so, fix some such $\al$, and let $Y=X\sm x$ for any $x\in X$.  Since $Y$ is independent, and $X$ is not, it follows that $Y$ is a basis of $A$.  Since $\al\in\G_X$, we have $Y\al=X\sm x\al$, and again this is a basis.  It follows from Lemma \ref{lem:auto} that the unique endomorphism $\ga$ of $A$ extending $\al\restr_Y$ is an automorphism.  But $\ga\restr_C=\id_C$ and $\ga\restr_Y=\al\restr_Y$, so $\ga$ and $\id_C\cup\al$ agree on $C\cup Y=A\sm x$.  Since $\ga$ and $\id_C\cup\al$ agree on the basis $Y$ of $A$, it follows that $\ga=\id_C\cup\al$.  This completes the proof that $\id_C\cup \al$ is an automorphism.

\pfitem{\ref{sub23}$\implies$\ref{sub21}}  Aiming to prove the contrapositive, we assume that \ref{sub21} does not hold; we must show that $\id_C\cup\al\not\in\Aut(A)$ for some $\al\in\G_X$.  Fix a basis $Y$ of $A$.  So $Y\sub X$, and by assumption $|X\sm Y|\geq2$.  The only automorphism extending $\id_Y$ is $\id_A=\id_C\cup\id_X$.  So we may take $\al\in\G_X$ to be any permutation that maps $Y$ identically, but permutes $X\sm Y$ non-trivially.
\epf

\begin{rem}
One may wonder how the conditions of Proposition \ref{prop:sub2} relate to the monoids $\End(A)$, $\PEnd(A)$ and $\PAut(A)$, and to what extent these have the same `freedom' as $\Aut(A)$, as discussed above.  We will treat $\PAut(A)$ in the next result.  The situation for the other two monoids is a little more complex, as (partial) endomorphisms can map elements of $X$ into~$C$, in the above notation.  But, for example, when $X$ is independent (so the conditions of Proposition~\ref{prop:sub1} are satisfied), a transformation $\al$ of $A$ belongs to $\End(A)$ if and only if $\al$ fixes~$C$ pointwise, so that
\[
\End(A) = \set{\al\in\T_A}{\al\restr_C=\id_C}.
\]
These submonoids of $\T_A$ have been studied extensively, and are typically denoted $\operatorname{Fix}(A,C)$; see for example \cite{HS2013}.  In the case that $C=\es$ we of course have $\End(A)=\T_A$.  If $|C|=1$, say with $C=\{c\}$, then $\End(A)\cong\PT_X$, via the standard trick of identifying a partial transformation $\al$ of $X$ with the transformation of $A=X\sqcup c$ mapping as $\al$ on $\dom(\al)$, and sending all elements of $A\sm\dom(\al)$ to $c$.
\end{rem}

Here is the promised version of Proposition \ref{prop:sub2}, concerning `maximum freedom' for \emph{partial} automorphisms.

\begin{prop}\label{prop:sub5}
Let $A$ be an independence algebra, and write $C=\CA$ and $X=A\sm C$.  Then the following are equivalent:
\ben
\item \label{sub51} $X$ is independent,
\item \label{sub52} $\PAut(A) = \set{\id_C\cup\al}{\al\in\I_X}$.
\een
When the above conditions hold, we have $\PAut(A)\cong\I_X$.
\end{prop}

\pf
The result is trivial if $\dim(A)=0$, as then $X=\es$, so we assume that $\dim(A)\not=0$.

\pfitem{\ref{sub51}$\implies$\ref{sub52}}  Since any element of $\PAut(A)$ is injective and maps $C$ identically, the forwards containment always holds.  Conversely, let $\al\in\I_X$, and write $B=\dom(\al)$ and $D=\im(\al)$.  We must show that $\id_C\cup\al\in\PAut(A)$; in fact, since $\id_C\cup\al$ is obviously injective, it is enough to check that it is a morphism $C\cup B\to C\cup D$.  Now, $C\cup B,C\cup D\leq A$ by Proposition \ref{prop:sub1}.  Since~$B$ is independent (as $X$ is), it is a basis for $C\cup B$, so there is therefore a unique morphism $C\cup B\to C\cup D$ extending~$\al:B\to D$, and this is of course $\id_C\cup\al$.

\pfitem{\ref{sub52}$\implies$\ref{sub51}}  By Proposition \ref{prop:sub1} it suffices to show that $B\leq A$ for any $C\sub B\sub A$.  But for any such $B$ we have $\id_B=\id_C\cup\id_{B\sm C}\in\PAut(A)$, so that $B=\dom(\id_B)\leq A$.
\epf

We have already noted that the conditions in Proposition \ref{prop:sub1} imply those in Proposition~\ref{prop:sub2}.  The next result concerns algebras that satisfy the latter but not the former.  An example of such an algebra is the two-dimensional vector space over the two-element field; another example is the one from \cite{FL1993} discussed in Remark \ref{rem:NS1}; see also Proposition \ref{prop:NS} below.

\begin{prop}\label{prop:sub3}
Let $A$ be an independence algebra with $\dim(A)\not=0$, write $C=\CA$ and $X=A\sm C$, and suppose $A$ does not satisfy the conditions of Proposition \ref{prop:sub1}.  Then the following are equivalent:
\ben
\item \label{sub31} $X\sm x$ is independent for some $x\in X$,
\item \label{sub32} $X\sm x$ is independent for all $x\in X$,
\item \label{sub33} $A\sm\{x,y\}\leq A$ for all $x,y\in X$ with $x\not=y$,
\item \label{sub34} $\Sub(A)=\set{B}{ C\sub B\sub A,\ |A\sm B|\not=1}$,
\item \label{sub35} $\Max(A)=\set{A\sm \{x,y\}}{x,y\in X,\ x\not=y}$.
\een
\end{prop}

\pf
\firstpfitem{\ref{sub31}$\iff$\ref{sub32}}  This was proved in Proposition \ref{prop:sub2}.

\pfitem{\ref{sub32}$\implies$\ref{sub34}}  Since $A$ does not satisfy the conditions of Proposition \ref{prop:sub1}, it has no subalgebras missing just one element.  This gives the forwards containment.

Conversely, suppose $C\sub B\sub A$ is such that $|A\sm B|\not=1$; we must show that $B\leq A$.  If $|A\sm B|=0$ then $B=A\leq A$, so suppose instead that $|A\sm B|\geq2$.
Aiming for a contradiction, suppose $B\not\leq A$, so that $x\in\la B\ra$ for some $x\in A\sm B$.  Since $|A\sm B|\geq2$, we may fix another element $y\in A\sm B$.
Since $x,y\not\in B$ and $C\sub B$, we have $x,y\in X$.  
Then $x\in\la B\ra=\la B\sm C\ra \sub \big\la X\sm\{x,y\}\big\ra$; but this implies that $X\sm y$ is not independent, and contradicts \ref{sub32}.

\pfitem{\ref{sub34}$\implies$\ref{sub35} and \ref{sub35}$\implies$\ref{sub33}}  These are clear.

\pfitem{\ref{sub33}$\implies$\ref{sub32}}  Aiming to prove the contrapositive, suppose $X\sm x$ is not independent for some $x\in X$, so that $y\in\big\la X\sm\{x,y\}\big\ra$ for some $y\in X\sm x$.  But then $y\in\big\la X\sm\{x,y\}\big\ra\sub\big\la A\sm\{x,y\}\big\ra$, so that $A\sm\{x,y\}\not\leq A$.
\epf

\begin{rem}\label{rem:sub3}
Comparing Propositions \ref{prop:sub1} and \ref{prop:sub3}, one might wonder why the following condition was not listed in the latter:
\bit
\item $A\sm\{x,y\}\leq A$ for \emph{some} $x,y\in X$ with $x\not=y$.
\eit
The reason is that this is strictly weaker than condition \ref{sub33} of Proposition~\ref{prop:sub3} in general.  

For example, let $A=\{a,b,c,d\}$ with a single unary operation mapping $a\leftrightarrow b$ and $c\leftrightarrow d$.  Then $A$ is an independence algebra (it is term-equivalent to the free two-dimensional $\G_2$-act, as in Example~\ref{eg:FXM}), and it is easy to see that the lattice $\Sub(A)$ is as follows:
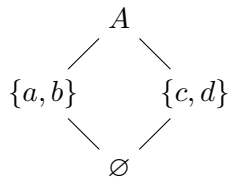
\[
\begin{tikzpicture}[scale=1]
\node (PE) at (0,2) {$A$};
\node (E) at (-1,1) {$\{a,b\}$};
\node (PA) at (1,1) {$\{c,d\}$};
\node (A) at (0,0) {$\es$};
\draw (PE)--(E)--(A)--(PA)--(PE);
\end{tikzpicture}
\]
In particular, $\{a,b\}$ and $\{c,d\}$ are the only two-element subalgebras of $A$.  This means that condition \ref{sub33} of Proposition \ref{prop:sub3} fails, even though the above weaker condition holds.  The conditions of Propositions \ref{prop:sub1} and \ref{prop:sub2} also fail for $A$.

Because of this, it follows from Proposition \ref{prop:sub2} that $\Aut(A)$ is a proper subgroup of $\G_A$ (note that $X=A\sm\C(A)=A$ since $A$ has no nullary operations), and we can also see this by direct computation.  Indeed, any basis of $A$ contains one element each from $\{a,b\}$ and $\{c,d\}$, and it follows that $\Aut(A)$ consists precisely of the following permutations, which are specified in two ways, by indicating the action on the basis $\{a,c\}$, and in standard cycle notation:
\begin{align*}
\id_A = \trans{a&c\\a&c} , & & \al_1& = \trans{a&c\\b&c} = (a,b), & \be_1 & = \trans{a&c\\c&a} = (a,c)(b,d),\\
 & & \al_2& = \trans{a&c\\a&d} = (c,d), & \be_2 & = \trans{a&c\\c&b} = (a,c,b,d),\\
 & & \al_3& = \trans{a&c\\b&d} = (a,b)(c,d), & \be_3 & = \trans{a&c\\d&b} = (a,d)(b,c),\\
 & & & & \be_4 & = \trans{a&c\\d&a} = (a,d,b,c).
\end{align*}
Consequently, $\Aut(A)$ is dihedral of order $8$, or equivalently a wreath product $\G_2\wr\G_2$ (see Chapter~\ref{chap:wreath}).
\end{rem}

We have already noted that the two-dimensional vector space over the two-element field satisfies the conditions of Proposition \ref{prop:sub2} but not those of Proposition~\ref{prop:sub1}.  The next result shows that this is unusual among such examples, as it is strong; note that $|X|=3$ for this example.

\begin{prop}\label{prop:NS}
Let $A$ be an independence algebra, write $C=\CA$ and $X=A\sm C$, and suppose $|X|\geq4$.  If $A$ satisfies the conditions of Proposition \ref{prop:sub2} but not those of Proposition~\ref{prop:sub1}, then~$A$ is not strong.
\end{prop}

\pf
Fix a partition $X=Y\sqcup Z$ where $|Y|,|Z|\geq2$.  Then $Y$ and $Z$ are both independent (as $X\sm x$ is independent for all $x\in X$).  By Proposition \ref{prop:sub3}, $C\cup Y=A\sm Z\leq A$, so $\la Y\ra=\la C\cup Y\ra=C\cup Y$, and similarly $\la Z\ra=C\cup Z$.  Certainly then $\la Y\ra\cap \la Z\ra=C$, yet $Y\cup Z=X$ is not independent.
\epf

\begin{rem}
By contrast, if $A$ satisfies the conditions of Proposition \ref{prop:sub1}, then $A$ is trivially strong, for then the independent sets are just the subsets of $X=A\sm\CA$.
\end{rem}

If the conditions of Proposition \ref{prop:sub1} hold for an independence algebra $A$, then it quickly follows that $A$ is equal to the union of any two distinct maximal subalgebras.  In Section \ref{sect:AutA} we will be concerned with the case in which this never happens.  In fact, we only need a weaker assumption, but the next result shows that these are equivalent when $A$ is strong.

\begin{lemma}\label{lem:sub4}
Let $A$ be an independence algebra, and consider the following two conditions:
\ben
\item \label{sub41} for any $B,C\in\Max(A)$, we have $B\cup C\not=A$,
\item \label{sub42} for any basis $X$ of $A$, and for any $x,y\in X$, we have $\la X\sm x\ra\cup\la X\sm y\ra\not=A$.
\een
Then \ref{sub41}$\implies$\ref{sub42}.  If $A$ is strong, then \ref{sub41}$\iff$\ref{sub42}
\end{lemma}

\pf
By Proposition \ref{prop:Max}, $\la X\sm x\ra$ is maximal for any basis $X$ of $A$ and any $x\in X$, so it is clear that \ref{sub41}$\implies$\ref{sub42}.

We now assume that $A$ is strong, and that \ref{sub42} holds.  Let $B,C\in\Max(A)$ be arbitrary; we must show that $B\cup C\not=A$.  This is clear if $B=C$, so suppose $B\not=C$.  Together with maximality it follows that $B\not\sub C$ and $C\not\sub B$, and also that $B\vee C=A$.  From the latter, we have $\codim(B\vee C)=0$, so it follows from Lemma \ref{lem:BC1}\ref{BC12} and Proposition \ref{prop:Max} that $\codim(B\cap C)=2$.  

Fix a basis $U$ for $B\cap C$, and extend this to bases $U\sqcup Y$ and $U\sqcup Z$ for $B$ and $C$, respectively.  Since $B\not=C$ we have $|Y|\geq1$ and $|Z|\geq1$.  By Lemma \ref{lem:strong}, $U\sqcup Y\sqcup Z$ is a basis of $B\vee C=A$.  But then $|Y|+|Z|=\codim(B\cap C)=2$, so in fact $|Y|=|Z|=1$, and we may write $Y=\{y\}$ and $Z=\{z\}$.  This all shows that $X=U\sqcup\{y,z\}$ is a basis of $A$.  We also have $B=\la U\sqcup y\ra=\la X\sm z\ra$, and similarly $C=\la X\sm y\ra$.  It then follows from \ref{sub42} that $B\cup C = \la X\sm z\ra\cup\la X\sm y\ra \not=A$.
\epf

\begin{rem}
The algebra $A$ from \cite{FL1993} discussed in Remark \ref{rem:NS1} satisfies condition \ref{sub42} of Lemma \ref{lem:sub4}, but not \ref{sub41}.  Indeed, each two-element subalgebra is maximal, and for example ${A=\{a_1,a_2\}\cup\{a_3,a_4\}}$.  However, any basis has the form $X=\{a_i,a_j,a_k\}$ for distinct $i,j,k$, and $\la X\sm x\ra\cup\la X\sm y\ra=(X\sm x)\cup(X\sm y)\sub X\not=A$ for any $x,y\in X$.
\end{rem}

\begin{rem}\label{rem:sub4}
We have already noted that the conditions in Proposition \ref{prop:sub1} and Lemma~\ref{lem:sub4} are mutually exclusive.  They are not, however, exhaustive.  For example, let $A$ be the four-element algebra from Remark \ref{rem:sub3}.  We have already observed that $A$ does not satisfy the conditions in Proposition \ref{prop:sub1}.  On the other hand, the maximal subalgebras of $A$ are precisely $B=\{a,b\}$ and $C=\{c,d\}$, and of course $B\cup C=A$.  In fact, $B=\la X\sm c\ra$ and $C=\la X\sm a\ra$ for the basis $X=\{a,c\}$.  Thus, $A$ does not satisfy either of the conditions of Lemma~\ref{lem:sub4}.
(This extends to arbitrary free $G$-acts of dimension $\geq2$ for non-trivial $G$.)
\end{rem}

\begin{rem}\label{rem:sub42}
Any vector space $A$ of dimension $\geq2$ satisfies the conditions of Lemma~\ref{lem:sub4}.  Indeed, if $X$ is a basis, and if $x,y\in X$ are distinct, then $x+y\not\in\la X\sm x\ra\cup\la X\sm y\ra$.  (To see this, note that if $x+y\in\la X\sm x\ra$, then since also $y\in X\sm x$, we would have $x=(x+y)-y\in\la X\sm x\ra$, contradicting independence of $X$.)

It follows that such a vector space $A$ cannot satisfy the conditions of Proposition \ref{prop:sub1}.  The only such~$A$ satisfying the conditions of Proposition \ref{prop:sub2} is two-dimensional over the two-element field, which we have already discussed.
\end{rem}

\sect{Finite-codimensional subalgebras}\label{sect:fca}

The previous section concerned maximal subalgebras of an independence algebra $A$: i.e., subalgebras of codimension $1$; cf.~Proposition \ref{prop:Max}.  We now consider the more general case of subalgebras of arbitrary finite codimension, the set of which we denote by
\[
\FSub(A) = \set{B\leq A}{\codim(B)<\aleph_0}.
\]
Our first result shows that $\FSub(A)$ is a sublattice of $\Sub(A)$ in the case that $A$ is finite-dimensional or strong, and Theorem \ref{thm:SubA} gives a presentation for $\FSub(A)$ in this case, considered as an $\cap$-semilattice.

\begin{lemma}\label{lem:FSub}
Let $A$ be an independence algebra, and let $B,C\in\FSub(A)$.  Then
\ben
\item \label{FSub1} $B\vee C\in\FSub(A)$,
\item \label{FSub2} $B\cap C\in\FSub(A)$ if $A$ is finite-dimensional or strong.
\een
\end{lemma}

\pf
\firstpfitem{\ref{FSub1}}  This follows immediately from Lemma \ref{lem:BCA}\ref{BCA1}, as $B\leq B\vee C\leq A$ and $\codim(B)<\aleph_0$.

\pfitem{\ref{FSub2}}  Since $\FSub(A)=\Sub(A)$ when $\dim(A)<\aleph_0$, we only need to consider the case that $A$ is strong.  But here Lemma \ref{lem:BC1}\ref{BC12} gives $\codim(B\cap C)\leq\codim(B)+\codim(C)<\aleph_0$.
\epf

Recall that $\Max(A)=\set{B\leq A}{\codim(B)=1}$ is the set of maximal subalgebras of $A$.

\begin{prop}\label{prop:MaxA}
If $A$ is a finite-dimensional or strong independence algebra, then $\Max(A)$ is the (unique) minimum monoid generating set for $\FSub(A)$, considered as an $\cap$-semilattice.  
\end{prop}

\pf
It follows immediately from Lemma \ref{lem:BCC} that $\Max(A)$ is a generating set.  

It remains to show that any generating set for $\FSub(A)$ contains $\Max(A)$, and for this it suffices to show that every element of $\Max(A)$ is an atom.  To do so, suppose $B=C\cap D$, with $B\in\Max(A)$ and $C,D\in\FSub(A)$; we must show that $B\in\{C,D\}$.  If $C=A$, then $B=A\cap D=D$.  Otherwise, $B=C\cap D\leq C<A$, so $B=C$ by maximality.
\epf

Now that we know $\FSub(A)$ contains a unique minimum generating set (when $A$ is finite-dimensional or strong), this is of course the generating set we wish to use in constructing a presentation.  We can do this in the case that $A$ is assumed to be strong.  We begin by defining an alphabet
\[
X_E = \bigset{x_B}{B\in\Max(A)}
\]
in one-one correspondence with $\Max(A)$.  (The reason for choosing the $E$ subscript will become clear in later sections.)  By Proposition \ref{prop:MaxA} we have a surmorphism
\[
\Phi_E:X_E^*\to\FSub(A):x_B\mt B.
\]
We define the set of relations $R_E=R_E^1\cup R_E^2\cup R_E^3$, where
\begin{align*}
R_E^1 &= \bigset{(x_B^2,x_B)}{B\in \Max(A)}, \\
R_E^2 &= \bigset{(x_Bx_C,x_Cx_B)}{B,C\in \Max(A)}, \\
R_E^3 &= \bigset{(x_Bx_C,x_Bx_D)}{B,C,D\in \Max(A),\ B\cap C=B\cap D}. 
\end{align*}
It is clear that $R_E\sub\ker(\Phi_E)$.

\begin{thm}\label{thm:SubA}
If $A$ is a strong independence algebra, then with the above notation, $\FSub(A)$ has presentation $\Mpres{X_E}{R_E}$ via $\Phi_E$.
\end{thm}

\pf
We begin with three technical lemmas, the first of which is rather simple.  The second may appear to be \emph{exceedingly} technical, but it turns out to be exactly what we need in the proof of the third, which leads to a swift conclusion to the proof of the theorem.

\begin{lemma}\label{lem:BC2}
If $B\leq A$ and $C\in \Max(A)$, then $\codim(B)\leq \codim(B\cap C)\leq \codim(B)+1$.
\end{lemma}

\pf
By Lemma \ref{lem:BC1}\ref{BC12} and Proposition \ref{prop:Max}, we have
\[
\codim(B\vee C) + \codim(B\cap C) = \codim(B)+1.
\]
Since $C\leq B\vee C\leq A$, Lemma \ref{lem:BCA}\ref{BCA1} gives $0\leq\codim(B\vee C)\leq 1$, and the result follows.
\epf

\begin{lemma}\label{lem:BC3}
If $B,C\in \Max(A)$ and $D\in\FSub(A)$ are such that $B\not=C$, $D\not\sub B$ and $D\cap C\sub B$, then there exists $B'\in \Max(A)$ such that $D\sub B'$ and $C\cap B=C\cap B'$.
\end{lemma}

\pf
We first claim that $D\cap B=D\cap C$.  To prove this, first note that $D\cap C\not=D$ (as $D\not\sub B$ but $D\cap C\sub B$), so Lemma \ref{lem:BC2} gives $\codim(D\cap C)=\codim(D)+1$.  Since also $D\cap B\not=D$ (as~$D\not\sub B$), Lemma \ref{lem:BC2} gives $\codim(D\cap B)=\codim(D)+1$ as well.  In particular, $\codim(D\cap C)=\codim(D\cap B)<\aleph_0$ (as $D\in\FSub(A)$).  By assumption we have $D\cap C\sub B$, so of course $D\cap C\sub D\cap B$.  It then follows from Lemma \ref{lem:BCA}\ref{BCA2} that $D\cap C=D\cap B$, as claimed.

We now fix a basis $X$ of $D\cap B(=D\cap C)$, and extend this to a basis $X\sqcup Y$ for $B\cap C$.  Since $\codim(B)=\codim(C)=1$ and $B\not=C$, Lemma \ref{lem:BC2} gives $\codim(B\cap C)=2$.  Since $D\not\sub B$, we may also fix some $z\in D\sm B(=D\sm C)$.  Since $z\not\in B$, certainly $z\not\in B\cap C=\la X\sqcup Y\ra$, so it follows from Lemma \ref{lem:Xy} that $X\sqcup Y\sqcup z$ is independent.  Let $B'=\la X\sqcup Y\sqcup z\ra$.  Since $B\cap C=\la X\sqcup Y\ra$ and $\codim(B\cap C)=2$, we have $B'<A$ and so $\codim(B')\geq1$.  Since $B\cap C\sub B'$, we have $\codim(B')\leq\codim(B\cap C)=2$.  If $\codim(B')=2$, then Lemma \ref{lem:BCA}\ref{BCA2} would give $B'=B\cap C$, contradicting $z\not\in B\cap C$.  So it follows that $\codim(B')=1$, so that $B'\in \Max(A)$.

Next let $D'=\la X\sqcup z\ra$.  Since $D\cap C=\la X\ra$, and since $X\sqcup z\sub D$, we have
\begin{equation}\label{eq:D'}
D\cap C \sub D' \sub D.
\end{equation}
Combined with $\codim(D\cap C) = \codim(D)+1$ (shown above), it follows that
\[
\codim(D)+1 \geq \codim(D') \geq \codim(D),
\]
and this of course implies that
\[
\codim(D') = \codim(D) \OR \codim(D') = \codim(D)+1 = \codim(D\cap C).
\]
In these two cases, and since $\codim(D)<\aleph_0$, Lemma \ref{lem:BCA}\ref{BCA2} and \eqref{eq:D'} give $D'=D$ or ${D'=D\cap C}$, respectively.  But the latter is impossible, since $D\cap C=\la X\ra$ and $D'=\la X\sqcup z\ra$, and since $X(\sqcup Y)\sqcup z$ is independent.  Thus, $D=D'$, and it follows that $D=\la X\sqcup z\ra\sub\la X\sqcup Y\sqcup z\ra=B'$.

It remains to show that $C\cap B=C\cap B'$, and for this it is enough to show that
\[
C\cap B\sub B' \AND C\cap B'\sub B.
\]
For the first, we have $C\cap B=\la X\sqcup Y\ra \sub B'$.  For the second, suppose to the contrary that there exists $w\in(C\cap B')\sm B$.  In particular we have $w\in B'=\la X\sqcup Y\sqcup z\ra$.  On the other hand, since $w\not\in B$, certainly $w\not\in B\cap C=\la X\sqcup Y\ra$.  It then follows from \ref{EP} that $z\in\la X\sqcup Y\sqcup w\ra\sub C$ (as $X\sqcup Y\sub B\cap C\sub C$ and $w\in C\cap B'\sub C$), and this contradicts $z\in D\sm C$.  
\epf

For the rest of the proof (of the theorem) we write~${\sim}=R_E^\sharp$.

\begin{lemma}\label{lem:wB}
If $w\in X_E^*$, and if $B\in \Max(A)$ is such that $w\Phi_E\sub B$, then $w\sim wx_B$.
\end{lemma}

\pf
Given the relations in $R_E^1\cup R_E^2$, it suffices to show that $w$ is $\sim$-equivalent to some word that involves the letter $x_B$.  

Write $w=x_{C_1}\cdots x_{C_k}$, where each $C_i\in \Max(A)$, and note that $w\Phi_E=C_1\cap\cdots\cap C_k$.  Since $w\Phi_E\sub B<A$, we must have $k\geq1$.  We now proceed by induction on $k$.  If $k=1$, then $C_1=w\Phi_E\sub B$, and it follows by maximality that $C_1=B$; thus, $w=x_B$ in this case.

Now suppose $k\geq2$, and let $u=x_{C_1}\cdots x_{C_{k-1}}$ and $D=u\Phi_E=C_1\cap\cdots\cap C_{k-1}\in\FSub(A)$.  If $D\sub B$, then by induction we have $u\sim ux_B$, and so $w=ux_{C_k}\sim ux_Bx_{C_k}$, and we are done.  So now we assume that $D\not\sub B$, and we also write $C=C_k$.  If $B=C$, then $w=ux_B$, and again we are done, so we also assume that $B\not=C$.  Since also $D\cap C=w\Phi_E\sub B$, the conditions of Lemma \ref{lem:BC3} are all satisfied, so by that lemma there exists $B'\in \Max(A)$ such that $D\sub B'$ and~$C\cap B=C\cap B'$.  By the latter,~$R_E^3$ contains the relation $(x_Cx_B,x_Cx_{B'})$.  Since $u$ has length~$k-1$, and since $u\Phi_E=D\sub B'$, it follows by induction that $u\sim ux_{B'}$.  But then 
\[
w=ux_C\sim  ux_{B'}x_C\sim  ux_Cx_{B'}\sim  ux_Cx_B,
\]
and the proof is complete.
\epf

We now return to the main proof.
We have already observed that $\Phi_E$ is surjective, and that $R_E\sub\ker(\Phi_E)$, so it remains to show that $\ker(\Phi_E)\sub R_E^\sharp$.  So suppose $(u,v)\in\ker(\Phi_E)$; we must show that $u\sim v$.  For $w\in X_E^*$, we have $w\Phi_E=A\iff w=\ew$ (the empty word), so it suffices to assume that $u$ and~$v$ are both non-empty.  Let $x_B\in X_E$ be any letter appearing in $v$, so that~$v\sim vx_B$ (by $R_E^1\cup R_E^2$).  Then $u\Phi_E=v\Phi_E=(vx_B)\Phi_E=v\Phi_E\cap B\sub B$, so it follows from Lemma \ref{lem:wB} that~$u\sim ux_B$.  Since this is the case for every letter appearing in $v$, it follows that $u\sim uv$.  By symmetry we also have $v\sim vu$, and since $uv\sim vu$ (by $R_E^2$), it follows that~$u\sim v$.
\epf

\begin{rem}\label{rem:NS4}
Lemma \ref{lem:BC2} can fail for non-strong independence algebras.  Indeed, consider the example from \cite{FL1993} discussed in Remark \ref{rem:NS1}.  Then with $B=\{a_1,a_2\}$ and $C=\{a_3,a_4\}$, both maximal, we have $\codim(B)=1$, yet $\codim(B\cap C)=3$.

Theorem \ref{thm:SubA} also does not (quite) hold for this algebra $A$.  According to GAP \cite{GAP,Semigroups}, ${\Mpres{X_E}{R_E}}$ has size $15$ in this case, while $\Sub(A)$ has size $12$; it contains all subsets of size~${\not=3}$.  We do, however, obtain a presentation for this $\Sub(A)$ by enlarging $R_E^3$ to
\[
\bigset{(x_Bx_C,x_{B'}x_{C'})}{B,C,B',C'\in\Max(A),\ B\cap C=B'\cap C'},
\]
as again confirmed by GAP.  We leave it as an open problem to determine whether a similar presentation exists for $\FSub(A)$ for an arbitrary (non-strong) independence algebra $A$.
\end{rem}

\sect{Generating (finitary) automorphisms}\label{sect:AutA}

The main purpose of this section is to prove results concerning (monoid) generating sets for the automorphism group of a finite-dimensional independence algebra; in fact, we prove more general results concerning \emph{finitary} automorphisms, which we define below.  

Let $A$ be an independence algebra.  As usual, the \emph{fix set} of a partial endomorphism ${\al\in\PEnd(A)}$ is the subalgebra
\[
\Fix(\al) = \set{x\in \dom(\al)}{x\al=x}.
\]
We also write
\[
\fix(\al)=\dim(\Fix(\al)) \AND \cofix(\al)=\codim(\Fix(\al)).
\]
Note that $\cofix(\al)$ is called the \emph{shift} of $\al$ in \cite{FL1993}, and denoted $s(\al)$.
We say $\al\in\PEnd(A)$ is \emph{finitary} if $\cofix(\al)<\aleph_0$, and we write
\[
\FPEnd(A) = \set{\al\in\PEnd(A)}{\cofix(\al)<\aleph_0}
\]
for the set of all finitary partial endomorphisms of $A$, and similarly for $\FEnd(A)$, $\FPAut(A)$ and $\FAut(A)$.  
Of course we have
\[
\FPEnd(A)=\PEnd(A) \iff \dim(A)<\aleph_0,
\]
and similarly with $\FEnd(A)$ and so on.
When $A$ is simply a set (with no operations), ${\FAut(A)=\FG_A}$ is the \emph{finitary symmetric group over~$A$}, consisting of all finitary permutations of~$A$: i.e., the permutations that move only finitely many elements of $A$.

\begin{lemma}\label{lem:FPEndA}
If $A$ is a finite-dimensional or strong independence algebra, then $\FPEnd(A)$ is a submonoid of $\PEnd(A)$.
\end{lemma}

\pf
The result is clear for $\dim(A)<\aleph_0$.  When $A$ is strong it follows from the identity
\[
\cofix(\al\be)\leq\cofix(\al)+\cofix(\be) \qquad\text{for $\al,\be\in\PEnd(A)$.}
\]
This was proved in \cite[Lemma 2.7]{FL1993} for (full) endomorphisms, and the proof there works unchanged for partial endomorphisms.
\epf

For any cardinal $\ka$, we define the set
\[
\Ga_\ka = \set{\al\in\Aut(A)}{\cofix(\al)=\ka} .
\]
For example, $\Ga_0=\{\id_A\}$.  If $\ka>\dim(A)$ then $\Ga_\ka=\es$.

\begin{thm}\label{thm:n-2}
Let $A$ be an independence algebra, and let $\Ga=\Ga_1\cup\Ga_2$.  Then for any ${\al\in\FAut(A)}$, we have $\al=\be_1\cdots\be_k$ for some $\be_1,\ldots,\be_k\in\Ga$ with ${\Fix(\al)\sub\Fix(\be_i)}$ for each~$i$.  Consequently, $\FAut(A)\sub\la\Ga\ra$, with equality if $A$ is finite-dimensional or strong.
\end{thm}

\pf
It suffices to prove the first claim, and for this we use induction on $d=\cofix(\al)$.  If $d=0$ then $\al=\id_A$, and we take $k=0$.  If $d=1$, then $\al\in\Ga$, and we take $k=1$ and $\be_1=\al$.  (We could include $d=2$ in the base case(s), but we prefer not to, as it allows us to only use certain elements of $\Ga_2$, which will then be used in the next proof.)  

We now assume that $d\geq 2$.  We will show that there exists $\be\in\Ga$ such that
\begin{equation}\label{eq:cofix}
\cofix(\al\be^{-1})\leq d-1 \AND \Fix(\al)\sub\Fix(\be).
\end{equation}
It follows from the latter that $\Fix(\al)\sub\Fix(\al\be^{-1})$.  By induction, it follows from the former that $\al\be^{-1}=\ga_1\cdots\ga_l$ for some $\ga_1,\ldots,\ga_l\in\Ga$ such that each $\Fix(\ga_i)\supseteq\Fix(\al\be^{-1})\supseteq\Fix(\al)$.  Since $\al=\ga_1\cdots\ga_l\cdot\be$, the inductive step will be complete.

For the rest of the proof, we fix a basis $X$ for $\Fix(\al)$.  Since $d>0$ we may also fix some $y\in A\sm\Fix(\al)$.  Since $y\not\in\Fix(\al)=\la X\ra$, and since $X$ is independent, it follows from Lemma \ref{lem:Xy} that $X\sqcup y$ is independent.

\pfcase1  Suppose first that $y\al\in\la X\sqcup y\ra$.  Note also that $y\al\not\in\Fix(\al)=\la X\ra$.  (Otherwise, $(y\al)\al=y\al$ would imply $y\al=y$, contradicting $y\not\in\Fix(\al)$.)  It then follows from~\ref{EP} that $y\in\la X\sqcup y\al\ra$.  Consequently, $\la X\sqcup y\ra=\la X\sqcup y\al\ra$.  Since $X\sqcup y$ is independent, we may extend it to a basis $X\sqcup y\sqcup Z$ of $A$, and we note that $|y\sqcup Z|=\cofix(\al)=d$.  So $Z$ is a $B$-basis of $A$, where $B=\la X\sqcup y\ra=\la X\sqcup y\al\ra$.  It then follows from Lemma \ref{lem:XY} that $X\sqcup y\al\sqcup Z$ is also a basis of $A$.  Now let $\be$ be the endomorphism of $A$ satisfying
\[
x\be=x \COMMA y\be=y\al \AND z\be=z \qquad\text{for all $x\in X$ and $z\in Z$.}
\]
Since $\be$ maps a basis bijectively to a basis, it is an automorphism.  And moreover, the maximal subalgebra $\la X\sqcup Z\ra$ is contained in $\Fix(\be)$, so it follows that $\cofix(\be)\leq1$.  Since ${y\be=y\al\not=y}$, we cannot have $\cofix(\be)=0$, so in fact $\be\in\Ga_1\sub\Ga$.  Certainly ${\Fix(\al)=\la X\ra\sub\Fix(\be)}$, and since $X\sqcup y\sub\Fix(\al\be^{-1})$, we also have $\cofix(\al\be^{-1})\leq |Z|=d-1$.  This completes the proof that~\eqref{eq:cofix} holds in this case.

\pfcase2  Now suppose $y\al\not\in\la X\sqcup y\ra$, and for convenience write $y'=y\al$.  It follows from Lemma~\ref{lem:Xy} that $X\sqcup\{y,y'\}$ is independent, so we may extend it to a basis $X\sqcup\{y,y'\}\sqcup Z$ of $A$.  This time we have $\big|\{y,y'\}\sqcup Z\big|=\cofix(\al)=d$, and we define $\be\in\Aut(A)$ by
\[
x\be=x \COMMA y\be=y' \COMMA y'\be = y \AND z\be=z \qquad\text{for all $x\in X$ and $z\in Z$.}
\]
Since $X\sqcup Z\sub\Fix(\be)$ but $y\not\in\Fix(\be)$, we have $1\leq\cofix(\be)\leq 2$, so $\be\in\Ga$, and again we have $\Fix(\al)\sub\Fix(\be)$.  Also, since $X\sqcup y\sub\Fix(\al\be^{-1})$ we have $\cofix(\al\be^{-1}) \leq |Z\sqcup y'|=d-1$, and again \eqref{eq:cofix} holds.
\epf

\begin{rem}\label{rem:yy'}
In the above proof, we noted that the automorphism $\be$ constructed during Case~2 satisfied $\cofix(\be)\leq 2$.  It is actually possible that $\cofix(\be)=1$, so that in fact $\be\in\Ga_1$.  For example, if $A$ is a vector space, then $X\sqcup\{y,y+y'\}\sqcup Z$ is also a basis of $A$, and $(y+y')\be=y+y'$.  It follows that $\FAut(A)=\la\Ga_1\ra$ if $A$ is a vector space.  
\end{rem}

Among other things, the next result shows that the situation just described is true of a wider class of algebras:

\newpage

\begin{thm}\label{thm:n-1}
Let $A$ be an independence algebra.
\ben
\item \label{n-11} If $A$ satisfies the conditions of Proposition \ref{prop:sub1}, then Theorem \ref{thm:n-2} holds with $\Ga_1\cup\Ga_2$ replaced by $\Ga_2$.
\item \label{n-12} If $A$ satisfies the conditions of Proposition \ref{prop:sub2}, but not those of Proposition \ref{prop:sub1}, then Theorem \ref{thm:n-2} holds with ${\Ga_1\cup\Ga_2}$ replaced by $\Ga_1$.
\item \label{n-13} If $A$ satisfies either of the conditions of Lemma \ref{lem:sub4}, then Theorem \ref{thm:n-2} holds with ${\Ga_1\cup\Ga_2}$ replaced by $\Ga_1$.
\een
\end{thm}

\pf
\firstpfitem{\ref{n-11}}  In this case, the subalgebras of codimension $1$ are of the form $A\sm x$, and if an automorphism fixes $A\sm x$ pointwise, it must also fix $x$.  So $\Ga_1=\es$, and the claim follows.

For later reference, it is worth discussing an alternative proof of this part.  We know from Proposition~\ref{prop:sub2} that ${\Aut(A)=\set{\id_C\cup\al}{\al\in\G_X}}$, where again we write $X=A\sm\CA$, and of course it follows that
\[
\FAut(A) = \set{\id_C\cup\al}{\al\in\FG_X}.
\]
It is well known that finite symmetric groups are generated by transpositions.  So too therefore is the finitary symmetric group $\FG_X$, and so $\FAut(A)$ is generated by automorphisms of the form $\be=\id_C\cup(x,y)$ for distinct $x,y\in X$ (again using standard cycle notation).  But for this $\be$ we have $\Fix(\be)=A\sm\{x,y\}=\big\la X\sm\{x,y\}\big\ra$.  Since $X$ is independent, so too is $X\sm\{x,y\}$, so $\Fix(\be)$ has codimension $2$.

\pfitem{\ref{n-12}}  This is similar to the second proof of \ref{n-11}.  Again we have ${\FAut(A)=\set{\id_C\cup\al}{\al\in\FG_X}}$, where $X=A\sm\CA$.  So $\FAut(A)$ is generated by automorphisms of the form $\be=\id_C\cup(x,y)$ for distinct $x,y\in X$, and we still have $\Fix(\be)=A\sm\{x,y\}=\big\la X\sm\{x,y\}\big\ra$.  Since $X$ is not independent, but $X\sm x$ and $X\sm y$ both are, it follows that $\Fix(\be)$ is maximal, and hence has codimension $1$.

\pfitem{\ref{n-13}}  Looking at the proof of Theorem \ref{thm:n-2}, it suffices to show that the automorphism $\be$ constructed during Case 2 can be expressed as $\be=\ga_1\ga_2\ga_3$, where each $\ga_i\in\Ga_1$ satisfies $\Fix(\al)\sub\Fix(\ga_i)$.  To do so, we keep the notation of that part of the proof, and it will be convenient to write $U=X\sqcup Z$.  So $U\sqcup\{y,y'\}$ is a basis of $A$, and $\be\in\Aut(A)$ acts via
\[
u\be=u \COMMA y\be=y' \AND y'\be=y \qquad\text{for all $u\in U$.}
\]
By Lemma \ref{lem:sub4}\ref{sub42}, we may choose some $w\in A$ such that $w\not\in\la U\sqcup y\ra\cup\la U\sqcup y'\ra$.  It follows quickly that $U\sqcup\{y,w\}$ and $U\sqcup\{y',w\}$ are both bases of $A$.  We then define $\ga_1,\ga_2,\ga_3\in\Ga_1$ by their action on the various bases by
\begin{align*}
u\ga_1 &= u, & u\ga_2 &= u, & u\ga_3 &= u,\\
y\ga_1 &= y, & y\ga_2 &= y', & y'\ga_3 &= y',\\
y'\ga_1 &= w, & w\ga_2 &= w, & w\ga_3 &= y,
\end{align*}
for each $u\in U$.  Verification that $\be=\ga_1\ga_2\ga_3$ is straightforward, and we have $\Fix(\al)=\la X\ra$, and~$X\sub U$.
\epf

\begin{rem}
As we observed in the proof, parts \ref{n-11} and \ref{n-12} of Theorem \ref{thm:n-1} boil down to the fact that finite symmetric groups are generated by transpositions.  The standard Coxeter presentation for a finite symmetric group involves (simple) transpositions \cite{Moore1897}; this is stated in Theorem \ref{thm:Gn} below.  We believe it would be interesting to look for presentations for $\FAut(A)$ in the cases that~$A$ does not satisfy the conditions of Propositions \ref{prop:sub1} or \ref{prop:sub2}, but this is beyond the scope of the current work.  To our knowledge, such presentations are not even known in the case that $A$ is a finite-dimensional vector space.

For infinite $\dim(A)$, generating sets for the full automorphism group $\Aut(A)$ can be very complicated, even when $A$ is simply a set (since $\G_A$ is uncountable).
\end{rem}

\begin{rem}
Since vector spaces of dimension $2\leq n<\aleph_0$ satisfy the conditions of Lemma~\ref{lem:sub4} (cf.~Remark \ref{rem:sub42}), Theorem \ref{thm:n-1}\ref{n-13} implies the well-known fact that general linear groups are generated by matrices whose fix spaces have codimension $1$.  Indeed, it is a basic fact of linear algebra that every invertible matrix is a product of `elementary row operation' matrices, and one can easily see that the fix space of each such matrix has codimension $1$.  This is clear for row operations of the form $R_i\to \lam R_i$ or $R_i\to R_i+R_j$; for a row swap $R_i\leftrightarrow R_j$, we use the argument in~Remark~\ref{rem:yy'}, with $y$ and $y'$ the $i$th and $j$th standard basis vectors.
\end{rem}

\begin{rem}
Consider again the algebra $A=\{a,b,c,d\}$ from Remark \ref{rem:sub3}.  Using the notation of that remark, we have
\[
\Ga_0=\{\id_A\} \COMMA \Ga_1=\{\al_1,\al_2\} \AND \Ga_2=\{\al_3,\be_1,\be_2,\be_3,\be_4\}.
\]
Then $A$ satisfies none of the conditions of Proposition \ref{prop:sub1}, Proposition~\ref{prop:sub2} or Lemma \ref{lem:sub4}, and one can check that
\[
\Aut(A) = \la\Ga_2\ra \qquad\text{but}\qquad \la\Ga_1\ra = \{\id_A,\al_1,\al_2,\al_3\ra \not= \Aut(A).
\]
\end{rem}

\sect{Action pairs and subsemigroups}\label{sect:AES}

Consider a partial endomorphism $\al:B\to A$, where $A$ is an independence algebra.  Since any basis for $B$ can be extended to a basis of $A$, it follows from \ref{FBP} that $\al$ can be extended to an endomorphism $\be:A\to A$.  Note then that $\al=\be\restr_B=\id_B\cdot\be$, so that every partial endomorphism is a restriction of an endomorphism.  This and more can be expressed in terms of (strong) action pairs, as we explore in this section.

In what follows, we continue to write
\[
\E_A = \set{\id_B}{B\leq A}
\]
for the semilattice of all partial identities of $A$.  We observed in Section \ref{sect:EndA} that $\E_A$ is isomorphic to the $\cap$-semilattice $\Sub(A)$.
We also write
\[
\Sing(\E_A) = \E_A\sm\{\id_A\} = \set{\id_B}{B<A},
\]
and we note that $\Sing(\E_A)$ is a subsemilattice of $\E_A$, regardless of the dimension of $A$.  We have $\Sing(\E_A)=\es$ if $\dim(A)=0$.

\begin{prop}\label{prop:ESA1}
For any independence algebra $A$, the following are all strong action pairs in~$\PEnd(A)$\emph{:}
\ben
\item \label{ESA11} $(\E_A,\End(A))$,
\item \label{ESA12} $(\Sing(\E_A),\End(A))$,
\item \label{ESA13} $(\E_A,\Sing(\End(A)))$ if $\dim(A)<\aleph_0$,
\item \label{ESA14} $(\Sing(\E_A),\Sing(\End(A)))$ if $\dim(A)<\aleph_0$,
\item \label{ESA15} $(\E_A,\Aut(A))$,
\item \label{ESA16} $(\Sing(\E_A),\Aut(A))$.  
\een
\end{prop}

\pf
This follows from Propositions \ref{prop:LR2} and \ref{prop:LRIA}, given \eqref{eq:EA} and \eqref{eq:TA}.
\epf

\begin{rem}
In fact the previous result holds for an arbitrary algebra $A$ (with the same proof), but note that parts \ref{ESA13} and \ref{ESA14} hold precisely when $\Sing(\End(A))$ is a subsemigroup of $\End(A)$.  (When $A$ is not an independence algebra, there might not be a meaningful notion of dimension.)
\end{rem}

Each of the pairs $(E,S)$ listed in Proposition \ref{prop:ESA1} leads as usual to the subsemigroup $ES$ of ${M=\PEnd(A)}$.  The next result identifies the set $ES$ for each pair.  Note that we do not assume $\dim(A)<\aleph_0$ in the third and fourth parts; even though the stated equality still holds when $\dim(A)\geq\aleph_0$, the subsets $S$ and $ES$ are not semigroups in this case.

\begin{thm}\label{thm:ESA2}
For any independence algebra $A$, we have:
\ben
\item \label{ESA21} $\E_A\cdot\End(A) = \PEnd(A)$,
\item \label{ESA22} $\Sing(\E_A)\cdot\End(A) = \PEnd(A)\sm\End(A)$,
\item \label{ESA23} $\E_A\cdot\Sing(\End(A)) \sub \Sing(\PEnd(A))$, with equality unless $\dim(A)=1$ and $A$ has no constants, 
\item \label{ESA24} $\Sing(\E_A)\cdot\Sing(\End(A)) \sub \PEnd(A)\sm\End(A)$, with equality unless $\dim(A)=1$ and $A$ has no constants, 
\item \label{ESA25} $\E_A\cdot\Aut(A) \sub \PAut(A)$, with equality if and only if $\dim(A)<\aleph_0$,
\item \label{ESA26} $\Sing(\E_A)\cdot\Aut(A) \sub \Sing(\PAut(A))$, with equality if and only if $\dim(A)<\aleph_0$.
\een
\end{thm}

\pf
The forwards inclusions are straightforward in each case (and trivial for \ref{ESA21}), given that
\[
\id_B\cdot\al=\al\restr_B \qquad\text{for all $\al\in\End(A)$ and $B\leq A$.}
\]
For cases \ref{ESA22}--\ref{ESA26}, we additionally rely on (combinations of) the following simple observations:
\bit
\item If $B\not=A$ (i.e., $\id_B\in\Sing(\E_A)$), then $\al\restr_B$ is not an endomorphism (and of course not an automorphism).
\item If $\al$ is not an automorphism (i.e., $\al\in\Sing(\End(A))$), then neither is $\al\restr_B$.  (Indeed, this is clear if $B=A$, and follows from the previous point if $B<A$.)
\item If $\al$ is injective, then so too is $\al\restr_B$.
\eit
We now work towards the reverse inclusions.  For this, we fix some ${\al\in\PEnd(A)}$, and write $B=\dom(\al)$.  Let $X$ be a basis of $B$, and extend this to a basis $X\sqcup Y$ of $A$.  Let $f:Y\to A$ be an arbitrary function, so that $\al\restr_X\sqcup f:X\sqcup Y\to A$ extends to an endomorphism $\be:A\to A$.  For any $x\in X$, we have $x\be=x\al$, so since $X$ is a basis it follows that $b\be=b\al$ for all $b\in B$.  Consequently, we have
\[
\al=\be\restr_B=\id_B\cdot\be \qquad\text{with $\id_B\in\E_A$ and $\be\in\End(A)$.}
\]
This completes the proof of part \ref{ESA21}.  For the other parts, we need to ensure that the additional constraints on $\al$ (by virtue of belonging to the appropriate right-hand set) allow us to define~$f$ in such a way that the appropriate conditions on~$\be$ (and $\id_B$) hold, or else explain why the backwards inclusions do not hold in the relevant cases.

\pfitem{\ref{ESA22}}  If $\al\in\PEnd(A)\sm\End(A)$, then $B\not=A$, so $\id_B\in\Sing(\E_A)$.

\pfitem{\ref{ESA23}}  First suppose $\dim(A)=1$ and $A$ has no constants.  Then by Remark \ref{rem:small}\ref{small7}, we have
\begin{equation}\label{eq:ES_es1}
\E_A\cdot\Sing(\End(A)) = \Sing(\E_A)\cdot\Sing(\End(A)) = \es 
\end{equation}
and
\begin{equation}\label{eq:ES_es2}
\Sing(\PEnd(A)) = \PEnd(A)\sm\End(A) = \{\es\},
\end{equation}
so we have strict inclusion in this case.

We now assume that $\dim(A)\not=1$ and/or $A$ has constants.  

Suppose $\al\in\Sing(\PEnd(A))$.  We must show that $\be$ can be chosen so that $\be\in\Sing(\End(A))$.  First, if $\al\in\End(A)$, then $\al\in\Sing(\End(A))$ and $B=A$, so that $\be=\al\in\Sing(\End(A))$.  So now suppose $\al\not\in\End(A)$.  In particular, $B<A$, which implies $Y\not=\es$, and so $\dim(A)\geq1$.  Since we have also assumed that $\dim(A)\not=1$ and/or $A$ has constants, it then follows from Lemma \ref{lem:Sing_es} that $\Sing(\End(A))\not=\es$.

If $X\not=\es$, then we fix some $x\in X$, and simply ensure that $yf=x\al$ for all $y\in Y$; in this case, $\al\restr_X\sqcup f$ is not injective, and so $\be\in\Sing(\End(A))$.

If $X=\es$, then $B=\CA$ and $\al=\id_B$ (since every partial endomorphism fixes $\CA$ pointwise).  So in fact, $\al=\id_B\cdot\be$ for \emph{any} $\be\in\End(A)$.  In particular, we can take $\be$ to be any element of $\Sing(\End(A))$, which we noted above is non-empty.

\pfitem{\ref{ESA24}}  If $\dim(A)=1$ and~$A$ has no constants, then strict inclusion again follows from~\eqref{eq:ES_es1} and~\eqref{eq:ES_es2}.  So now suppose $\dim(A)\not=1$ and/or $A$ has constants, and let ${\al\in\PEnd(A)\sm\End(A)}$.  As in \ref{ESA22}, we have $\id_B\in\Sing(\E_A)$.  As in \ref{ESA23}, we can choose~$\be$ so that $\be\in\Sing(\End(A))$.

\pfitem{\ref{ESA25}}  First suppose $\dim(A)\geq\aleph_0$, and fix any basis $Z$ for $A$.  Let $z\in Z$, and let $\be$ be an endomorphism extending some bijection $Z\to Z\sm z$.  Then in fact $\be\in\Sing(\PAut(A))$.  But if $\be=\id_B\cdot\ga$ for some $B\leq A$ and $\ga\in\Aut(A)$, then from $\be=\ga\restr_B$ it follows that $B=\dom(\be)=A$, and so $\be=\ga\in\Aut(A)$, a contradiction.  

Now suppose $\dim(A)<\aleph_0$.  Since $\al\in\PAut(A)$, it follows that $X\al$ is a basis for $\im(\al)$ and $|X\al|=|X|$.  We then extend $X\al$ to a basis $X\al\sqcup W$ for $A$, noting that $|W|=|Y|$, and let $f:Y\to W$ be any bijection.  Clearly then $\be\in\Aut(A)$.

\pfitem{\ref{ESA26}}  The proof is essentially the same as \ref{ESA25}.  The same $\be\in\Sing(\PAut(A))$ deals with the case $\dim(A)\geq\aleph_0$.  When $\dim(A)<\aleph_0$, we define $f$ as above, and note that $\al\not\in\Aut(A)$ implies~$B<A$.
\epf

\begin{rem}\label{rem:not_proper_A}
One might wonder if any of the pairs $(E,S)$ in Proposition \ref{prop:ESA1} are proper, in the sense of Definition \ref{defn:proper}.  By Proposition \ref{prop:lr} this is equivalent to the left restriction semigroup $ES$ being proper, in the sense of Definition \ref{defn:lr_proper2}.  But it is easy to see that this is never the case (apart from trivially small exceptions).  Indeed, $\id_C$ is a left zero element of $\PEnd(A)$, where $C=\C(A)$, and we have $\id_C\in E=P(ES)$ for each of the pairs $(E,S)$ listed in the proposition.  It follows immediately that $\si=\si_{ES}=\nabla_S$ (the universal relation) in each case.  The equivalence~\eqref{eq:lr_proper2} becomes
\[
\al=\be \iff \dom(\al)=\dom(\be) \qquad\text{for all $\al,\be\in ES$,}
\]
and this clearly does not hold (apart from trivially small exceptions).

On the other hand, it follows from Theorem \ref{thm:cover} (cf.~Theorem \ref{thm:McA}) that the left restriction semigroup $ES$ is covered by the proper left restriction semigroup
\[
\bigset{(\id_B,\al)\in E\rtimes S}{\id_B=\id_B\cdot\al^+} = \bigset{(\id_B,\al)\in E\rtimes S}{B\sub\dom(\al)}.
\]
\end{rem}

\begin{rem}
Part \ref{ESA22} of Theorem \ref{thm:ESA2} is a special case of the right ideal $\Pf T=M\sm T$ discussed in Remark \ref{rem:PS}; cf.~\eqref{eq:EA} and \eqref{eq:TA}.
\end{rem}

By Proposition \ref{prop:ESA1}, $(\E_A,\Aut(A))$ and $(\Sing(\E_A),\Aut(A))$ are strong action pairs in $\PEnd(A)$, for arbitrary $A$.  Theorem \ref{thm:ESA2} characterises the resulting subsemigroups $\E_A\cdot\Aut(A)$ and $\Sing(\E_A)\cdot\Aut(A)$ in the case that $\dim(A)$ is finite.  For the infinite case, we need the notion of \emph{codimension} from Definition \ref{defn:codim}:

\begin{thm}\label{thm:ESA3}
For any independence algebra $A$, we have
\ben
\item \label{ESA31} $\E_A\cdot\Aut(A) = \bigset{\al\in\PAut(A)}{\codim(\dom(\al))=\codim(\im(\al))}$,
\item \label{ESA32} $\Sing(\E_A)\cdot\Aut(A) = \bigset{\al\in\PAut(A)}{\codim(\dom(\al))=\codim(\im(\al))\geq1}$.
\een
\end{thm}

\pf
Since
\[
\E_A\cdot\Aut(A) = \big({\Sing}(\E_A)\cdot\Aut(A)\big) \sqcup \Aut(A),
\]
and since
\[
\Aut(A) = \set{\al\in\PAut(A)}{\codim(\dom(\al))=\codim(\im(\al))=0},
\]
it suffices to prove~\ref{ESA31}.

Beginning with the forwards inclusion, let $\al\in\E_A\cdot\Aut(A)$, so that $\al=\id_B\cdot\be=\be\restr_B$ for some $B\leq A$ and $\be\in\Aut(A)$.  Fix some $B$-basis $X$ of $A$, and also write $C=\im(\al)=B\be$ and $Y=X\be$.  We first claim that $Y$ is a $C$-basis of $A$.  To see this, first note that since $\be$ is an automorphism, we have $A=\la B\cup X\ra = \la (B\cup X)\be\ra = \la C\cup Y\ra$.  It remains (for the claim) to show that $Y$ is $C$-independent, and this is also straightforward.  If $y\in\big\la C\cup (Y\sm y)\big\ra$ for some $y\in Y$, then applying the automorphism $\be^{-1}$, we have $y\be^{-1}\in\big\la B\cup (X\sm y\be^{-1})\big\ra$, with $y\be^{-1}\in X$, contradicting $B$-independence of $X$.  With the claim established, and since $\be$ is injective, we then have
\[
\codim(\dom(\al)) = \codim(B) = |X| = |X\be| = |Y| = \codim(C) = \codim(\im(\al)).
\]

Conversely, suppose $\al\in\PAut(A)$ is such that $\codim(B)=\codim(C)$, where we have written $B=\dom(\al)$ and $C=\im(\al)$.  Let $X$ be a basis for $B$, and let $Y$ and~$Z$ be $B$- and $C$-bases for~$A$, respectively.  Since $\al:B\to C$ is an isomorphism, $X\al$ is a basis of $C$, and it follows from Lemma~\ref{lem:XY} that $X\cup Y$ and $X\al\cup Z$ are bases for $A$.  
Since ${|Y|=\codim(B)=\codim(C)=|Z|}$, we may fix a bijection $\phi:Y\to Z$.  Since $X\cup Y$ is a basis of $A$, we may define $\be\in\End(A)$ to be the endomorphism extending $\al\restr_X\cup\phi:X\cup Y\to X\al\cup Z$.  Since $\al\restr_X\cup\phi$ is a bijection between bases, it follows that $\be\in\Aut(A)$.  For any $x\in X$ we have $x\be=x\al$, so since $B=\la X\ra$ it follows that $b\be=b\al$ for all $b\in B$.  We then have $\al=\be\restr_B=\id_B\cdot\be\in\E_A\cdot\Aut(A)$.
\epf

\begin{rem}
Theorem \ref{thm:ESA3} holds for arbitrary $\dim(A)$.  When $\dim(A)<\aleph_0$, we have
\[
\dim(B)=\dim(C)\implies\codim(B)=\codim(C) \qquad\text{for $B,C\leq A$},
\]
so parts \ref{ESA25} and \ref{ESA26} of Theorem \ref{thm:ESA2} follow from Theorem \ref{thm:ESA3}.
\end{rem}

\begin{rem}\label{rem:FA}
Note that $\E_A$ and $\Aut(A)$ respectively consist of all idempotents and units of the inverse monoid $\PAut(A)$.  It follows that $\E_A\cdot\Aut(A)$ is the largest factorisable inverse submonoid of $\PAut(A)$.  Theorem \ref{thm:ESA3}\ref{ESA31} is therefore a generalisation of the main result of \cite{CH1974}, which concerns the case that $A$ is a set, and says that the largest factorisable inverse submonoid of the symmetric inverse monoid $\I_A$ is
\[
\F_A = \bigset{\al\in\I_A}{|A\sm{\dom(\al)}|=|A\sm{\im(\al)}|}.
\]
(In this case, an arbitrary subset $B\sub A$ is a subalgebra, and $\codim(B)=|A\sm B|$.)
\end{rem}

\newpage

\begin{rem}\label{rem:PTX}
It is worth considering Proposition \ref{prop:ESA1} and Theorems \ref{thm:ESA2} and \ref{thm:ESA3} in the special case that $A=X$ is simply a set with no operations.  As explained above, we then have ${\PEnd(A)=\PT_X}$, and so on.  We then obtain the following strong action pairs in $M=\PT_X$:
\ben
\item $(\E_X,\T_X)$,
\item $(\Sing(\E_X),\T_X)$,
\item $(\E_X,\Sing(\T_X))$ if $|X|<\aleph_0$,
\item $(\Sing(\E_X),\Sing(\T_X))$ if $|X|<\aleph_0$,
\item $(\E_X,\G_X)$,
\item $(\Sing(\E_X),\G_X)$.  
\een
These action pairs give rise, respectively, to the following product semigroups:
\ben
\item $\E_X\cdot\T_X = \PT_X$,
\item $\Sing(\E_X)\cdot\T_X = \PT_X\sm\T_X$,
\item $\E_X\cdot\Sing(\T_X) = \Sing(\PT_X)$ if $1\not=|X|<\aleph_0$,
\item $\Sing(\E_X)\cdot\Sing(\T_X) = \PT_X\sm\T_X$ if $1\not=|X|<\aleph_0$,
\item $\E_X\cdot\G_X = \F_X $,
\item $\Sing(\E_X)\cdot\G_X = \Sing(\F_X)$.  
\een
The results of the current chapter all have interpretations for these subsemigroups of~$\PT_X$.  We will not explicitly state these, however, as they are more readily deduced as special cases of the results of Chapter~\ref{chap:wreath}.  
\end{rem}

\sect{Congruence conditions}\label{sect:Acong}

The action pairs $(E,S)$ listed in Proposition~\ref{prop:ESA1} lead to the subsemigroups $ES$ of ${M=\PEnd(A)}$, which are characterised in Theorems \ref{thm:ESA2} and~\ref{thm:ESA3}.  In each case we have $ES\cong(E\rtimes S)/\th$, as in Proposition \ref{prop:A1}.  Section \ref{sect:cong} gave a number of conditions under which relatively small generating sets exist for the congruence~$\th$; such conditions also featured in a number of results of Chapter~\ref{chap:pres} concerning presentations for $ES$.  The purpose of the current section is to identify which of these conditions hold for the pairs in Proposition~\ref{prop:ESA1}.

\begin{lemma}\label{lem:Athgen1}
Let $(E,S)$ be one of the (strong) action pairs in $\PEnd(A)$ listed in Proposition~\ref{prop:ESA1}, with the exception of $(E,S)=(\Sing(\E_A),\End(A))$ in the case that $\dim(A)\geq\aleph_0$.  Then $(E,S)$ satisfies at least one of the conditions listed in Lemma \ref{lem:w}.
\end{lemma}

\pf
The pairs in parts \ref{ESA11}, \ref{ESA13} or \ref{ESA15} of Proposition \ref{prop:ESA1} satisfy \emph{all} the conditions of Lemma~\ref{lem:w}, since $\E_A$ is a submonoid of $\PEnd(A)$; cf.~Remark \ref{rem:w}.  

The other pairs (not including the stated exception) satisfy condition \ref{w8}, and hence also~\ref{w10}.  To prove this, we need to show that every element of $ES$ has a right identity from $E$.  So let $\al\in ES$, and put $B=\im(\al)$.  Since we clearly have $\al=\al\cdot\id_B$, it remains to show that $B\not=A$.  For this, we consider the three cases (as listed in Proposition \ref{prop:ESA1}) separately.

\pfitem{\ref{ESA12}}  Here Theorem \ref{thm:ESA2} gives $ES=\PEnd(A)\sm\End(A)$, and $\dim(A)$ is finite (because the infinite case is excluded in the lemma).  Since $\al$ is therefore not an endomorphism, $\dom(\al)$ is a proper subalgebra of $A$, and hence $\dim(\dom(\al))<\dim(A)$, as the latter is finite.  It follows that $\dim(B)=\dim(\im(\al))\leq\dim(\dom(\al))<\dim(A)$, and so $B\not=A$.

\pfitem{\ref{ESA14}}  The proof is essentially the same as for \ref{ESA12}, except that finiteness of $\dim(A)$ is built into the assumption of this case in Proposition \ref{prop:ESA1}.  

\pfitem{\ref{ESA16}}  Here we use Theorem \ref{thm:ESA3}\ref{ESA32}, which tells us that $\codim(B)\geq1$, whence $B\not=A$.
\epf

\begin{rem}
When $\dim(A)\geq\aleph_0$, $(E,S)=(\Sing(\E_A),\End(A))$ does not satisfy any of the conditions of Lemma \ref{lem:w}.  Indeed, in this case Theorem \ref{thm:ESA2}\ref{ESA22} gives $ES=\PEnd(A)\sm\End(A)$, and if $\al$ is any partial endomorphism with $\dom(\al)\not=A=\im(\al)$, then no element of ${E=\Sing(\E_A)}$ is a right identity for $\al$.  This shows that condition \ref{w8} of Lemma \ref{lem:w} does not hold, and hence neither do any of the other conditions.  (No such $\al$ exists if $\dim(A)<\aleph_0$.)
\end{rem}

\begin{rem}
In the above proof we showed that condition \ref{w8} of Lemma \ref{lem:w} was satisfied for the cases in which $E=\Sing(\E_A)$, apart from the excluded case; condition \ref{w10} therefore holds as well.  On the other hand, condition~\ref{w7} is never satisfied (in these cases), apart from trivially small exceptions.  Examining the flow of implications in Lemma \ref{lem:w}, it follows that none of \ref{w1}, \ref{w3}--\ref{w7} or~\ref{w9} hold.  Condition \ref{w2} is satisfied for ${(E,S)=(\Sing(\E_A),\Sing(\End(A))}$ with $\dim(A)<\aleph_0$, but not for the other two pairs.
\end{rem}

It follows from Lemma \ref{lem:Athgen1} that Lemma \ref{lem:Om3} applies to each of the action pairs listed in Proposition~\ref{prop:ESA1}, with the exception of $(\Sing(\E_A),\End(A))$ in the case that $\dim(A)\geq\aleph_0$.  More powerful results from Section \ref{sect:cong}---such as Lemmas~\ref{lem:Om4}, \ref{lem:Om5} and~\ref{lem:Om6}, and also Theorems~\ref{thm:Msimp} and~\ref{thm:simp} from Chapter~\ref{chap:pres}---involve pairs $(E,S)$ with~$E$ a (commutative) monoid, satisfying various conditions on joins of the right congruences $\th_e$ (or~$\Th_e$) or sub(semi)groups $\S_e$, for $e\in E$.  The pairs $(E,S)$ from Proposition~\ref{prop:ESA1} with $E=\Sing(\E_A)$ are therefore excluded.  Thus, the remainder of Section \ref{sect:Acong} concerns pairs $(E,S)$ where:
\bit
\item $E=\E_A=\set{\id_B}{B\leq A}$, and
\item $S$ is one of $\End(A)$, $\Sing(\End(A))$ or $\Aut(A)$, where we must additionally assume that ${\dim(A)<\aleph_0}$ in the second case.
\eit
Our main goal here is to determine which of the above-mentioned conditions hold for these pairs.

Recall that for $B\leq A$, we have the right congruence $\th_{\id_B}$  on~$S$, defined by
\begin{equation}\label{eq:thB}
\th_{\id_B} = \bigset{(\al,\be)\in S\times S}{\id_B\cdot\al=\id_B\cdot\be} = \bigset{(\al,\be)\in S\times S}{\al\restr_B=\be\restr_B}.
\end{equation}
(Since this depends on $S$, we will be careful to specify which of the above choices of $S$ statements apply to.)  To simplify notation, we write $\th_B$ for $\th_{\id_B}$.  It is clear that (for any $S$)
\begin{equation}\label{eq:BC}
B\sub C \implies \th_B\supseteq\th_C \qquad\text{for all $B,C\leq A$.}
\end{equation}
It is also worth noting that $\al\restr_B=\be\restr_B$ (i.e., $(\al,\be)\in\th_B$) if and only if $\al$ and $\be$ agree on a basis of $B$.

The next two results apply to \emph{strong} independence algebras (cf.~Definition \ref{defn:strong}).  The first shows that Lemma \ref{lem:Om4}\ref{Om42} applies to the pair $(\E_A,\End(A))$ for strong $A$:

\begin{prop}\label{prop:Ejoins}
If $S=\End(A)$ for a strong independence algebra $A$, then $\th_B\vee\th_C=\th_{B\cap C}$ for all $B,C\leq A$.
\end{prop}

\pf
By \eqref{eq:BC}, we have $\th_B,\th_C\sub\th_{B\cap C}$, and hence $\th_B\vee\th_C\sub\th_{B\cap C}$.  Thus, it remains to show that ${\th_{B\cap C}\sub\th_B\vee\th_C}$.  To do so, fix some $(\al,\be)\in\th_{B\cap C}$, meaning that
\[
\al,\be\in\End(A) \AND \al\restr_{B\cap C}=\be\restr_{B\cap C}.
\]
Fix a basis $X$ for $B\cap C$, and extend this to bases $X\sqcup Y$ and $X\sqcup Z$ for $B$ and $C$, respectively.  By Lemma \ref{lem:strong}, $X\sqcup Y\sqcup Z$ is independent, so there exists an endomorphism $\ga\in\End(A)$ extending the map
\[
X\sqcup Y\sqcup Z \to A:\begin{cases}
x\mt x\al=x\be &\text{for $x\in X$,}\\
y\mt y\al  &\text{for $y\in Y$,}\\
z\mt z\be  &\text{for $z\in Z$.}
\end{cases}
\]
Since $\ga$ agrees with $\al$ on the basis $X\sqcup Y$ of $B$ we have $(\al,\ga)\in\th_B$.  Similarly, $(\ga,\be)\in\th_C$, so it follows that $(\al,\be)\in\th_B\vee\th_C$, as required.
\epf

The situation for the pair $(\E_A,\Sing(\End(A)))$ is almost the same as for $(\E_A,\End(A))$, with one exception (when $\dim(A)=2$ and $A$ has no constants).  The next result concerns the right congruences $\Th_B=\Th_{\id_B}$ on $\Sing(\End(A))\cup\{\id_A\}$, defined in \eqref{eq:The}.  It is again clear that $B\sub C \implies \Th_B\supseteq\Th_C$ for all $B,C\leq A$; cf.~\eqref{eq:BC}.

\begin{prop}\label{prop:SEjoins}
If $S=\Sing(\End(A))$ for a finite-dimensional strong independence algebra~$A$, and if $[\dim(A)\not=2$ or $A$ has constants$]$, then $\Th_B\vee\Th_C=\Th_{B\cap C}$ for all $B,C\leq A$.
\end{prop}

\pf
Again we must show that ${\Th_{B\cap C}\sub\Th_B\vee\Th_C}$.  This is clear if $B$ and $C$ are comparable in the inclusion order, so suppose instead that $B$ and $C$ are incomparable, and let $(\al,\be)\in\Th_{B\cap C}$.  Fix a basis $X$ for $B\cap C$, and extend this to bases $X\sqcup Y$ and $X\sqcup Z$ for $B$ and $C$, respectively.  So again $X\sqcup Y\sqcup Z$ is independent.  While it is possible for $X$ to be empty, both $Y$ and $Z$ are non-empty (by incomparability).  We now consider separate cases.  In each one, we define three singular endomorphisms $\ga_1,\ga_2,\ga_3\in\Sing(\End(A))$ such that
\begin{equation}\label{eq:AB}
\al \mr\Th_B \ga_1 \mr\Th_C \ga_2 \mr\Th_B \ga_3 \mr\Th_C \be,
\end{equation}
which leads to $(\al,\be)\in\Th_B\vee\Th_C$.  The verification that $\ga_1,\ga_2,\ga_3$ are indeed singular, and that~\eqref{eq:AB} holds, is left to the reader.

\pfcase1  If $X\not=\es$, then we fix some $x_0\in X$, and we let $\ga_1,\ga_2,\ga_3$ be any endomorphisms such that for $x\in X$, $y\in Y$ and $z\in Z$, 
\begin{align*}
x\ga_1 &= x\al=x\be, & x\ga_2 &= x\al=x\be, & x\ga_3 &= x\al=x\be,\\
y\ga_1 &= y\al, & y\ga_2 &= x_0\al, & y\ga_3 &= x_0\al,\\
z\ga_1 &= x_0\al, & z\ga_2 &= x_0\al, & z\ga_3 &= z\be.
\end{align*}
\pfcasens2  If $X=\es$ and $A$ has a constant $c$, then we let $\ga_1,\ga_2,\ga_3$ be any endomorphisms such that for $y\in Y$ and $z\in Z$,
\begin{align*}
y\ga_1 &= y\al, & y\ga_2 &= c, & y\ga_3 &= c,\\
z\ga_1 &= c, & z\ga_2 &= c, & z\ga_3 &= z\be.
\end{align*}
\pfcasens3  Finally, suppose $X=\es$ and $A$ has no constants.  By the assumption in the statement of the lemma, we have $\dim(A)\not=2$, and so since $\dim(A)\geq|Y|+|Z|\geq2$, we must have $\dim(A)\geq3$.  By symmetry we may assume that $|Y|\geq|Z|$, and we fix some $y_0\in Y$.  We also extend $Y\sqcup Z$ to a basis $Y\sqcup Z\sqcup W$ for $A$, noting that $W$ could be empty.  We then let $\ga_1,\ga_2,\ga_3$ be any endomorphisms such that for $y\in Y$, $z\in Z$ and $w\in W$,
\begin{align*}
y\ga_1 &= y\al, & y\ga_2 &= y_0\al, & y\ga_3 &= y_0\al,\\
z\ga_1 &= y_0\al, & z\ga_2 &= y_0\al, & z\ga_3 &= z\be,\\
w\ga_1 &= y_0\al, & w\ga_2 &= y_0\al, & w\ga_3 &= y_0\al.
\end{align*}
(Singularity of $\ga_1$ and $\ga_2$ is clear, and also of $\ga_3$ when $W\not=\es$.  If $W=\es$, then ${3\leq\dim(A)=|Y|+|Z|}$ and $|Y|\geq|Z|$ force $|Y|\geq2$, and singularity of $\ga_3$ quickly follows.)
\epf

Proposition \ref{prop:SEjoins} excluded the case in which $\dim(A)=2$ and $A$ has no constants.  We will soon show, in Proposition \ref{prop:2D}, that a weaker (but equally powerful) condition holds in this case; cf.~Remark \ref{rem:2D}.  We begin with a lemma.

\begin{lemma}\label{lem:2D}
Suppose $A$ is a two-dimensional independence algebra with no constants.  Then for any $\al\in\Sing(\End(A))\cup\{\id_A\}$ we have $(\al,\id_A)\in\Th_B\vee\Th_C$ for some $B,C\in\Max(A)$.
\end{lemma}

\pf
Only the case $\al\not=\id_A$ needs proof.
Fix a basis $\{x,y\}$ for $A$, and write $u=x\al$ and $v=y\al$.  Since $\al$ is singular, and since~$A$ has no constants, we have $\rank(\al)=1$, and it quickly follows that $\im(\al)=\la u\ra=\la v\ra$.
Let $B_1=\la x\ra$, $B_2=\la y\ra$ and $C=\la u\ra=\la v\ra$; each is maximal, by Proposition \ref{prop:Max}.  Since $C\not=A=\la x,y\ra$, we may assume by symmetry that $x\not\in C=\la u\ra$.  It follows from Lemma \ref{lem:Xy} that $\{x,u\}$ is independent, and hence a basis for $A$.  In this case we define $\be\in\Sing(\End(A))$ with $x\be=u\be=u$.  Then $\al \mr\Th_{B_1} \be \mr\Th_C \id_A$, so that $(\al,\id_A)\in\th_{B_1}\vee\th_C$.
\epf

\begin{prop}\label{prop:2D}
If $S=\Sing(\End(A))$ for a two-dimensional independence algebra $A$ with no constants, then for any $B<A$ we have
\[
\Th_B = \bigvee_{C\in\Max(A), \atop B\sub C} \Th_C.
\]
\end{prop}

\pf
Fix $B<A$, and write $\Si = \set{C\in\Max(A)}{B\sub C}$.  By the assumptions on $A$, either $B=\es(=\CA)$ or else $B\in\Max(A)$.  In the latter case we have $\Si=\{B\}$, and the result is trivial.  So now we assume that $B=\es$, and here we have $\Si=\Max(A)$.  Since $\Th_B=\Th_\es=\nabla_{S^1}$, we must show that $\bigvee_{C\in\Max(A)}\Th_C=\nabla_{S^1}$.  For this it suffices to show that any $\al\in S=\Sing(\End(A))$ is $\si$-related to $\id_A$, where $\si=\bigvee_{C\in\Max(A)}\Th_C$, and this follows immediately from Lemma~\ref{lem:2D}.
\epf

\begin{rem}\label{rem:2D}
Proposition \ref{prop:SEjoins} feeds into Theorem \ref{thm:simp}\ref{simp2}, and shows that (in the notation of Theorems \ref{thm:ES} and \ref{thm:simp}) relations $R_2$ may be replaced by $R_2''$ for the pair $(U,S)=(\E_A,\Sing(\End(A)))$.  Note that $R_2''$ is defined with respect to a fixed generating set $V=X_U\phi_U$ for $U=\E_A$, coming from a presentation $U\cong\Mpres{X_U}{R_U}$.  When $A$ is finite-dimensional, Proposition \ref{prop:Max} (cf.~Theorem~\ref{thm:SubA}) tells us that the most obvious generating set~$V$ for ${U=\E_A\cong(\Sub(A),\cap)}$ is 
\[
V = \bigset{\id_B}{B\in\Max(A)}.
\]
Proposition \ref{prop:SEjoins} excluded the case in which $A$ is two-dimensional and has no constants.  However, Proposition \ref{prop:2D} shows that Theorem \ref{thm:simp}\ref{simp1} applies in this case with respect to the same set $V=\set{\id_B}{B\in\Max(A)}$, so that relations~$R_2$ may be replaced by $R_2'$ in the notation of that theorem.  It follows that in fact $R_2'=R_2''$ in this case.
\end{rem}

\begin{rem}\label{rem:2D2}
Even though Proposition \ref{prop:SEjoins} excluded the case in which $\dim(A)=2$ and $A$ has no constants, the conclusion of the proposition may still hold in this case.  For example, suppose $A=\{1,2\}$ with no operations, so that $S=\Sing(\End(A))=\Sing(\T_2)=\{\al,\be\}$, where $\al=\trans{1&2\\1&1}$ and $\be=\trans{1&2\\2&2}$.   Recall that Proposition \ref{prop:SEjoins} only needs proof when $B$ and $C$ are incomparable, and up to symmetry this is only the case when $B=\{1\}$ and $C=\{2\}$.  But here we have $\al \mr\Th_B\id_A$ and $\be\mr\Th_C\id_A$, and it quickly follows that $\Th_B\vee\Th_C=\nabla_{S^1}=\Th_\es=\Th_{B\cap C}$.
\end{rem}

Keeping Remark \ref{rem:2D} in mind, in order to give presentations for $\PEnd(A)$ and $\Sing(\PEnd(A))$, we will need to have generating sets for the right congruences $\th_B$ and $\Th_B$ for $B\in\Max(A)$.  It turns out that these are extremely simple:

\begin{lemma}\label{lem:AthB1}
Let $A$ be an independence algebra with $\dim(A)\geq1$, and additionally assume that $\dim(A)\geq2$ if $A$ has no constants.  Then for any ${B\in\Max(A)}$, the right congruence $\th_B$ on $\End(A)$ is generated by the pair~$(\al,\id_A)$ for any idempotent $\al\in\End(A)$ with $\im(\al)=B$.
\end{lemma}

\pf
Throughout the proof we write $\si$ for the right congruence on $\End(A)$ generated by~$(\al,\id_A)$.

Since $\al$ maps $\im(\al)=B$ identically (as it is an idempotent), we have $\id_B\cdot\al=\id_B=\id_B\cdot\id_A$.  This shows that $(\al,\id_A)\in\th_B$, and so $\si\sub\th_B$.

For the reverse inclusion, fix some $(\be,\ga)\in\th_B$.  For any $x\in A$ we have $x\al\in B$ (as $\im(\al)=B$), so since $\be\restr_B=\ga\restr_B$ it follows that $x(\al\be)=(x\al)\be=(x\al)\ga=x(\al\ga)$, so that $\al\be=\al\ga$.  But also, since $\id_A\mr\si\al$, and since $\si$ is a right congruence, we have $\be=\id_A\cdot\be\mr\si\al\be=\al\ga\mr\si\id_A\cdot\ga=\ga$.
\epf

The next result has essentially the same proof:

\begin{lemma}\label{lem:AthB2}
Let $A$ be an independence algebra with $\dim(A)\geq1$, and additionally assume that $\dim(A)\geq2$ if $A$ has no constants.  Then for any ${B\in\Max(A)}$, the right congruence $\Th_B$ on $\Sing(\End(A))\cup\{\id_A\}$ is generated by the pair~$(\al,\id_A)$ for any idempotent $\al\in\End(A)$ with $\im(\al)=B$.  \epfres
\end{lemma}

\begin{rem}
Lemmas \ref{lem:AthB1} and \ref{lem:AthB2} excluded the case in which $\dim(A)=1$ and $A$ has no constants.  But in this case $\PEnd(A)=\Aut(A)\cup\{\es\}$ and $\Sing(\PEnd(A))=\{\es\}$, as in Remark~\ref{rem:small}\ref{small7}.  So it is trivial to derive a presentation for $\PEnd(A)$ from a presentation for $\Aut(A)$, and of course a presentation for $\Sing(\PEnd(A))$ is trivial.
\end{rem}

We now move on to the pair $(\E_A,\Aut(A))$, and here we do not need to assume $A$ is strong.  Since $\Aut(A)$ is a group, Lemma \ref{lem:Om5} applies, and shows that a generating set for the congruence~$\th$ may be constructed using (generating sets for) the subgroups
\begin{equation}\label{eq:SB}
\S_{\id_B} = \set{\al\in\Aut(A)}{\id_B=\id_B\cdot\al} = \set{\al\in\Aut(A)}{\al\restr_B=\id_B} \qquad\text{for $B\leq A$.}
\end{equation}
For simplicity again, we will abbreviate $\S_{\id_B}$ to $\S_B$.  Note that we have the alternative formulation in terms of fix sets:
\[
\S_B = \bigset{\al\in\Aut(A)}{B\sub\Fix(\al)} \qquad\text{for $B\leq A$.}
\]
Lemma \ref{lem:Om6} gives conditions on these subgroups under which the generating set for~$\th$ from Lemma \ref{lem:Om5} can be reduced further.  The simpler of the two conditions is in Lemma \ref{lem:Om6}\ref{Om62}, and says (in this case) that ${\S_B\vee\S_C=\S_{B\cap C}}$ for all $B,C\leq A$.  Unfortunately, this does not hold in general, even in the simple case that $A$ is a set of size $\geq2$ (with no operations), where $\Aut(A)$ is the symmetric group~$\G_A$.  For example, if $A=B\sqcup C$, with $B$ and $C$ both non-empty, then each element of~$\S_B$ and~$\S_C$ fixes $B$ and $C$ setwise; so too therefore does each element of $\S_B\vee\S_C$.  However, $\S_{B\cap C}=\S_\es$ is all of $\G_A$.
Despite this, we will now show that the weaker assumption of Lemma~\ref{lem:Om6}\ref{Om61} does hold for the pair $(\E_A,\Aut(A))$ when $\dim(A)=n<\aleph_0$, with respect to the set
\[
V = \set{\id_B}{B\leq A,\ n-2\leq\dim(B)\leq n-1}.
\]

\begin{prop}\label{prop:Ajoins}
Let $A$ be an independence algebra with $\dim(A)=n<\aleph_0$, and let
\[
Q=\set{B\leq A}{n-2\leq\dim(B)\leq n-1}.
\]
Then for any~$B<A$ we have
\[
\S_B = \bigvee_{C\in Q, \atop B\sub C} \S_C.
\]
\end{prop}

\pf
As explained in Remark \ref{rem:Om6}, we only need to demonstrate the forwards inclusion.  So let $\al\in\S_B$, and fix a factorisation $\al=\be_1\cdots\be_k$, for $\be_1,\ldots,\be_k\in\Ga_1\cup\Ga_2$, as in Theorem \ref{thm:n-2}.  Also write $C_i=\Fix(\be_i)\in Q$ for each $i$, so that each $C_i\supseteq\Fix(\al)\supseteq B$.  But then
\[
\al=\be_1\cdots\be_k\in\S_{C_1}\vee\cdots\vee\S_{C_k} \sub \bigvee_{C\in Q, \atop B\sub C} \S_C.  \qedhere
\]
\epf

\begin{rem}\label{rem:Qn-2n-1}
As in Theorem \ref{thm:n-1}, the set $Q$ in Proposition \ref{prop:Ajoins} can be replaced by
\[
\set{B\leq A}{\dim(B)=n-2} \OR \set{B\leq A}{\dim(B)=n-1},
\]
as appropriate, if the various conditions of Propositions \ref{prop:sub1} and \ref{prop:sub2}, or Lemma \ref{lem:sub4} hold.
\end{rem}

\sect{Presentations}\label{sect:EA}

We conclude Chapter \ref{chap:IA} by discussing presentations for the semigroups
\[
T_1=\PEnd(A) \COMMA T_2=\Sing(\PEnd(A)) \AND T_3=\PAut(A),
\]
where $A$ is a finite-dimensional strong independence algebra.
As we have seen, each arises as a product $T_i=ES_i$, where $E=\E_A$ and 
\[
S_1=\End(A) \COMMA S_2=\Sing(\End(A)) \AND S_3=\Aut(A).
\]
We have already noted that presentations are not known for the $S_i$ in general.  However, it is possible to describe presentations for each $T_i$ modulo a presentation for the corresponding $S_i$.

For example, suppose $S=\End(A)$ has presentation $\Mpres{X_S}{R_S}$ via $\phi_S:X_S^*\to \End(A)$.  We know from Theorem \ref{thm:SubA} that $E=\E_A(\cong\Sub(A)=\FSub(A))$ has presentation $\Mpres{X_E}{R_E}$ via $\phi_E:X_E^*\to E:x_B\mt \id_B$.  It then follows from Theorem \ref{thm:ESmon0} that $ES=\PEnd(A)$ has presentation
\[
\Mpres{X_E\cup X_S}{R_E\cup R_S\cup R_1\cup R_\Om},
\]
where the additional sets of relations $R_1$ and~$R_\Om$ are as follows.  For the former (and making the respective notational substitutions $z$ and $x_B$ for the letters $x\in X_S$ and $y\in X_E(=X_U)$ in the original definition of $R_1$ in \eqref{eq:R1}), we have
\[
R_1 = \bigset{(z x_B , {}^zx_B\cdot z)}{z\in X_S, \ B\in\Max(A)},
\]
where ${}^zx_B$ is some word over $X_E$ mapping to ${}^{z\phi_S}\id_B=\id_{B(z\phi_S)^{-1}}$.
As explained in Remark \ref{rem:ESmon0} (see \eqref{eq:ROm3}, by Proposition~\ref{prop:Ejoins} and Lemma \ref{lem:AthB1} (and excluding trivial small cases), we can take
\[
R_\Om = \bigset{(x_B,x_B N_S(\al_B))}{B\in\Max(A)},
\]
where for each $B\in\Max(A)$, $\al_B\in\End(A)$ is some fixed idempotent with image $B$, and where ${N_S:\End(A)\to X_S^*}$ is a normal form function.

The situation for $S=\Sing(\End(A))$ is very similar.  We first note that $(E,S)$ satisfies Assumption \ref{ass:US}, with $E$ commutative.  We may therefore apply Theorem \ref{thm:simp}, Propositions~\ref{prop:SEjoins} and~\ref{prop:2D}, and Remark \ref{rem:2D}.  The latter explains whether to use part \ref{simp1} or \ref{simp2} of Theorem \ref{thm:simp}, but in either case the resulting presentation for $\Sing(\PEnd(A))$ has the form
\[
\Spres{X_E\cup X_S}{R_E\cup R_S\cup R_1\cup R_2''},
\]
where this time $\Spres{X_S}{R_S}$ is a presentation for $S=\Sing(\End(A))$.  The sets $R_1$ and $R_2''$ here have the same form as $R_1$ and $R_\Om$ above (for $S=\End(A)$).  Note that ${}^zx_B\in X_E^*$ might be empty, but $N_S(\al_B)\in X_S^+$ is always non-empty.

In the case that $S=\Aut(A)$, Theorem \ref{thm:ESmon0} and Remark \ref{rem:ESmon0} apply, drawing from Lemma~\ref{lem:Om6}\ref{Om61} and Proposition~\ref{prop:Ajoins}.  These yield the presentation 
\[
\Mpres{X_E\cup X_S}{R_E\cup R_S\cup R_1\cup R_\Om},
\]
where $R_1$ is yet again as above, and
\begin{align}
\nonumber R_\Om &= \bigset{(x_B,x_B\cdot N_S(\al))}{B\in\Max(A),\ \al\in\Ga_B}\\
\label{eq:ROm} & \quad \cup \bigset{(x_Bx_C,x_Bx_C\cdot N_S(\al))}{B,C\in\Max(A),\ B\not=C,\ \al\in\Ga_{B\cap C}},
\end{align}
where each $\Ga_B$ is a generating set for $\S_B=\S_{\id_B}$.  Such generating sets can be deduced from Theorems~\ref{thm:n-2} and \ref{thm:n-1}.  As in Remark \ref{rem:Qn-2n-1}, we only need one of the two sets in the union \eqref{eq:ROm} if any of the conditions of Theorem \ref{thm:n-1} hold.

\chap{Wreath products and free acts}\label{chap:wreath}

This final chapter concerns a number of (transformational) wreath products $M\wr S$, where $M$ is an arbitrary monoid, and~$S$ is a subsemigroup of some partial transformation semigroup $\PT_X$.  The standard way to define the wreath product $M\wr\T_X$ is as a semidirect product $M^X\rtimes\T_X$, where~$M^X$ is the direct product of $|X|$ copies of $M$, and where $\T_X$ acts on the coordinates of tuples from~$M^X$; the full definitions are given below.  (This distinguishes the \emph{transformational} wreath products from the \emph{unrestricted} wreath products considered implicitly in Section \ref{sect:embed1}.  The unrestricted wreath product of $M$ with $\T_X$ is the semidirect product $M^{\T_X}\rtimes\T_X$.)

There are (at least) two standard ways to generalise this construction of $M\wr\T_X$ in order to define~$M\wr\PT_X$: either as a suitable subsemigroup or quotient of a semidirect product $M_0^X\rtimes\PT_X$, where $M_0$ is $M$ with a zero adjoined.  These definitions are of course interchangeable, but among other things we will see in this chapter that action pairs provide a natural mechanism for passing between the two viewpoints.  When $M$ is a group, $M\wr\PT_X$ is isomorphic to the partial endomorphism monoid of a free $M$-act of rank $|X|$.  Since free group-acts are (strong) independence algebras, it follows that the results of this chapter and those of Chapter \ref{chap:IA} have a common specialisation.  We will sometimes be able to use results from Chapter \ref{chap:IA} to assist in proofs in the current chapter. 

We begin in Section \ref{sect:Wprelim} with definitions and basic results, including that $M\wr\PT_X$ is a left restriction monoid, and we comment on the above connection with free group-acts. The rest of Chapter \ref{chap:wreath} then proceeds in two largely parallel strands.  In Sections~\ref{sect:WAPI} and \ref{sect:WAPII} we identify two families of action pairs in $M\wr\PT_X$, leading to natural subsemigroups such as $M\wr\T_X$, $M\wr\I_X$, $M\wr(\PT_X\sm\T_X)$, $M\wr\Sing(\PT_X)$ for finite $X$, and so on.  Sections~\ref{sect:WcongI} and~\ref{sect:WcongII} explore `congruence conditions' that feed into general results from Chapters~\ref{chap:AP} and~\ref{chap:pres}.  In Sections \ref{sect:WpresI} and~\ref{sect:WpresII} we apply these general results to obtain presentations for several of our wreath products; see Theorems~\ref{thm:MwrSingPTn},~\ref{thm:MwrPTn},~\ref{thm:MwrGn},~\ref{thm:MwrTn} and~\ref{thm:MwrIn}.  
Taking $M=\{1\}$, these theorems reduce to well-known presentations for various (partial) transformation semigroups \cite{JEptnsn2,Aizenstat1958,Moore1897,Popova1961,Lavers1997,East2007,East2007b,EL2004}.

Throughout this chapter we will again often identify one-element subsets of $X$ with their unique elements, and so use abbreviations such as $X\sm x\equiv X\sm\{x\}$.

\sect{Preliminaries}\label{sect:Wprelim}

Let $M$ be a monoid with identity $1$, and $X$ an arbitrary set.  Let $M_0=M\sqcup\{0\}$, where $0$ is a symbol not belonging to $M$, acting as an adjoined zero element, even if $M$ already had a zero element.  

We write~$M_0^X$ for the set of all $X$-tuples $\ba=(a_x)_{x\in X}$ over $M_0$; we usually abbreviate such a tuple to $\ba=(a_x)$.  As usual, $M_0^X$ is a monoid under the componentwise product.  The identity of $M_0^X$ is $\bone=(1)_{x\in X}$, the $X$-tuple with all entries equal to $1$.  Before we define our wreath products, we fix some basic notation.

The \emph{support} of $\ba=(a_x)\in M_0^X$ is the set $\supp(\ba)=\set{x\in X}{a_x\not=0}$.  It is worth noting that
\begin{equation}\label{eq:suppab}
\supp(\ba\bb)=\supp(\ba)\cap\supp(\bb) \qquad\text{for all $\ba,\bb\in M_0^X$.}
\end{equation}
For $B\sub X$, we write $\bone_B\in M_0^X$ for the \emph{indicator function} of $B$, defined by
\begin{align}
\nonumber
\bone_B = (b_x) &\WHERE b_x = \begin{cases}
\mathrlap{1}\hphantom{a_{x\al}} &\text{if $x\in B$,}\\
\mathrlap{0}\hphantom{a_{x\al}} &\text{if $x\in X\sm B$.}
\end{cases}
\intertext{For $\ba=(a_x)\in M_0^X$, and for $B\sub X$, we write $\ba\restr_B\in M_0^X$ for the \emph{restriction of $\ba$ to $B$}, defined by}
\nonumber
\ba\restr_B = \ba\cdot\bone_B = \bone_B\cdot\ba = (b_x) &\WHERE b_x = \begin{cases}
\mathrlap{a_x}\hphantom{a_{x\al}} &\text{if $x\in B$,}\\
\mathrlap{0}\hphantom{a_{x\al}} &\text{if $x\in X\sm B$.}
\end{cases}
\longintertext{Note that $\supp(\ba\restr_B)=\supp(\ba)\cap B$.  \endgraf
The partial transformation monoid $\PT_X$ has a left action on $M_0^X$ by semigroup morphisms, defined for $\al\in\PT_X$ and $\ba=(a_x)\in M_0^X$ by}
\label{eq:ala}
{}^\al\ba = (b_x) &\WHERE b_x = \begin{cases}
a_{x\al} &\text{if $x\in\dom(\al)$,}\\
0 &\text{otherwise.}
\end{cases}
\end{align}
It is worth noting that there are two ways for $b_x$ to equal $0$ in \eqref{eq:ala}; either $x\not\in\dom(\al)$ or else $x\in\dom(\al)$ and $a_{x\al}=0$.  This means that
\begin{equation}\label{eq:suppala}
\supp({}^\al\ba)=\supp(\ba)\al^{-1} = \bigset{x\in\dom(\al)}{x\al\in\supp(\ba)} \qquad\text{for all $\ba\in M_0^X$ and $\al\in\PT_X$.}
\end{equation}
The action from \eqref{eq:ala} has a natural diagrammatic interpretation, as shown in Figure \ref{fig:act} with $X=\{1,\ldots,6\}$ and $\al=\trans{1&2&3&4&5&6\\2&-&3&2&6&6}$.

\begin{figure}[h]
\begin{center}
\scalebox{.9}{
\begin{tikzpicture}[scale=.95]
\tikzstyle{vertex}=[circle,draw=black, fill=white, inner sep = 0.06cm]
\draw[-{latex}](7.5,1)--(9.5,1);
\begin{scope}[shift={(0,0)}]	
\uvs{1,...,6}
\lvs{1,...,6}
\stlines{1/2,3/3,4/2,5/6,6/6}
\foreach \x in {1,...,6} {\draw[dotted] (\x,0)--(\x,-1);}
\draw[|-|] (.5,0)--(.5,2);
\node[left] () at (.5,1) {$\al$};
\end{scope}
\begin{scope}[shift={(0,-3)}]	
\foreach \x in {1,...,6} {\node[vertex] () at (\x,2){$a_{\x}$};}
\node[left] () at (0,2) {$\ba$};
\draw[-{latex}](0,2)--(.5,2);
\end{scope}
\begin{scope}[shift={(10,0)}]	
\uvs{1,...,6}
\lvs{1,...,6}
\stlines{1/2,3/3,4/2,5/6,6/6}
\foreach \x in {1,...,6} {\draw[dotted] (\x,2)--(\x,3);}
\draw[|-|] (6.5,0)--(6.5,2);
\node[right] () at (6.5,1) {$\al$};
\end{scope}
\begin{scope}[shift={(10,1)}]	
\foreach \x/\y in {1/2,3/3,4/2,5/6,6/6} {\node[vertex] () at (\x,2){$a_{\y}$};}
\node[vertex,minimum size=6mm] () at (2,2) {$0$};
\node[right] () at (7,2) {${}^\al\ba$};
\draw[-{latex}](7,2)--(6.5,2);
\end{scope}
\end{tikzpicture}
}
\caption{The action of $\al\in\PT_X$ on $\ba\in M_0^X$ from \eqref{eq:ala}, with $X=\{1,\ldots,6\}$.}
\label{fig:act}
\end{center}
\end{figure}

It is clear that the action in \eqref{eq:ala} is monoidal: i.e., that ${}^{\id_X}\ba=\ba$ for all $\ba\in M_0^X$.  It is not, however, by monoid morphisms, as for any $\al\in\PT_X$, we have ${}^\al\bone=\bone_{\dom(\al)}$.  It follows from this that
\begin{equation}\label{eq:albone}
{}^\al\bone=\bone \iff \al\in\T_X.
\end{equation}

For any subsemigroup $S\leq\PT_X$, the action in \eqref{eq:ala} leads as usual to the semidirect product
\[
M_0^X\rtimes S = \bigset{(\ba,\al)}{\ba\in M_0^X,\ \al\in S} \qquad\text{with operation}\qquad (\ba,\al)\cdot(\bb,\be) = (\ba\cdot{}^\al\bb,\al\be).
\]
If $S$ is a submonoid of $\PT_X$, then $(\bone,\id_X)$ is a left identity for $M_0^X\rtimes S$, as the action is monoidal, but it is only a right identity when the action of $S$ on $M_0^X$ is by monoid morphisms; by \eqref{eq:albone}, this is equivalent to having $S\sub\T_X$.  This all follows from Lemma \ref{lem:USmon}, or by examining products of the form
\begin{equation}\label{eq:abone}
(\ba,\al)\cdot(\bone,\id_X) = (\ba\cdot{}^\al\bone,\al) = (\ba\cdot\bone_{\dom(\al)},\al) = (\ba\restr_{\dom(\al)},\al).
\end{equation}
In any case, when $S$ is a submonoid of $\PT_X$, the local monoid of $M_0^X\rtimes S$ with identity $(\bone,\id_X)$ is $(\bone,\id_X)\cdot(M_0^X\rtimes S)\cdot(\bone,\id_X) = (M_0^X\rtimes S)\cdot(\bone,\id_X)$.  Using~\eqref{eq:abone} it is easy to see that this is
\begin{equation}\label{eq:MSbone}
(M_0^X\rtimes S)\cdot(\bone,\id_X) = \bigset{(\ba,\al)\in M_0^X\rtimes S}{\supp(\ba)\sub\dom(\al)}.
\end{equation}
This leads us naturally to the following:

\begin{defn}\label{defn:wreath}
For any set $X$, any subsemigroup $S$ of $\PT_X$, and any monoid $M$, the \emph{(transformational) wreath product} $M\wr S$ is the set
\[
M\wr S = \bigset{(\ba,\al)\in M_0^X\rtimes S}{\supp(\ba)=\dom(\al)}.
\]
It is routine to check that $M\wr S$ is a subsemigroup of $M_0^X\rtimes S$.  Indeed, given $(\ba,\al),(\bb,\be)\in M\wr S$, we have $(\ba,\al)\cdot(\bb,\be)=(\ba\cdot{}^\al\bb,\al\be)$, and \eqref{eq:suppab} and \eqref{eq:suppala} give
\[
\supp(\ba\cdot{}^\al\bb) = \supp(\ba)\cap\supp(\bb)\al^{-1} = \dom(\al)\cap\dom(\be)\al^{-1} = \dom(\al\be).
\]
\end{defn}

Wreath products of this kind (and similar) have been studied by numerous authors.  See for example \cite{DGY2015,KKM2000,KM1980a,KM1980b,KM1980c,Brookes2020,Skornjakov1979,CE2022,MS2021,Lima_thesis}.  The introduction to \cite{KM1980a} discusses some of the early history of the idea, going back to the work of Specht \cite{Specht1933}.  

\begin{rem}\label{rem:wreath}
If $S\sub\T_X$, then every element $(\ba,\al)\in M\wr S$ satisfies ${\supp(\ba)=\dom(\al)=X}$: i.e., $\ba\in M^X$.  In this case, we have $M\wr S=M^X\rtimes S$.

Given \eqref{eq:abone}, it is clear that $M\wr S$ is a monoid whenever $S$ is a submonoid of $\PT_X$.
\end{rem}

The elements of $M\wr\PT_X$, and their products, have a natural diagrammatic representation, which will be useful in all the calculations to follow.  
Figure \ref{fig:MwrPT6} gives an example with ${X=\{1,\ldots,6\}}$, and
\[
\al = \trans{1&2&3&4&5&6\\2&-&3&2&6&6} \AND \be = \trans{1&2&3&4&5&6\\1&1&-&4&-&4}.
\]
When drawing an element $(\ba,\al)$ of $M\wr\PT_X$, we omit the label $a_x$ of upper vertex $x$ if this label is $0$.  We will also often omit the label if it is $1$.  One can tell whether an omitted label is $0$ or $1$ by the non/existence of an edge at that vertex.  See Figures \ref{fig:Baa}--\ref{fig:lx}.

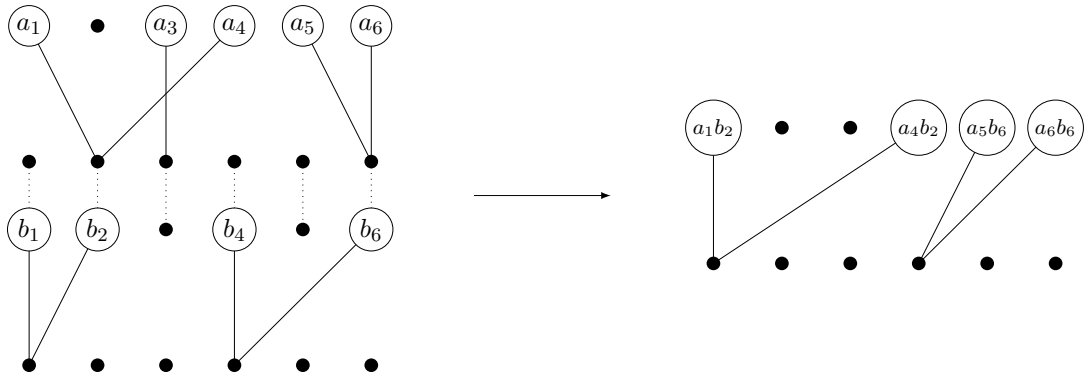
\begin{figure}[h]
\begin{center}
\scalebox{.9}{
\begin{tikzpicture}[scale=1]
\tikzstyle{vertex}=[circle,draw=black, fill=white, inner sep = 0.06cm]
\draw[-{latex}](7.5,-.5)--(9.5,-.5);
\begin{scope}[shift={(0,0)}]	
\lvs{1,...,6}
\stlines{1/2,3/3,4/2,5/6,6/6}
\foreach \x in {1,3,4,5,6} {\node[vertex] () at (\x,2){$a_{\x}$};}
\uvs{2}
\foreach \x in {1,...,6} {\draw[dotted] (\x,0)--(\x,-1);}
\end{scope}
\begin{scope}[shift={(0,-3)}]	
\lvs{1,...,6}
\stlines{1/1,2/1,4/4,6/4}
\foreach \x in {1,2,4,6} {\node[vertex] () at (\x,2){$b_{\x}$};}
\uvs{3,5}
\end{scope}
\begin{scope}[shift={(10,-1.5)}]	
\lvs{1,...,6}
\stlines{1/1,4/1,5/4,6/4}
\foreach \x/\y in {1/2,4/2,5/6,6/6} {\node[vertex] () at (\x,2){{\footnotesize $a_{\x}b_{\y}$}};}
\uvs{2,3}
\end{scope}
\end{tikzpicture}
}
\caption{Elements of $M\wr\PT_X$ (left) and their product (right), with $X=\{1,\ldots,6\}$.}
\label{fig:MwrPT6}
\end{center}
\end{figure}

An important special case arises when $M$ is a group, when wreath products can be viewed as (partial) endomorphism monoids.

\begin{eg}\label{eg:FXM}
For an arbitrary monoid $M$, the class of $M$-acts forms a variety.  Consequently, free $M$-acts exist, and a number of equivalent characterisations exist; see for example \cite{KKM2000}.  Roughly speaking, the free $M$-act of rank $\rho$ can be thought of as $\rho$ disjoint copies of $M$, with the action coinciding with multiplication in $M$.  More formally, given a set $X$, the \emph{free (left) $M$-act over~$X$} is the algebra~$F_X(M)$ with:
\bit
\item underlying set $M\times X$, and 
\item a unary operation $f_a$ for each $a\in M$, defined by $f_a(b,x)=(ab,x)$.
\eit
It is known \cite{Urbanik1966} that $F_X(M)$ is an independence algebra (of dimension $|X|$) if and only if $M$ is a group.  In this case, we have
\begin{align*}
\PEnd(\FXM) &\cong M\wr\PT_X , & \PAut(\FXM) &\cong M\wr\I_X , \\
\End(\FXM) &\cong M\wr\T_X , & \Aut(\FXM) &\cong M\wr\G_X .
\end{align*}
For example, the isomorphism $\PAut(\FXM) \cong M\wr\I_X$ (with $X$ finite) was demonstrated in \cite[Theorem~2.2]{Brookes2020}, and the argument there works virtually unchanged for the others (for any $X$).  Thus, the group case of any result proved in this chapter can also be viewed as a special case of a result from Chapter \ref{chap:IA} (when $A$ is a free group-act).

When $M$ is not a group, we do not have the above isomorphisms.  Indeed, the proof in the group case relies on the fact that the subalgebras of $\FXM$ are all of the form $\FYM=M\times Y$ for some $Y\sub X$, but this is no longer true when $M$ is not a group.  For example, $I\times X$ is a subalgebra of $\FXM$ for any left ideal $I\sub M$.  More generally, one can show that the subalgebras of $\FXM$ are precisely the subsets of the form $\bigcup_{y\in Y} (I_y\times\{y\})$, where $Y$ is an arbitrary subset of $X$, and where each $I_y$ ($y\in Y$) is a left ideal of $M$.  (When $M$ is a group, each $I_y=M$, as groups have no proper non-empty left ideals, and then the subalgebra in question is just $M\times Y=\FYM$.)
\end{eg}

Returning now to the more general discussion of wreath products, we fix a monoid $M$, a set~$X$, and a subsemigroup $S\leq\PT_X$.
It is important to note that $M\wr S$ contains a natural copy of~$S$, namely:
\[
S \equiv \bigset{(\bone_{\dom(\al)},\al)}{\al\in S}.
\]
In this way, $M\wr\PT_X$ contains the semilattice
\[
\E_X = \set{\id_B}{B\sub X} \equiv \bigset{(\bone_B,\id_B)}{B\sub X}
\]
of partial identities.  As in Chapter \ref{chap:IA}, we also write $\Sing(\E_X)=\E_X\sm\{\id_X\}$.  

It is easy to see that for all $(\ba,\al)\in M\wr\PT_X$ and $B\sub X$, we have
\begin{equation}\label{eq:Baa}
\id_B\cdot(\ba,\al)=(\ba\restr_B,\al\restr_B) \AND (\ba,\al)\cdot\id_B = \id_{B\al^{-1}}\cdot(\ba,\al) .
\end{equation}
Figures \ref{fig:Baa} and \ref{fig:aaB} illustrate these identities, with $X=\{1,\ldots,6\}$, $B=\{1,2,4,5\}$, and $\al=\trans{1&2&3&4&5&6\\2&-&3&2&6&6}$.  

\begin{figure}[h]
\begin{center}
\scalebox{.9}{
\begin{tikzpicture}[scale=1]
\tikzstyle{vertex}=[circle,draw=black, fill=white, inner sep = 0.06cm]
\draw[-{latex}](7.5,-.5)--(9.5,-.5);
\begin{scope}[shift={(0,0)}]	
\uvs{1,...,6}
\lvs{1,...,6}
\stlines{1/1,2/2,4/4,5/5}
\foreach \x in {1,...,6} {\draw[dotted] (\x,0)--(\x,-1);}
\end{scope}
\begin{scope}[shift={(0,-3)}]	
\lvs{1,...,6}
\stlines{1/2,3/3,4/2,5/6,6/6}
\foreach \x in {1,3,4,5,6} {\node[vertex] () at (\x,2){$a_{\x}$};}
\uvs{2}
\end{scope}
\begin{scope}[shift={(10,-1.5)}]	
\lvs{1,...,6}
\stlines{1/2,4/2,5/6}
\foreach \x in {1,4,5} {\node[vertex] () at (\x,2){$a_{\x}$};}
\uvs{2,3,6}
\end{scope}
\end{tikzpicture}
}
\caption{An example of $\id_B\cdot(\ba,\al)=(\ba\restr_B,\al\restr_B)$ in $M\wr\PT_X$, with $X=\{1,\ldots,6\}$.}
\label{fig:Baa}
\end{center}
\end{figure}

\begin{figure}[h]
\begin{center}
\scalebox{.9}{
\begin{tikzpicture}[scale=1]
\tikzstyle{vertex}=[circle,draw=black, fill=white, inner sep = 0.06cm]
\draw[{latex}-{latex}](7.5,-.5)--(9.5,-.5);
\begin{scope}[shift={(0,0)}]	
\lvs{1,...,6}
\stlines{1/2,3/3,4/2,5/6,6/6}
\foreach \x in {1,3,4,5,6} {\node[vertex] () at (\x,2){$a_{\x}$};}
\uvs{2}
\foreach \x in {1,...,6} {\draw[dotted] (\x,0)--(\x,-1);}
\end{scope}
\begin{scope}[shift={(0,-3)}]	
\uvs{1,...,6}
\lvs{1,...,6}
\stlines{1/1,2/2,4/4,5/5}
\end{scope}
\begin{scope}[shift={(10,0)}]	
\uvs{1,...,6}
\lvs{1,...,6}
\stlines{1/1,4/4}
\foreach \x in {1,...,6} {\draw[dotted] (\x,0)--(\x,-1);}
\end{scope}
\begin{scope}[shift={(10,-3)}]	
\lvs{1,...,6}
\stlines{1/2,3/3,4/2,5/6,6/6}
\foreach \x in {1,3,4,5,6} {\node[vertex] () at (\x,2){$a_{\x}$};}
\uvs{2}
\end{scope}
\begin{scope}[shift={(5,-9)}]	
\lvs{1,...,6}
\stlines{1/2,4/2}
\foreach \x in {1,4} {\node[vertex] () at (\x,2){$a_{\x}$};}
\uvs{2,3,5,6}
\draw[-{latex}](0,5)--(2,3);
\draw[-{latex}](7,5)--(5,3);
\end{scope}
\end{tikzpicture}
}
\caption{An example of $(\ba,\al)\cdot\id_B = \id_{B\al^{-1}}\cdot(\ba,\al)$ in $M\wr\PT_X$, with $X=\{1,\ldots,6\}$.}
\label{fig:aaB}
\end{center}
\end{figure}

The next result is crucial in all that follows.  As far as we are aware, it has not been stated explicitly in the literature.  Of course it follows from Proposition \ref{prop:LRIA} in the special case that~$M$ is a group, when $M\wr\PT_X\cong\PEnd(\FXM)$ is the partial endomorphism monoid of the (independence) algebra $\FXM$.

\begin{prop}\label{prop:LRMX}
For any set $X$ and monoid $M$, the wreath product $M\wr\PT_X$ is a left restriction monoid under the unary operation $(\ba,\al)^+=\al^+=\id_{\dom(\al)}\equiv(\bone_{\dom(\al)},\id_{\dom(\al)})$.
\end{prop}

\pf
The identities \ref{L1} and \ref{L2} are clear.  For \ref{L3} and \ref{L4}, we use \eqref{eq:Baa}.  We just give the details for \ref{L4}, as \ref{L3} is easier.  So consider elements $x=(\ba,\al)$ and $y=(\bb,\be)$ from $M\wr\PT_X$, and write $B=\dom(\be)$.  Then
\[
xy^+ = (\ba,\al)\cdot\id_B \AND (xy)^+x = \id_D\cdot(\ba,\al) \WHERE D=\dom(\al\be) = B\al^{-1},
\]
so $xy^+=(xy)^+x$ follows from the second identity in \eqref{eq:Baa}.
\epf

In the notation of Section \ref{sect:LRAP}, we have
\begin{align}
\label{eq:EXTX1} P(M\wr\PT_X)  &= \bigset{(\ba,\al)^+}{(\ba,\al)\in M\wr\PT_X} = \E_X \\
\label{eq:EXTX2} \text{and}\qquad T(M\wr\PT_X) &= \bigset{(\ba,\al)\in M\wr\PT_X}{(\ba,\al)^+=\id_X} = M\wr\T_X.
\end{align}

Next we note that we may also identify $M_0^X$ with a submonoid of $M\wr\PT_X$:
\[
M_0^X \equiv \bigset{(\ba,\id_{\supp(\ba)})}{\ba\in M_0^X}.
\]
For $\ba\in M_0^X$ and $\al\in\PT_X$, we have
\begin{align}
\nonumber \ba \cdot \al &= (\ba\restr_B,\al\restr_C) &&\text{where $B=\dom(\al)$ and $C=\supp(\ba)$}\\
\label{eq:aa} &= (\ba\restr_D,\al\restr_D) &&\text{where $D=\dom(\al)\cap\supp(\ba)$.}
\end{align}
This is illustrated in Figure \ref{fig:aa}, with $X=\{1,\ldots,6\}$, $C=\{1,2,4,5\}$ and $\al=\trans{1&2&3&4&5&6\\2&-&3&2&6&6}$.  
In particular, if $(\ba,\al)\in M\wr\PT_X$, then $\supp(\ba)=\dom(\al)$, and so $(\ba,\al)=\ba\cdot\al$.  It is also important to note that
\begin{equation}\label{eq:alba}
\al\cdot\ba = {}^\al\ba\cdot\al \qquad\text{for all $\ba\in M_0^X$ and $\al\in\PT_X$.}
\end{equation}
Indeed, this can be verified by directly showing that both sides evaluate to $(({}^\al\ba)\restr_{\dom(\al)},\al\restr_{\supp({}^\al\ba)})$, keeping in mind the identifications $\al\equiv(\bone_{\dom(\al)},\al)$, etc.  It is also easy to see diagrammatically; cf.~Figure \ref{fig:act}.

\begin{figure}[h]
\begin{center}
\scalebox{.9}{
\begin{tikzpicture}[scale=1]
\tikzstyle{vertex}=[circle,draw=black, fill=white, inner sep = 0.06cm]
\draw[-{latex}](7.5,-.5)--(9.5,-.5);
\begin{scope}[shift={(0,0)}]	
\lvs{1,...,6}
\stlines{1/1,2/2,4/4,5/5}
\foreach \x in {1,2,4,5} {\node[vertex] () at (\x,2){$a_{\x}$};}
\uvs{3,6}
\foreach \x in {1,...,6} {\draw[dotted] (\x,0)--(\x,-1);}
\end{scope}
\begin{scope}[shift={(0,-3)}]	
\uvs{1,...,6}
\lvs{1,...,6}
\stlines{1/2,3/3,4/2,5/6,6/6}
\end{scope}
\begin{scope}[shift={(10,-1.5)}]	
\lvs{1,...,6}
\stlines{1/2,4/2,5/6}
\foreach \x in {1,4,5} {\node[vertex] () at (\x,2){$a_{\x}$};}
\uvs{2,3,6}
\end{scope}
\end{tikzpicture}
}
\caption{An example of $\ba\cdot\al = (\ba\restr_{\dom(\al)},\al\restr_{\supp(\ba)})$ in $M\wr\PT_X$, with $X=\{1,\ldots,6\}$.}
\label{fig:aa}
\end{center}
\end{figure}

The remainder of this chapter will be split into two largely parallel strands, each dealing with one family of action pairs:
\bit
\item Sections \ref{sect:WAPI}--\ref{sect:WpresI} concern pairs of the form $(E,S)$, where $E=\E_X$ or $\Sing(\E_X)$, and where $S\leq M\wr\T_X$.
\item Sections \ref{sect:WAPII}--\ref{sect:WpresII} concern pairs of the form $(U,S)$, where $U=M_0^X$ or $M^X$, and where $S\leq\PT_X$.  
\eit
The pairs of the first type are all strong, as are some (but not all) of the second type.  No pair of the first type is proper (apart from trivially small exceptions), but some of the second type are proper.

\sect{Action pairs and subsemigroups I}\label{sect:WAPI}

Here is the first family of action pairs in $M\wr\PT_X$:

\begin{prop}\label{prop:WAP1}
For any set $X$ and monoid $M$, the following are all strong action pairs in~${M\wr\PT_X}$\emph{:}
\ben
\item \label{ESW11} $(\E_X,M\wr\T_X)$,
\item \label{ESW12} $(\Sing(\E_X),M\wr\T_X)$,
\item \label{ESW13} $(\E_X,M\wr\Sing(\T_X))$ if $|X|<\aleph_0$,
\item \label{ESW14} $(\Sing(\E_X),M\wr\Sing(\T_X))$ if $|X|<\aleph_0$,
\item \label{ESW15} $(\E_X,M\wr\G_X)$,
\item \label{ESW16} $(\Sing(\E_X),M\wr\G_X)$.  
\een
\end{prop}

\pf
This follows from Propositions \ref{prop:LR2} and \ref{prop:LRMX}, together with \eqref{eq:EXTX1} and \eqref{eq:EXTX2}.
\epf

The next result characterises the semigroups $ES$ arising from the pairs $(E,S)$ in Proposition~\ref{prop:WAP1}.  Again, parts \ref{ESW23} and \ref{ESW24} do not assume that $X$ is finite, even though $\Sing(\T_X)=\T_X\sm\G_X$ and $\Sing(\PT_X)=\PT_X\sm\G_X$ are not semigroups for infinite $X$; in these cases, $M\wr\Sing(\T_X)$ simply refers to the relevant \emph{subset} of $M\wr\T_X$, and similarly for $M\wr\Sing(\PT_X)$.  (The $|X|=1$ cases are excluded in these parts because the left-hand sets are empty, and the right-hand sets contain the empty map.)
For parts \ref{ESW25} and \ref{ESW26}, we again write
\[
\F_X = \bigset{\al\in\I_X}{|X\sm\dom(\al)|=|X\sm\im(\al)|}.
\]
As discussed in Remark~\ref{rem:FA}, $\F_X$ is the largest factorisable inverse submonoid of $\I_X$ \cite{CH1974}.  Note that $\Sing(\F_X)=\F_X\sm\G_X$ is always an (almost-factorisable inverse) subsemigroup of $\F_X$, even when $X$ is infinite; cf.~Example~\ref{eg:FIM}.

\begin{prop}\label{prop:WAP2}
For any set $X$ and monoid $M$, we have:
\ben
\item \label{ESW21} $\E_X\cdot(M\wr\T_X) = M\wr\PT_X$,
\item \label{ESW22} $\Sing(\E_X)\cdot(M\wr\T_X) = M\wr(\PT_X\sm\T_X)$,
\item \label{ESW23} $\E_X\cdot(M\wr\Sing(\T_X)) = M\wr\Sing(\PT_X)$ if $|X|\not=1$, 
\item \label{ESW24} $\Sing(\E_X)\cdot(M\wr\Sing(\T_X)) = M\wr(\PT_X\sm\T_X)$ if $|X|\not=1$, 
\item \label{ESW25} $\E_X\cdot(M\wr\G_X) =M\wr\F_X\sub M\wr\I_X$, with equality if and only if $|X|<\aleph_0$,
\item \label{ESW26} $\Sing(\E_X)\cdot(M\wr\G_X)=M\wr\Sing(\F_X) \sub M\wr\Sing(\I_X)$, with equality if and only if $|X|<\aleph_0$.
\een
\end{prop}

\pf
The proof is straightforward, and in each case we can use the relevant part of Theorems~\ref{thm:ESA2} and~\ref{thm:ESA3}, applied to $\PEnd(A)=\PT_X$, where $A=X$ is viewed as an independence algebra with no operations.  

For example, in~\ref{ESW23}, the main work is in showing the backwards inclusion.  For this, consider some $(\ba,\al)$ from $M\wr\Sing(\PT_X)$, where $|X|\not=1$.  By Theorem \ref{thm:ESA2}\ref{ESA23} we have $\al=\id_B\cdot\be$ for some $\be\in\Sing(\T_X)$, where $B=\dom(\al)$.  And then for any $\bb\in M^X$ with $\ba=\bb\restr_B$, we have $(\ba,\al)=\id_B\cdot(\bb,\be)$, with $(\bb,\be)\in M\wr\Sing(\T_X)$.
 \epf

\begin{rem}\label{rem:not_proper_W}
As in Remark \ref{rem:not_proper_A}, and again excluding trivially small exceptions, none of the action pairs $(E,S)$ listed in Proposition \ref{prop:WAP1} are proper.  Hence (cf.~Proposition \ref{prop:lr}), none of the left restriction semigroups $ES$ in Proposition \ref{prop:WAP2} are proper.  But by Theorem \ref{thm:cover} each~$ES$ is covered by the proper left restriction semigroup
\[
\bigset{(\id_B,(\ba,\al))\in E\rtimes S}{\id_B=\id_B\cdot(\ba,\al)^+} = \bigset{(\id_B,(\ba,\al))\in E\rtimes S}{B\sub\dom(\al)}.
\]
\end{rem}

\sect{Congruence conditions I}\label{sect:WcongI}

In Section \ref{sect:WpresI} we will use the results of Chapter \ref{chap:pres} to give presentations for the semigroups arising from (some of) the pairs from Proposition \ref{prop:WAP1}.  We first examine which of the simplifying conditions on the right congruences hold for these pairs.

\begin{lemma}\label{lem:Wthgen1}
Let $(E,S)$ be one of the (strong) action pairs in $M\wr\PT_X$ listed in Proposition~\ref{prop:WAP1}, with the exception of $(E,S)=(\Sing(\E_X),M\wr\T_X)$ in the case that $|X|\geq\aleph_0$.  Then $(E,S)$ satisfies at least one of the conditions listed in Lemma \ref{lem:w}.  
\end{lemma}

\pf
This is proved in similar fashion to Lemma \ref{lem:Athgen1}.  Again, the cases in which $E=\E_X$ are clear, and for those involving $E=\Sing(\E_X)$ it reduces to checking (using Proposition \ref{prop:WAP2}) that any element $(\ba,\al)$ of the semigroup $ES$ satisfies $\im(\al)\not=X$.
\epf

\newpage

As in Section \ref{sect:Acong}, we now restrict our attention to the pairs $(E,S)$ from Proposition \ref{prop:WAP1} in which 
\bit
\item $E=\E_X$, and
\item $S$ is one of $M\wr\T_X$, $M\wr\Sing(\T_X)$ or $M\wr\G_X$, where we must additionally assume that $|X|<\aleph_0$ in the second case.
\eit
We are particularly concerned with determining which of the stronger congruence conditions considered in Section~\ref{sect:cong} and Chapter \ref{chap:pres} are satisfied by these pairs.  Recall that these conditions were stated in terms of the right congruences on $S$ or $S^1$ defined, for each $B\sub X$, by
\begin{align*}
\th_{\id_B} &= \bigset{((\ba,\al),(\bb,\be))\in S\times S}{\id_B\cdot(\ba,\al)=\id_B\cdot(\bb,\be)} \\
&= \bigset{((\ba,\al),(\bb,\be))\in S\times S}{\ba\restr_B=\bb\restr_B,\ \al\restr_B=\be\restr_B}\\
\intertext{for $S=M\wr\T_X$ and $S=M\wr\G_X$, and}
\Th_{\id_B} &= \bigset{((\ba,\al),(\bb,\be))\in S^1\times S^1}{\id_B\cdot(\ba,\al)=\id_B\cdot(\bb,\be)} \\
&= \bigset{((\ba,\al),(\bb,\be))\in S^1\times S^1}{\ba\restr_B=\bb\restr_B,\ \al\restr_B=\be\restr_B} 
\end{align*}
for $S=M\wr\Sing(\T_X)$, where $S^1=(M\wr\Sing(\T_X))\cup\{\id_X\}$.
Again we will abbreviate these to~$\th_B$ and $\Th_B$.  Note that we do not use the sub(semi)groups $\S_B$ in the case of $S=M\wr\G_X$, since this is not a group in general, even though $\G_X$ is.

\begin{lemma}\label{lem:Wjoins1}
If $S=M\wr\T_X$ for any $X$, then $\th_B\vee\th_C=\th_{B\cap C}$ for any $B,C\sub X$.
\end{lemma}

\pf
This follows quickly from Proposition \ref{prop:Ejoins}.  
We need to show that any pair $((\ba,\al),(\bb,\be))$ from $\th_{B\cap C}$ belongs to $\th_B\vee\th_C$, so fix some such pair.  Since $\al\restr_{B\cap C}=\be\restr_{B\cap C}$, the pair $(\al,\be)$ belongs to the `$\th_{B\cap C}$ relation of $\T_X$'.  So we let $\ga\in\T_X$ be as in the proof of Proposition \ref{prop:Ejoins}, so that $\ga\restr_B=\al\restr_B$ and $\ga\restr_C=\be\restr_C$.  Since $\ba\restr_{B\cap C}=\bb\restr_{B\cap C}$, we may also define
\begin{equation}\label{eq:bc}
\bc = (c_x) \WHERE c_x = \begin{cases}
a_x &\text{if $x\in B$,}\\
b_x &\text{if $x\in C$,}\\
1 &\text{otherwise,}
\end{cases}
\end{equation}
and then $(\ba,\al) \mr\th_B (\bc,\ga) \mr\th_C (\bb,\be)$.
\epf

\begin{lemma}\label{lem:Wjoins2}
If $S=M\wr\Sing(\T_X)$ for $|X|<\aleph_0$, then $\Th_B\vee\Th_C=\Th_{B\cap C}$ for any $B,C\sub X$.
\end{lemma}

\pf
This is similar to the proof of Lemma \ref{lem:Wjoins1}, with a little care taken for the $|X|=2$ case (cf.~Proposition \ref{prop:SEjoins} and Remark \ref{rem:2D2}).

Beginning with the case $|X|\not=2$, we fix some $((\ba,\al),(\bb,\be))$ from $\Th_{B\cap C}$, aiming to show it belongs to $\Th_B\vee\Th_C$.  We let $\ga_1,\ga_2,\ga_3$ be as in the proof of Proposition \ref{prop:SEjoins} (Case 1 or 3), and take $\bc$ as in~\eqref{eq:bc}, and we have $(\ba,\al) \mr\Th_B (\bc,\ga_1) \mr\Th_C (\bc,\ga_2) \mr\Th_B (\bc,\ga_3) \mr\Th_C (\bb,\be)$.

For the $|X|=2$ case, we assume without loss of generality that $X=\{1,2\}$.  We need only show that $\Th_B\vee\Th_C=\Th_{B\cap C}$ when $B$ and $C$ are incomparable, and by symmetry the only case to consider is $B=\{1\}$ and $C=\{2\}$.  Since $B\cap C=\es$ and $\Th_\es$ is the universal relation, we must show that every element of $M\wr\Sing(\T_X)$ is $(\Th_B\vee\Th_C)$-related to $\id_X$.  So fix some such $(\ba,\al)\in M\wr\Sing(\T_X)$.  By symmetry, we may assume that $\al=\trans{1&2\\1&1}$.  Write $\ba=(a_1,a_2)$, and set $\bb=(1,a_2)$.  Then $(\ba,\al) \mr\Th_C (\bb,\al) \mr\Th_B \id_X$.
\epf

Lemmas \ref{lem:Wjoins1} and \ref{lem:Wjoins2} allow us to simplify presentations for $M\wr\PT_X$ and $M\wr\Sing(\PT_X)$ arising from Theorems \ref{thm:ESmon0} and \ref{thm:ES}; see Remark \ref{rem:ESmon0} and Theorem \ref{thm:simp}\ref{simp2}.
To utilise these results, we need generators for the right congruences $\th_B$ and $\Th_B$ in the case that $|X\sm B|=1$, and for these we first recall the definition of certain idempotents from $\T_X$.  Given distinct $x,y\in X$, we denote by $\ve_{xy}\in\T_X$ the unique idempotent with image $X\sm y$ and mapping $y\mt x$.  That is,
\begin{equation}\label{eq:vexy}
z\ve_{xy} = \begin{cases}
z &\text{if $z\in X\sm y$,}\\
x &\text{if $z=y$.}
\end{cases}
\end{equation}
See Figure \ref{fig:vexyab} for $\ve_{24}$ and $\ve_{42}$, with $X=\{1,\ldots,6\}$, but ignore the vertex labels.

\begin{lemma}\label{lem:thB1}
If $S=M\wr\T_X$ for any $X$, and if $B=X\sm y$ for some $y\in X$, then $\th_B$ is generated as a right congruence by any pair $(\ve_{xy},\id_X)$ with $x\in B$.
\end{lemma}

\pf
Let $\si$ be the right congruence generated by the stated pair.  Since $\id_B=\id_B\cdot\ve_{xy}$, we have $\si\sub\th_B$.  For the reverse inclusion, let $((\ba,\al),(\bb,\be))\in\th_B$.  So $\ba\restr_B=\bb\restr_B$ and $\al\restr_B=\be\restr_B$, and from these it quickly follows that $\ve_{xy}\cdot(\ba,\al) = \ve_{xy}\cdot(\bb,\be)$.  Since $\id_X \mr\si \ve_{xy}$, and since $\si$ is a right congruence, we have
\[
(\ba,\al)=\id_X\cdot(\ba,\al) \mr\si \ve_{xy}\cdot(\ba,\al) = \ve_{xy}\cdot(\bb,\be) \mr\si \id_X\cdot(\bb,\be) = (\bb,\be).  \qedhere
\]
\epf

The proof of the next result is essentially identical.

\begin{lemma}\label{lem:thB2}
If $S=M\wr\Sing(\T_X)$ for $|X|<\aleph_0$, and if $B=X\sm y$ for some $y\in X$, then $\Th_B$ is generated as a right congruence by any pair $(\ve_{xy},\id_X)$ with $x\in B$.  \epfres
\end{lemma}

\begin{rem}\label{rem:thB}
In Lemmas \ref{lem:thB1} and \ref{lem:thB2}, since the right congruences $\th_B$ and $\Th_B$ are generated by the pair $(\ve_{xy},\id_X)$ for \emph{any} $x\in B=X\sm y$, they are also generated by the set $\bigset{(\ve_{xy},\id_X)}{x\in B}$ of \emph{all} such pairs.
\end{rem}

We now move on to the case of $S=M\wr\G_X$, for which we must assume $X$ is finite.  It will be convenient to first give generators for the right congruences $\th_B$ for arbitrary $B\sub X$.  To describe these, we begin with some notation.  For $x\in X$ and $a\in M$, we write $\bt_{x;a}\in M^X$ for the $X$-tuple whose $x$th coordinate is $a$, and all other entries are $1$.  For distinct $x,y\in X$, we write $\tau_{xy}\in\G_X$ for the transposition that interchanges $x$ and $y$.  

For the next lemma and its proof, recall that we identify a tuple $\ba\in M^X$ with the element $(\ba,\id_X)$ of $M\wr\T_X$.  The statement gives a generating set $\Om_B$ for the right congruence $\th_B$, and $\Om_B$ is a union of two sets of pairs.  A typical pair from the first has the form $(\bt_{x;a},\id_X)$, and as just mentioned this is shorthand for the pair $(\bt_{x;a},\id_X) \equiv ((\bt_{x;a},\id_X),(\bone,\id_X))$.  Similar comments apply to the second set, where we identify $(\tau_{xy},\id_X)\equiv((\bone,\tau_{xy}),(\bone,\id_X))$.

\begin{lemma}\label{lem:thB3}
Suppose $S=M\wr\G_X$ for $|X|<\aleph_0$, and let $\Ga$ be a (monoid) generating set for~$M$.  Then for any $B\sub X$, $\th_B$ is generated as a right congruence by the set
\[
\Om_B = \bigset{(\bt_{x;a},\id_X)}{x\in X\sm B,\ a\in\Ga} \cup \bigset{(\tau_{xy},\id_X)}{x,y\in X\sm B,\ x\not=y}.
\]
\end{lemma}

\pf
Let $\si$ be the right congruence generated by $\Om_B$.  It is easy to check that $\Om_B\sub\th_B$, so we have $\si\sub\th_B$.  It remains to show that $\th_B\sub\si$.

First we claim that
\begin{equation}\label{eq:xia}
\bt_{x;a} \mr\si \id_X \qquad\text{for all $x\in X\sm B$ and $a\in M$.}
\end{equation}
To prove this, write $a=c_1\cdots c_k$, where $c_1,\ldots,c_k\in\Ga$.  Since $(\bt_{x;c_i},\id_X)\in\Om_B\sub\si$ for all $i$, it follows from Lemma \ref{lem:RC} that $\bt_{x;a} = \bt_{x;c_1}\cdots\bt_{x;c_k} \mr\si \id_X$, and \eqref{eq:xia} is proved.

Now let $((\ba,\al),(\bb,\be))\in\th_B$, so that $\ba\restr_B=\bb\restr_B$ and $\al\restr_B=\be\restr_B$.  We must show that $(\ba,\al) \mr\si(\bb,\be)$.
We begin by defining
\[
\bc\in M^X \WHERE c_x = \begin{cases}
a_x=b_x &\text{if $x\in B$,}\\
1 &\text{if $x\in X\sm B$.}
\end{cases}
\]
Without loss of generality, we may assume that $X\sm B=\{1,\ldots,k\}$ for some $k$, and we note that~\eqref{eq:xia} gives $\bt_{i;a_{i}} \mr\si \id_X$ for all $1\leq i\leq k$.  It follows from Lemma \ref{lem:RC} that $\bt_{1;a_{1}}\cdots\bt_{k;a_{k}}\mr\si\id_X$.  Since $\si$ is a right congruence, it follows from this that $\ba=\bt_{1;a_{1}}\cdots\bt_{k;a_{k}}\cdot\bc \mr\si \id_X\cdot\bc = \bc$, and then in turn that $(\ba,\al) = \ba\cdot\al \mr\si \bc\cdot\al$.  

A symmetrical calculation gives $(\bb,\be)\mr\si\bc\cdot\be$, so we can complete the proof by showing that~${\bc\cdot\al\mr\si\bc\cdot\be}$.  Since $\al\restr_B=\be\restr_B$, we have ${}^{\al^{-1}}\bc={}^{\be^{-1}}\bc$ (the action is defined in \eqref{eq:ala}), so if we write $\bd$ for this tuple, then we have $\bc\cdot\al=\al\cdot\bd$ and $\bc\cdot\be=\be\cdot\bd$.  (Indeed, for the former we have $\al\cdot\bd = {}^\al\bd\cdot\al = {}^\al({}^{\al^{-1}}\bc)\cdot\al = \bc\cdot\al$, and the latter is similar.)  Since $\si$ is a right congruence, it therefore suffices to show that $\al\mr\si\be$.  
In fact, since $\al=\al\be^{-1}\cdot\be$ and $\be=\id_X\cdot\be$, it is enough to show that $\al\be^{-1}\mr\si\id_X$.  

From $\al\restr_B=\be\restr_B$, we see that $\al\be^{-1}$ fixes $B$ pointwise, so we may write~$\al\be^{-1}=\tau_{y_1z_1}\cdots\tau_{y_lz_l}$ as a product of transpositions, where $y_i,z_i\in X\sm B$ for each $i$.  Since $(\tau_{y_iz_i},\id_X)\in\Om_B\sub\si$ for all $i$, it follows from Lemma \ref{lem:RC} that $\al\be^{-1}=\tau_{y_1z_1}\cdots\tau_{y_lz_l} \mr\si \id_X$, as required.
\epf

In Lemma \ref{lem:Wjoins4} we show that the right congruences $\th_B$ on $S=M\wr\G_X$ satisfy condition \ref{Om41} of Lemma \ref{lem:Om4}.  In order to do this, it is first convenient to prove the next lemma, which shows that the stronger condition \ref{Om42} of the same lemma \emph{almost} holds.

\begin{lemma}\label{lem:Wjoins3}
If $S=M\wr\G_X$ for $|X|<\aleph_0$, then $\th_B\vee\th_C=\th_{B\cap C}$ for any $B,C\sub X$ with $(X\sm B)\cap (X\sm C)\not=\es$.
\end{lemma}

\pf
Write $\si=\th_B\vee\th_C$; as usual, it is enough to show that $\th_{B\cap C}\sub\si$.  For this, it suffices to show that the generating set $\Om_{B\cap C}$ is contained in $\si$ (cf.~Lemma \ref{lem:thB3}).  Much of $\Om_{B\cap C}$ is in fact contained in~$\Om_B\cup\Om_C$ (keeping $X\sm(B\cap C)=(X\sm B)\cup(X\sm C)$ in mind).  Up to symmetry, any pair belonging to $\Om_{B\cap C}\sm(\Om_B\cup\Om_C)$ has the form~$(\tau_{xy};\id_X)$, for some $x\in X\sm B$ and $y\in X\sm C$.  By assumption we may fix some $z\in(X\sm B)\cap (X\sm C)$, and we note that $(\tau_{xz};\id_X)\in\Om_B\sub\si$ and $(\tau_{yz};\id_X)\in\Om_C\sub\si$.  It then follows from Lemma \ref{lem:RC} that $\tau_{xy} = \tau_{xz}\tau_{yz}\tau_{xz} \mr\si \id_X$, which completes the proof.
\epf

\begin{lemma}\label{lem:Wjoins4}
Let $S=M\wr\G_X$ for $|X|=n<\aleph_0$, and let $Q = \set{B\sub X}{n-2\leq|B|\leq n-1}$.  Then for any $B\sub X$, we have
\[
\th_B = \bigvee_{C\in Q,\atop B\sub C} \th_C.
\]
\end{lemma}

\pf
If $B=X$, then $\th_B=\De_S$, the stated join is empty, and the result is trivial.  If $|B|=n-1$ then we take $k=1$ and $C_1=B$.  

Now suppose $|B|\leq n-2$, and write $X\sm B=\{x_1,\ldots,x_l\}$, noting that $l\geq2$.  For $1\leq i<l$ let $C_i=X\sm\{x_i,x_{i+1}\}\in Q$, so that $B=C_1\cap\cdots\cap C_{l-1}$.  Iterating Lemma \ref{lem:Wjoins3}, we obtain
\[
\th_B = \th_{C_1\cap\cdots\cap C_{l-1}} = \th_{C_1}\vee\cdots\vee\th_{C_{l-1}} \sub \bigvee_{C\in Q,\atop B\sub C} \th_C.
\]
By \eqref{eq:thv1k}, this completes the proof.
\epf

\sect{Presentations I}\label{sect:WpresI}

It is now possible to give explicit presentations for the three semigroups
\bit
\item $M\wr\PT_X = \E_X\cdot(M\wr\T_X)$,
\item $M\wr\I_X = \E_X\cdot(M\wr\G_X)$, and
\item $M\wr\Sing(\PT_X)=\E_X\cdot(M\wr\Sing(\T_X))$,
\eit
where $M$ is an arbitrary monoid, and where $X$ is finite.  Without loss of generality we may assume that $X=\bn=\{1,\ldots,n\}$ for some fixed integer $n$, and to avoid trivialities we assume that $n\geq2$.  As usual, we denote the semigroups involved by $\E_n$, $\T_n$, and so on.

In what follows, we just give the full details for $M\wr\Sing(\PT_n)$, and indicate how to proceed for $M\wr\PT_X$ and $M\wr\I_X$ in Remark \ref{rem:others}; these other two (and more) will be treated fully in Section \ref{sect:WpresII}.  

To deal with $M\wr\Sing(\PT_n)=\E_n\cdot(M\wr\Sing(\T_n))$, we need presentations for $\E_n$ and ${M\wr\Sing(\T_n)}$.
The following is well known (see for example \cite[p.~115]{MMT1987}), and is also a special case of Theorem~\ref{thm:SubA}.  

\begin{thm}\label{thm:En}
For $n\geq0$, the semilattice $\E_n=\set{\id_A}{A\sub\bn}$ has presentation $\Mpres{X_E}{R_E}$ via~$\phi_E:X_E^*\to \E_n$, where:
\bit
\item $X_E=\{t_1,\ldots,t_n\}$,
\item $t_i\phi_E=\id_{\bn\sm i}$ for all $i$,
\item $R_E=\bigset{(t_i^2,t_i)}{i\in\bn}\cup\bigset{(t_it_j,t_jt_i)}{i,j\in\bn}$.  \epfres
\eit
\end{thm}

Next we recall the presentation for $M\wr\Sing(\T_n)$ from \cite{FADEG2019}.  For distinct $i,j\in\bn$, and for $a,b\in M$, we define
\[
\ve_{ij;ab} = (\bc,\ve_{ij}) \in M\wr\Sing(\T_n) \WHERE c_k = \begin{cases}
a &\text{if $k=i$,}\\
b &\text{if $k=j$,}\\
1 &\text{if $k\in \bn\sm\{i,j\}$.}
\end{cases}
\]
See Figure \ref{fig:vexyab} for two examples with $n=6$.

\begin{figure}[h]
\begin{center}
\scalebox{.9}{
\begin{tikzpicture}[scale=1]
\tikzstyle{vertex}=[circle,draw=black, fill=white, inner sep = 0.06cm]
\uvs{1,...,6}
\lvs{1,...,6}
\foreach \x in {1,2,3,5,6} {\stline\x\x}
\stline42
\foreach \x/\y in {2/a,4/b} {\node[vertex] () at (\x,2){$\y$};}
\begin{scope}[shift={(9,0)}]
\uvs{1,...,6}
\lvs{1,...,6}
\foreach \x in {1,3,4,5,6} {\stline\x\x}
\stline24
\foreach \x/\y in {2/b,4/a} {\node[vertex] () at (\x,2){$\y$};}
\end{scope}
\end{tikzpicture}
}
\caption{The generators $\ve_{24;ab}$ (left) and $\ve_{42;ab}$ (right) from $M\wr\Sing(\T_6)$.}
\label{fig:vexyab}
\end{center}
\end{figure}
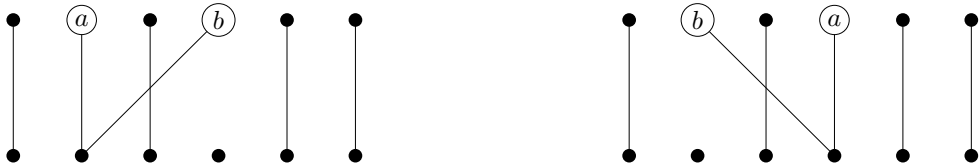

The following is Theorem 5.2 of \cite{FADEG2019}.  Theorem 5.9 of the same paper gives a presentation in terms of the smaller (all-idempotent) generating set 
\[
\set{\ve_{ij;1b}}{i,j\in\bn,\ i\not=j,\ b\in M}
\]
in the special case that the underlying monoid~$M$ has the property that its left ideals are totally ordered under inclusion.  In particular this property holds for groups, as they have no proper non-empty left ideals.  In the next statement, we abbreviate~$e_{ij;11}$ to~$e_{ij}$ in the last two relations.  We display relations as equations, so writing $u=v$ instead of $(u,v)$, for readability; we also do this for several subsequent presentations.  (Note that ${\Sing(\T_n)=\es}$ if~$n<2$.)

\newpage

\begin{thm}\label{thm:MwrSingTn} 
For $n\geq2$, the wreath product $M\wr\Sing(\T_n)$ has presentation $\Spres{X_S}{R_S}$ via $\phi_S:X_S^+\to M\wr\Sing(\T_n)$, where:
\bit
\item $X_S = \set{e_{ij;ab}}{i,j\in\bn,\ i\not=j,\ a,b\in M}$,
\item $e_{ij;ab}\phi_S=\ve_{ij;ab}$ for all $i,j,a,b$,
\item $R_S$ consists of the following relations, with ${a,b,c,d\in M}$ arbitrary, and $i,j,k,l\in\bn$ distinct, in each:
\begin{align*}
e_{ij;ab}e_{ij;cd} = e_{ij;ac,bc} &= e_{ji;ba}e_{ij;dc}, \\
e_{ij;ab}e_{kl;cd} &= e_{kl;cd}e_{ij;ab}, \\
e_{ik;ab}e_{jk;1c} &= e_{ik;ab}, \\
e_{ik;ab}e_{jk;c1} &= e_{ki;ba}e_{ji;c1}e_{ik;11}, \\
e_{ik;aa}e_{jk;b1} &= e_{ik;11}e_{jk;b1}e_{ik;a1}, \\
e_{ij;ab}e_{ik;cd} = e_{ik;ac,d}e_{ij;1,bc} &= e_{jk;bc,d}e_{ij;ac,1}, \\
e_{ij;c,ad}e_{ik;1,bd} = e_{ik;c,bd}e_{ij;1,ad} &= e_{jk;ab}e_{ij;cd}, \\
e_{ki}e_{ij}e_{jk} &= e_{ik}e_{kj}e_{ji}e_{ik}, \\
e_{ki}e_{ij}e_{jk}e_{kl} &= e_{ik}e_{kl}e_{li}e_{ij}e_{jl} . \epfreseq  
\end{align*}
\eit
\end{thm}

Note that in the above relations, products in $M$ sometimes appear in subscripts of letters from $X_S$, and additional commas are included to avoid any ambiguity.  For example, in the first relation, the two subscripts from $M$ appearing in $e_{ij;ac,bc}$ are $ac$ and $bc$.

Here then is the promised presentation for $M\wr\Sing(\PT_n)$.

\begin{thm}\label{thm:MwrSingPTn} 
For $n\geq2$, and with the notation of Theorems \ref{thm:En} and \ref{thm:MwrSingTn}, the wreath product $M\wr\Sing(\PT_n)$ has presentation $\Spres{X_E\cup X_S}{R_E\cup R_S\cup R}$ via
\[
\Phi:(X_E\cup X_S)^+\to M\wr\Sing(\PT_n) : \begin{cases}
e_{ij;ab} \mt \ve_{ij;ab}, \\
t_i \mt \id_{\bn\sm i} ,
\end{cases}
\]
where $R$ consists of the following relations, with $a,b\in M$ arbitrary, and $i,j,k\in\bn$ with $i,j$ distinct, in each:
\begin{align}
\label{eq:MPTX1} e_{ij;ab}t_k &= \begin{cases}
e_{ij;ab} &\text{if $k=j$,}\\
t_it_je_{ij;ab} &\text{if $k=i$,}\\
t_ke_{ij;ab} &\text{otherwise,}
\end{cases}\\
\label{eq:MPTX2} t_je_{ij} &= t_j.
\end{align}
\end{thm}

\pf
This is obtained from Theorems \ref{thm:simp}\ref{simp2}, \ref{thm:En} and \ref{thm:MwrSingTn}, and Lemmas \ref{lem:Wjoins2} and \ref{lem:thB2}.  Explicitly, these results lead to the presentation $\Spres{X_E\cup X_S}{R_E\cup R_S\cup R_1\cup R_2''}$ via~$\Phi$, where 
\begin{align}
\label{eq:MPTX3} R_1 &= \bigset{(e_{ij;ab} t_k , {}^{e_{ij;ab}}t_k\cdot e_{ij;ab})}{i,j,k\in\bn,\ i\not=j,\ a,b\in M} \\
\label{eq:MPTX4} \AND R_2'' &= \bigset{(t_j N(\al),t_jN(\be))}{j\in\bn,\ (\al,\be)\in\Om_j},
\end{align}
where
\bit
\item ${}^{e_{ij;ab}}t_k$ is a (possibly empty) word over $X_E$ mapping to ${}^{\ve_{ij;ab}}\id_{\bn\sm k}$,
\item $N:M\wr\Sing(\T_n)\to X_S^+$ is a normal form function, and additionally $N(\id_n)=\ew$, and
\item $\Om_j$ generates $\Th_{\bn\sm j}=\Th_{\id_{\bn\sm j}}$ as a right congruence.
\eit
(Note that $R_3$ from Theorem \ref{thm:simp}\ref{simp2} is empty because the action pair $(\E_n,\Sing(\T_n))$ is strong.)

By direct computation, we can take
\[
{}^{e_{ij;ab}}t_k = \begin{cases}
\ew &\text{if $k=j$,}\\
t_it_j &\text{if $k=i$,}\\
t_k &\text{otherwise,}
\end{cases}
\]
and then $R_1$ is precisely the set of relations in \eqref{eq:MPTX1}.

By Remark \ref{rem:thB}, we can take $\Om_j=\bigset{(\ve_{ij},\id_X)}{i\in\bn\sm j}$, and then $R_2''$ is precisely the set of relations in \eqref{eq:MPTX2}.
\epf

\begin{rem}
The $k=i$ case of \eqref{eq:MPTX1} says $e_{ij;ab}t_i = t_it_je_{ij;ab}$.  One can show that (in the presence of the other relations) this is equivalent the slightly simpler~${e_{ij;ab}t_i = t_it_j}$.  

Taking $M$ to be the trivial monoid, $M\wr\Sing(\PT_n) \equiv\Sing(\PT_n)$, and the above presentation reduces to the main result of \cite{JEptnsn2}.
\end{rem}

\begin{rem}\label{rem:others}
One can adapt the above method to deal with
\[
M\wr\PT_n = \E_n\cdot(M\wr\T_n) \AND M\wr\I_n = \E_n\cdot(M\wr\G_n).
\]
Here, in addition to the presentation for $\E_n$ from Theorem \ref{thm:En}, we need presentations for ${M\wr\T_n}$ and $M\wr\G_n$.  Since the latter are (monoid) semidirect products $M\wr\T_n=M^n\rtimes\T_n$ and ${M\wr\G_n=M^n\rtimes\G_n}$, one could obtain such presentations from Lavers' Theorem \ref{thm:Lavers}.  However, since we will deal with such wreath products in Section \ref{sect:WpresII}, we do not provide the details here; see Theorems \ref{thm:MwrPTn} and~\ref{thm:MwrGn}--\ref{thm:MwrIn}.
\end{rem}

\sect{Action pairs and subsemigroups II}\label{sect:WAPII}

We now move on to the second family of action pairs in $M\wr\PT_X$.  Recall that $M\wr\PT_X$ contains natural copies of $M_0^X$ and $\PT_X$, via the identifications
\[
\ba \equiv (\ba,\id_{\supp(\ba)}) \AND \al \equiv (\bone_{\dom(\al)},\al) \qquad\text{for $\ba\in M_0^X$ and $\al\in\PT_X$.}
\]

\begin{prop}\label{prop:WAP3}
For any set $X$ and monoid $M$, the following are all action pairs in $M\wr\PT_X$, with respect to the action \eqref{eq:ala}\emph{:}
\ben\bmc3
\item \label{WAP31} $(M_0^X,\PT_X)$, 
\item \label{WAP32} $(M_0^X,\I_X)$, 
\item \label{WAP33} $(M_0^X,\T_X)$, 
\item \label{WAP34} $(M_0^X,\G_X)$, 
\item \label{WAP35} $(M^X,\T_X)$, 
\item \label{WAP36} $(M^X,\G_X)$.
\emc\een
For finite $X$, the following are also action pairs:
\ben\bmc2\addtocounter{enumi}{6}
\item \label{WAP37} $(M_0^X,\Sing(\PT_X))$, 
\item \label{WAP39} $(M_0^X,\Sing(\I_X))$, 
\item \label{WAP38} $(M_0^X,\Sing(\T_X))$, 
\item \label{WAP310} $(M^X,\Sing(\T_X))$.
\emc\een
\end{prop}

\pf
Given Lemma \ref{lem:US0}, it suffices to show that $(M_0^X,\PT_X)$ is an action pair.  
For~\ref{A1}, we have already noted in \eqref{eq:alba} that $\al\cdot\ba={}^\al\ba\cdot\al$ for all $\al\in\PT_X$ and $\ba\in M_0^X$.  For \ref{A2} we need to show that
\[
\ba\cdot\al = \bb\cdot\be \implies \ba\cdot\al^+ = \bb\cdot\be^+ \qquad\text{for all $\ba,\bb\in M_0^X$ and $\al,\be\in\PT_X$,}
\]
where $\al^+={}^\al\bone = \bone_{\dom(\al)} \equiv \id_{\dom(\al)}$.  By \eqref{eq:aa}, we have
\begin{align*}
\ba\cdot\al &= (\ba\restr_{\dom(\al)},\al\restr_{\supp(\ba)}) , & \ba\cdot\al^+ &= (\ba\restr_{\dom(\al^+)},\al^+\restr_{\supp(\ba)}) , \\
\bb\cdot\be &= (\bb\restr_{\dom(\be)},\be\restr_{\supp(\bb)}) , & \bb\cdot\be^+ &= (\bb\restr_{\dom(\be^+)},\be^+\restr_{\supp(\bb)}).
\end{align*}
So the required implication quickly follows from:
\bit
\item the identities $\dom(\al^+)=\dom(\al)$ and $\dom(\al^+\restr_{\supp(\ba)})=\dom(\al\restr_{\supp(\ba)})$,
\item the analogous identities for $\bb$ and $\be$, and
\item the fact that $\id_C=\id_D \iff C=D$ for $C,D\sub X$, noting that $\al^+\restr_{\supp(\ba)}$ and $\be^+\restr_{\supp(\bb)}$ are both partial identities.  \qedhere
\eit
\epf

\begin{rem}\label{rem:WAP3}
Pairs \ref{WAP35}, \ref{WAP36} and \ref{WAP310} from Proposition \ref{prop:WAP3} have the form $(M^X,S)$.  Any such pair~$(M^X,S)$ necessarily satisfies $S\sub\T_X$, since $M^X$ is not closed under the action of elements from~${\PT_X\sm\T_X}$.  Any such pair also has the \emph{right-}uniqueness property.  Indeed, it quickly follows from~\eqref{eq:aa} that
\[
\ba\cdot\al=\bb\cdot\be\implies \al=\be \qquad\text{for all $\ba,\bb\in M^X$ and $\al,\be\in\PT_X$.}
\]
But this is to be expected, since for $S\leq\T_X$ we have $M\wr S = M^X\rtimes S$; cf.~Remark \ref{rem:wreath}.
\end{rem}

We now identify the strong and proper pairs from Proposition \ref{prop:WAP3}.

\begin{prop}\label{prop:WAP5}
Of the action pairs listed in Proposition \ref{prop:WAP3}, and apart from trivially small exceptions:
\ben
\item \label{WAP51} the strong pairs are precisely \ref{WAP33}--\ref{WAP36}, \ref{WAP38} and \ref{WAP310},
\item \label{WAP52} the proper pairs are precisely \ref{WAP35}, \ref{WAP36} and \ref{WAP310}.  
\een
\end{prop}

\pf
\firstpfitem{\ref{WAP51}}  The stated pairs are precisely those of the form $(U,S)$ with $S\sub\T_X$, so this part follows immediately from Lemma \ref{lem:SAPAP}\ref{SAPAP2} and \eqref{eq:albone}.

\pfitem{\ref{WAP52}}  The stated pairs are precisely those of the form $(M^X,S)$ with $S\sub\T_X$.  For each of these, we have $P(U,S)=\{\id_X\equiv\bone\}$, so that $\si=\De_S$ is the trivial relation.  The required equivalence~\eqref{eq:proper} then becomes 
\[
\ba\cdot\al = \bb\cdot\be  \qquad\Leftrightarrow\qquad \ba = \bb  \ANd \al=\be \qquad\qquad\text{for all $\ba,\bb\in M^X$ and $\al,\be\in S$,}
\]
and this again follows immediately from \eqref{eq:aa}.

The pairs not listed are all of the form $(M_0^X,S)$ with $S\leq\PT_X$.  To show that these are not proper we consider separate cases according to whether or not $S\sub\T_X$.  We must show that the following fails for some $\ba,\bb\in M_0^X$ and $\al,\be\in S$:
\begin{equation}\label{eq:fail}
\ba\cdot\al = \bb\cdot\be  \qquad\Leftrightarrow\qquad \ba\cdot\al^+ = \bb\cdot\be^+  \ANd \al\mr\si\be.
\end{equation}

\pfcase1  If $S\sub\T_X$, then as above we have $\si=\De_S$, and $\al^+=\bone$ for all $\al\in S$.  In this case, the forwards implication in~\eqref{eq:fail} fails for $\ba=\bb=\bzero$ and distinct $\al,\be\in S$, where $\bzero=(0)_{x\in X}\in M_0^X$ is the all-zero $X$-tuple.

\pfcase2  Now suppose $S\not\sub\T_X$.  By the form of the pairs in question, we have $\es=\id_\es\in S$, and so $\bzero = \es^+\in P$, so that $\si=\nabla_S$ is the universal relation in this case.  Here the backwards implication in \eqref{eq:fail} fails for $\ba=\bb=\bone$ and any $\al,\be\in S$ with $\dom(\al)=\dom(\be)$ but $\al\not=\be$.
\epf

\begin{rem}
As suggested by Case 1 in the above proof, one can make a general statement about proper action pairs in a certain special case concerning zero elements.  Indeed, suppose $(U,S)$ is a proper action pair in a monoid $M$ with zero $0$, and suppose additionally that $0\in U$.  Since $0\cdot s=0\cdot t$ (and $0\cdot s^+=0\cdot t^+$) for all $s,t\in S$, it follows from \eqref{eq:proper} that $s\mr\si t$ for all $s,t\in S$, meaning that $\si=\nabla_S$ for such pairs.  

On the other hand, if $(U,S)$ is a weak action pair with $\si=\nabla_S$, then it follows from Definition~\ref{defn:proper} that $(U,S)$ is proper precisely when
\[
us=vt \iff us^+=vt^+ \qquad\text{for all $u,v\in U^1$ and $s,t\in S$.}
\]
This is of course an if-and-only-if version of the defining property \ref{A2} from Definition \ref{defn:AP}.
\end{rem}

It is a routine matter to determine the subsemigroups $US$ of $M\wr\PT_X$ corresponding to the action pairs $(U,S)$ in Proposition~\ref{prop:WAP3}.  In the next statement, we again do not assume $X$ is finite in the parts involving singular subsets/subsemigroups, and we again refer to the subsemigroup $\F_X\leq\I_X$ from Remark \ref{rem:FA}; recall that $\F_X=\I_X$ when $X$ is finite.

\begin{prop}\label{prop:WAP4}
For any set $X$ and monoid $M$, we have:
\ben\bmc2
\item \label{WAP41} $M_0^X \cdot \PT_X = M\wr\PT_X$,
\item \label{WAP42} $M_0^X \cdot \I_X = M\wr\I_X$,
\item \label{WAP43} $M_0^X \cdot \T_X = M\wr\PT_X$,
\item \label{WAP44} $M_0^X \cdot \G_X = M\wr\F_X$,
\item \label{WAP45} $M^X \cdot \T_X = M\wr\T_X$,
\item \label{WAP46} $M^X \cdot \G_X = M\wr\G_X$,
\item \label{WAP47} $M_0^X \cdot \Sing(\PT_X) = M\wr\Sing(\PT_X)$,
\item \label{WAP48} $M_0^X \cdot \Sing(\I_X) = M\wr\Sing(\I_X)$,
\item \label{WAP49} $M_0^X \cdot \Sing(\T_X) = M\wr\Sing(\PT_X)$,
\item \label{WAP410} $M^X \cdot \Sing(\T_X) = M\wr\Sing(\T_X)$.  \epfres
\emc\een
\end{prop}

\sect{Congruence conditions II}\label{sect:WcongII}

We now examine the congruence conditions for the pairs from Proposition \ref{prop:WAP3}.  Each such pair is of the form $(U,S)$, where $S\leq\PT_X$, and $U$ is either $M_0^X$ or $M^X$, the latter only when $S\sub\T_X$.  Since~$U$ is a monoid in each case, all such pairs satisfy all the conditions of Lemma \ref{lem:w}, and hence Lemma \ref{lem:Om3} applies to each pair, which then feeds into Proposition \ref{prop:ES_sd}.  

The stronger congruence conditions involve the behaviour of the right congruences and sub(semi)groups
\begin{equation}
\label{eq:tha} \th_\ba = \bigset{(\al,\be)\in S\times S}{\ba\cdot\al=\ba\cdot\be} 
\ANd \S_\ba = \set{\al\in S}{\ba\cdot\al=\ba} \qquad \text{for $\ba\in M_0^X$,}
\end{equation}
or the related right congruences $\Th_\ba$ on $S^1$, as appropriate.  Such conditions then feed into one of Theorems~\ref{thm:ESmon0}, \ref{thm:ESmon}, \ref{thm:Msimp}, \ref{thm:ES} or \ref{thm:simp}, as applicable to the pair $(U,S)$.  Theorem~\ref{thm:ESmon} applies when~$U$ and~$S$ are both submonoids, which obviously excludes the pairs with $S\sub\Sing(\PT_X)$.  Theorem~\ref{thm:simp} (and indeed even Theorem \ref{thm:ES}) is not generally applicable to \emph{any} of the pairs, as $M_0^X\sm\{\bone\}$ and $M^X\sm\{\bone\}$ are typically not subsemigroups.  Thus, for the remainder of the current section, we will only be concerned with the first six pairs from Proposition~\ref{prop:WAP3}.  
By Proposition~\ref{prop:WAP4}, the four excluded pairs produce the following three product semigroups:
\[
M\wr\Sing(\T_X) \COMMA M\wr\Sing(\PT_X) \AND M\wr\Sing(\I_X).
\]
Presentations for the first two (for finite $X$) are stated above in Theorems \ref{thm:MwrSingTn} and \ref{thm:MwrSingPTn}.  For the third (again for finite $X$), see \cite{CE2022}.

We have so far restricted our attention to the first six pairs $(U,S)$ from Proposition \ref{prop:WAP3}.  In fact, since we are primarily concerned with finding presentations for the corresponding subsemigroup~$US$ (for finite $X$), and since the same subsemigroup can arise from different pairs, as seen in Proposition~\ref{prop:WAP4} (and recalling that $\F_X=\I_X$ when $X$ is finite), we can further restrict our attention to the following four pairs:
\bena\bmc4
\item \label{Wa} $(M_0^X,\T_X)$, 
\item \label{Wb} $(M_0^X,\G_X)$, 
\item \label{Wc} $(M^X,\T_X)$, 
\item \label{Wd} $(M^X,\G_X)$.
\emc\een
As in Proposition \ref{prop:WAP4}, these give rise to the following subsemigroups (respectively), where we additionally assume $X$ is finite in \ref{Wb}:
\bena\bmc4
\item $M\wr\PT_X$, 
\item $M\wr\I_X$, 
\item $M\wr\T_X$, 
\item $M\wr\G_X$.
\emc\een
Since the pairs \ref{Wc} and \ref{Wd} have the right-uniqueness property (cf.~Remark \ref{rem:WAP3}), the congruence~$\th$ is trivial, as are the right congruences $\th_\ba$ and sub(semi)groups $\S_\ba$.  So in fact, we need only focus on pairs \ref{Wa} and \ref{Wb} in the current section, though some of the results we prove are more general.

\begin{lemma}\label{lem:tha}
Let $\ba\in M_0^X$, and write $B=\supp(\ba)$.  Then for any subsemigroup $S\leq\PT_X$ we have
\[
\th_\ba = \th_{\bone_B} = \bigset{(\al,\be)\in S\times S}{\al\restr_B = \be\restr_B}.
\]
\end{lemma}

\pf
Since $\al\restr_B=\bone_B\cdot\al$ for all $\al\in\PT_X$, we have $\th_{\bone_B} = \bigset{(\al,\be)\in S\times S}{\al\restr_B = \be\restr_B}$ by definition.

Next, let $\al,\be\in S$.  Then since $\ba=\ba\cdot\bone_B$, we have
\[
(\al,\be)\in\th_{\bone_B} \implies \bone_B\cdot\al = \bone_B\cdot\be \implies \ba\cdot\al = \ba\cdot \bone_B\cdot\al = \ba\cdot\bone_B\cdot\be = \ba\cdot\be \implies (\al,\be)\in\th_\ba.
\]
The converse follows quickly from \eqref{eq:aa} and the definition of $\th_\ba$.
\epf

\begin{lemma}\label{lem:thab}
If $S=\T_X$ or $\PT_X$, then
\[
\th_{\ba\bb} = \th_\ba\vee\th_\bb \qquad\text{for all $\ba,\bb\in M_0^X$.}
\]
\end{lemma}

\pf
Throughout the proof we write $B=\supp(\ba)$ and $C=\supp(\bb)$, so $\supp(\ba\bb)=B\cap C$.  By Lemma \ref{lem:tha}, we need to show that
\[
\th_{\bone_{B\cap C}} = \th_{\bone_B}\vee\th_{\bone_C}.
\]
Since $\bone_B$ and $\bone_C$ commute (and $\bone_B\bone_C=\bone_{B\cap C}$), it follows from \eqref{eq:thv1k} that we only need to show the forwards inclusion.
To do so, suppose $(\al,\be)\in\th_{\bone_{B\cap C}}$, so that $\al\restr_{B\cap C}=\be\restr_{B\cap C}$.  We then let $\ga\in\PT_X$ be such that $\ga\restr_B=\al\restr_B$ and $\ga\restr_C=\be\restr_C$; in the case $S=\T_X$ we also define $x\ga=x$ for all $x\in X\sm(B\cup C)$.  We then have $\ga\in S$, and $\al\mr\th_B\ga\mr\th_C\be$.
\epf

\begin{rem}
The previous result excluded the cases $S=\I_X$ and $S=\G_X$ since injectivity of~$\al$ and $\be$ does not imply that the $\ga$ constructed in the proof is injective.
\end{rem}

Lemma \ref{lem:thab} concerns the congruence condition from Lemma \ref{lem:Om4}\ref{Om42}.  However, that result can only be applied if $U$ is commutative, and since $U=M_0^X$ or $M^X$, this is only the case if~$M$ is commutative, and we do not need or want to assume this is the case.  It turns out that while part~\ref{Om42} of Lemma~\ref{lem:Om4} does not apply in general, part \ref{Om41} always does.  While this might seem less desirable, it is in fact an improvement in this case, since we can choose the relevant subset $V\sub U$ to be rather a lot smaller than a generating set for $U$, at least when $X$ is finite, as we now explain.

For $x\in X$, we write $\bv_x = \bone_{X\sm x} \equiv \id_{X\sm x}$, and we write 
\[
V = \set{\bv_x}{x\in X}.
\]
When $X$ is finite we have
\[
\la V\ra = \set{\bone_B}{B\sub X} \equiv \set{\id_B}{B\sub X} = \E_X.
\]
Again writing $\pre$ for the relation in \eqref{eq:pre}, note that
\[
\bv_x\pre\ba \iff x\not\in\supp(\ba) \qquad\text{for all $\ba\in M_0^X$ and $x\in X$.}
\]
The following therefore verifies the conditions of Lemma \ref{lem:Om4}\ref{Om41} with respect to this $V$.

\begin{lemma}\label{lem:thav}
If $S=\PT_X$ or $\T_X$ for $|X|<\aleph_0$, and if $\ba\in M_0^X$, then
\[
\th_\ba = \bigvee_{x\in X\sm B} \th_{\bv_x} \WHERE B=\supp(\ba).
\]
\end{lemma}

\pf
By Lemma \ref{lem:tha} we have $\th_\ba = \th_{\bone_B}$.  Since $\bone_B = \prod_{x\in X\sm B}\bv_x$, it follows from Lemma \ref{lem:thab} that $\th_{\bone_B}=\bigvee_{x\in X\sm B}\th_{\bv_x}$.
\epf

The proof of the next result is similar to (but easier than) that of Lemma \ref{lem:thB1}.  For this result we do not need to assume that $X$ is finite.

\begin{lemma}\label{lem:thvy}
If $S=\T_X$, and if $y\in X$, then $\th_{\bv_y}$ is generated as a right congruence by any pair $(\ve_{xy},\id_X)$ with $x\in X\sm y$.  \epfres
\end{lemma}

We now have all the information we need for the pair $(M_0^X,\T_X)$, so we now turn our attention to $(M_0^X,\G_X)$ for finite $X$.  For $\ba\in M_0^X$ and $\al\in\G_X$, we have $\al\in\S_\ba \iff(\al,\id_X)\in\th_\ba$, so by Lemma \ref{lem:tha} we have 
\[
\S_\ba = \set{\al\in\G_X}{\al\restr_B=\id_B}  \WHERE B=\supp(\ba).
\]
The next result immediately follows, where again $\tau_{xy}$ denotes the transposition $(x,y)$.

\begin{lemma}\label{lem:Sagens}
If $S=\G_X$ for $|X|<\aleph_0$, and if $\ba\in M_0^X$, then $\S_\ba$ is generated by the set $\set{\tau_{xy}}{x,y\in X\sm B,\ x\not=y}$.  \epfres
\end{lemma}

As in the proof of Lemma \ref{lem:Wjoins3}, it is then easy to show that
\begin{equation}\label{eq:thab}
\S_{\ba\bb} = \S_\ba\vee\S_\bb \qquad\text{for any $\ba,\bb\in M_0^X$ with $(X\sm\supp(\ba))\cap(X\sm\supp(\bb))\not=\es$.}
\end{equation}
For distinct $x,y\in X$, we write $\bv_{xy} = \bv_x\bv_y = \bone_{X\sm\{x,y\}} \equiv \id_{X\sm\{x,y\}}$.  The next result shows that Lemma \ref{lem:Om6}\ref{Om61} applies to the pair $(M_0^X,\G_X)$, with respect to the set
\[
V = \set{\bv_{xy}}{x,y\in X,\ x\not=y}.
\]

\begin{lemma}\label{lem:thav2}
If $S=\G_X$ for $|X|<\aleph_0$, and if $\ba\in M_0^X$, then
\[
\S_\ba = \bigvee_{x,y\in X\sm B\atop x\not=y} \S_{\bv_{xy}} \WHERE B=\supp(\ba).  
\]
\end{lemma}

\pf
If $|X\sm B|\leq1$, then $\S_\ba=\{\id_X\}$, and the stated join is empty.  If $|X\sm B|\geq2$, then we argue as in the proof of Lemma \ref{lem:Wjoins4}, applying \eqref{eq:thab} in place of Lemma \ref{lem:Wjoins3}.
\epf

\sect{Presentations II}\label{sect:WpresII}

We now assume that $X=\bn=\{1,\ldots,n\}$ for some integer $n\geq2$, with the goal of giving presentations for the wreath products 
\[
M\wr\PT_n = M_0^n\cdot\T_n \COMMA
M\wr\I_n = M_0^n\cdot\G_n \COMMA
M\wr\T_n = M^n\cdot\T_n \AND
M\wr\G_n = M^n\cdot\G_n.
\]
The wreath products $M\wr\Sing(\T_n)$ and $M\wr\Sing(\I_n)$ were treated in \cite{FADEG2019,CE2022}, and $M\wr\Sing(\PT_n)$ in Theorem \ref{thm:MwrSingPTn}.
For the remaining wreath products listed above, we need presentations for $M^n$, $M_0^n$, $\G_n$ and $\T_n$.  For~$\G_n$ we have the following classical result, in which we continue to use the notation $\tau_{xy}$ for the transposition $(x,y)\in\G_n$:

\newpage

\begin{thm}[Moore \cite{Moore1897}]\label{thm:Gn}
For $n\geq2$, the symmetric group $\G_n$ has presentation ${\Mpres{X_G}{R_G}}$ via $\phi_G:X_G^*\to\G_n$, where
\bit
\item $X_G=\{s_1,\ldots,s_{n-1}\}$,
\item $s_i\phi_G = \tau_{i,i+1}$ for all $i$, and
\item $R_G$ consists of the following relations:
\[
s_i^2 = \ew \text{ for all $i$}\COMMA
s_is_j = s_js_i \text{ if $|i-j|>1$}\COMMA
s_is_js_i = s_js_is_j \text{ if $|i-j|=1$.}  \epfreseq
\]
\eit
\end{thm}

The first presentation for $\T_n$ was given by A{\u\i}zen{\v{s}}tat \cite{Aizenstat1958}, but it is convenient here to use the following formulation from \cite{JEptc}; see also \cite{Lavers1997}.  The statement (and following ones) uses the maps~$\ve_{ij}$ from \eqref{eq:vexy}.

\begin{thm}[A{\u\i}zen{\v{s}}tat \cite{Aizenstat1958}]\label{thm:Tn}
For $n\geq2$, the full transformation semigroup $\T_n$ has presentation ${\Mpres{X_T}{R_T}}$ via $\phi_T:X_T^*\to\T_n$, where
\bit
\item $X_T=\{s_1,\ldots,s_{n-1}\}\cup\{\lam_1,\ldots,\lam_{n-1}\}\cup\{\rho_1,\ldots,\rho_{n-1}\}$,
\item $s_i\phi_T = \tau_{i,i+1}$, $\lam_i\phi_T = \ve_{i,i+1}$ and $\rho_i\phi_T = \ve_{i+1,i}$ for all~$i$, 
\item $R_T$ consists of the following relations:
\begin{align*}
& s_{i}^2 = \ew_n \COMMA \lam_{i} = \lam_{i}^2 = \rho_{i}\lam_{i} = s_{i}\lam_{i} = \rho_{i}s_{i} \COMMA \rho_{i} = \rho_{i}^2 = \lam_{i}\rho_{i} = s_{i}\rho_{i}  = \lam_{i}s_{i}, \\
& \lam_{i}\lam_{i+1}=\lam_{i}s_{i+1} \COMMA \rho_{i+1}\rho_{i}=\rho_{i+1}s_{i} \COMMA \lam_{i}\rho_{i+1} = \lam_{i} \COMMA \rho_{i+1}\lam_{i} = \rho_{i+1}, \\
& \lam_{i+1}\lam_{i} = \lam_{i}\lam_{i+1}\lam_{i} = \lam_{i+1}\lam_{i}\lam_{i+1} \COMMA \rho_{i}\rho_{i+1} = \rho_{i}\rho_{i+1}\rho_{i} = \rho_{i+1}\rho_{i}\rho_{i+1}, \\
& \lam_{i+1}s_{i}=s_{i}s_{i+1}\lam_{i}\lam_{i+1} \COMMA \rho_{i}s_{i+1}=s_{i+1}s_{i}\rho_{i+1}\rho_{i}, \\
& s_{i}s_{j} = s_{j}s_{i} \COMMA \lam_{i}\lam_{j} = \lam_{j}\lam_{i} \COMMA \rho_{i}\rho_{j} = \rho_{j}\rho_{i}, &&\hspace{-3.5cm}\text{if $|i-j|>1$,} \\
& s_{i}\lam_{j}=\lam_{j}s_{i} \COMMA s_{i}\rho_{j}=\rho_{j}s_{i}, &&\hspace{-3.5cm}\text{if $|i-j|>1$,} \\
& s_{i}s_{j}s_{i} = s_{j}s_{i}s_{j} ,&&\hspace{-3.5cm}\text{if $|i-j|=1$,} \\
& \lam_{i}\rho_{j} = \rho_{j}\lam_{i},  &&\hspace{-3.5cm}\text{if $j\not=i,i+1$.}   \epfreseq
\end{align*}
\eit
\end{thm}

The generators for $\G_n$ and $\T_n$ from Theorems \ref{thm:Gn} and \ref{thm:Tn} are pictured in Figure \ref{fig:gensGnTn}.

\begin{figure}[h]
\begin{center}
\scalebox{.9}{
\begin{tikzpicture}[scale=.7]
\tikzstyle{vertex}=[circle,draw=black, fill=white, inner sep = 0.06cm]
\begin{scope}[shift={(0,0)}]
\uvs{1,3,4,5,6,8}
\lvs{1,3,4,5,6,8}
\udotted{1.5}{2.5}
\udotted{6.5}{7.5}
\ddotted{1.5}{2.5}
\ddotted{6.5}{7.5}
\foreach \x in {1,3,6,8} {\stline\x\x}
\stline45
\stline54
\foreach \x/\y in {1/1,4/i,8/n} {\node () at (\x,2.4){\footnotesize $\y$};}
\end{scope}
\begin{scope}[shift={(-7,4)}]
\uvs{1,3,4,5,6,8}
\lvs{1,3,4,5,6,8}
\udotted{1.5}{2.5}
\udotted{6.5}{7.5}
\ddotted{1.5}{2.5}
\ddotted{6.5}{7.5}
\foreach \x in {1,3,6,8} {\stline\x\x}
\stline44
\stline54
\foreach \x/\y in {1/1,4/i,8/n} {\node () at (\x,2.4){\footnotesize $\y$};}
\end{scope}
\begin{scope}[shift={(7,4)}]
\uvs{1,3,4,5,6,8}
\lvs{1,3,4,5,6,8}
\udotted{1.5}{2.5}
\udotted{6.5}{7.5}
\ddotted{1.5}{2.5}
\ddotted{6.5}{7.5}
\foreach \x in {1,3,6,8} {\stline\x\x}
\stline45
\stline55
\foreach \x/\y in {1/1,4/i,8/n} {\node () at (\x,2.4){\footnotesize $\y$};}
\end{scope}
\end{tikzpicture}
}
\caption{The generators $\lam_i\phi_T$ (left), $s_i\phi_G=s_i\phi_T$ (middle) and $\rho_i\phi_T$ (right), for $1\leq i\leq n-1$.}
\label{fig:gensGnTn}
\end{center}
\end{figure}
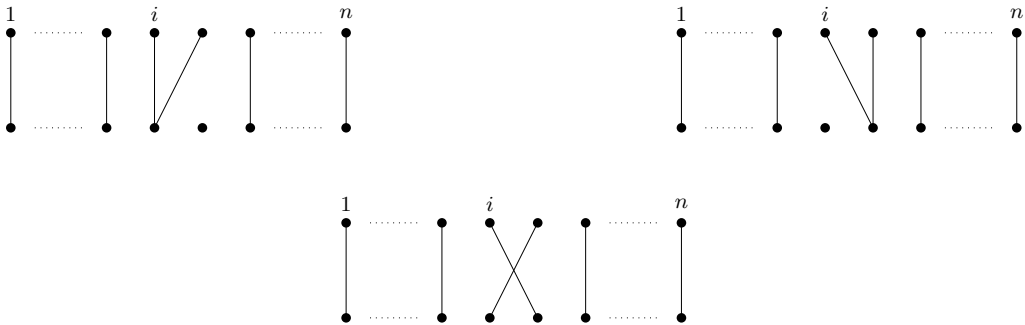

Presentations for $M^n$ and $M_0^n$ are easy to obtain (either directly, or by iterating Theorem~\ref{thm:Lavers} in the special case of direct products), so we fix notation for these as follows.  For $a\in M$ and~$i\in\bn$, we denote by $a^{(i)}\in M^n$ the $n$-tuple over $M$ whose $i$th entry is $a$, and all other entries are~$1$.  (This tuple was denoted $\bt_{i;a}$ in Section \ref{sect:WcongI}, but the $a^{(i)}$ notation is more convenient for our present purposes.)  If $\ba=(a_1,\ldots,a_n)\in M^n$, then $\ba=a_1^{(1)}\cdots a_n^{(n)}$.

We assume the monoid $M$ has presentation
\[
\Mpres{X_M}{R_M} \qquad\text{via}\qquad \phi_M:X_M^*\to M.
\]
For $w\in X_M^*$, we write $\ol w=w\phi_M$.  
For each $i\in\bn$, we define an alphabet $X_M^{(i)} = \set{x^{(i)}}{x\in X_M}$ in one-one correspondence with~$X_M$.  For a word $w=x_1\cdots x_k\in X_M^*$, we write $w^{(i)}=x_1^{(i)}\cdots x_k^{(i)}$, and we also define ${R_M^{(i)} = \bigset{(w^{(i)},v^{(i)})}{(w,v)\in R_M}}$.  We then define the alphabet
\[
X_{M^n}=X_M^{(1)}\cup\cdots\cup X_M^{(n)},
\]
and the set of relations
\[
R_{M^n}=R_M^{(1)}\cup\cdots\cup R_M^{(n)} \cup \bigset{(x^{(i)}y^{(j)},y^{(j)}x^{(i)}}{x,y\in X_M,\ i,j\in\bn,\ i\not=j}.
\]

\begin{thm}\label{thm:Mn}
With the above notation, $M^n$ has presentation $\Mpres{X_{M^n}}{R_{M^n}}$ via
\[
\phi_{M^n}:X_{M^n}^*\to M^n: x^{(i)} \mt \ol x^{(i)}.  \epfreseq
\]
\end{thm}

A presentation for $M_0^n$ can be derived as a special case of Theorem \ref{thm:Mn}, beginning instead with a presentation for $M_0$.  However, it is convenient to give a more explicit formulation, making use of the fact that $M$ is a subsemigroup of $M_0$ by definition.  For this, we additionally use the notation of Theorem \ref{thm:En}, in particular the alphabet $X_E$ and the relations $R_E$.  We define the alphabet
\[
X_{M_0^n} = X_{M^n} \cup X_E,
\]
and set of relations
\begin{align*}
R_{M_0^n} = R_{M^n} \cup R_E &\cup \bigset{(t_ix^{(j)},x^{(j)}t_i)}{x\in X_M,\ i,j\in\bn,\ i\not=j} \\
{} & \cup \bigset{(t_ix^{(i)},t_i),(x^{(i)}t_i,t_i)}{x\in X_M,\ i\in\bn}.
\end{align*}

\begin{thm}\label{thm:M0n}
With the above notation, $M_0^n$ has presentation $\Mpres{X_{M_0^n}}{R_{M_0^n}}$ via
\[
\phi_{M_0^n}:X_{M_0^n}^*\to M_0^n: x^{(i)} \mt \ol x^{(i)},\ t_i \mt \bone_{\bn\sm i}.  \epfreseq
\]
\end{thm}

The generators for $M^n$ and $M_0^n$ from Theorems \ref{thm:Mn} and \ref{thm:M0n} are pictured in Figure \ref{fig:xiti}, each identified with a pair from $M\wr\PT_n$ in the usual way.

\begin{figure}[h]
\begin{center}
\scalebox{.9}{
\begin{tikzpicture}[scale=1]
\tikzstyle{vertex}=[circle,draw=black, fill=white, inner sep = 0.06cm]
\uvs{1,3,4,5,7}
\lvs{1,3,4,5,7}
\foreach \x in {1,3,4,5,7} {\stline\x\x}
\foreach \x/\y in {4/\ol x} {\node[vertex] () at (\x,2){$\y$};}
\udotted{1.5}{2.5}
\udotted{5.5}{6.5}
\ddotted{1.5}{2.5}
\ddotted{5.5}{6.5}
\foreach \x/\y in {1/1,4/i,7/n} {\node () at (\x,2.5){\footnotesize $\y$};}
\begin{scope}[shift={(10,0)}]
\uvs{1,3,4,5,7}
\lvs{1,3,4,5,7}
\foreach \x in {1,3,5,7} {\stline\x\x}
\udotted{1.5}{2.5}
\udotted{5.5}{6.5}
\ddotted{1.5}{2.5}
\ddotted{5.5}{6.5}
\foreach \x/\y in {1/1,4/i,7/n} {\node () at (\x,2.5){\footnotesize $\y$};}
\end{scope}
\end{tikzpicture}
}
\caption{The generators $\ol x^{(i)}$ (left) and $\bone_{\bn\sm i}\equiv\id_{\bn\sm i}$ (right).}
\label{fig:xiti}
\end{center}
\end{figure}

We can now give the main results of this section, which are presentations for the wreath products $M\wr S$, where $S$ is one of $\G_n$, $\T_n$, $\I_n$ or $\PT_n$.  It is convenient to begin with the largest of these, $M\wr\PT_n$.

\newpage

\begin{thm}\label{thm:MwrPTn} 
For $n\geq2$, the wreath product $M\wr\PT_n$ has presentation
\[
\Mpres{X_{M_0^n}\cup X_T}{R_{M_0^n}\cup R_T\cup R}
\]
via 
\[
\Psi:(X_{M_0^n}\cup X_T)^* \to M\wr\PT_n: \begin{cases}
x^{(i)} \mt \ol x^{(i)} ,\\ 
t_i \mt \id_{\bn\sm i} ,\\ 
s_i \mt \tau_{i,i+1} ,\\ 
\lam_i \mt \ve_{i,i+1} ,\\ 
\rho_i \mt \ve_{i+1,i} , 
\end{cases}
\]
where $R$ consists of the following relations, with $x\in X_M$, $1\leq i\leq n-1$ and $1\leq k\leq n$ in each:
\begin{align}
\label{eq:sx}
s_i x^{(k)} &= \begin{cases}
\mathrlap{x^{(i+1)}s_i}\hphantom{XXXXXX} &\text{if $k=i$,}\\
x^{(i)}s_i &\text{if $k=i+1$,}\\
x^{(k)}s_i &\text{otherwise,}
\end{cases}
\\
\label{eq:lx} \lam_i x^{(k)} &= \begin{cases}
\mathrlap{x^{(i)}x^{(i+1)}\lam_i}\hphantom{XXXXXX} &\text{if $k=i$,}\\
\lam_i &\text{if $k=i+1$,}\\
x^{(k)}\lam_i &\text{otherwise,}
\end{cases}
\\
\label{eq:rx} \rho_i x^{(k)} &= \begin{cases}
\mathrlap{\rho_i}\hphantom{XXXXXX} &\text{if $k=i$,}\\
x^{(i)}x^{(i+1)}\rho_i &\text{if $k=i+1$,}\\
x^{(k)}\rho_i &\text{otherwise,}
\end{cases}
\\
\label{eq:st}
s_i t_k &= \begin{cases}
\mathrlap{t_{i+1}s_i}\hphantom{XXXXXX} &\text{if $k=i$,}\\
t_is_i &\text{if $k=i+1$,}\\
t_ks_i &\text{otherwise,}
\end{cases}\\
\label{eq:lt} \lam_i t_k &= \begin{cases}
\mathrlap{t_it_{i+1}\lam_i}\hphantom{XXXXXX} &\text{if $k=i$,}\\
\lam_i &\text{if $k=i+1$,}\\
t_k\lam_i &\text{otherwise,}
\end{cases}
\\
\label{eq:rt} \rho_i t_k &= \begin{cases}
\mathrlap{\rho_i}\hphantom{XXXXXX} &\text{if $k=i$,}\\
t_it_{i+1}\rho_i &\text{if $k=i+1$,}\\
t_k\rho_i &\text{otherwise,}
\end{cases}
\\
\label{eq:tr} t_i\rho_i &= t_i.
\end{align}
\end{thm}

\pf
By Proposition \ref{prop:WAP3} $(M_0^n,\T_n)$ is a strong action pair, and by Proposition \ref{prop:WAP4} we have $M\wr\PT_n = M_0^n\cdot\T_n$.  Theorem \ref{thm:ESmon0} therefore applies, and tells us that $M\wr\PT_n$ has presentation 
\[
\Mpres{X_{M_0^n}\cup X_T}{R_{M_0^n}\cup R_T\cup R_1\cup R_\Om}
\]
via $\Psi$, where:
\bit
\item $R_1 = \bigset{(xy,{}^xy\cdot x)}{x\in X_T,\ y\in X_{M_0^n}}$, and
\item $R_\Om = \bigset{(N(\ba),N(\bb))}{(\ba,\bb)\in\Om}$ for any generating set $\Om$ for the congruence $\th$ on $M_0^n\rtimes\T_n$.
\eit
In $R_1$, ${}^xy$ is a word over $X_{M_0^n}$ mapping to ${}^{\ol x}\ol y$.  So we must find a suitable set of such words ${}^xy$, for each combination of $x=s_i$, $\lam_i$ or $\rho_i$, and $y=x^{(k)}$ or $t_k$.  This has been carried out, and the resulting relations are precisely \eqref{eq:sx}--\eqref{eq:rt}, as can be easily checked.  We give a sample calculation in Figure \ref{fig:lx} for the $k=i,i+1$ cases of \eqref{eq:lx}.

For the relations $R_\Om$, we begin by defining the tuples $\bv_i = \bone_{\bn\sm i} \equiv \id_{\bn\sm i} = t_i\Psi$ for each $i\in\bn$.
By Lemmas \ref{lem:thav} and \ref{lem:Om4}\ref{Om41}, we have $\th=\Om^\sharp$ for
\[
\Om = \bigset{((\bv_i,\al),(\bv_i,\be))}{i\in\bn,\ (\al,\be)\in\Om_{\bv_i}},
\]
where each $\Om_{\bv_i}$ generates $\th_{\bv_i}$ as a right congruence.  By Lemma \ref{lem:thvy} we can take $\Om_{\bv_i}$ to be the set consisting of the single pair $(\id_\bn,\ve_{i+1,i})$.  That is, we can take
\[
\Om = \bigset{((\bv_i,\ve_{i+1,i}),(\bv_i,\id_\bn))}{i\in\bn}.
\]
The resulting set $R_\Om$ therefore consists precisely of the relations \eqref{eq:tr}.
\epf

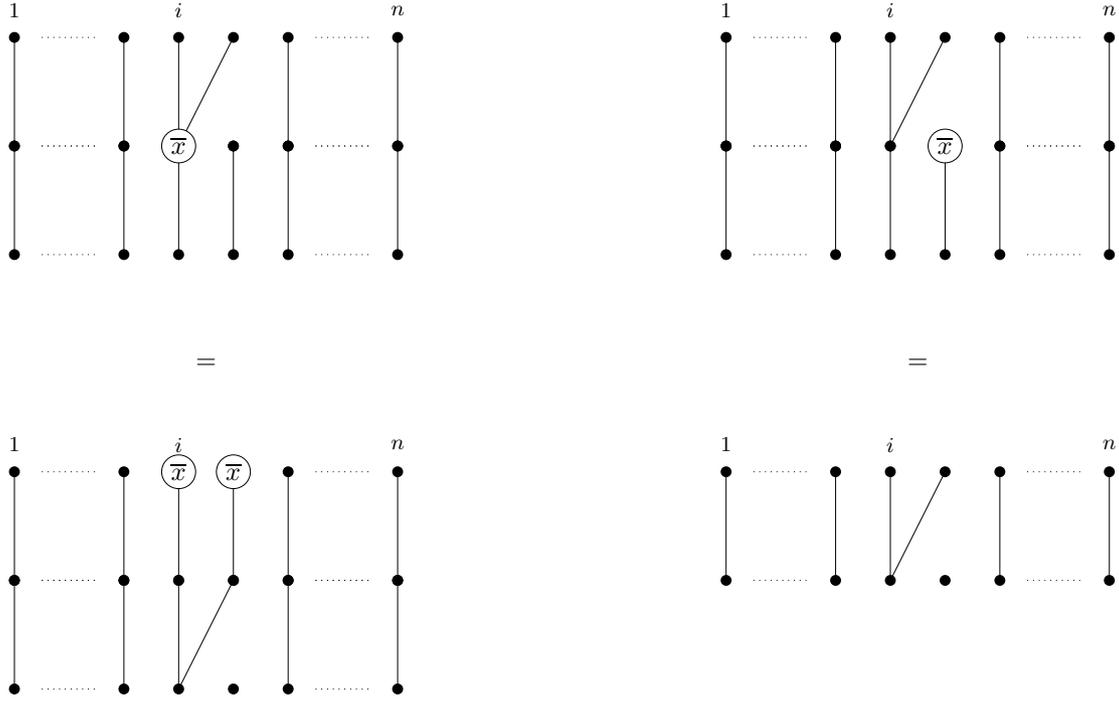
\begin{figure}[h]
\begin{center}
\scalebox{.9}{
\begin{tikzpicture}[scale=.8]
\tikzstyle{vertex}=[circle,draw=black, fill=white, inner sep = 0.06cm]
\begin{scope}
\uvs{1,3,4,5,6,8}
\lvs{1,3,4,5,6,8}
\foreach \x in {1,3,4,6,8} {\stline\x\x}
\stline54
\udotted{1.5}{2.5}
\udotted{6.5}{7.5}
\ddotted{1.5}{2.5}
\ddotted{6.5}{7.5}
\foreach \x/\y in {1/1,4/i,8/n} {\node () at (\x,2.5){\footnotesize $\y$};}
\end{scope}
\begin{scope}[shift={(0,-2)}]
\uvs{1,3,4,5,6,8}
\lvs{1,3,4,5,6,8}
\foreach \x in {1,3,4,5,6,8} {\stline\x\x}
\foreach \x/\y in {4/\ol x} {\node[vertex] () at (\x,2){$\y$};}
\udotted{1.5}{2.5}
\udotted{6.5}{7.5}
\ddotted{1.5}{2.5}
\ddotted{6.5}{7.5}
\node () at (4.5,-2) {$=$};
\end{scope}
\begin{scope}[shift={(0,-8)}]
\uvs{1,3,4,5,6,8}
\lvs{1,3,4,5,6,8}
\foreach \x in {1,3,4,5,6,8} {\stline\x\x}
\foreach \x/\y in {4/\ol x,5/\ol x} {\node[vertex] () at (\x,2){$\y$};}
\udotted{1.5}{2.5}
\udotted{6.5}{7.5}
\ddotted{1.5}{2.5}
\ddotted{6.5}{7.5}
\foreach \x/\y in {1/1,4/i,8/n} {\node () at (\x,2.5){\footnotesize $\y$};}
\end{scope}
\begin{scope}[shift={(0,-10)}]
\uvs{1,3,4,5,6,8}
\lvs{1,3,4,5,6,8}
\foreach \x in {1,3,4,6,8} {\stline\x\x}
\stline54
\udotted{1.5}{2.5}
\udotted{6.5}{7.5}
\ddotted{1.5}{2.5}
\ddotted{6.5}{7.5}
\end{scope}
\begin{scope}[shift={(13,0)}]
\uvs{1,3,4,5,6,8}
\lvs{1,3,4,5,6,8}
\foreach \x in {1,3,4,6,8} {\stline\x\x}
\stline54
\udotted{1.5}{2.5}
\udotted{6.5}{7.5}
\ddotted{1.5}{2.5}
\ddotted{6.5}{7.5}
\foreach \x/\y in {1/1,4/i,8/n} {\node () at (\x,2.5){\footnotesize $\y$};}
\end{scope}
\begin{scope}[shift={(13,-2)}]
\uvs{1,3,4,5,6,8}
\lvs{1,3,4,5,6,8}
\foreach \x in {1,3,4,5,6,8} {\stline\x\x}
\foreach \x/\y in {5/\ol x} {\node[vertex] () at (\x,2){$\y$};}
\udotted{1.5}{2.5}
\udotted{6.5}{7.5}
\ddotted{1.5}{2.5}
\ddotted{6.5}{7.5}
\node () at (4.5,-2) {$=$};
\end{scope}
\begin{scope}[shift={(13,-8)}]
\uvs{1,3,4,5,6,8}
\lvs{1,3,4,5,6,8}
\foreach \x in {1,3,4,6,8} {\stline\x\x}
\stline54
\udotted{1.5}{2.5}
\udotted{6.5}{7.5}
\ddotted{1.5}{2.5}
\ddotted{6.5}{7.5}
\foreach \x/\y in {1/1,4/i,8/n} {\node () at (\x,2.5){\footnotesize $\y$};}
\end{scope}
\end{tikzpicture}
}
\caption{The $k=i$ (left) and $k=i+1$ (right) cases of relation \eqref{eq:lx}.}
\label{fig:lx}
\end{center}
\end{figure}

\begin{rem}
The relations in Theorem \ref{thm:MwrPTn} contain a fair few redundancies.  For example, the following relations could all be deleted:
\ben
\item \label{red1} the $k=i+1$ (or $k=i$) cases of \eqref{eq:sx} and \eqref{eq:st},
\item \label{red2} all of \eqref{eq:lx} and \eqref{eq:lt}, or alternatively all of \eqref{eq:rx} and \eqref{eq:rt}.
\een
Indeed, for \ref{red1}, we use $R_G$ and the $k=i$ case of \eqref{eq:st} to transform $s_it_{i+1}$ into $t_is_{i+1}$ as follows:
\[
s_it_{i+1} \to s_it_{i+1}s_is_i \to s_is_it_is_i \to t_is_i.
\]
(A similar calculation works for $s_ix^{(i+1)}\to x^{(i)}s_i$.)  As an example calculation for \ref{red2}, we transform $\lam_it_i$ into $t_it_{i+1}\lam_i$ using $R_T$, \eqref{eq:st} and \eqref{eq:rt} as follows:
\[
\lam_it_i \to \rho_is_it_i \to \rho_it_{i+1}s_i \to t_it_{i+1}\rho_is_i \to t_it_{i+1}\lam_i.
\]
\end{rem}

Each of the remaining results will involve a restriction of the surmorphism $\Psi$ from Theorem~\ref{thm:MwrPTn}, and a subset of the relations \eqref{eq:sx}--\eqref{eq:tr}.  The proofs of the next two results are exactly as for Theorem \ref{thm:MwrPTn}, but with a little less to check.

\begin{thm}\label{thm:MwrGn} 
For $n\geq2$, the wreath product $M\wr\G_n$ has presentation
\[
\Mpres{X_{M^n}\cup X_G}{R_{M^n}\cup R_G\cup R}
\]
via $\Psi\restr_{(X_{M^n}\cup X_G)^*}$, where $R$ consists of the relations \eqref{eq:sx}.  \epfres
\end{thm}

\begin{thm}\label{thm:MwrTn} 
For $n\geq2$, the wreath product $M\wr\T_n$ has presentation
\[
\Mpres{X_{M^n}\cup X_T}{R_{M^n}\cup R_T\cup R}
\]
via $\Psi\restr_{(X_{M^n}\cup X_T)^*}$, where $R$ consists of the relations \eqref{eq:sx}--\eqref{eq:rx}.  \epfres
\end{thm}

For the final result, we need to work just a little harder.

\begin{thm}\label{thm:MwrIn} 
For $n\geq2$, the wreath product $M\wr\I_n$ has presentation
\[
\Mpres{X_{M_0^n}\cup X_G}{R_{M_0^n}\cup R_G\cup R}
\]
via $\Psi\restr_{(X_{M_0^n}\cup X_G)^*}$, where $R$ consists of \eqref{eq:sx}, \eqref{eq:st}, and the following relations:
\begin{equation}
\label{eq:ts} t_it_{i+1}s_i = t_it_{i+1}. 
\end{equation}
\end{thm}

\pf
As in the proof of Theorem \ref{thm:MwrPTn}, we quickly obtain the presentation
\[
\Mpres{X_{M_0^n}\cup X_G}{R_{M_0^n}\cup R_G\cup R_1\cup R_\Om},
\]
where:
\bit
\item $R_1$ consists of relations \eqref{eq:sx} and \eqref{eq:st}, and
\item $R_\Om = \bigset{(N(\ba),N(\bb))}{(\ba,\bb)\in\Om}$ for any generating set $\Om$ for the congruence $\th$ on $M_0^n\rtimes\G_n$.
\eit
This time, Lemmas \ref{lem:Sagens}, \ref{lem:thav2} and \ref{lem:Om6}\ref{Om61} tell us that $\th=\Om^\sharp$ for
\[
\Om = \bigset{((\bv_{ij},\id_\bn),(\bv_{ij},\tau_{ij}))}{\oijn},
\]
where $\bv_{ij}=\bv_i\bv_j=\bone_{\bn\sm\{i,j\}}\equiv\id_{\bn\sm\{i,j\}} = (t_it_j)\Psi$ for each $\oijn$, and $\tau_{ij}$ is the transposition $(i,j)$.  Noting that $\tau_{ij}=((s_{j-1}\cdots s_{i+1})s_i(s_{i+1}\cdots s_{j-1}))\phi_G$, the resulting set $R_\Om$ then consists of the relations 
\begin{equation}\label{eq:titj}
t_it_j = t_it_j \cdot (s_{j-1}\cdots s_{i+1})s_i(s_{i+1}\cdots s_{j-1}) \qquad\text{for each $\oijn$,}
\end{equation}
and we note that \eqref{eq:ts} contains only those with $j=i+1$.  But it is easy to see that \eqref{eq:st} and~\eqref{eq:ts} together imply \eqref{eq:titj}.  Indeed, writing $\sim$ for the congruence generated by the claimed set of relations $R_{M_0^n}\cup R_G\cup R$, it is easy to see that
\begin{equation}\label{eq:titjs}
t_i ( s_{j-1}\cdots s_{i+1}) \sim (s_{j-1}\cdots s_{i+1} ) t_i \AND t_j ( s_{j-1}\cdots s_{i+1}) \sim (s_{j-1}\cdots s_{i+1}) t_{i+1},
\end{equation}
and then
\begin{align*}
t_it_j  (s_{j-1}\cdots s_{i+1})s_i(s_{i+1}\cdots s_{j-1}) &\sim  (s_{j-1}\cdots s_{i+1})t_it_{i+1} s_i(s_{i+1}\cdots s_{j-1}) &&\text{by \eqref{eq:titjs}}\\
&\sim  (s_{j-1}\cdots s_{i+1})t_it_{i+1} (s_{i+1}\cdots s_{j-1})  &&\text{by \eqref{eq:ts}}\\
&\sim t_it_j  (s_{j-1}\cdots s_{i+1}) (s_{i+1}\cdots s_{j-1})  &&\text{by \eqref{eq:titjs}}\\
&\sim t_it_j  &&\text{by $R_G$.}  \qedhere
\end{align*}
\epf

\begin{rem}
Relations \eqref{eq:ts} are not included in the presentation from Theorem \ref{thm:MwrPTn}, even though these relations obviously hold in the monoid $M\wr\PT_n$.  We leave it as an exercise to show that \eqref{eq:ts} follows from the relations listed in Theorem \ref{thm:MwrPTn}.
\end{rem}

\footnotesize
\def\bibspacing{-1.1pt}
\bibliography{biblio}

\def\cprime{$'$}
\begin{thebibliography}{100}

\bibitem{Aizenstat1958}
A.~J. A{\u\i}zen{\v{s}}tat.
\newblock Defining relations of finite symmetric semigroups (in {R}ussian).
\newblock {\em Mat. Sb. N.S.}, 45 (87):261--280, 1958.

\bibitem{Aizenstat1962}
A.~J. A{\u\i}zen{\v{s}}tat.
\newblock The defining relations of the endomorphism semigroup of a finite
  linearly ordered set.
\newblock {\em Sibirsk. Mat. \u Z.}, 3:161--169, 1962.

\bibitem{APW1992}
J.~Almeida, J.-E. Pin, and P.~Weil.
\newblock Semigroups whose idempotents form a subsemigroup.
\newblock {\em Math. Proc. Cambridge Philos. Soc.}, 111(2):241--253, 1992.

\bibitem{Artin1925}
E.~Artin.
\newblock Theorie der {Z}\"{o}pfe.
\newblock {\em Abh. Math. Sem. Univ. Hamburg}, 4(1):47--72, 1925.

\bibitem{Artin1947}
E.~Artin.
\newblock Theory of braids.
\newblock {\em Ann. of Math. (2)}, 48:101--126, 1947.

\bibitem{Ault1974}
J.~E. Ault.
\newblock Semigroups with midunits.
\newblock {\em Trans. Amer. Math. Soc.}, 190:375--384, 1974.

\bibitem{BP1985}
J.~Berstel and D.~Perrin.
\newblock {\em Theory of codes}, volume 117 of {\em Pure and Applied
  Mathematics}.
\newblock Academic Press, Inc., Orlando, FL, 1985.

\bibitem{Billhardt1997}
B.~Billhardt.
\newblock Extensions of semilattices by left type-{$A$} semigroups.
\newblock {\em Glasgow Math. J.}, 39(1):7--16, 1997.

\bibitem{BGG2010}
M.~J.~J. Branco, G.~M.~S. Gomes, and V.~Gould.
\newblock Extensions and covers for semigroups whose idempotents form a left
  regular band.
\newblock {\em Semigroup Forum}, 81(1):51--70, 2010.

\bibitem{BGG2011}
M.~J.~J. Branco, G.~M.~S. Gomes, and V.~Gould.
\newblock Left adequate and left {E}hresmann monoids.
\newblock {\em Internat. J. Algebra Comput.}, 21(7):1259--1284, 2011.

\bibitem{Brookes2020}
M.~D. G.~K. Brookes.
\newblock Congruences on the partial automorphism monoid of a free group
  action.
\newblock {\em Internat. J. Algebra Comput.}, 31(6):1147--1176, 2021.

\bibitem{BS1981}
S.~Burris and H.~P. Sankappanavar.
\newblock {\em A course in universal algebra}, volume~78 of {\em Graduate Texts
  in Mathematics}.
\newblock Springer-Verlag, New York-Berlin, 1981.

\bibitem{CFP1996}
J.~W. Cannon, W.~J. Floyd, and W.~R. Parry.
\newblock Introductory notes on {R}ichard {T}hompson's groups.
\newblock {\em Enseign. Math. (2)}, 42(3-4):215--256, 1996.

\bibitem{CH1974}
S.~Y. Chen and S.~C. Hsieh.
\newblock Factorizable inverse semigroups.
\newblock {\em Semigroup Forum}, 8(4):283--297, 1974.

\bibitem{CE2022}
C.~Clark and J.~East.
\newblock Presentations for wreath products involving symmetric inverse monoids
  and categories.
\newblock {\em Preprint}, 2022, {\tt arXiv:2204.06992}.

\bibitem{CPbook}
A.~H. Clifford and G.~B. Preston.
\newblock {\em The algebraic theory of semigroups. {V}ol. {I}}.
\newblock Mathematical Surveys, No. 7. American Mathematical Society,
  Providence, R.I., 1961.

\bibitem{CPbook2}
A.~H. Clifford and G.~B. Preston.
\newblock {\em The algebraic theory of semigroups. {V}ol. {II}}.
\newblock Mathematical Surveys, No. 7. American Mathematical Society,
  Providence, R.I., 1967.

\bibitem{DGY2015}
Y.~Dandan, I.~Dolinka, and V.~Gould.
\newblock Free idempotent generated semigroups and endomorphism monoids of free
  {$G$}-acts.
\newblock {\em J. Algebra}, 429:133--176, 2015.

\bibitem{DE2018}
I.~Dolinka and J.~East.
\newblock Semigroups of rectangular matrices under a sandwich operation.
\newblock {\em Semigroup Forum}, 96(2):253--300, 2018.

\bibitem{Sandwiches1}
I.~Dolinka, I.~{\DJ}ur{\dj}ev, J.~East, P.~Honyam, K.~Sangkhanan, J.~Sanwong,
  and W.~Sommanee.
\newblock Sandwich semigroups in locally small categories {I}: foundations.
\newblock {\em Algebra Universalis}, 79(3):Art. 75, 35 pp, 2018.

\bibitem{DR2009}
E.~R. Dombi and N.~Ru{\v{s}}kuc.
\newblock On generators and presentations of semidirect products in inverse
  semigroups.
\newblock {\em Bull. Aust. Math. Soc.}, 79(3):353--365, 2009.

\bibitem{Dubreil1941}
P.~Dubreil.
\newblock Contribution \`a la th\'{e}orie des demi-groupes.
\newblock {\em M\'{e}m. Acad. Sci. Inst. France (2)}, 63(3):52, 1941.

\bibitem{Dyck1882}
W.~Dyck.
\newblock Gruppentheoretische {S}tudien.
\newblock {\em Math. Ann.}, 20(1):1--44, 1882.

\bibitem{EEF2005}
D.~Easdown, J.~East, and D.~G. FitzGerald.
\newblock Presentations of factorizable inverse monoids.
\newblock {\em Acta Sci. Math. (Szeged)}, 71(3-4):509--520, 2005.

\bibitem{EL2004}
D.~Easdown and T.~G. Lavers.
\newblock The inverse braid monoid.
\newblock {\em Adv. Math.}, 186(2):438--455, 2004.

\bibitem{JEinsn}
J.~East.
\newblock A presentation of the singular part of the symmetric inverse monoid.
\newblock {\em Comm. Algebra}, 34(5):1671--1689, 2006.

\bibitem{East2007b}
J.~East.
\newblock Braids and partial permutations.
\newblock {\em Adv. Math.}, 213(1):440--461, 2007.

\bibitem{East2007}
J.~East.
\newblock Vines and partial transformations.
\newblock {\em Adv. Math.}, 216(2):787--810, 2007.

\bibitem{JEtnsn}
J.~East.
\newblock A presentation for the singular part of the full transformation
  semigroup.
\newblock {\em Semigroup Forum}, 81(2):357--379, 2010.

\bibitem{JEptnsn}
J.~East.
\newblock Presentations for singular subsemigroups of the partial
  transformation semigroup.
\newblock {\em Internat. J. Algebra Comput.}, 20(1):1--25, 2010.

\bibitem{JEtnsn2}
J.~East.
\newblock Defining relations for idempotent generators in finite full
  transformation semigroups.
\newblock {\em Semigroup Forum}, 86(3):451--485, 2013.

\bibitem{JEptnsn2}
J.~East.
\newblock Defining relations for idempotent generators in finite partial
  transformation semigroups.
\newblock {\em Semigroup Forum}, 89(1):72--76, 2014.

\bibitem{JEinsn2}
J.~East.
\newblock A symmetrical presentation for the singular part of the symmetric
  inverse monoid.
\newblock {\em Algebra Universalis}, 74(3-4):207--228, 2015.

\bibitem{JE2019a}
J.~East.
\newblock Idempotents and one-sided units in infinite partial {B}rauer monoids.
\newblock {\em J. Algebra}, 534:427--482, 2019.

\bibitem{JE2020a}
J.~East.
\newblock Idempotents and one-sided units: lattice invariants and a semigroup
  of functors on the category of monoids.
\newblock {\em J. Algebra}, 560:1219--1252, 2020.

\bibitem{JEptc}
J.~East.
\newblock Presentations for tensor categories.
\newblock {\em Preprint}, 2020, {\tt arXiv:2005.01953}.

\bibitem{East2021}
J.~East.
\newblock Presentations for {$\mathbb P^K$}.
\newblock {\em Monatsh. Math.}, to appear, {\tt arXiv:2105.06127}.

\bibitem{EF2012}
J.~East and D.~G. FitzGerald.
\newblock The semigroup generated by the idempotents of a partition monoid.
\newblock {\em J. Algebra}, 372:108--133, 2012.

\bibitem{EF2010}
B.~Everitt and J.~Fountain.
\newblock Partial symmetry, reflection monoids and {C}oxeter groups.
\newblock {\em Adv. Math.}, 223(5):1782--1814, 2010.

\bibitem{EF2013}
B.~Everitt and J.~Fountain.
\newblock Partial mirror symmetry, lattice presentations and algebraic monoids.
\newblock {\em Proc. Lond. Math. Soc. (3)}, 107(2):414--450, 2013.

\bibitem{FADEG2019}
Y.-Y. Feng, A.~Al-Aadhami, I.~Dolinka, J.~East, and V.~Gould.
\newblock Presentations for singular wreath products.
\newblock {\em J. Pure Appl. Algebra}, 223(12):5106--5146, 2019.

\bibitem{FitzGerald2010}
D.~G. FitzGerald.
\newblock Factorizable inverse monoids.
\newblock {\em Semigroup Forum}, 80(3):484--509, 2010.

\bibitem{Fountain1977}
J.~Fountain.
\newblock A class of right {PP} monoids.
\newblock {\em Quart. J. Math. Oxford Ser. (2)}, 28(111):285--300, 1977.

\bibitem{Fountain1979}
J.~Fountain.
\newblock Adequate semigroups.
\newblock {\em Proc. Edinburgh Math. Soc. (2)}, 22(2):113--125, 1979.

\bibitem{Fountain1991b}
J.~Fountain.
\newblock Free right type {$A$} semigroups.
\newblock {\em Glasgow Math. J.}, 33(2):135--148, 1991.

\bibitem{FG1993}
J.~Fountain and G.~M.~S. Gomes.
\newblock Proper left type-{$A$} monoids revisited.
\newblock {\em Glasgow Math. J.}, 35(3):293--306, 1993.

\bibitem{FGG2009}
J.~Fountain, G.~M.~S. Gomes, and V.~Gould.
\newblock The free ample monoid.
\newblock {\em Internat. J. Algebra Comput.}, 19(4):527--554, 2009.

\bibitem{FL1993}
J.~Fountain and A.~Lewin.
\newblock Products of idempotent endomorphisms of an independence algebra of
  infinite rank.
\newblock {\em Math. Proc. Cambridge Philos. Soc.}, 114(2):303--319, 1993.

\bibitem{FPW2004}
J.~Fountain, J.-E. Pin, and P.~Weil.
\newblock Covers for monoids.
\newblock {\em J. Algebra}, 271(2):529--586, 2004.

\bibitem{GAP}
The GAP~Group.
\newblock {\em {GAP -- Groups, Algorithms, and Programming}}.

\bibitem{GG1999}
G.~M.~S. Gomes and V.~Gould.
\newblock Proper weakly left ample semigroups.
\newblock {\em Internat. J. Algebra Comput.}, 9(6):721--739, 1999.

\bibitem{GG2000}
G.~M.~S. Gomes and V.~Gould.
\newblock Graph expansions of unipotent monoids.
\newblock {\em Comm. Algebra}, 28(1):447--463, 2000.

\bibitem{GG2001}
G.~M.~S. Gomes and V.~Gould.
\newblock Fundamental {E}hresmann semigroups.
\newblock {\em Semigroup Forum}, 63(1):11--33, 2001.

\bibitem{Gould1995}
V.~Gould.
\newblock Independence algebras.
\newblock {\em Algebra Universalis}, 33(3):294--318, 1995.

\bibitem{Gould1996}
V.~Gould.
\newblock Graph expansions of right cancellative monoids.
\newblock {\em Internat. J. Algebra Comput.}, 6(6):713--733, 1996.

\bibitem{Gould_notes}
V.~Gould.
\newblock Notes on restriction semigroups and related structures.
\newblock \verb|https://www-users.york.ac.uk/~varg1/restriction.pdf|, 2010.

\bibitem{GS2013}
V.~Gould and M.~B. Szendrei.
\newblock Proper restriction semigroups---semidirect products and
  {$W$}-products.
\newblock {\em Acta Math. Hungar.}, 141(1-2):36--57, 2013.

\bibitem{Green1951}
J.~A. Green.
\newblock On the structure of semigroups.
\newblock {\em Ann. of Math. (2)}, 54:163--172, 1951.

\bibitem{Hamilton1856}
W.~R. Hamilton.
\newblock Memorandum respecting a new system of roots of unity.
\newblock {\em Philos. Mag.}, 12(81):446--446, 1856.

\bibitem{Hickey1983}
J.~B. Hickey.
\newblock Semigroups under a sandwich operation.
\newblock {\em Proc. Edinburgh Math. Soc. (2)}, 26(3):371--382, 1983.

\bibitem{Hickey1986}
J.~B. Hickey.
\newblock On variants of a semigroup.
\newblock {\em Bull. Austral. Math. Soc.}, 34(3):447--459, 1986.

\bibitem{Higgins1992}
P.~M. Higgins.
\newblock {\em Techniques of semigroup theory}.
\newblock Oxford Science Publications. The Clarendon Press, Oxford University
  Press, New York, 1992.

\bibitem{HS2013}
P.~Honyam and J.~Sanwong.
\newblock Semigroups of transformations with fixed sets.
\newblock {\em Quaest. Math.}, 36(1):79--92, 2013.

\bibitem{Howie1995}
J.~M. Howie.
\newblock {\em Fundamentals of semigroup theory}, volume~12 of {\em London
  Mathematical Society Monographs. New Series}.
\newblock The Clarendon Press, Oxford University Press, New York, 1995.
\newblock Oxford Science Publications.

\bibitem{HR1994}
J.~M. Howie and N.~Ru{\v{s}}kuc.
\newblock Constructions and presentations for monoids.
\newblock {\em Comm. Algebra}, 22(15):6209--6224, 1994.

\bibitem{Humphreys1990}
J.~E. Humphreys.
\newblock {\em Reflection groups and {C}oxeter groups}, volume~29 of {\em
  Cambridge Studies in Advanced Mathematics}.
\newblock Cambridge University Press, Cambridge, 1990.

\bibitem{Jones1989}
P.~R. Jones.
\newblock Exchange properties and basis properties for closure operators.
\newblock {\em Colloq. Math.}, 57(1):29--33, 1989.

\bibitem{Kambites2011}
M.~Kambites.
\newblock Free adequate semigroups.
\newblock {\em J. Aust. Math. Soc.}, 91(3):365--390, 2011.

\bibitem{KKM2000}
M.~Kilp, U.~Knauer, and A.~V. Mikhalev.
\newblock {\em Monoids, acts and categories}, volume~29 of {\em de Gruyter
  Expositions in Mathematics}.
\newblock Walter de Gruyter \& Co., Berlin, 2000.
\newblock With applications to wreath products and graphs, A handbook for
  students and researchers.

\bibitem{KM1980a}
U.~Knauer and A.~Mikhalev.
\newblock Endomorphism monoids of free acts and {$0$}-wreath products of
  monoids. {I}. {A}nnihilator properties.
\newblock {\em Semigroup Forum}, 19(2):177--187, 1980.

\bibitem{KM1980b}
U.~Knauer and A.~Mikhalev.
\newblock Endomorphism monoids of free acts and {$0$}-wreath products of
  monoids. {II}. {R}egularity.
\newblock {\em Semigroup Forum}, 19(3):189--198, 1980.

\bibitem{KM1980c}
U.~Knauer and A.~Mikhalev.
\newblock Endomorphism monoids of free acts and {$0$}-wreath products of
  monoids. {III}. {S}tandard involution and continuous endomorphisms.
\newblock {\em Semigroup Forum}, 19(4):355--369, 1980.

\bibitem{Lavers1997}
T.~G. Lavers.
\newblock The theory of vines.
\newblock {\em Comm. Algebra}, 25(4):1257--1284, 1997.

\bibitem{Lavers1998}
T.~G. Lavers.
\newblock Presentations of general products of monoids.
\newblock {\em J. Algebra}, 204(2):733--741, 1998.

\bibitem{Lawson1986}
M.~V. Lawson.
\newblock The structure of type {$A$} semigroups.
\newblock {\em Quart. J. Math. Oxford Ser. (2)}, 37(147):279--298, 1986.

\bibitem{Lawson1994}
M.~V. Lawson.
\newblock Almost factorisable inverse semigroups.
\newblock {\em Glasgow Math. J.}, 36(1):97--111, 1994.

\bibitem{Lawson1998}
M.~V. Lawson.
\newblock {\em Inverse semigroups}.
\newblock World Scientific Publishing Co., Inc., River Edge, NJ, 1998.
\newblock The theory of partial symmetries.

\bibitem{LM2007}
M.~V. Lawson and S.~W. Margolis.
\newblock In {M}c{A}lister's footsteps: a random ramble around the
  {$P$}-theorem.
\newblock In {\em Semigroups and formal languages}, pages 145--163. World Sci.
  Publ., Hackensack, NJ, 2007.

\bibitem{Lima_thesis}
L.~Lima.
\newblock {\em The Local Automorphism Monoid of an Independence Algebra}.
\newblock PhD thesis, University of York, 1994.

\bibitem{Marczewski1959}
E.~Marczewski.
\newblock Independence in algebras of sets and {B}oolean algebras.
\newblock {\em Fund. Math.}, 48:135--145, 1959/60.

\bibitem{MSS2014}
S.~Margolis, F.~Saliola, and B.~Steinberg.
\newblock Semigroups embeddable in hyperplane face monoids.
\newblock {\em Semigroup Forum}, 89(1):236--248, 2014.

\bibitem{MSS2015}
S.~Margolis, F.~Saliola, and B.~Steinberg.
\newblock Combinatorial topology and the global dimension of algebras arising
  in combinatorics.
\newblock {\em J. Eur. Math. Soc. (JEMS)}, 17(12):3037--3080, 2015.

\bibitem{MP1987}
S.~W. Margolis and J.-E. Pin.
\newblock Inverse semigroups and extensions of groups by semilattices.
\newblock {\em J. Algebra}, 110(2):277--297, 1987.

\bibitem{Markov1947}
A.~Markov.
\newblock The impossibility of certain algorithms in the theory of associative
  systems. {II}.
\newblock {\em Doklady Akad. Nauk SSSR (N.S.)}, 58:353--356, 1947.

\bibitem{MS2021}
V.~Mazorchuk and S.~Srivastava.
\newblock Jucys-murphy elements and grothendieck groups for generalized rook
  monoids.
\newblock {\em Preprint}, 2021, {\tt arXiv:2104.13632}.

\bibitem{McAlister1974}
D.~B. McAlister.
\newblock Groups, semilattices and inverse semigroups.
\newblock {\em Trans. Amer. Math. Soc.}, 192:227--244, 1974.

\bibitem{McAlister1974b}
D.~B. McAlister.
\newblock Groups, semilattices and inverse semigroups {II}.
\newblock {\em Trans. Amer. Math. Soc.}, 196:351--370, 1974.

\bibitem{McAlister1976}
D.~B. McAlister.
\newblock Some covering and embedding theorems for inverse semigroups.
\newblock {\em J. Austral. Math. Soc. Ser. A}, 22(2):188--211, 1976.

\bibitem{MMT1987}
R.~N. McKenzie, G.~F. McNulty, and W.~F. Taylor.
\newblock {\em Algebras, lattices, varieties. {V}ol. {I}}.
\newblock The Wadsworth \& Brooks/Cole Mathematics Series. Wadsworth \&
  Brooks/Cole Advanced Books \& Software, Monterey, CA, 1987.

\bibitem{Meldrum1995}
J.~D.~P. Meldrum.
\newblock {\em Wreath products of groups and semigroups}, volume~74 of {\em
  Pitman Monographs and Surveys in Pure and Applied Mathematics}.
\newblock Longman, Harlow, 1995.

\bibitem{Semigroups}
J.~D. Mitchell et~al.
\newblock {\em Semigroups - GAP package}.

\bibitem{Moore1897}
E.~H. Moore.
\newblock Concerning the abstract groups of order $k!$ and $\tfrac{1}{2}k!$
  holohedrically isomorphic with the symmetric and the alternating
  substitution-groups on $k$ letters.
\newblock {\em Proc. London Math. Soc.}, 28(1):357--366, 1897.

\bibitem{Munn1970}
W.~D. Munn.
\newblock Fundamental inverse semigroups.
\newblock {\em Quart. J. Math. Oxford Ser. (2)}, 21:157--170, 1970.

\bibitem{Munn1974}
W.~D. Munn.
\newblock Free inverse semigroups.
\newblock {\em Proc. London Math. Soc. (3)}, 29:385--404, 1974.

\bibitem{Munn1976}
W.~D. Munn.
\newblock A note on {$E$}-unitary inverse semigroups.
\newblock {\em Bull. London Math. Soc.}, 8(1):71--76, 1976.

\bibitem{Narkiewicz1961}
W.~Narkiewicz.
\newblock Independence in a certain class of abstract algebras.
\newblock {\em Fund. Math.}, 50:333--340, 1961/62.

\bibitem{NB2021}
C.-F. Nyberg-Brodda.
\newblock The {B} {B} {N}ewman spelling theorem.
\newblock {\em Br. J. Hist. Math.}, 36(2):117--131, 2021.

\bibitem{OCarroll1976}
L.~O'Carroll.
\newblock Embedding theorems for proper inverse semigroups.
\newblock {\em J. Algebra}, 42(1):26--40, 1976.

\bibitem{PR1979}
M.~Petrich and N.~R. Reilly.
\newblock A representation of {$E$}-unitary inverse semigroups.
\newblock {\em Quart. J. Math. Oxford Ser. (2)}, 30(119):339--350, 1979.

\bibitem{Popova1961}
L.~M. Popova.
\newblock Defining relations in some semigroups of partial transformations of a
  finite set (in {R}ussian).
\newblock {\em Uchenye Zap. Leningrad Gos. Ped. Inst.}, 218:191--212, 1961.

\bibitem{Popova1962}
L.~M. Popova.
\newblock Defining relations of a semigroup of partial endomorphisms of a
  finite linearly ordered set (in {R}ussian).
\newblock {\em Leningrad. Gos. Ped. Inst. U\v cen. Zap.}, 238:78--88, 1962.

\bibitem{Post1947}
E.~L. Post.
\newblock Recursive unsolvability of a problem of {T}hue.
\newblock {\em J. Symbolic Logic}, 12:1--11, 1947.

\bibitem{Preston1986}
G.~B. Preston.
\newblock Semidirect products of semigroups.
\newblock {\em Proc. Roy. Soc. Edinburgh Sect. A}, 102(1-2):91--102, 1986.

\bibitem{RSbook}
J.~Rhodes and B.~Steinberg.
\newblock {\em The {$q$}-theory of finite semigroups}.
\newblock Springer Monographs in Mathematics. Springer, New York, 2009.

\bibitem{RRT2003}
E.~F. Robertson, N.~Ru{\v{s}}kuc, and M.~R. Thomson.
\newblock Finite generation and presentability of wreath products of monoids.
\newblock {\em J. Algebra}, 266(2):382--392, 2003.

\bibitem{RRW1998}
E.~F. Robertson, N.~Ru\v{s}kuc, and J.~Wiegold.
\newblock Generators and relations of direct products of semigroups.
\newblock {\em Trans. Amer. Math. Soc.}, 350(7):2665--2685, 1998.

\bibitem{Ruskuc1999}
N.~Ru\v{s}kuc.
\newblock Presentations for subgroups of monoids.
\newblock {\em J. Algebra}, 220(1):365--380, 1999.

\bibitem{Scheiblich1973}
H.~E. Scheiblich.
\newblock Free inverse semigroups.
\newblock {\em Proc. Amer. Math. Soc.}, 38:1--7, 1973.

\bibitem{Skornjakov1979}
L.~A. Skornjakov.
\newblock Regularity of the wreath product of monoids.
\newblock {\em Semigroup Forum}, 18(1):83--86, 1979.

\bibitem{Specht1933}
W.~Specht.
\newblock Eine {V}erallgemeinerung der {P}ermutationsgruppen.
\newblock {\em Math. Z.}, 37(1):321--341, 1933.

\bibitem{Steinberg2003}
B.~Steinberg.
\newblock Mc{A}lister's {$P$}-theorem via {S}ch\"{u}tzenberger graphs.
\newblock {\em Comm. Algebra}, 31(9):4387--4392, 2003.

\bibitem{Szendrei1980}
M.~B. Szendrei.
\newblock On a pullback diagram for orthodox semigroups.
\newblock {\em Semigroup Forum}, 20(1):1--10, 1980.

\bibitem{Szendrei1982}
M.~B. Szendrei.
\newblock Correction and supplement to: ``{O}n a pullback diagram for orthodox
  semigroups'' [{S}emigroup {F}orum {\bf 20} (1980), no. 1, 1--10; {MR}
  81m:20085].
\newblock {\em Semigroup Forum}, 25(3-4):311--324, 1982.

\bibitem{Szendrei2013}
M.~B. Szendrei.
\newblock Embedding into almost left factorizable restriction semigroups.
\newblock {\em Comm. Algebra}, 41(4):1458--1483, 2013.

\bibitem{Trotter1995}
P.~G. Trotter.
\newblock Covers for regular semigroups and an application to complexity.
\newblock {\em J. Pure Appl. Algebra}, 105(3):319--328, 1995.

\bibitem{Turing1950}
A.~M. Turing.
\newblock The word problem in semi-groups with cancellation.
\newblock {\em Ann. of Math. (2)}, 52:491--505, 1950.

\bibitem{Urbanik1966}
K.~Urbanik.
\newblock Linear independence in abstract algebras.
\newblock {\em Colloq. Math.}, 14:233--255, 1966.

\bibitem{Yamada1955}
M.~Yamada.
\newblock A note on middle unitary semigroups.
\newblock {\em K\=odai Math. Sem. Rep.}, 7:49--52, 1955.

\bibitem{Zenab2018}
R.~Zenab.
\newblock Algebraic properties of {Z}appa-{S}z\'{e}p products of semigroups and
  monoids.
\newblock {\em Semigroup Forum}, 96(2):316--332, 2018.

\end{thebibliography}
\bibliographystyle{abbrv}

\end{document}